\documentclass[11pt,twoside]{article}

\usepackage{fullpage}

\usepackage{epsf}
\usepackage{fancyhdr}
\usepackage{graphics}
\usepackage{graphicx}
\usepackage{psfrag}
\usepackage{microtype}
\usepackage{subfigure}
\usepackage{algorithmic}
\usepackage{color,xcolor}
\usepackage[linesnumbered,ruled]{algorithm2e}
\DontPrintSemicolon
\usepackage{color}

\newenvironment{carlist}
 {\begin{list}{$\bullet$}
 {\setlength{\topsep}{0in} \setlength{\partopsep}{0in}
  \setlength{\parsep}{0in} \setlength{\itemsep}{\parskip}
  \setlength{\leftmargin}{0.07in} \setlength{\rightmargin}{0.08in}
  \setlength{\listparindent}{0in} \setlength{\labelwidth}{0.08in}
  \setlength{\labelsep}{0.1in} \setlength{\itemindent}{0in}}}
 {\end{list}}

\newcommand{\bcar}{\begin{carlist}}
\newcommand{\ecar}{\end{carlist}}

\usepackage{amsthm}
\usepackage{amsfonts}
\usepackage{amsmath}
\usepackage{amssymb,bbm}
\usepackage{comment}
\def\ball{\mathbb{B}}

\newcommand{\stepsize}{\eta}

\newcommand{\real}{\ensuremath{\mathbb{R}}}

\newcommand{\order}[1]{\ensuremath{\mathcal{O}\parenth{#1}}}

\newcommand{\thetahat}{\ensuremath{\widehat{\theta}}}

\newcommand{\thetabar}{\ensuremath{\bar{\theta}}}

\newcommand{\Xspace}{\ensuremath{\mathbb{X}}}

\newcommand{\parenth}[1]{\left( #1 \right)}

\newcommand{\abss}[1]{\left| #1 \right |}

\newcommand{\law}{\ensuremath{\mathcal{L}}}

\newcommand{\sphere}{\ensuremath{\mathbb{S}}}

\newcommand{\mydefn}{\ensuremath{:=}}

\newcommand{\defn}{:=}

\newcommand{\matsnorm}[2]{|\!|\!| #1 | \! | \!|_{{#2}}}
\newcommand{\vecnorm}[2]{\| #1\|_{#2}}

\newcommand{\opnorm}[1]{\ensuremath{\matsnorm{#1}{\tiny{\mbox{op}}}}}

\newcommand{\inprod}[2]{\ensuremath{\langle #1 , \, #2 \rangle}}

\newcommand{\kull}[2]{\ensuremath{D_{\text{KL}}(#1\; \| \; #2)}}

\newcommand{\Exs}{\ensuremath{{\mathbb{E}}}}
\newcommand{\Prob}{\ensuremath{{\mathbb{P}}}}

\DeclareMathOperator{\var}{var}
\DeclareMathOperator{\cov}{cov}
\DeclareMathOperator{\trace}{trace}

\newtheoremstyle{named}{}{}{\itshape}{}{\bfseries}{.}{.5em}{\thmnote{#3's }#1}
\theoremstyle{named}

\theoremstyle{plain}

\newtheorem{theorem}{Theorem}
\newtheorem{proposition}{Proposition}

\newtheorem{lemma}{Lemma}

\newtheorem{corollary}{Corollary}

\newlength{\widebarargwidth}
\newlength{\widebarargheight}
\newlength{\widebarargdepth}

\makeatletter
\long\def\@makecaption#1#2{
        \vskip 0.8ex
        \setbox\@tempboxa\hbox{\small {\bf #1:} #2}
        \parindent 1.5em  
        \dimen0=\hsize
        \advance\dimen0 by -3em
        \ifdim \wd\@tempboxa >\dimen0
                \hbox to \hsize{
                        \parindent 0em
                        \hfil
                        \parbox{\dimen0}{\def\baselinestretch{0.96}\small
                                {\bf #1.} #2
                                }
                        \hfil}
        \else \hbox to \hsize{\hfil \box\@tempboxa \hfil}
        \fi
        }
\makeatother

\long\def\comment#1{}
\definecolor{battleshipgrey}{rgb}{0.52, 0.52, 0.51}
\definecolor{darkgray}{rgb}{0.66, 0.66, 0.66}
\definecolor{darkgreen}{rgb}{0.0, 0.2, 0.13}
\definecolor{darkspringgreen}{rgb}{0.09, 0.45, 0.27}
\definecolor{dukeblue}{rgb}{0.0, 0.0, 0.61}
\definecolor{olivedrab7}{rgb}{0.24, 0.2, 0.12}
\definecolor{darkblue}{rgb}{0.0, 0.0, 0.55}
\definecolor{darkscarlet}{rgb}{0.34, 0.01, 0.1}
\definecolor{candyapplered}{rgb}{1.0, 0.03, 0.0}
\definecolor{ao(english)}{rgb}{0.0, 0.5, 0.0}
\definecolor{applegreen}{rgb}{0.55, 0.71, 0.0}

\usepackage{mathrsfs}

\usepackage{url}
\usepackage[colorlinks,linkcolor=magenta,citecolor=blue, pagebackref=true,backref=true]{hyperref}
\renewcommand*{\backrefalt}[4]{%
    \ifcase #1 \footnotesize{(Not cited.)}%
    \or        \footnotesize{(Cited on page~#2.)}%
    \else      \footnotesize{(Cited on pages~#2.)}%
    \fi}

\usepackage{nicefrac}
\usepackage{comment}

\usepackage{chngpage}

 \usepackage{tabularx}%

\usepackage{enumitem}
\usepackage{booktabs}
\usepackage{caption}

\usepackage{bm,bbm}
\usepackage{mathtools}

\usepackage{fullpage}

\newcommand{\myiid}{\ensuremath{\mathrm{i.i.d.}}}

\newcommand{\Bmat}{B} \newcommand{\coordinate}{e}
\newcommand{\Wass}{\mathcal{W}}

\newcommand{\filtration}{\mathcal{F}}

\newcommand{\testMat}{Q}

\newcommand{\numcold}{{n_c}}

\newcommand{\quadvar}[2]{[ #1 ]_{#2}} 

\newcommand{\varbound}{\bar{\sigma}}

\newcommand{\statespace}{\mathbb{X}}
\newcommand{\goodendex}{\ensuremath{\clubsuit}}

\newcommand{\vectorize}{\mathrm{vec}}

\newcommand{\Ltwospace}{\mathbb{L}^2}

\newtheorem{example}{Example} \newcommand{\metric}{\rho}

\newcommand{\pmax}{\bar{p}}

\newcommand{\support}{\mathrm{supp}}

\newcommand{\discount}{\gamma}

\newcommand{\lyap}{\Phi}
\newtheorem{assumption}{Assumption}
\newcommand{\numobs}{\ensuremath{n}} \newcommand{\contraction}{\gamma}

\newcommand{\usedim}{\ensuremath{d}}
\newcommand{\SigStar}{\ensuremath{\Sigma^*}}

 \newcommand{\reward}{r}

\newcommand{\totalvariation}{d_{\mathrm{TV}}}

\newcommand{\sigmaA}{\sigma_\Lmat} \newcommand{\sigmab}{\sigma_b}

 \newcommand{\lspace}{\LinSpace}
\newcommand{\LinSpace}{\mathbb{S}} 

\newcommand{\projecttolin}{\Pi_{\lspace}}

\newcommand{\Lbar}{{\bar{\Lmat}}} \newcommand{\martingale}{\Psi}

\newcommand{\matLip}{\sigmaA} \newcommand{\veclip}{\sigmab}
\newcommand{\bbar}{{\bar{\bvec}}}

\newcommand{\complex}{\mathbb{C}}

\newcommand{\transition}{P}
\newcommand{\stationary}{\xi}

\newcommand{\greenFunc}{\boldsymbol{g}}

\newcommand{\Event}{\mathscr{E}} \newcommand{\Term}{T}
 
\newcommand{\chisqdiv}[2]{\chi^2 \left( #1 ~||~ #2 \right)}
\long\def\comment#1{} \newcommand{\maxdiv}[2]{D_{\infty} \left( #1
  ~||~ #2 \right)} \long\def\comment#1{}

 \newcommand{\state}{s}

\newcommand{\substationary}{}
\newcommand{\statnorm}[1]{\vecnorm{#1}{\substationary}}
\newcommand{\statinprod}[2]{\inprod{#1}{#2}_{\substationary}}

\newcommand{\Amat}{\ensuremath{A}} \newcommand{\bvec}{\ensuremath{b}}

\newcommand{\mixingtime}{t_{\mathrm{mix}}} \newcommand{\valuefunc}{V}
\newcommand{\valuestar}{\valuefunc^*}
\newcommand{\valuebar}{\bar{\valuefunc}}
\newcommand{\valuehat}{\widehat{\valuefunc}}
 
\newcommand{\slidingwindow}{\omega}

\newcommand{\Lmat}{\ensuremath{L}} 

\newcommand{\multnoiseMG}{Z} \newcommand{\addnoiseMG}{\zeta}
\newcommand{\multnoiseMarkov}{N} \newcommand{\addnoiseMarkov}{\nu}
\newcommand{\conmax}{\ensuremath{\contraction_{\tiny{\operatorname{max}}}}}

\newcommand{\numburn}{\ensuremath{{\numobs_0}}}

\newcommand{\SpecMat}{\ensuremath{M}}

\newcommand{\neighborhood}{\mathfrak{N}}

\newcommand{\lammax}{\ensuremath{\lambda_{\mbox{\tiny{max}}}}}

\newcommand{\crosscorrelation}{\mu}
\newcommand{\bigquantity}{\mathfrak{H}}
\newcommand{\crosstermdebias}{\lambda} \newcommand{\perturb}{h}

\newcommand{\Htil}{\ensuremath{\widetilde{H}}}

\newcommand{\statnullspace}{\mathbb{H}}
\newcommand{\greenOp}{\mathcal{A}}
\newcommand{\transitionOp}{\mathcal{P}}

\newcommand{\Lmap}{\ensuremath{\boldsymbol{L}}}
\newcommand{\bmap}{\ensuremath{\boldsymbol{b}}}

\newcommand{\myorder}{\ensuremath{\mathcal{O}}}

\newcommand{\AmatInd}[1]{\ensuremath{A_{#1}}}
\newcommand{\AstarInd}[1]{\ensuremath{A^*_{#1}}}

\newcommand{\RandOp}{\ensuremath{\mathcal{T}}}
\newcommand{\RandOpBar}{\ensuremath{\bar{\RandOp}}}

\newcommand{\CovOp}{\ensuremath{\Sigma_{0}}}
\newcommand{\CrossOp}{\ensuremath{\Sigma_1}}

\newcommand{\Pstar}{\ensuremath{P_*}}
\newcommand{\Rstar}{\ensuremath{R_*}}
\newcommand{\Qstar}{\ensuremath{Q_*}}

\newcommand{\vardim}{\ensuremath{m}}

\begin{document}

\begin{center}
{\bf{\LARGE{Optimal and instance-dependent guarantees for Markovian
      linear stochastic approximation}}}

\vspace*{.2in} {\large{
 \begin{tabular}{cc}
  Wenlong Mou$^{ \diamond}$ & Ashwin Pananjady$^{\star}$ \\ Martin
  J. Wainwright$^{\diamond, \dagger}$ &Peter L.  Bartlett$^{\diamond,
    \dagger}$
 \end{tabular}

}

\vspace*{.2in}

{\small{
\begin{tabular}{c}
Department of Electrical Engineering and Computer Sciences$^\diamond$
\\
Department of Statistics$^\dagger$ \\ UC Berkeley
\end{tabular}
  
\vspace*{.1in}

\begin{tabular}{c}
Schools of Industrial \& Systems Engineering, and \\
Electrical \& Computer Engineering$^\star$\\
Georgia Tech
 \end{tabular}
  } }}
\end{center}

\begin{abstract}
We study stochastic approximation procedures for approximately solving
a $d$-dimensional linear fixed point equation based on observing a
trajectory of length $n$ from an ergodic Markov chain.  We first
exhibit a non-asymptotic bound of the order $t_{\mathrm{mix}}
\tfrac{d}{n}$ on the squared error of the last iterate of a standard
scheme, where $t_{\mathrm{mix}}$ is a mixing time.  We then prove a
non-asymptotic instance-dependent bound on a suitably averaged
sequence of iterates, with a leading term that matches the local
asymptotic minimax limit, including sharp dependence on the parameters
$(d, t_{\mathrm{mix}})$ in the higher order terms. We complement these
upper bounds with a non-asymptotic minimax lower bound that
establishes the instance-optimality of the averaged SA estimator. We
derive corollaries of these results for policy evaluation with Markov
noise---covering the TD($\lambda$) family of algorithms for all
$\lambda \in [0, 1)$---and linear autoregressive models. Our
  instance-dependent characterizations open the door to the design of
  fine-grained model selection procedures for hyperparameter tuning
  (e.g., choosing the value of $\lambda$ when running the TD($\lambda$)
  algorithm).
\end{abstract}


\section{Introduction}
\label{sec:intro}

Linear $Z$-estimation problems---in which we are interested in
computing the fixed point of a linear system of equations---arise in
many application domains, including reinforcement learning and
approximate dynamic
programming~\cite{bertsekas2019reinforcement,szepesvari2010algorithms},
stochastic control and
filtering~\cite{benveniste2012adaptive,borkar2009stochastic,kushner2003stochastic},
and time-series analysis~\cite{hamilton2020time}. In many of these
applications, the data-generating mechanism is modeled using an
underlying Markov chain.  The resulting dependency among the
observations presents challenges for algorithm design as well as
statistical analysis.  In this paper, our goal is to provide an
instance-dependent statistical analysis---one that captures the
difficulty of the particular $Z$-estimation problem at hand---and to
develop computationally efficient algorithms that match these
fundamental limits.

A linear $Z$-estimation problem in $\real^\usedim$ is specified by a
fixed point equation of the form
\begin{align}
\label{EqnFixedPoint} 
\theta = \Lbar \theta + \bbar,
\end{align}
where the matrix $\Lbar \in \real^{\usedim \times \usedim}$ and the
vector $\bbar \in \real^\usedim$ are parameters of the problem.  In
settings of interest in this paper, the problem parameters $(\Lbar,
\bbar)$ are unknown, and we observe only a sequence $(\Lmat_t,
\bvec_t)_{t \geq 1}$ of noisy observations, generated according to a
Markov process in the following manner.  The Markov process generates
a sequence $(\state_t)_{t \geq 0}$ of states taking values in some
underlying state space $\statespace$.  This chain is assumed to be
ergodic, with a unique stationary distribution $\stationary$. The
observed pair $(\Lmat_{t + 1}, \bvec_{t + 1})$ at each time $t$
depends on the current state $\state_t$, and moreover, their
expectations under the stationary distribution $\stationary$ are equal
to their population-level counterparts $(\Lbar, \bbar)$.

This general formulation includes a number of special cases of
interest.  In the simplest setting, at each time $t$, we observe a
matrix-vector pair of the form $L_{t+1} = \Lmap(\state_t)$ and
$b_{t+1} = \bmap(\state_t)$, where $\Lmap: \statespace \rightarrow
\real^{\usedim \times \usedim}$ and $\bmap: \statespace \rightarrow
\real^\usedim$ are deterministic mappings such that
\begin{subequations}
\label{eq:markov-obs-model}  
\begin{align}
\Exs_\stationary \big[ \Lmap (\state) \big] = \Lbar, \quad \mbox{and}
\quad \Exs_\stationary \big[ \bmap (\state) \big] =
\bbar.\label{EqnUnbiasedStationary}
\end{align}
Many applications involve additional sources of randomness beyond that
naturally associated with the Markov chain itself.  In order to
accommodate this possibility, we can consider observations of the form
\begin{align}
L_{t+1} = \Lmap_{t + 1}(\state_t), \quad \mbox{and} \quad b_{t+1} =
\bmap_{t+1}(\state_t).
\end{align}
Here the mappings $\Lmap_{t+1}$ and $\bmap_{t+1}$ are now allowed to
be $\myiid$ random, independent of $\state_t$, but are required to be
related to the deterministic mappings $\Lmap$ and $\bmap$ via the
relation
\begin{align}
\label{EqnConditionalObservations}  
  \Exs \big[ \Lmap_{t + 1} (\state) \big] = \Lmap (\state), \quad
  &\Exs \big[ \bmap_{t + 1} (\state) \big] = \bmap (\state), \quad
  \mbox{for all $\state \in \statespace$.}
\end{align}
\end{subequations}
By the tower property of conditional expectation, for a stationary
Markov chain, equations~\eqref{EqnUnbiasedStationary}
and~\eqref{EqnConditionalObservations} imply that
$\Lmap_{t+1}(\state_t)$ and $\bmap_{t+1}(\state_t)$ are unbiased
estimates of $\Lbar$ and $\bbar$, respectively.\footnote{However,
equation~\eqref{EqnConditionalObservations} does not require the
observations to be conditionally unbiased.} The random operator
observed at each iteration is therefore given by $\theta \mapsto
\Lmap_{t + 1} (\state_t) \theta + \bmap_{t + 1} (\state_t)$. This is a
natural generalization of ``random field
noise''~\cite{dieuleveut2020bridging,yu2020analysis} to the Markovian
setting: instead of observing $\mathrm{i.i.d.}$ random fields at
iteration, we observe random functional of a Markov chain's states.

Stochastic approximation (SA) methods, dating back to the seminal work
of Robbins and Monro~\cite{robbins1951stochastic}, are standard
iterative procedures for using data to approximately compute $\theta$.
These algorithms proceed in a streaming fashion: upon receiving each
data point, an incremental update is made and the (averaged or) final
iterate is returned in a single pass. In this way, each iteration of
stochastic approximation incurs only mild computational and storage
costs. Given these attractive computational properties, it is natural
to ask if there are SA methods that also enjoy optimal statistical
performance. To motivate the SA updates, we could start by considering
a stochastic version of the fixed-point iteration to solve
equation~\eqref{EqnFixedPoint}
\begin{align*}
  \theta_{t + 1} = \Lmat_{t + 1} \theta_t + \bvec_{t + 1}.
\end{align*}
With the randomness of the observations $(\Lmat_{t + 1}, \bvec_{t +
  1})$, the iterates will fluctuate at a constant order, and may not
converge. In order to stabilize the stochastic fixed-point iteration,
a stepsize $\stepsize \in (0, 1)$ may be introduced, leading to the
canonical SA updates.

In this paper, we analyze the SA procedure based on the updates
\begin{subequations}
\label{eq:lsa}
\begin{align}
\label{eq:lsa-iterates}  
  \theta_{t + 1} & \mydefn (1 - \stepsize) \theta_t + \stepsize
  (\Lmat_{t + 1} \theta_t + \bvec_{t + 1}), \quad \mbox{for } t =
  0,1,\ldots \\
\label{eq:lsa-average}  
 \thetahat_\numobs & \mydefn \frac{1}{\numobs - \numburn}\sum_{t =
   \numburn}^{\numobs - 1} \theta_t \quad \mbox{for $\numobs =
   \numburn + 1, \numburn + 2, \ldots$.}
\end{align}
\end{subequations}
Equation~\eqref{eq:lsa-iterates} describes a standard stochastic
approximation update with constant stepsize $\stepsize > 0$, whereas
equation~\eqref{eq:lsa-average} corresponds to an application of the
Polyak--Ruppert averaging
procedure~\cite{polyak1992acceleration,ruppert1988efficient} to the
iterates, with burn-in period $\numburn$. When each matrix observation
$\Lmat_{t + 1}$ has a constant rank independent of the dimension
$\usedim$---as is the case for temporal difference learning methods in
reinforcement learning (see Section~\ref{subsec:examples})---the SA
method~\eqref{eq:lsa} can be implemented with $\order{\usedim}$
computational and storage cost per iteration.

There is an extensive body of past work on stochastic approximation
methods with Markov data. Here we provide an overview of the
literature most germane to our contributions, and defer a more
detailed review to
Appendix~\ref{subsec:additional-related-works}. Asymptotic convergence
of SA procedures with Markovian data can be established using either
the ODE method~\cite{borkar2009stochastic} or the Poisson equation
method~\cite{benveniste2012adaptive}. Tsitsiklis and Van
Roy~\cite{tsitsiklis1997analysis} analyze the asymptotic convergence
of SA in the specific context of temporal difference methods in
reinforcement learning.  Although asymptotic guarantees provide
helpful guidance, it is often most useful to have non-asymptotic
guarantees that account both for limited sample size and scale of
modern problems, and for these reasons, non-asymptotic analysis of
Markovian SA procedures has attracted much recent attention.

Assuming a mixing time bound on the Markov chain, a projected variant
of linear SA was analyzed in the paper~\cite{bhandari2018finite},
leading to non-asymptotic rates that are near-optimal in their
dependence on the sample size $n$. Srikant and
Ying~\cite{srikant2019finite} analyzed the standard SA scheme without
the projection step used in the paper~\cite{bhandari2018finite}, and
obtained the same convergence rate in both mean-squared error and
higher moments.  Under an appropriate Lyapunov function assumption on
the Markov chain, Durmus et al.~\cite{durmus2021stability} proved
finite-time bounds for linear SA using stability properties of random
matrix products. Variants and special cases of SA procedures with
Markov data have also been studied, including two-time-scale
algorithms~\cite{kaledin2020finite}, gradient-based optimization under
Markov data~\cite{doan2020finite}, and estimation in auto-regressive
models~\cite{bresler2020least,jain2021streaming}. \\

\noindent Despite this encouraging progress to date, two important
questions still remain open, and are the focus of this paper:
\bcar
\item \textbf{Sample complexity with optimal dimension dependence:}
  The primary goal of non-asymptotic analysis is to provide guarantees
  on the estimation error that have an explicit dependence on the
  problem at hand, and that hold true for a reasonable range of values
  of the sample size $\numobs$. For instance, suppose the linear
  $Z$-estimation problem in $\real^d$ is driven by an underlying
  Markov chain of mixing time $\mixingtime$. Then under natural noise
  assumptions, one should expect an effective sample size of the order
  $\numobs/\mixingtime$, so that the mean-squared error should scale
  as $\order{\mixingtime \usedim / \numobs}$, with this being the
  dominant term whenever $\numobs \gtrsim \mixingtime \usedim$. Such
  an error bound is particularly important for sieve estimators, where
  the problem dimension $\usedim$ is adaptively chosen based on the
  sample size $\numobs$. However, existing analyses of linear SA do
  not provide such tight dimension-dependence.  Using the notation of
  Eq~\eqref{eq:lsa-iterates}, the estimation error bounds in the
  papers~\cite{srikant2019finite,bhandari2018finite} rely on a uniform
  upper bound on the operator norm of the stochastic matrix $\Lmat_{t
    + 1} (\state_t)$; this quantity scales linearly with dimension
  $\usedim$ in many applications.  Consequently, the resulting bounds
  on the MSE have a sub-optimal dependence on dimension, which is
  unsatisfactory for problems with growing dimensions. Similarly, the
  bounds in the
  papers~\cite{durmus2021stability,kotsalis2020simple,chen2021lyapunov}
  also exhibit a sub-optimal dependence on dimension. To the best of
  our knowledge, the question of whether linear SA succeeds under the
  minimal conditions on sample size---in particular, with $\numobs$
  mildly larger than $\mixingtime \cdot \usedim$---remains open.
\item \textbf{Instance-dependent optimality:} While many estimators
  may exhibit near-optimal statistical performance in the globally
  minimax (i.e., worst-case) sense, some of them perform significantly
  better than others when applied to practical problem instances. This
  phenomenon motivates the study of local (i.e., instance-dependent)
  performance in the non-asymptotic regime. Such results have recently
  been established for linear $Z$-estimation in the
  i.i.d. setting~\cite{pananjady2020instance,li2020breaking,khamaru2020temporal,mou2020optimal}.
  The latter two papers listed provide non-asymptotic analogs of
  classical theory on local asymptotic minimaxity
  (c.f.~\cite{van2000asymptotic}),
  which establishes lower bounds by looking at the worst-case family
  of instances in a local neighborhood of a given problem. In the
  Markov setting, two questions naturally arise: (1)~What does it mean
  for an estimator to be locally optimal in a non-asymptotic sense?
  (2)~Does the linear SA estimator~\eqref{eq:lsa} match the local
  lower bound for every problem instance?  \ecar
  %


\subsection{Contributions and organization}

The primary goal of this paper is to resolve these challenges, and
provide a sharp analysis of (averaged) linear SA algorithms.
Arguably, our results are not merely of theoretical interest; they
also provide important guidance for practice, such as in choosing
algorithm parameters including the burn-in period and stepsize.  In
more detail:
\bcar
\item We perform a fine-grained analysis of linear SA and produce an
  upper bound on its statistical error that explicitly tracks the
  dependence on problem-specific complexity as well as step-size.
  Furthermore, our bound holds true provided $\numobs \gtrsim
  \mixingtime \cdot \usedim$, establishing that the algorithm does
  indeed attain a sharp sample complexity guarantee with optimal
  dimension dependence.

\item In a complementary direction to our upper bounds, we show a
  local minimax lower bound with an appropriately defined notion of
  local neighborhood of Markov chains. This lower bound certifies the
  statistical optimality of the linear SA estimator, again in an
  instance-dependent sense.

\item We derive consequences of our general analysis for temporal
  difference methods in reinforcement learning, demonstrating a key
  problem-dependent quantity in matching upper and lower bounds.
  \ecar

One technical aspect of our analysis is noteworthy.  En route to
establishing bounds with sharp dimension dependence, we introduce a
careful ``bootstrapping'' argument: starting with a loose bound, we
progressively refine it via the repeated application of certain
self-bounding inequalities. We suspect that this method may be of
independent interest in providing sharp analyses of other stochastic
approximation methods.

The remainder of this paper is organized as follows.  We complete this
section by introducing notation to be used throughout the paper, and
then providing a more detailed discussion of related work.  In
Section~\ref{sec:problem-setup}, we provide the basic problem set-up,
discuss the underlying assumptions, and give some illustrative
examples.  Section~\ref{sec:main-results} is devoted to the
presentation of our main results, which include upper bounds on the
estimation error of stochastic approximation procedures, along with
local minimax lower bounds that apply to any estimator.  In
Section~\ref{sec:consequences}, we develop some consequences of these
results for specific models, including policy evaluation in
reinforcement learning and estimation in autoregressive models.
Sections~\ref{subsec:proof-prop-lsa-iterate} and
\ref{subsec:proof-main-prj} are devoted to the proofs of
Proposition~\ref{prop:lsa-markov-iterate-bound} and
Theorem~\ref{thm:markov-main}, respectively. We conclude with a
discussion in Section~\ref{sec:discussion}. The proof of
Theorem~\ref{thm:local-minimax} and some auxiliary results, as well as
some corollaries, are postponed to the appendix.

\paragraph{Notation:} We let $(\statespace, \metric)$ denote a metric
space. For any $x \in \statespace$, we use $\delta_x$ to denote the
distribution that places all its mass on $\{x\}$. Given a random
variable $X$, we use the notation $\law (X)$ to denote its probability
distribution.  For a pair $(\pi, \mu)$ of probability distributions on
$\statespace$, let $\Gamma(\pi, \mu)$ denote the space of all possible
couplings of $\mu$ and $\pi$. For any $p \geq 1$, the Wasserstein-$p$
distance between $\pi$ and $\mu$ is given by
\begin{align}
  \label{EqnWasserstein}
\mathcal{W}_{p, \metric} (\pi, \mu) & \defn \Big \{
\inf_{\gamma\in\Gamma(\pi, \mu)} \int_{\statespace \times \statespace}
\metric (x, y)^p d \gamma(x,y) \Big \}^{1/p},
\end{align}
and the total variation distance between $\pi$ and $\mu$ by
\begin{align*}
    \totalvariation (\pi, \mu) \mydefn \sup_{A \subseteq \statespace}
    \abss{\pi (A) - \mu (A)}.
\end{align*}

Our analysis also involves various other divergences between
probability measures. For any pair of probability distributions $P$
and $Q$ on the same space, we use $P \ll Q$ to denote the fact that
$P$ is absolute continuous with respect to $Q$, and use
$\frac{dP}{dQ}$ to indicate the Radon--Nikodym derivative. Given $P
\ll Q$, we define:
\begin{align*}
\mbox{KL Divergence:} &\quad \kull{P}{Q} \mydefn \Exs_P \big[ \log
  \tfrac{d P}{d Q} (X) \big],\\ \mbox{$\chi^2$ divergence:} &\quad
\chisqdiv{P}{Q} \mydefn \Exs_P \big[\tfrac{d P}{d Q}(X)
  - 1 \big],\\ \mbox{Max divergence:} &\quad D_{\infty} (P || Q)
\mydefn \sup_{x \in \mathrm{supp} (Q)} \abss{\log \tfrac{d P}{d Q} (x)
}.
\end{align*}

Given any matrix \mbox{$\Amat = (a_{ij}) \in \real^{n \times m}$,} its
vectorization is obtained by concatenating its
columns---viz. $\vectorize (\Amat) \mydefn [\begin{smallmatrix} a_{11 }&
  a_{21}& \cdots & a_{n1} & a_{12}& \cdots & a_{n2} & \cdots& a_{1m}
  & \cdots& a_{nm} \end{smallmatrix}]^\top \in \real^{n m}$.  We use
$\{e_j\}_{j=1}^\usedim$ to denote the standard basis vectors in the
Euclidean space $\real^\usedim$, i.e., $e_j$ is a vector with a $1$ in
the $j$-th coordinate and zeros elsewhere.  For two matrices $A \in
\real^{d_1 \times d_2}$ and $B \in \real^{d_3 \times d_4}$, we use $A
\otimes B$ to denote their Kronecker product, a $d_1 d_3 \times d_2
d_4$ real matrix. For symmetric matrices $A, B \in \real^{d \times
  d}$, the notation $A \preceq B$ means that $B - A$ is a positive
semi-definite matrix, whereas $A \prec B$ indicates that $B - A$ is
positive definite. We use $\lammax (A)$ and $\lambda_{\min} (A)$ to
denote the largest and smallest eigenvalue of the matrix $A$,
respectively. We use the following notation for matrix norms: for any
matrix $A \in \real^{d_1 \times d_2}$, we use the notation
$\opnorm{A}, ~\matsnorm{A}{F}$ and $\matsnorm{A}{\mathrm{nuc}}$ to
denote its operator norm, Frobenius norm and nuclear norm,
respectively.

Finally, throughout the paper, we use $\filtration_t \mydefn \sigma
\big( ( \bvec_i, \Lmat_i, \state_i )_{i \leq t} \big)$ to denote the
natural filtration induced by the Markovian observations.


\section{Problem set-up}
\label{sec:problem-setup}

Recall from our earlier set-up (cf. equation~\eqref{EqnFixedPoint})
that we are interested in solving a fixed point equation of the form
$\theta = \Lbar \theta + \bbar$, based on noisy observations of the
pair $(\Lbar, \bbar)$, as defined by the Markov observation
model~\eqref{eq:markov-obs-model}.  We require that the matrix $\Lbar$
satisfies the conditions
\begin{align}
\label{eq:kappa-opnorm}  
\kappa \mydefn \frac{1}{2} \lammax \big(\Lbar + \Lbar^\top \big) < 1,
\quad \mbox{and} \quad \opnorm{\Lbar} \leq \conmax.
\end{align}
This condition is used throughout the paper.


\subsection{Assumptions}

We now introduce and discuss the remaining four assumptions that
underlie our analysis.

\subsubsection{Conditions on Markov chain}

We first describe the conditions imposed on the underlying Markov
chain in our observation model.  Let $\{\state_t \}_{t \geq 0}$ denote
a trajectory drawn from a Markov chain with transition kernel
$\transition$.  We assume that this chain has a unique stationary
distribution $\stationary$, and impose the following mixing condition
in Wasserstein-$1$ distance:
\begin{assumption}
\label{assume-markov-mixing}
There exists a natural number $\mixingtime$ and a universal constant
$c_0 > 0$ such that for any $x, y \in \statespace$, we have 
\begin{align}
\label{eq:assume-mixing-tmix}      
\Wass_{1, \metric} (\delta_x \transition^{\mixingtime}, \delta_y
\transition^{\mixingtime}) \stackrel{(a)}{\leq} \frac{1}{2} \metric
(x, y), \quad \mbox{and} \quad \Wass_{1, \metric} (\delta_x
\transition^{t}, \delta_y \transition^{t}) \stackrel{(b)}{\leq} c_0
\metric (x, y)
\end{align}
for all $t = 1, 2, \ldots$.
\end{assumption}
\noindent It is known that such a condition implies rapid mixing (see~\cite{levin2017markov}, Section 4.5). For most parts of the paper, we assume that the chain is initialized with a
sample $\state_0 \sim \stationary$ from the stationary distribution.
Given that our mixing time bound guarantees exponential decay of the
Wasserstein distance, this condition is mild: it can be removed by
waiting $\order{ \mixingtime }$ iterations for the process to mix. By making this intuition rigorous, we will also present a slightly weaker error bound under arbitrary initial distribution (see Corollary~\ref{cor:non-stationary-initial}).\\

\subsubsection{Tail conditions on noise}

In our observation model, the ``noise'' terms correspond to the
differences $\Lmap_{t + 1}(\state_t) - \Lmap(\state_t)$ and
$\Lmap(\state_t) - \Lbar$, along with analogous quantities for the
vector $\bvec$.  Our second assumption imposes conditions on these noise
variables. We consider separate conditions on these martingale (i.e., $\Lmap_{t + 1} (s_t) - \Lmap (s_t)$ and $\bmap_{t + 1} (s_t) - \bmap (s_t)$) and Markov (i.e., $\Lmap (s_t)- \Lbar$ and $\bmap (s_t) - \bbar$) parts of the noise.
\begin{assumption}
\label{assume:noise-moments}
There exists an even integer $\pmax \in [2, + \infty]$ and non-negative constants $\sigmaA$ and $\sigmab$, such that for any positive even
integer $p \leq \pmax$, scalar $t \geq 0$, vector $u \in
\sphere^{\usedim-1}$, and index $j \in \{1, \ldots, \usedim\}$, we
have
\begin{subequations}
\begin{align}
  \Exs \big[ \inprod{\coordinate_j}{ \big( \Lmap_{t + 1}(\state_t) -
      \Lmap(\state_t) \big) u}^p \mid \filtration_t \big] &\leq p!
  \sigmaA^p, \quad \mbox{and}\\
   \Exs_{\state \sim \stationary}
  \big[\Exs \big[ \inprod{\coordinate_j}{\bmap_{t + 1} (\state) -
        \bmap (\state)}^p \mid \state\big] \big] &\leq p!  \sigmab^p,
\end{align}
as well as
\begin{align}
  \Exs_{\state \sim \stationary} \big[ \inprod{\coordinate_j}{ \big(
      \Lmap (\state) - \Lbar \big) u}^p \big] \leq p!  \sigmaA^p,
  \quad \mbox{and} \quad \Exs_{\state \sim \stationary} \big[
    \inprod{\coordinate_j}{\bmap (\state_t) - \bbar}^p \big] \leq p!
  \sigmab^p.
\end{align}
\end{subequations}
\end{assumption}
Note that this assumption is mildest for $\pmax = 2$, and strongest
for $\pmax = \infty$.  In the latter case, when $\pmax = \infty$, the
assumption requires $\Lmat_{t + 1}$ and $\bvec_{t + 1}$ to be
sub-exponential random variables in the standard coordinate directions
(since $\log(p!) \le p \log (p/2)$ by concavity of the log
function). This condition covers, for instance, the case where
$\Lmat_{t + 1}$ is the outer product of sub-Gaussian random vectors,
as in temporal difference learning methods. In addition to
accommodating this case, Assumption~\ref{assume:noise-moments} also
covers the heavier-tailed setting in which only finitely many moments
exist. In particular, when $\pmax = 2$, the second moment assumption
coincides with the assumption made in the paper~\cite{mou2020optimal}.

An important quantity in our analysis is the \emph{effective noise
level} given by
\begin{align}
    \varbound \mydefn \sup_{p \in [2, \pmax]} \sup_{j \in [d]} \sup_{t \geq 0}~ p^{-1} \big( \Exs \big[
      \inprod{\coordinate_j}{(\Lmap_{t + 1} (s_t) - \Lbar) \thetabar +
        (\bmap_{t + 1} (s_t) - \bbar)}^{p} \big] \big)^{1/p}.\label{eq:defn-effective-noise}
\end{align}
Note that under Assumption~\ref{assume:noise-moments}, we have the
upper bound $\varbound \leq \sigmaA \vecnorm{\thetabar}{2} + \sigmab$.

\subsubsection{Metric space conditions}

For most of our analysis, we impose the following condition:
\begin{assumption}
\label{assume:stationary-tail}
The metric space $(\statespace, \metric)$ has diameter at most one.
\end{assumption}
\noindent Note that our assumption of unit diameter is arbitrary;
boundedness suffices.  In order to accommodate the general case, it
suffices to rescale the parameters $\matLip$ and $\veclip$.

When applying our theory to unbounded spaces (e.g., $\statespace =
\real^d$), we use a truncation argument to show that there is an event
over a reduced state space on which this condition holds with
probability tending exponentially to $1$. (See Appendix~\ref{SecTruncate} for the details of this argument.) To unify the notation, we always assume the distance to be of constant order with hih probability, which results in constant diameter of the truncated space. In high-dimensional Euclidean spaces, the distance between two generic random vectors can easily become dimension-dependent. In such cases, we rescale the space to make it a constant. The rescaling could lead to dimension-dependent Lipschitz constants, which is captured in Assumption~\ref{assume-lip-mapping} to follow.

\subsubsection{Lipschitz condition}

Finally, we place a Lipschitz assumption---under the metric
$\metric$---on the mapping from the metric space $\statespace$ to the
stochastic operators. Given the Markov chain setup in the metric space
$(\statespace, \metric)$, it is tempting to assume a dimension-free
Lipschitz bound on the mappings $(\Lmap_{t}, \bmap_t)$. However, such
Lipschitz constants typically depend on dimension for practical
problems.  Concretely, view the $\Lbar$-scale parameters $(\kappa,
\conmax)$ as constants and assume that the observations $\Lmap_{t + 1}
(\state_t)$ each have rank at most $r$. We then have
\begin{align}
\label{eq:constant-rank-obs-has-large-opnorm}  
    \Exs \big[ \opnorm{\Lmap_{t + 1} (\state_t)} \big] \geq \frac{
      \Exs \big[ \matsnorm{\Lmap_{t + 1} (\state_t)}{\mathrm{nuc}}
        \big]}{r} \geq \frac{ \mathrm{trace} \big( \Exs \big[ \Lmap_{t
          + 1} (\state_t) \big] \big)}{r} = \frac{\mathrm{trace} (
      \Lbar)}{r}.
\end{align}
Note that the term $\mathrm{trace} ( \Lbar)$ typically scales as
$\Theta (d)$, even in the ``easy case'' when $\Lbar$ is a constant
multiple of identity matrix.

Consequently, the Lipschitz constant for the mapping $\Lmap_t
: \statespace \rightarrow \real^{d \times d}$ grows at least linearly
in dimension $d$.  On the other hand, as a $d$-dimensional standard
Gaussian random variable has norm $\sqrt{d} -
\widetilde{\mathcal{O}}(1)$ with high probability, it is natural to
assume the Lipschitz constant for the vector-valued mapping
$\bmap_t: \statespace \rightarrow \real^d$ to be of order at least
$\Omega (\sqrt{d})$. We therefore make the following assumption:

\begin{assumption}
\label{assume-lip-mapping}
There exist constants $\matLip, \veclip > 0$ such that, almost surely
for any $x, y \in \statespace$, we have
\begin{align}
\opnorm{\Lmap_{t} (x) - \Lmap_t (y)} \leq \sigmaA \usedim \cdot
\metric (x, y) \quad \mbox{and} \quad \vecnorm{\bmap_{t} (x) - \bmap_t
  (y)}{2} \leq \sigmab \sqrt{\usedim} \cdot \metric (x, y)
\end{align}
for all $t = 1, 2, \ldots$.
\end{assumption}
Note that in Assumption~\ref{assume-lip-mapping}, we have explicitly
scaled the RHS of the inequalities with factors that depend on the
problem dimension $\usedim$, so that the pair $(\matLip, \veclip)$
should indeed be viewed as dimension-free.  It is also worth noting
that the notation $(\matLip, \veclip)$ is overloaded, since we can
take the maximum of the bounds in Assumptions
\ref{assume:noise-moments} and~\ref{assume-lip-mapping}.  As shown in
Appendix~\ref{SecTruncate}, for certain natural problem classes,
Assumption~\ref{assume:noise-moments} indeed implies
Assumption~\ref{assume-lip-mapping} with discrete metric, up to
logarithmic factors.


\subsection{Some illustrative examples}
\label{subsec:examples}

Our assumptions cover a broad range of ergodic Markov chains, and the
fixed-point equation~\eqref{EqnFixedPoint} associated with their
stationary distribution naturally arises from several problems. In
this section, we describe a few concrete examples of our general
setup. We first discuss the class of Markov chains satisfying our
assumptions, and then describe the linear $Z$-estimators associated
with such problems.


\subsubsection{Examples of Markov chains}

By varying our choice of the metric $\metric$, we recover several
important classes of Markov chains that satisfy
Assumptions~\ref{assume-markov-mixing}
and~\ref{assume:stationary-tail}.
\bcar
\item Consider a Markov chain defined on a countable state space
  $\statespace$, and consider the discrete metric $\metric (x, y)
  \mydefn \bm{1}_{x \neq y}$.  In this context,
  Assumption~\ref{assume-markov-mixing} corresponds to mixing time
  bound in total variation---viz.
\begin{align*}
  \totalvariation (\delta_x \transition^{\mixingtime}, \delta_y
  \transition^{\mixingtime}) \leq \tfrac{1}{2} \qquad \mbox{for all
    pairs $x, y \in \statespace$.}
\end{align*}
This mixing condition is satisfied for some finite $\mixingtime$ when
the Markov chain is irreducible, aperiodic and positive
recurrent. Moreover, this metric space has unit diameter, so that
Assumption~\ref{assume:stationary-tail} holds as well.

\item As another example, consider the state space $\statespace =
  \ball (0, 1) \subseteq \real^d$ equipped with the Euclidean metric
  $\metric(x, y) = \vecnorm{x - y}{2}$.  We can define a Markov chain
  on this space via the random evolution $X_{k + 1} = \RandOp_{k + 1}
  (X_k)$, where the random non-linear operators $\{ \RandOp_k \}_{k
    \geq 1} \subseteq \statespace^{\statespace}$ are drawn $\myiid$
  from some distribution.  We assume that the expected operator
  $\bar{\RandOp} \mydefn \Exs[\RandOp_1]$ satisfies the contraction
  condition $\vecnorm{\RandOpBar(x) - \RandOpBar(y)}{2} \leq \gamma
  \vecnorm{x - y}{2}$ with some $\gamma < 1$. Assuming the stochastic
  operator $\RandOp$ to be Lipschitz and to satisfy a second moment
  bound, this dynamical system satisfies the Wasserstein contraction
  condition under the Euclidean metric.
\ecar



\subsubsection{Examples of linear $Z$-estimators}

We now describe some interesting examples of linear $Z$-estimators, to
which we will return in later sections.

\begin{example}[Approximate policy evaluation]
  \label{example:td0}
  \upshape

We begin by considering the temporal difference (TD) algorithm for
approximate estimation of value functions.  This problem arises in the
context of Markov reward processes (MRPs), which are Markov chains
that are augmented with a reward function $\reward: \statespace
\rightarrow \real$.  A trajectory from a Markov reward process is a
sequence $\{ (\state_t, R_t) \}_{t \geq 0}$, where $\{ \state_t \}_{t
  \geq 0}$ is the Markov trajectory of states, and $R_t$ is a random
reward, corresponding to a conditionally unbiased estimate (given
$\state_t$) of the reward function value $\reward(\state_t)$.  Given a
discount factor \mbox{$\discount \in [0,1)$,} the expected discount
  reward defines the \emph{value function} \mbox{$\valuestar(\state) =
    \Exs \big[ \sum_{t = 0}^\infty \discount^t R_t \mid \state_0 =
      \state \big]$.}

This value function is connected to linear $Z$-estimators via the
Bellman principle. Let $\transition$ denote the transition operator of
the Markov chain, and let $\stationary$ denote the stationary
distribution.  Note that the $\transition$ maps the space $\Ltwospace
(\statespace, \stationary)$ to itself.  With this notation, the value
function $\valuestar$ is known to be the unique fixed point of the
\emph{Bellman evaluation equation}
\begin{align}
\label{eq:bellman}  
  \valuefunc & = \discount \transition \valuefunc + \reward.
\end{align}
In general, this equation is non-trivial to solve, especially given a
limited trajectory length.  In practice, it is standard to compute
approximate solutions using linear basis expansions
~\cite{bradtke1996linear,tsitsiklis1997analysis}, and this approach
underlies the family of TD algorithms with linear function approximation.

Let $\{\phi_j \}_{j=1}^\usedim$ be a collection of linearly
independent real-valued functions defined on the state space, and
consider the linear subspace $\LinSpace$ of all functions of the form
$V_\theta(\state) = \sum_{j=1}^\usedim \theta_j \phi_j(\state)$.  This
subspace defines the \emph{projected Bellman equation}
\begin{align}
\label{eq:projected-fixed-point}  
\valuebar = \projecttolin \big( \discount \transition \valuebar + r
\big),
 \end{align}
where $\projecttolin$ is the orthogonal projection operator under
$\Ltwospace (\statespace, \stationary)$.

By definition, the projected fixed point $\valuebar$ can be written in
the form $\valuebar(\state) = \sum_{j=1}^\usedim \thetabar_j
\phi_j(\state)$ for some vector $\thetabar \in \real^\usedim$. Define
the vector-valued mapping $\phi = [\phi_j]_{j = 1}^\usedim$, some
simple calculations show that this parameter vector must satisfy the
linear system
\begin{align}
  \label{eq:lstd-in-low-dimensional-space}
  \CovOp \thetabar = \discount \CrossOp \thetabar + \Exs_{s \sim
    \stationary} \big[ R_0(\state) \phi(\state) \big],
\end{align}
where $\CovOp = \Exs_{\state \sim \stationary} \big[ \phi(\state)
  \phi(\state)^\top \big]$ is the second-moment matrix of
$\phi(\state)$ under the stationary distribution, and $\CrossOp =
\Exs[\phi(\state) \phi(\state^+)^\top]$ is the cross-moment operator
of the Markov chain.  In defining this cross-moment, the expectation
is taken over $\state \sim \stationary$ and $\state^+ \sim
\transition(\state, \cdot)$.

This problem can be viewed within our framework by considering a
Markov chain on the augmented state space \mbox{$\slidingwindow_t =
  (\state_t, \state_{t + 1})$.}
Equation~\eqref{eq:lstd-in-low-dimensional-space} defines a fixed
point equation under the stationary distribution of this Markov chain.
Define the minimum and maximum eigenvalues $\mu \mydefn
\lambda_{\min}(\CovOp)$ and $\beta \mydefn \lambda_{\max}(\CovOp)$,
along with the observation functions
\begin{align}
  \bmap_{t + 1} (\slidingwindow_t) = \tfrac{1}{\beta} R_t (s_t)
  \phi(s_t), \quad \mbox{and} \quad
  \Lmap_{t + 1}(\slidingwindow_t) &= I_d - \tfrac{1}{\beta} \big [
    \phi(s_t) \phi(s_{t})^\top - \discount \phi(s_t)
    \phi(s_{t+1})^\top \big ].\label{eq:fit-td0-into-our-setup}
\end{align}
With these choices, the stochastic approximation
procedure~\eqref{eq:lsa} is the widely used TD$(0)$ algorithm. On the
other hand, for a stationary Markov chain $(s_t)_{t \in \mathbb{Z}}$,
the fixed-point equation $\thetabar = \Exs \left[ \Lmap_{t + 1}
  (\slidingwindow_t) \right] \cdot \thetabar + \Exs \left[ \bmap_{t +
    1} (\slidingwindow_t) \right]$ is equivalent to
equation~\eqref{eq:lstd-in-low-dimensional-space}. Note that though
the expression for the mappings $\bmap_{t + 1}$ and $\Lmap_{t + 1}$
depends on unknown parameter $\beta$, they can be absorbed into the
stepsize choice, and the algorithm works well without such knowledge.

Typically, the Euclidean norm $\vecnorm{\phi(s)}{2}$ of the feature
vectors scales as $\sqrt{\usedim}$, and under the stationary
distribution $\stationary$, the variance of any coordinate of
$\phi(\state)$ is of constant order.  Under these conditions, the
cross-moment matrix $\CrossOp$ has operator norm of constant order.
On the other hand, as for the random observations, we have the
scalings $\opnorm{\Lmat_{t + 1}} = \myorder(\usedim)$ and
$\vecnorm{\bvec_{t + 1}}{2} = \myorder(\sqrt{d})$, so that
Assumptions~\ref{assume:noise-moments} and~\ref{assume-lip-mapping}
are satisfied.
\hfill $\goodendex$
\end{example}

In the context of TD, it is natural to consider a \emph{sieve
estimator}.  Given a collection of basis functions $\{\phi_j
\}_{j=1}^\infty$, we can define the nested family $\LinSpace_1
\subset \LinSpace_2 \subset \cdots$, where $\LinSpace_\usedim$ denotes
the span of the sub-collection $\{\phi_j\}_{j=1}^\usedim$.  Here the
choice of the sieve parameter $\usedim$ is key: larger values reduce
the approximation error at the expense of increasing the estimation
error.  We discuss how this can be done in
Section~\ref{sec:consequences}. \\

Another extension of the TD$(0)$ algorithm---one that becomes feasible
under the Markovian observation model---is the TD$(\lambda)$ family of
procedures.  A fundamental question is how well the solution of the
projected fixed-point equation~\eqref{eq:projected-fixed-point}
approximates the true value function $V^*$. Prior work by a subset of 
the current authors~\cite{mou2020optimal} analyzes this quantity, and provides
matching upper and lower bounds in the $\myiid$ setting.  However, the
Markovian observation model actually allows this approximation error
to be reduced, albeit at the cost of increased estimation error, as
discussed in our next example.

\begin{example}[Policy evaluation with TD$(\lambda)$] \upshape
  \label{example:td-lambda}
  The family of TD$(\lambda)$ algorithms is motivated by the following
  observation: since the value function $\valuestar$ is the fixed
  point of equation~\eqref{eq:bellman}, it is also the fixed point of
  the composition of itself. Concretely, for any $k \geq 1$, we have:
  \begin{align*}
      \valuestar = (\discount \transition)^k \valuestar + \sum_{j =
        0}^{k - 1} (\discount \transition)^j \reward.
  \end{align*}
  For any $\lambda \in [0, 1)$, we take the weighted average of the
    above (infinite) collection of equations using
    exponentially-decaying weight $(1, \lambda, \lambda^2, \cdots)$,
    and obtain the following equation.
\begin{subequations}
  \begin{align}
\label{eq:population-td-lambda-unprojected}        
     \valuefunc = (1 - \lambda) \sum_{k = 0}^{\infty} \lambda^k
     (\discount \transition)^{k+1} \valuefunc + \sum_{k = 0}^{\infty}
     \lambda^k (\discount \transition)^k \reward.
  \end{align}
  The solution $\valuestar$ to the equation~\eqref{eq:bellman} also
  solves equation~\eqref{eq:population-td-lambda-unprojected}.
  
Following the same route as TD$(0)$, for a given subspace $\LinSpace$
of functions, we seek a solution $\valuebar^{(\lambda)}$ to the
projected fixed equation equation
\begin{align}
  \label{eq:population-td-lambda}  
\valuebar^{(\lambda)} = (1 - \lambda) \sum_{k = 0}^{\infty} \lambda^k
\projecttolin (\discount \transition)^{k+1} \valuebar^{(\lambda)} +
\sum_{k = 0}^{\infty} \lambda^k \projecttolin (\discount
\transition)^k \reward,
\end{align}
\end{subequations}
in which the operator $\transition$ has been replaced by the
projection $\projecttolin \transition$.  Although the fixed points of
equation~\eqref{eq:population-td-lambda-unprojected} and the Bellman
equation~\eqref{eq:bellman} coincide, the projected
version~\eqref{eq:population-td-lambda} has a different set of fixed
points.

Since the value function $\valuebar^{(\lambda)}$ lies in the linear
space $\LinSpace$, it has a representation of the form
$\valuebar^{(\lambda)}(\state) = \sum_{j=1}^\usedim
\thetabar^{(\lambda)}_j \phi_j(\state)$ for some coefficient vector
$\thetabar^{(\lambda)} \in \real^\usedim$.  From
equation~\eqref{eq:population-td-lambda}, this vector must satisfy a
linear system of the form
\begin{align}
  \label{EqnTDLambdaLow}
\big[ \sum_{k=0}^{\infty} (\lambda \discount)^k \Sigma_k \big]
\thetabar^{(\lambda)} & = \big[ \sum_{k=0}^{\infty} (\lambda
  \discount)^k \discount \Sigma_{k+1} \big] \thetabar^{(\lambda)} +
\sum_{k=0}^{\infty} (\lambda \discount)^k \Exs \big[ R_0 (s_0)
  \phi(s_{-k}) \big],
\end{align}
where $\{s_k\}_{k = - \infty}^{\infty}$ is a stationary Markov chain
following the transition kernel $\transition$, and we define $\Sigma_k
= \Exs[\phi(s_{-k}) \phi(s_0)^\top]$ for each integer $k$.  As it
should, when we set $\lambda = 0$, equation~\eqref{EqnTDLambdaLow}
reduces to the TD$(0)$ update from
equation~\eqref{eq:lstd-in-low-dimensional-space}.

In order to use stochastic approximation methods to solve this
equation, we consider an augmented Markov process $(\state_{t + 1},
\state_{t}, g_t)_{t \in \mathbb{Z}}$ in the space $\statespace^2
\times \real^\usedim$, which evolves as
\begin{subequations}
\begin{align}
    \state_{t + 1} \sim \transition (\state_t, \cdot), \quad
    \mbox{and} \quad g_{t} = \phi (\state_{t}) + \discount \lambda
    g_{t - 1}.
\end{align}
If feature vectors $\phi (s_t)$ lie in a compact set almost surely, we
have $g_t = \sum_{k = 0}^{+ \infty} (\discount \lambda)^k \phi (s_{t -
  k})$. Let $\widetilde{\stationary}$ be the stationary distribution
of this augmented Markov chain.\footnote{Such a stationary
distribution exists and is unique under suitable assumptions. See
Section~\ref{subsec:consequence-td-lambda} for details.}  In terms of
an element $\slidingwindow = (s, s^+, g)$ drawn according this
stationary distribution, the fixed-point
equation~\eqref{eq:population-td-lambda} admits the succinct
representation
\begin{align}
  \label{eq:TD-lambda-succinct}
\Exs_{\widetilde{\stationary}} \big[ g \phi (s)^\top \big]
\thetabar^{(\lambda)} = \discount \Exs_{\widetilde{\stationary}} \big[
  g \phi (s^+)^\top \big] \thetabar^{(\lambda)} +
\Exs_{\widetilde{\stationary}} \big[ R_0 (\state) g \big].
\end{align}
By choosing the observation functions
\begin{align}
 \Lmap_{t + 1}(\slidingwindow_t) = I_d - \nu \cdot \big( g_t \phi
 (s_t)^\top - \discount g_t \phi (s_{t + 1})^\top \big), \quad
 \bmap_{t + 1}(\slidingwindow_t) = \nu \cdot R_t (s_t) \phi
 (s_t),\label{eq:fit-tdlambda-into-our-setup}
\end{align}
\end{subequations}
for a scalar $\nu > 0$, this algorithm is a special case of our
general set-up. In particular, by substituting the infinite-sum
expression for the random variable $g_t$ into
equation~\eqref{eq:TD-lambda-succinct}, we obtain the projected linear
equation~\eqref{EqnTDLambdaLow} under the low-dimensional
representation. See Section~\ref{sec:consequences} for a more detailed
verification of the assumptions needed to apply our main results for
this problem.  \hfill \goodendex
\end{example}


For our last example, we turn to a different class of problems
involving vector autoregressive (VAR) models for time
series~\cite{Lut05}.
\begin{example}[Parameter estimation in autoregressive models]
  \label{example:autoregressive}
  \upshape
An $\vardim$-dimensional VAR model of order $k$ describes the
evolution of a random vector $X_t$ as a $k^{th}$-order Markov process.
The model is specified by a collection of $\vardim \times \vardim$
matrices $\{\AstarInd{j} \}_{j=1}^k$, and the random vector evolves
according to the recursion
\begin{align}
\label{eq:autoregressive-process}  
X_{t + 1} = \sum_{j=1}^k \AstarInd{j} X_{t-j+1} + \varepsilon_{t + 1},
\end{align}
where the noise sequence $\big(\varepsilon_{t}\big)_{t \geq 0}$ is
$\myiid$ and zero-mean and supported on a bounded set.

Considering the $(k + 1)$-fold tuple $\slidingwindow_t = (X_{t + 1},
X_t, \cdots, X_{t -k + 1})$, the process $\big(\slidingwindow_t
\big)_{t \geq 0}$ is Markovian. Under appropriate stability
assumptions on the model parameter, the process mixes rapidly under
the $(k + 1) \vardim$-dimensional Euclidean metric. Let
$\widetilde{\stationary}$ denote its stationary distribution, and
suppose for convenience that the chain is observed at stationarity.
    
In order to estimate the model parameters, we consider the following
set of Yule--Walker estimation equations:
\begin{align}
\label{eq:yule-walker-pop}  
  \Exs \big[ X_{t + 1} X_{t - \ell}^\top \big] = \AstarInd{1} \Exs
  \big[ X_t X_{t - \ell}^\top \big] + \AstarInd{2} \Exs \big[ X_{t - 1}
    X_{t - \ell}^\top \big] + \cdots + \AstarInd{k} \Exs \big[ X_{t - k
      + 1} X_{t - \ell}^\top\big],
\end{align}
for $\ell = 0,1, \cdots, k - 1$.
    
These equations form a $k \vardim^2$-dimensional linear system for
estimating $k \vardim^2$-dimensional parameters. Note that the
parameters live in the space of matrix sequences, and so we slightly
abuse our notation for simplicity: $\Lmat$ denotes a linear operator
from $\real^{k \times \vardim \times \vardim}$ to itself, and $\bvec$
is an element in $\real^{k \times \vardim \times \vardim}$.  At the
sample level, for any collection $\AmatInd{} \defn \{ \AmatInd{j}
\}_{j=1}^k \in \real^{k \times \vardim \times \vardim}$ of system
matrices, the stochastic observations are given by
\begin{align*}
  \big[\bmap_{t + 1}(\slidingwindow_t)\big]_{\ell} & = \nu \: X_{t +
    1} X_{t - \ell}^\top \quad \mbox{for $\ell = 0, 1, \ldots, k-1$,
    and} \\
  \big(\Lmap_{t + 1} (\slidingwindow_t) \big) [\AmatInd{}]_{\ell} & =
  \AmatInd{\ell} - \nu \; \sum_{j = 0}^{k - 1} \AmatInd{j} X_{t - j}
  X_{t - \ell}^\top, \quad \mbox{for $\ell = 0, 1, \ldots, k-1$.}
\end{align*}
Once again, the parameter $\nu$ is a scaling constant needed to fit
into the fixed-point equation framework, and is absorbed into the
stepsize choice of the algorithm. \hfill $\goodendex$
\end{example}


\section{Main results}
\label{sec:main-results}

We now turn to the statement of our main results, beginning with our
upper bounds in Section~\ref{SecUpper}, followed by lower bounds in
Section~\ref{SecLower}.

\subsection{Instance-dependent upper bounds}
\label{SecUpper}

In this section, we begin by stating some upper bounds
(Theorem~\ref{thm:markov-main}) on the behavior of the Polyak--Ruppert
averaged SA scheme~\eqref{eq:lsa-average}. These bounds are
instance-dependent, in the sense that they are specified in terms of
an explicit function of the operator $\Lbar$ and the fixed point
$\thetabar$.  We then state a second result
(Proposition~\ref{prop:lsa-markov-iterate-bound}) on the non-averaged
iterates, which plays a key role in proving
Theorem~\ref{thm:markov-main}.

\subsubsection{Instance-dependent bounds on the averaged iterates}
\label{subsubsec:main-upper-bound}

For any state $\state \in \statespace$, define the functions
\begin{align*}
\varepsilon_{\mathrm{MG}} (\state) \mydefn (\bmap_1 (\state) - \bmap
(\state) ) + (\Lmap_1 (\state) - \Lmap (\state)) \thetabar, \quad
\mbox{and} \quad \varepsilon_{\mathrm{Mkv}} (\state) \mydefn \bmap
(\state) + \Lmap (\state) \thetabar - \thetabar.
\end{align*}
Note that for a fixed state $\state$, the quantity
$\varepsilon_{\mathrm{MG}}(\state)$ depends on the random variables
$\bmap_1(\state)$ and $\Lmap_1(\state)$, and so is a random vector,
whereas by contrast, the quantity $\varepsilon_{\mathrm{Mkv}}(\state)$
is deterministic. Letting $(\widetilde{s}_t)_{t = - \infty}^{\infty}$
be a stationary Markov chain under the transition kernel
$\transition$, we then define the matrices
  \begin{align}
\label{eq:def-sigstar}    
    \SigStar_{\mathrm{MG}} \mydefn \Exs_{\stationary} \big[ \cov \big(
      \varepsilon_{\mathrm{MG}} (\state) \mid \state\big) \big], \quad
    \mbox{and} \quad
\SigStar_{\mathrm{Mkv}} \mydefn \sum_{t = - \infty}^{\infty} \Exs
\big[\varepsilon_{\mathrm{Mkv}} (\widetilde{\state}_t)
  \varepsilon_{\mathrm{Mkv}} (\widetilde{\state}_0)^\top \big].
\end{align}
Overall, the performance of our algorithm depends on the \emph{matrix
sum} $\SigStar \mydefn \SigStar_{\mathrm{MG}} +
\SigStar_{\mathrm{Mkv}}$, as well as the effective noise
variance $\varbound^2$ defined in equation~\eqref{eq:defn-effective-noise}.  In terms of these quantities, we have the following
guarantee: \\

\begin{theorem}
\label{thm:markov-main}
Under
Assumptions~\ref{assume-markov-mixing}--\ref{assume:stationary-tail},
suppose that we set the stepsize $\stepsize$ and burn-in parameter
$\numburn$ as $\stepsize = \big(c (\sigmaA^2 d + \conmax^2) (1 -
\kappa) \numobs^2 \mixingtime \big)^{-1/3}$ and $\numburn =
\frac{1}{2} \numobs$, where $c$ is a suitably chosen universal
constant. There exist universal constants $c_1, c_2 > 0$, such that for any sample size $\numobs$ satisfying
$\frac{\numobs}{\log^2 \numobs} \geq \frac{2 \mixingtime (\sigmaA^2
  \usedim + \conmax^2)}{(1 - \kappa)^2} \log (c_0 d)$, the
Polyak--Ruppert estimate~\eqref{eq:lsa-average} has MSE bounded as
\begin{align}
\label{EqnMainBound}  
\Exs \big[ \vecnorm{\thetahat_\numobs - \thetabar}{2}^2 \big] & \leq
\tfrac{c_1}{\numobs} \mathrm{Tr} \big( (I - \Lbar)^{-1}
(\SigStar_{\mathrm{MG}} + \SigStar_{\mathrm{Mkv}})(I - \Lbar)^{-\top}
\big) + c_2 \big( \tfrac{ \sigmaA^2 \usedim \mixingtime}{(1 -
  \kappa)^2 \numobs} \big)^{4/3} \varbound^2 \log^2 \numobs.
\end{align}
\end{theorem}
\noindent See Section~\ref{subsec:proof-main-prj} for the proof of
this theorem. \\

A few remarks are in order. First, and as shown in the next section,
the first term $\numobs^{-1} \mathrm{Tr} \big( (I - \Lbar)^{-1}
\SigStar (I - \Lbar)^{-1} \big)$ is optimal for the Markovian
stochastic approximation problem in an instance-dependent sense.  This
term appears in existing central limit results for Markovian
stochastic approximation~\cite{fort2015central}, and our bound
captures this dependence in a non-asymptotic manner up to a universal
constant. It is worth noting that when the Markov chain is uniformly
geometrically ergodic, a central limit theorem for the averaged
iterate $\thetahat_\numobs$ directly follows from classical Markovian
CLT (see~\cite{meyn2012markov}, Chapter 17).

The first term in the bound~\ref{EqnMainBound} can always be further upper bounded\footnote{This can be easily seen from exponential decay of the correlation; in particular, see equation~\eqref{EqnCoconutBall} in the
proof of the theorem.} by $c_1
\tfrac{\varbound^2}{(1 - \kappa)^2 \numobs} \mixingtime \usedim \cdot \log^2
(c_0 d)$.  On the other hand, disregarding dependence on
$(\sigmaA, \sigmab)$ and logarithmic factors in the sample size, the
second term in the bound scales as $\myorder \big( \big(
\frac{\mixingtime \usedim}{(1 - \kappa)^2 \numobs} \big)^{4/3} \big)$.
Consequently, up to polylogarithmic factors, we have
\begin{align}
  \label{eq:worst-case-bd}
  \Exs \big[ \vecnorm{\thetahat_\numobs - \thetabar}{2}^2 \big]
  \lesssim \frac{\varbound^2 \mixingtime \usedim}{(1 - \kappa)^2
    \numobs} .
  \end{align}
Thus, at least in a worst-case sense, the second term is always
dominated by the first term, and our instance-dependent analysis also recovers a worst-case optimal statistical rate for linear $Z$-estimation with Markovian data. It is also worth noting that the second term in equation~\eqref{EqnMainBound} decays with sample size at at $\numobs^{- 4/3}$ rate, faster than the $O(\numobs^{-1})$ leading-order term. For sufficiently large sample size $\numobs$, this term is dominated by the first term, and the behavior of the estimator $\thetahat_\numobs$ is governed by the instance-optimal quantity. It should be noted that $\numobs^{-4/3}$-rate of the second-order term---indicating how fast the exact instance-optimal behavior kicks in---may not be optimal. Indeed, it decays more slowly than the $\numobs^{-2}$ second-order asymptotic efficiency in regular parametric models~\cite{rao1962first,ghosh1980second}, and we conjecture that such a second-order term is unavoidable for stochastic approximation. That being said, the sub-optimality is only a second-order phenomenon, and the main message of Theorem~\ref{thm:markov-main} is unaffected: with a reasonable sample size, the Polyak--Ruppert estimator is instance-optimal, up to constant factors.

Note that Theorem~\ref{thm:markov-main} involves inexplicit universal constants $(c, c_1, c_2)$. Since our theory focuses on optimal instance-dependent quantities and sample complexities up to universal constant factors, we do not try to optimize these constants. Our proof gives an upper bound with $c_1 = 16$.~\footnote{The universal constants $c_2$ in the high-order term depend on the constant pre-factor in Proposition~\ref{prop:lsa-markov-iterate-bound} to follow. We will not track their explicit values for simplicity.} That being said, for the tail-averaging procedure described above, we can make the constant $c_1$ arbitrarily close to 2, as we are using the latter half of the data in Polyak--Ruppert averaging. With a more careful choice of the burn-in period, e.g., $\numburn \asymp \frac{\log \numobs}{\stepsize (1 - \kappa)}$, the constant $c_1$ in equation~\eqref{EqnMainBound} can be made arbitrarily close to 1. The proof is straightforward---the leading-order term in equation~\eqref{EqnMainBound} comes from the variance of sum of a functional on the Markov chain state space, which can be computed directly.

We note that Theorem~\ref{thm:markov-main} makes two types of tail
assumptions on the random observations:
Assumption~\ref{assume:noise-moments} with $\pmax = 2$ requires
dimension-free second moment bounds in any coordinate direction,
whereas the Lipschitz condition (Assumption~\ref{assume-lip-mapping}) together with Assumption~\ref{assume:stationary-tail}
(boundedness of the domain) imply a (dimension-dependent)
uniform upper bound on the noise. The two assumptions play very
different roles when the dimension dependence is taken into account. As we
will see in Proposition~\ref{prop:sieve-td-results}, such assumptions
are naturally satisfied in the context of sieve estimators, for which
dimension $d$ of the problem is selected adaptively based on sample
size $\numobs$.

Finally, we also note that the requirement on the sample size
$\numobs$ is nearly optimal, since we require $n = \widetilde{\Omega}
\big( \frac{\mixingtime \usedim}{(1 - \kappa)^2} \big)$ to make the
estimation error~\eqref{eq:worst-case-bd} less than a constant (by seeing $\sigmaA$ and $\conmax$ as constants). Up to
an additional $\order{\mixingtime}$ factor, the sample size
requirement in Theorem~\ref{thm:markov-main} also matches that of
linear stochastic approximation in the $\myiid$
setting~\cite{lakshminarayanan2018linear,mou2020linear,mou2020optimal}.
This additional $\order{\mixingtime}$ factor is unavoidable, which can
be seen from the following reduction from the Markov to the $\myiid$ setting. Consider a problem instance in the $\myiid$ setup, given by a probability distribution $\Prob$ over $\real^{d \times d} \times \real^d$. Defining the state $(\Lmat_{t}, \bvec_t)$, consider a lazy
Markov chain that remains at the same state with probability $1 -
\frac{1}{\mixingtime}$, and jumps to an independent state drawn from
$\Prob$ with probability
$\frac{1}{\mixingtime}$. A Markov trajectory of size $\numobs$ in this
lazy Markov chain is approximately equivalent to $\order{\numobs /
  \mixingtime}$ samples in the $\myiid$ model, and results in a
multiplicative blow-up of $\order{\mixingtime}$ in the sample
complexity requirement for the Markov case.

\paragraph{Starting from non-stationary $\state_0$:} Note that Theorem~\ref{thm:markov-main} is shown under an initial state satisfying $\state_0 \sim \stationary$. Such an assumption may not be always available in practice. However, in Corollary~\ref{cor:non-stationary-initial} to follow, we will show that under minor modification, our conclusion easily extends to non-stationary initial distributions.

For non-stationary initial distributions, we wait for a cold-start period $\numcold \mydefn \numobs / 4$ to start running the stochastic approximation iterate~\eqref{eq:lsa-iterates}, i.e., we take the stepsize sequence as
\begin{align*}
  \stepsize_t \mydefn \begin{cases}
    0 & t \in \{0,1, \cdots, \numcold - 1\},\\
    \stepsize & t \geq \numcold.
  \end{cases}
\end{align*}
The rest of the SA procedure, including the average step~\eqref{eq:lsa-average}, remains the same as before. For notational simplicity, we let $\theta_0 = 0$. 
\begin{corollary}\label{cor:non-stationary-initial}
  Under the same setup and parameter choice as in Theorem~\ref{thm:markov-main}, assume furthermore that Assumption~\ref{assume-markov-mixing} is satisfied with $c_0 = 1$ and set $\numcold = \numobs / 4$. There exists an event $\Event$, such that $\Prob (\Event) \geq 1 - e^{- \numobs / (16 \mixingtime)}$ and
  \begin{multline*}
    \Exs \big[ \vecnorm{\thetahat_\numobs - \thetabar}{2}^2 \bm{1}_{\Event} \big] \leq
\tfrac{c'}{\numobs} \mathrm{Tr} \big( (I - \Lbar)^{-1}
(\SigStar_{\mathrm{MG}} + \SigStar_{\mathrm{Mkv}})(I - \Lbar)^{-\top}
\big) + c' \big( \tfrac{ \sigmaA^2 \usedim \mixingtime}{(1 -
  \kappa)^2 \numobs} \big)^{4/3} \varbound^2 \log^2 \numobs\\
  + (\sigmab^2 + \sigmaA^2 \vecnorm{\thetabar}{2}^2)\cdot \exp \left( - \frac{\numobs}{32 \mixingtime} \right).
  \end{multline*}
\end{corollary}
\noindent See Section~\ref{app:subsec-proof-non-stationary-initial} for the proof of
this corollary. A few remarks are in order. Compared to the MSE bound in Theorem~\ref{thm:markov-main}, the bound in Corollary~\ref{cor:non-stationary-initial} exhibits two differences: we need to exclude an extreme event $\Event^c$ that occurs with exponentially small probability, and the right-hand-side of the bound involves an additional, exponentially decaying term. Roughly speaking, the high-probability event $\Event$ corresponds to the Markov chain states being close to a coupled chain. In the regime $\numobs \gg \mixingtime$ (which is implied by the sample size requirement in Theorem~\ref{thm:markov-main}), both the probability of the extreme event and the additional term are very small, and the guarantees under a non-stationary initial distribution behave qualitatively similar to the stationary case. For technical reasons, Corollary~\ref{cor:non-stationary-initial} requires a slightly stronger condition on the Markov chain---the transition kernel needs to be non-expansive under the metric $\metric$, that is, we require that $c_0 = 1$ in Assumption~\ref{assume-markov-mixing}. However, we note that such a non-expansive property holds for a wide range of applications: it is automatically satisfied in the case of $\metric (x, y) = \bm{1}_{x \neq y}$ and $\Wass_{1, \metric} = \totalvariation$. For general metric spaces, the mixing time bound in Assumption~\ref{assume-markov-mixing}(a) is usually established by showing the transition kernel $\transition$ is a contraction, i.e., $c_0 < 1$, which implies non-expansiveness (see e.g.~\cite{bou2020coupling}).
Finally, we note that the higher moment bounds for the last iterate established in Proposition~\ref{prop:lsa-markov-iterate-bound} can also be extended to the case of non-stationary initial distributions, yielding similar results.

\subsubsection{Bounds on the non-averaged iterates}

The proof of Theorem~\ref{thm:markov-main} involves first analyzing
the non-averaged iterates.  Since the upper bound established in this
step is of independent interest, we state and discuss it here:
\begin{proposition}
\label{prop:lsa-markov-iterate-bound}
Under
Assumptions~\ref{assume-markov-mixing}---\ref{assume:stationary-tail},
there are universal positive constants $(c_0, c_1)$ such that for any
integer $p \in \{1\} \cup [\log \numobs, \pmax / 2]$, scalar $\tau
\geq 2 p \mixingtime \log (c_0 d / \stepsize)$, and positive stepsize
\mbox{$\stepsize \in \big(0, \tfrac{1 - \kappa}{2 c p^3 ( \sigmaA^2 d
    + \conmax^2 ) \tau} \big]$,} we have
\begin{align}
  \big( \Exs \vecnorm{\theta_t - \thetabar}{2}^{2p} \big)^{1/p} \leq
  e^{- \frac{1}{2} \stepsize (1 - \kappa) t} \big( \Exs
  \vecnorm{\theta_0 - \thetabar}{2}^{2p} \big)^{1/p} + \frac{c p^3
    \stepsize}{1 - \kappa} \varbound^2 \tau d
\end{align}
\mbox{for all $t = 1, \ldots, \numobs$.}
\end{proposition}
\noindent See Section~\ref{subsec:proof-prop-lsa-iterate} for the
proof of this proposition.

Note that the guarantees on the unaveraged iterates in
Proposition~\ref{prop:lsa-markov-iterate-bound}---unlike those of
Theorem~\ref{thm:markov-main} for the averaged iterates---do not match
the optimal instance-dependent behavior.  This is to be expected,
since at least asymptotically, the unaveraged sequence converges to a
Gaussian random vector with covariance specified by the solution of a
Riccati equation.  (For details, see Section 4.5.3 of the book~\cite{benveniste2012adaptive}).  This covariance term need not
match the optimal statistical error.

On the other hand, by choosing $\stepsize \asymp \frac{\log
  \numobs}{(1 - \kappa) \numobs}$, the bound in
Proposition~\ref{prop:lsa-markov-iterate-bound} matches the worst-case
bound in equation~\eqref{eq:worst-case-bd}, up to log factors. We also
note that in Proposition~\ref{prop:lsa-markov-iterate-bound}, the
exponent $p$ can take values in two ranges: regardless of the value of
$\pmax \in [2, \infty]$, one can always take $p = 1$ and obtain an
upper bound on the mean-squared error $\Exs \big[ \vecnorm{\theta_t -
    \thetabar}{2}^{2} \big]$. This bound only requires
Assumption~\ref{assume:noise-moments} to hold true with $\pmax \geq
2$, which covers many important examples (see
Section~\ref{sec:consequences}). On the other hand, when
Assumption~\ref{assume:noise-moments} is satisfied with $\pmax \geq 2
\log \numobs$ and a stronger moment assumption is imposed, one can
obtain a $p$-th moment bound for any $p \geq [2\log \numobs,
  \pmax]$. This bound can be readily converted into a high-probability
bound for the last iterate of stochastic approximation. It is worth
noting that we study these two cases separately, using slightly
different proof techniques.

Let us now make some comparisons between
Proposition~\ref{prop:lsa-markov-iterate-bound} and existing results
on the unaveraged forms of Markovian stochastic approximation.  As we
have noted in our examples, in many cases, the quantities $(\sigmaA,
\sigmab, \varbound)$ do not depend on the dimension, in which case the
error bound in Proposition~\ref{prop:lsa-markov-iterate-bound} grows
linearly with dimension $\usedim$. In comparison, in
terms of our notation, the error bounds in the
papers~\cite{bhandari2018finite,srikant2019finite} both exhibit
quadratic dependency on the quantity $\tfrac{\max_{s \in \statespace}
  \opnorm{\Lmap_t (s)}}{1 - \kappa}$.  As we noted previously in
equation~\eqref{eq:constant-rank-obs-has-large-opnorm}, this quantity
scales linearly in dimension when the observations have a constant
rank (independent of dimension), so that (even after optimal parameter
tuning), the bounds from these parameters scale at least
proportionally to $\frac{d^2}{\numobs}$.  This scaling should be
contrasted with the $\myorder(\usedim/\numobs)$ guarantees from our
bounds. On a complementary note, the analysis in~\cite{durmus2021stability} involves a different mixing
assumption, and so is not directly comparable to our results.
However, it is worth noting that their bound $\vecnorm{\theta_t -
  \thetabar}{2}$ also has an explicit $\order {d / \sqrt{\numobs}}$
term (cf. equation (32) in their paper), showing that the MSE bound
grows quadratically with dimension.


\subsection{Local minimax lower bounds}
\label{SecLower}

Thus far, we established instance-dependent upper bounds for the averaged
SA scheme with Markov noise.  It is natural to wonder whether these
bounds can be improved.  Answering this question requires the
development of local minimax lower bounds, which we describe in this
section.

\subsubsection{Set-up and local neighborhoods}

We begin with the set-up and the definition of local neighborhoods for
our lower bounds.  Let $\transition$ be an irreducible Markov
transition kernel on a finite state space $\statespace$ with
associated stationary measure $\stationary_\transition$.  Consider the
solution $\thetabar(\transition)$ to the fixed-point equation
\begin{align}
\label{eq:fixed-point-in-lower-bound-setup}  
\thetabar(\transition) = \Exs_{\stationary_\transition} \big[ \Lmap
  (\state) \big] \cdot \thetabar (\transition) +
\Exs_{\stationary_\transition} \big[ \bmap (\state) \big].
\end{align}
where the maps $\bmap$ and $\Lmap$ are known to the estimator, whereas
the Markov transition kernel is unknown.  For some fixed
$\transition_0$ with stationary measure $\stationary_0$, we would like
to lower bound the number of observations required to estimate
$\thetabar(\transition_0)$ to a given accuracy.  In order to obtain
such a lower bound, we consider the fixed point
problem~\eqref{eq:fixed-point-in-lower-bound-setup} over a local
neighborhood\footnote{Doing so is necessary to rule out trivial
estimators, and the possibility of super-efficiency.} of the pair
$(\transition_0, \stationary_0)$.  We assume that the estimator is
based on a Markov trajectory $\{\state_t \}_{t=0}^\numobs$, with
initial state $\state_0$ drawn according to the original\footnote{ In
our construction, both kernels $\transition_0$ and $\transition$ are
rapidly mixing and their stationary measure are sufficiently close in
TV distance that the choice of initial distribution does not affect
the result.  Drawing $\state_0 \sim \stationary_0$ is made for
theoretical convenience.}  stationary distribution $\stationary_0$,
and successive states evolving according to the transition kernel
$\transition$.

In order to quantify the complexity of estimation localized around the
Markov transition kernel $\transition_0$, we define the following two
notions of local neighborhood:
\begin{subequations}
\begin{align}
\neighborhood_{\mathrm{Prob}} \big( \transition_0, \varepsilon \big)
&\mydefn \big\{ \transition: \sum_{x \in \statespace} \stationary_0
(x) \cdot \chisqdiv{ \transition (x, \cdot) }{\transition_0 (x,
  \cdot)} \leq \varepsilon^2 \big\}, \\
\neighborhood_{\mathrm{Est}} \big( \transition_0, \varepsilon \big)&
\mydefn \big\{ \transition: \vecnorm{\thetabar (\transition) -
  \thetabar (\transition_0)}{2} \leq \varepsilon \big\}.
\end{align}
\end{subequations}
The two notions of neighborhood focus on different types of locality
restrictions on the model class: the local problem class
$\neighborhood_{\mathrm{Prob}}$ contains all the Markov transition
kernels that are ``globally close'' to a given kernel $\transition_0$,
measured by a weighted $\chi^2$ divergence. It is worth noting that
this weighted $\chi^2$ divergence has an operational
interpretation. Suppose we draw $x \sim \stationary_0$, and then draw
the next state $y \sim \transition_0 (x, \cdot)$ accordingly the
original Markov kernel $\transition_0$, as well as $y' \sim
\transition (x, \cdot)$ under the kernel $\transition$. Then the
weighted $\chi^2$ divergence is the $\chi^2$ divergence between the
joint laws of $(x, y)$ and $(x, y')$.

On the other hand, the local class $\neighborhood_{\mathrm{Est}}$
contains Markov transition kernels $\transition$ such that the
solution $\thetabar (\transition)$ to the fixed-point
equation~\eqref{eq:fixed-point-in-lower-bound-setup} lies in a local
neighborhood of the given solution $\thetabar (\transition_0)$,
measured by the Euclidean distance. This problem class captures the
complexity specifically for solving the fixed-point equation, without
the need to estimate the entire transition kernel. In particular, it
is easy to construct a Markov kernel $\transition$ such that the
solution $\thetabar (\transition)$ is very close to $\thetabar
(\transition_0)$, but the distance between the transition kernels
$\transition$ and $\transition_0$ (e.g. measured in weighted $\chi^2$
divergence) is arbitrarily large.


\subsubsection{Instance-dependent lower bound}

Our lower bound is proved on the smallest worst-case risk attainable
over the intersection of $\neighborhood_{\mathrm{Prob}}$ and
$\neighborhood_{\mathrm{Est}}$.  We use the shorthand notation
$\Lbar^{(0)} \mydefn \Exs_{\stationary_0} \big[ \Lmap (\state)
  \big]$. Also recall the covariance matrix $\SigStar_{\mathrm{Mkv}}$ defined in
equation~\eqref{eq:def-sigstar}, for a stationary trajectory
$(\widetilde{\state}_t)_{t \in \mathbb{Z}}$ under the transition
kernel $\transition_0$.  Our bound depends on the \emph{local radius}
\begin{align}
\label{eq:radius-in-lower-bound}  
  \varepsilon_\numobs = \numobs^{-1/2} \, \sqrt{\trace \big( (I -
    \Lbar^{(0)})^{-1} \SigStar_{\mathrm{Mkv}} (I -
    \Lbar^{(0)})^{-\top} \big)},
    \end{align}
which is the contribution of Markovian noise to the upper bound stated
in Theorem~\ref{thm:markov-main}.

We are now ready to state our lower bound.  Recall that we have
assumed that the kernel $\transition_0$ is irreducible and aperiodic.
We also assume the mixing condition
(Assumption~\ref{assume-markov-mixing}) holds with the discrete metric
$\metric(x, y) = \bm{1}_{\{x \neq y\}}$ and mixing time $\mixingtime$,
and that $\mathrm{supp} \big( \transition_0(\state, \cdot)\big) \geq
2$ for all $\state \in \statespace$.
\begin{theorem}
\label{thm:local-minimax}
Under the assumptions stated above, there exist universal positive
constants $(c, c_1, c_2)$ such that for any sample size $\numobs$
lower bounded as
\begin{subequations}
\begin{align}
\numobs \geq \tfrac{c \mixingtime^2 \sigmaA^2 \usedim^2 \log^2
  \usedim}{(1 - \kappa)^2}, \quad \mbox{and} \quad \numobs^2
\varepsilon_\numobs^2 \geq \tfrac{2c (1 + \sigmaA^2) \varbound^2
  \mixingtime^4 \usedim^2}{(1 - \kappa)^4 } \log^6 \big(
\tfrac{\usedim}{\min \limits_\state \stationary_0(\state)} \big),
\end{align}
we have the minimax lower bound
\begin{align}
  \inf_{\thetahat_\numobs} ~\sup_{\transition \in \neighborhood'} \Exs
  \big[ \vecnorm{\thetahat_\numobs - \thetabar (\transition)}{2}^2
    \big] \geq c_2 \varepsilon_\numobs^2,
\end{align}
\end{subequations}
where $\neighborhood' \defn \neighborhood_{\mathrm{Prob}}
\big(\transition_0, c_1 \sqrt{\tfrac{\usedim}{\numobs}} \big) \cap
\neighborhood_{\mathrm{Est}} (\transition_0, c_1
\varepsilon_\numobs)$.
\end{theorem}
\noindent See Appendix~\ref{subsec:proof-lower-bound} for the proof of
this theorem.\\

A few remarks are in order. First, note that the minimax lower bound
is with respect to the problem class $\neighborhood_{\mathrm{Prob}}
\big(\transition_0, c_1 \sqrt{\tfrac{d}{n}} \big) \cap
\neighborhood_{\mathrm{Est}} (\transition_0, c_1
\varepsilon_\numobs)$, which requires both the transition kernel
$\transition$ and the solution $\thetabar (\transition)$ to be close
to the given problem instance $(\transition_0, \thetabar
(\transition_0))$. The size of the weighted $\chi^2$ neighborhood
scales with the standard parametric rate $\sqrt{d/n}$, as desired in
such problems.  On the other hand, the size of the neighborhood around
$\thetabar(P_0)$ is proportional to the local radius
$\varepsilon_\numobs$ that appears in the lower bound.  Operationally,
this result indicates that even if the estimator knows in advance that
$\thetabar (\transition)$ lies in the ball $\ball (\thetabar
(\transition_0), c_1\varepsilon_\numobs)$, one cannot do much better
than simply outputting an arbitrary point in this ball without looking
at the data. Let $\neighborhood_{\mathrm{global}}$ be the set of
Markov chain fixed-point equation problem instances satisfying
Assumptions~\ref{assume-markov-mixing}--\ref{assume-lip-mapping}. Following
the discussion in Section~\ref{subsubsec:main-upper-bound}, a
worst-case lazy Markov chain trajectory of length $\numobs$ yields an
effective $\mathrm{i.i.d.}$ sample size of order $\numobs /
\mixingtime$. Note that in the $\mathrm{i.i.d.}$ settings, the
fixed-point equation problem considered in this paper covers linear
regression (see~\cite{mou2020optimal}). Following the well-known
minimax lower bound for linear regression
(see~\cite{wainwright2019high}, Section 15.3), we can obtain a global
minimax lower bound using this reduction:
\begin{align*}
    \inf_{\thetahat_\numobs} ~\sup_{\transition \in
      \neighborhood_{\mathrm{global}}} \Exs \big[
      \vecnorm{\thetahat_\numobs - \thetabar (\transition)}{2}^2 \big]
    \gtrsim \frac{\mixingtime \usedim}{(1 - \kappa)^2\numobs},
\end{align*}
which is also achieved by Theorem~\ref{thm:markov-main}. Compared to
this global minimax lower bound, the local minimax formulation in
Theorem~\ref{thm:local-minimax} provides a more fine-grained
characterization of the minimax risk landscape across different
problem instances: there could be many different estimators that
achieve the global minimax lower bound (for example,
Proposition~\ref{prop:lsa-markov-iterate-bound} shows that the last
iterate is near-optimal up to logarithmic factors), but the
Polyak--Rupert averaged estimator is optimally adaptive to the
complexity associated to any problem instance, characterized by the
quantity $\varepsilon_\numobs^2$.

Second, it should be noted that quantity $\varepsilon_\numobs^2$
matches (up to a constant factor) the optimal mean-squared error given
by the local asymptotic minimax
theorem~\cite{van2000asymptotic,greenwood1995efficiency}.  In contrast
to such asymptotic theory, however, Theorem~\ref{thm:local-minimax}
applies when $\numobs$ is finite, and does not impose any regularity
assumptions on the estimator.  Furthermore, the radius
$\varepsilon_\numobs$ that is used to define the local neighborhood
$\neighborhood_{\mathrm{Est}} (\transition_0, \varepsilon_\numobs)$ is
optimal in the following sense. On the one hand, since the plug-in
estimator is asymptotically normal~\cite{greenwood1995efficiency}, for
any decreasing sequence $\varepsilon_\numobs'$ such that
$\varepsilon_\numobs' > \varepsilon_\numobs$ and $\varepsilon_\numobs'
\rightarrow 0^+$, the minimax risk within the neighborhood
$\neighborhood_{\mathrm{Est}} (\transition_0, \varepsilon_\numobs')$
behaves asymptotically as $\varepsilon_\numobs^2$ up to constant
factors. On the other hand, for any decreasing sequence
$\varepsilon_\numobs'$ such that $\varepsilon_\numobs' <
\varepsilon_\numobs$, the minimax risk in the neighborhood
$\neighborhood_{\mathrm{Est}} (\transition_0, \varepsilon_\numobs')$
is at most $\varepsilon_\numobs'$. In the latter case, the
neighborhood is so small that it provides more information than the
data provides.

Third, note that Theorem~\ref{thm:local-minimax} involves inexplicit
universal constants $(c_1, c_2)$. We do not optimize these constants
in our proof, and our proof gives a bound with $c_1 = 1$ and $c_2 =
\frac{1}{4 (5 + \pi)}$. Note that the local asymptotic minimax lower
bound for this problem~\cite{greenwood1995efficiency} implies that
\begin{align*}
  \limsup_{c_1 \rightarrow + \infty} ~\lim_{\numobs \rightarrow +
    \infty} ~\inf_{\thetahat_\numobs} ~\sup_{\transition \in
    \neighborhood_{\mathrm{Prob}} (\transition_0, c_1 /
    \sqrt{\numobs})} \Exs \big[ \vecnorm{\thetahat_\numobs - \thetabar
      (\transition)}{2}^2 \big] \geq \varepsilon_\numobs^2,
\end{align*}
which suggests that the constants can be sharpened.  Indeed, using
more careful arguments, the non-asymptotic lower bound exhibits a
similar nature: in Theorem~\ref{thm:local-minimax}, if we take the
constant $c_1$ in the size of the neighborhood sufficiently large, the
pre-factor $c_2$ in the minimax lower bound can be made close to $1$,
exactly matching the asymptotic lower bound and the refined upper
bound (see discussion following Theorem~\ref{thm:markov-main} for
details). In doing so, we can re-scale the prior distribution in the
proof of Theorem~\ref{thm:local-minimax} with $c_1$, and the
leading-order term in the Bayesian Cram\'{e}r--Rao lower bound will
come with a pre-factor $1$. For brevity, we do not dive into details
of this argument.

Fourth, Theorem~\ref{thm:local-minimax} matches the Markov noise term
in Theorem~\ref{thm:markov-main}, establishing its optimality when the
martingale part of the noise vanishes, i.e., $\Lmap_{t} (s) = \Lmap
(s)$ and $\bmap_{t} (s) = \bmap (s)$. The lower bound does not capture
the martingale part of the noise because we assume that the functions
$\Lmat: \statespace \rightarrow \real^{d \times d}$ and
$\bvec: \statespace \rightarrow \real^d$ are known to the
estimator. In the setting where these functions are also observed only
through noisy $\myiid$ data $(\Lmat_{t}, \bvec_t)$, Theorem 3 in the
paper~\cite{mou2020optimal} implies a lower bound of the form $c_2
\numobs^{-1} \trace \big( (I - \Lbar^{(0)})^{-1}
\SigStar_{\mathrm{MG}} (I - \Lbar^{(0)})^{-\top} \big)$. Combining it
with Theorem~\ref{thm:local-minimax} implies a minimax lower bound
involving the term $c_2' \numobs^{-1} \trace \big( (I -
\Lbar^{(0)})^{-1} (\SigStar_{\mathrm{Mkv}} + \SigStar_{\mathrm{MG}})
(I - \Lbar^{(0)})^{-\top} \big)$ in a properly defined local
neighborhood, thus establishing the optimality of
Theorem~\ref{thm:markov-main}. At the same time, we note that
Theorem~\ref{thm:local-minimax} requires the sample size to be at
least $\mixingtime^2 d^2$, which is more stringent than the $\order
{\mixingtime d}$ requirement in the upper bound. While
Theorem~\ref{thm:markov-main} holds true with a linear sample-size $n
= \order{d}$, it is only shown to be instance-optimal for larger
$\numobs = \Omega(d^2)$. This mismatch is due to the fact that small
perturbations of the Markov transition kernel in certain directions
can destroy its fast mixing property. That being said,
Theorem~\ref{thm:local-minimax} is still a finite-sample result, with
polynomial dependency on the quantities $\big(\mixingtime, \usedim,
\frac{1}{1 - \kappa} \big)$, and poly-logarithmic dependency on the
quantity $\min_s \xi_0(s)$.


\section{Some consequences for specific problems}
\label{sec:consequences}

In this section, we specialize our analysis to the examples described
in Section~\ref{subsec:examples}, namely approximate policy evaluation
using TD algorithms, and estimation in autoregressive time series
models.  By verifying the conditions needed to apply
Theorem~\ref{thm:markov-main} and
Proposition~\ref{prop:lsa-markov-iterate-bound}, we obtain some more
concrete corollaries of our general theory.


\subsection{TD$(0)$ method}

Recall the TD$(0)$ algorithm for policy evaluation, as previously
described in Example~\ref{example:td0}. We are interested in estimating the solution $\valuestar$ of the Bellman equation~\eqref{eq:bellman} when an approximation scheme is employed using the basis functions $(\phi_j)_{j = 1}^{\usedim}$. Using the shorthand
$\inprod{\theta}{\phi(\state)} = \sum_{j=1}^\usedim \theta_j
\phi_j(\state)$ for the Euclidean inner product in $\real^\usedim$, with observation model $(\Lmap_{t + 1} (\slidingwindow_t), \bmap_{t + 1} (\slidingwindow_t))$ defined in equation~\eqref{eq:fit-td0-into-our-setup}, the averaged SA procedure~\eqref{eq:lsa} is given by:
\begin{align}
\label{eq:td-0-iterates} 
    \theta_{t + 1} \stackrel{(a)}{=} \theta_t - \stepsize \big \{
    \inprod{\phi(s_t) - \discount \phi(s_{t + 1})}{\theta_t} - R_{t+1} (s_t)
    \big \} \phi(s_t), \quad \mbox{and} \quad
    \thetahat_\numobs \stackrel{(b)}{=} \tfrac{1}{\numobs - \numburn}
    \sum_{t = \numburn}^{\numobs - 1} \theta_t.
\end{align}
To be clear, the update~\eqref{eq:td-0-iterates}(a) is the standard
TD$(0)$ algorithm with stepsize $\stepsize$, whereas the addition of
the averaging step~\eqref{eq:td-0-iterates}(b) yields the
Polyak--Ruppert averaged version of the scheme. Note that we re-scale the stepsize $\stepsize$ by a factor of $\beta$ for notational convenience. In the following
subsections, we derive corollaries of our general theory for the
averaged scheme under different mixing conditions on the underlying
Markov chain.


\subsubsection{Markov chains with mixing in total variation distance}

We first assume that the Markov chain satisfies a mixing condition
(cf.  Assumption~\ref{assume-markov-mixing}) in the discrete metric:
i.e., after $\mixingtime$ steps, we have $\totalvariation
(\delta_{\state} \transition^{\mixingtime}, \delta_{\state'}
\transition^{\mixingtime}) \leq \tfrac{1}{2}$ \mbox{for any pair
  $\state, \state' \in \statespace$.}  Let $\stationary$ denote the
stationary distribution of the Markov chain that generates the
trajectory $\{ \state_t \}_{t \geq 0}$, and let $\transition$ denote
its transition kernel.  Note that the augmented state vector
\mbox{$\slidingwindow_t = (\state_t, \state_{t + 1})$} evolves
according to a Markov process with mixing time $\mixingtime +
1$. Moreover, the stationary distribution of the pair $\slidingwindow
= (s, s^+)$ has the form $s \sim \stationary, ~ s^+ \sim \transition
(\cdot \mid s)$.  We denote the stationary covariance of the feature vectors as
$\Bmat \mydefn \Exs_{s \sim \stationary} \big[ \phi(\state)
  \phi(\state)^\top \big]$, and also define the minimum and maximum
eigenvalues $\mu \mydefn \lambda_{\min} (\Bmat)$ and $\beta \mydefn
\lammax(\Bmat)$.  We assume that
\begin{subequations}
\begin{align}
\label{eq:td0-bounded-feature-vec}  
 \vecnorm{\Bmat^{-1/2} \phi (s)}{2} \stackrel{(a)}{\leq} \varsigma
 \sqrt{d} \quad \mbox{and} \quad |R_t (s)| \stackrel{(b)}{\leq}
 \varsigma \quad \mbox{for all $s \in \statespace$, and} \\
\label{eq:td0-fourth-moment} 
\Exs_{\stationary} \big[ \inprod{\Bmat^{-1/2} \phi (s)}{u}^4 \big]
\leq \varsigma^4 \quad \mbox{for all $u \in \sphere^{d - 1}$.}
\end{align}
\end{subequations}
In order to state our result, we define the following quantities:
\begin{align*}
\SpecMat & \mydefn \discount \Bmat^{-1/2} \cdot \Exs_{s \sim
  \stationary, s^+ \sim \transition (s, \cdot)} \big[ \phi(s) \phi
  (s^+)^\top \big] \cdot \Bmat^{-1/2},\\
\varepsilon_{\mathrm{Mkv}}(s, s^+) & \mydefn \Bmat^{-1/2} \big( \phi
(s)^\top \thetabar - \discount \phi(s^+)^\top \thetabar - r (s) \big)
\phi (s),\\
 \varepsilon_{\mathrm{MG}} (s) &\mydefn \Bmat^{-1/2} (R
(s) - r (s)) \phi (s) 
\end{align*}
We also define the following covariance matrices according to equation~\eqref{eq:def-sigstar}:
\begin{align*}
\SigStar_{\mathrm{Mkv}} & \mydefn \sum_{t = - \infty}^{\infty} \Exs
\big[\varepsilon_{\mathrm{Mkv}} (\state_t, \state_{t + 1})
  \varepsilon_{\mathrm{Mkv}} (\state_0, \state_1)^\top \big], \\
\SigStar_{\mathrm{MG}} &\mydefn \Exs_{s \sim \stationary} \big[ \Exs
  \big[ \varepsilon_{\mathrm{MG}} (s) \varepsilon_{\mathrm{MG}}
    (s)^\top \mid s \big] \big].
\end{align*}
Finally, we define the quantity
\begin{align}
  \varbound^2 \mydefn \varsigma^2 \cdot \sqrt{\Exs \big[
      \big(\phi(s_t)^\top \thetabar - \discount \phi (s_{t + 1})
      \thetabar - R_t (s_t) \big)^4
      \big]}, \label{eq:varbound-in-td}
\end{align}
and let $\kappa \mydefn \frac{1}{2 }\lammax (\SpecMat +
\SpecMat^\top)$. It is easy to see that $\kappa \leq \discount <
1$. Assuming that $\mu > 0$, we are then ready to state our main
result for the TD$(0)$ method.
\begin{corollary}
\label{cor:TD-0-finite-state-markov}
Under the setup above, take the stepsize $\stepsize$ and burn-in
period $\numburn$ as
\begin{align}
\label{eq:td0-stepsize-burnin-req}  
\stepsize = \tfrac{1}{c \beta ( (\varsigma^4 + 1) d (1 - \kappa)
  \numobs^2 \mixingtime)^{1/3}}, \quad \mbox{and} \quad \numburn =
\tfrac{1}{2} \numobs,
\end{align}
and suppose that $\frac{\numobs}{\log^3 \numobs} \geq \frac{2
  \mixingtime (\varsigma^4 + 1) d\beta^2}{(1 - \kappa)^2 \mu^2}$. The
estimator $\valuehat_\numobs \mydefn \thetahat_\numobs \phi$ obtained
from the Polyak--Ruppert procedure~\eqref{eq:td-0-iterates} satisfies the bound
\begin{multline}
\label{eq:td0-final-bound}  
 \Exs \big[ \vecnorm{\valuehat_\numobs - \valuebar}{\Ltwospace
     (\statespace, \stationary)}^2 \big] \leq \tfrac{c}{\numobs}
 \mathrm{Tr} \big \{ (I_d - M)^{-1} (\SigStar_{\mathrm{Mkv}} +
 \SigStar_{\mathrm{MG}}) (I_d - M)^{- \top} \big \}\\
  + c \big(
 \tfrac{\beta^2 \varsigma^4 d \mixingtime }{\mu^2 (1 - \kappa)^2
   \numobs} \big)^{4/3} \varbound^2 \log^2 \numobs,
\end{multline}
where $\valuebar$ is the solution to the projected fixed-point
equation~\eqref{eq:projected-fixed-point} and $c > 0$ is a universal
constant.
\end{corollary}
\noindent See Appendix~\ref{subsec:proof-cor-td-0-finite-state-markov}
for the proof of this corollary.\\

A few remarks are in order. First, we measure the estimation error in
the canonical $\vecnorm{\cdot}{\Ltwospace (\statespace, \stationary)}$
norm, instead of the Euclidean distance in $\real^d$.  Consequently,
the proof of this corollary actually uses a generalized version of
Theorem~\ref{thm:markov-main} proved for weighted $\ell^2$ norms. On
the other hand, we note that the error
bound~\eqref{eq:td0-final-bound} is with respect to the solution
$\valuebar$ to the projected fixed-point equation. In the
well-specified case where $\valuestar \in \LinSpace$, this solution
coincides with the value function $\valuestar$. In general, the
approximation error needs to be taken into account, and this was the focus
of our prior paper~\cite{mou2020optimal}.  In conjunction with this
result, Corollary~\ref{cor:TD-0-finite-state-markov} implies the error
bound
\begin{multline}
  \Exs \big[ \vecnorm{\valuehat_\numobs - \valuestar}{\Ltwospace
      (\statespace, \stationary)}^2 \big] \\
      \leq c \big[1 + \lammax
    \big( (I_d - \SpecMat)^{-1} (\discount^2 I_d - \SpecMat
    \SpecMat^\top) (I_d - \SpecMat)^{- \top} \big) \big]
  \inf_{\valuefunc \in \LinSpace} \vecnorm{\valuefunc -
    \valuestar}{\Ltwospace (\statespace, \stationary)}^2 \\
  \label{eq:td0-approx-factor-tradeoff}
   + \tfrac{c}{\numobs} \mathrm{Tr} \big \{ (I_d - M)^{-1}
  (\SigStar_{\mathrm{Mkv}} + \SigStar_{\mathrm{MG}}) ) (I_d - M)^{-
    \top} \big \} + c \big( \tfrac{\beta^2 \varsigma^4 d \mixingtime
  }{\mu^2 (1 - \kappa)^2 \numobs} \big)^{4/3} \varbound^2 \log^2 \numobs.
\end{multline}
In Section~\ref{subsec:consequence-td-lambda} to follow, we provide a
general recipe to trade off approximation and estimation errors to
choose the value of $\lambda$ in the class of TD$(\lambda)$
algorithms. Before that, we discuss two extensions of
Corollary~\ref{cor:TD-0-finite-state-markov}.

\subsubsection{Markov chains with mixing in Wasserstein metric}
Note that for Corollary~\ref{cor:TD-0-finite-state-markov}, the mixing
time condition is imposed with total variation distance. When the
state space $\statespace$ is continuous, e.g., the set
$\statespace$ is a subset of $\real^m$, mixing in Wasserstein distance
could capture the geometry of the underlying metric better. In this
section, we extend our analysis to such settings, highlighting the
dimension dependency in the sample complexity.

Concretely, we consider a Markov chain $(\state_t)_{t \geq 0}$ on a
compact domain $\statespace \subseteq \real^m$, and a feature mapping
$\phi: \statespace \rightarrow \real^d$. We assume that the Markov
chain admits a unique stationary measure $\stationary$, and the mixing
time assumption holds in Wasserstein-1 distance, so that $\Wass_1
\big( \delta_x \transition^{\mixingtime}, \delta_y
\transition^{\mixingtime} \big) \leq \tfrac{1}{2} \vecnorm{x - y}{2}$
for all $x, y \in \statespace$.  For the sake of normalization, we
assume that $\statespace \subseteq \ball (0, 1)$ and $\phi(0) = 0$. On
the feature mapping $\phi$, we assume the following:
\begin{subequations}
\label{eq:continuous-space-td-assumptions}
\begin{align}
 \exists
    \mu, \beta > 0, &\quad \mu I_d \preceq \Bmat \mydefn \Exs_{s \sim
      \stationary} \big[ \phi (s) \phi (s)^\top \big] \preceq \beta
    I_d,\label{eq:continuous-space-td-well-conditioned}\\ 
    \forall x, y \in \statespace, &\quad \vecnorm{\Bmat^{-1/2} \big(\phi (x) - \phi
      (y) \big)}{2} \leq \varsigma \sqrt{d} \vecnorm{x -
      y}{2},\label{eq:continuous-space-td-lip}\\
      \forall u \in
    \sphere^{d - 1}, & \quad \Exs_{s \sim \stationary} \big[
      \inprod{u}{\Bmat^{-1/2}\phi (s)}^4 \big] \leq
    \varsigma^4, \label{eq:continuous-space-td-fourth-moment}\\
    \forall s,
    s' \in \statespace, ~ t \geq 1, &\quad |R_t (s) - R_t (s')| \leq
    \varsigma \vecnorm{s - s'}{2},~ |R_t (s)| \leq \varsigma \quad
    \mbox{ a.s.}
\end{align}
\end{subequations}
Here, we regard the parameters $(\varsigma, \mu, \beta)$ as
dimension-independent positive constants. Since the state space
$\statespace$ has diameter bounded by $2$, the feature mapping $\phi$
satisfying equation~\eqref{eq:continuous-space-td-well-conditioned}
necessarily has Lipschitz constant of order $\order {\sqrt{d}}$. For a simple example, 
take the state $x$ itself as the feature
vector (after appropriate re-scaling), which corresponds to the case
of $m = d$ and $\phi(x) = \sqrt{d} \cdot x$. \\

\noindent With this set-up, we have the following guarantee:
\begin{corollary}
\label{prop:td-0-continuous-statespace}
Assuming the conditions in
equation~\eqref{eq:continuous-space-td-assumptions}, taking stepsize
and burn-in period as equation~\eqref{eq:td0-stepsize-burnin-req}, for
the Polyak--Ruppert averaged stochastic approximation procedure~\eqref{eq:td-0-iterates}, the
bound~\eqref{eq:td0-final-bound} holds.
\end{corollary}
\noindent See Appendix~\ref{subsubsec:proof-prop-td-0-cts-statespace}
for the proof.\\

Corollary~\ref{prop:td-0-continuous-statespace} shows that the same
instance-dependent bound holds true for a continuous state space
setting. Such a bound is useful for many applications; for example, in
the case of quadratic value functions on a subset of $\real^m$, the
feature mapping takes the form $x \mapsto \phi(x) \defn m x x^T$, so
that the dimension $d = m^2$.  Assuming that the process $(s_t)_{t
  \geq 0}$ is supported in a unit ball $\ball (0, 1)$ and has
well-conditioned stationary covariance, it is easy to verify that
Assumptions~\eqref{eq:continuous-space-td-assumptions} are satisfied
with dimension-free constants $(\varsigma, \mu, \beta)$. This example
is particularly useful for policy evaluation in Linear Quadratic
Regulators (LQR), and more generally for other stochastic dynamical
systems.

\subsubsection{Analysis of a sieve estimator}
\label{subsubsec:sieve-consequence}

The optimal dimension dependency in Theorem~\ref{thm:markov-main}
allows us to obtain optimal estimators for various classes of
non-parametric problems, in which the dimension is a parameter to be
chosen.  In particular, sieve methods are a class of non-parametric
estimators based on nested sequences of finite-dimensional
approximations.  In this section, we analyze the behavior of a
stochastic approximation sieve estimator in the Markovian setting. The
optimal dimension dependence in our theorem recovers the minimax
optimal rates for estimation, while our instance-dependent bounds help
in capturing more refined structure in the problem instance.

Concretely, assuming that the Hilbert space $\Ltwospace (\statespace,
\stationary)$ is separable, let $(\phi_j)_{j = 1}^{\infty}$ be a set
of (not necessarily orthogonal) basis functions. We consider the case
where the mixing condition holds true with total variation
distance\footnote{By following the approach in the previous
subsection, the analysis can also be extended to the case of mixing in
Wasserstein distance.}. The following assumptions are imposed on the
basis functions:
\begin{subequations}\label{eq:sieve-td-assumptions}
  \begin{align}
      \forall j \in \mathbb{N}^+, &\quad \sup_{x \in \statespace}
      |\phi_j (x)| \leq \varsigma,\\ \forall d \in \mathbb{N}^+,
      &\quad \mu I_d \leq \big[ \Exs_{s \sim \stationary} \big( \phi_j
        (s) \phi_\ell (s) \big) \big]_{j, \ell \in [d]} \leq \beta
      I_d,\\ \forall t \geq 1, &\quad \sup_{x \in \statespace} |R_t
      (x)| \leq \varsigma.
  \end{align}
\end{subequations}
The first assumption is standard in nonparametric regression, and
satisfied by many useful basis functions such as the Fourier basis and
Walsh-Hadamard basis. The second assumption relaxes the orthogonality
requirement on the bases, by only requiring the Gram matrix to be
well-conditioned.

We define the noise level $\varbound$ using the second moment:
\begin{align}
    \varbound^2 \mydefn \varsigma^2 \cdot \Exs \big[ \big(
        \valuebar (s_t) - \discount \valuebar (s_{t + 1}) - R_t
        (s_t) \big)^2 \big].
\end{align}

Once again, we run the averaged stochastic approximation
procedure~\eqref{eq:td-0-iterates} on this problem. A crucial point of
departure from the parametric models discussed above is that the
number of basis functions $d_\numobs$ in sieve estimators is chosen
based on the problem structure and sample size. Let $\LinSpace (d_n)
\mydefn \mathrm{span} (\phi_1, \phi_2, \cdots, \phi_{d_n})$ denote the
subspace spanned by the first $d_n$ basis functions. The following
result is a direct corollary of our theorem, and covers the case of
fixed $d_n$; we discuss the trade-off between approximation and
estimation error in the choice of $d_\numobs$ presently.

\begin{corollary}
\label{prop:sieve-td-results}
Assuming the conditions in equation~\eqref{eq:sieve-td-assumptions},
take the stepsize and burn-in period as in
equation~\eqref{eq:td0-stepsize-burnin-req}. Assuming that $\mu,
\beta, \varsigma \asymp 1$, the Polyak--Ruppert averaged stochastic
approximation procedure~\eqref{eq:td-0-iterates} satisfies the
bound~\eqref{eq:td0-final-bound} with $d = d_n$.
\end{corollary}
\noindent See Appendix~\ref{subsubsec:proof-prop-sieve-td-results} for
the proof. \\

Recall that by taking into account the approximation error, the error
for estimating the true value function $\valuestar$ takes the
following form:
\begin{multline*}
\Exs \big[ \vecnorm{\valuehat_\numobs - \valuestar}{\Ltwospace
    (\statespace, \stationary)}^2 \big] \\
    \leq c \big[1 + \lammax
  \big( (I - \SpecMat)^{-1} (\discount^2 I_d - \SpecMat \SpecMat^\top)
  (I - \SpecMat)^{- \top} \big) \big] \inf_{\valuefunc \in \LinSpace}
\vecnorm{\valuefunc - \valuestar}{\Ltwospace (\statespace,
  \stationary)}^2 \\  + \tfrac{c}{\numobs} \mathrm{Tr}
\big( (I - M)^{-1} (\SigStar_{\mathrm{Mkv}} + \SigStar_{\mathrm{MG}})
(I - M)^{- \top} \big) + c \big( \tfrac{\varbound^2 \mixingtime d_n}{
  (1 - \kappa)^2 \numobs} \big)^{4/3} \log^2 \numobs.
\end{multline*}

Let $\{\psi_j \}_{j = 1}^{+\infty}$ be an orthonormal basis of
$\Ltwospace (\statespace, \stationary)$ such that $\mathrm{span}
(\psi_1, \cdots, \psi_d) = \mathrm{span} (\phi_1, \cdots, \phi_d)$ for
any $d \geq 1$.  (For instance, one can let $\{\psi_j\}_{j = 1}^{
  +\infty}$ be the Gram-Schmidt orthonormalization of the original
basis functions). Given a non-increasing sequence $\{\alpha_j\}_{j =
  1}^{\infty}$ of positive reals such that $\lim_{j \rightarrow +
  \infty} \alpha_j = 0$, we first let $\mathcal{H}_0$ be a linear
subspace of $\Ltwospace (\statespace, \stationary)$, consisting of all
the finite linear combination of basis vectors $\{\psi_j\}_{j =
  1}^{+\infty}$, equipped with the following inner product:
  \begin{align*}
      \forall u, v \in \mathcal{H}_0, \quad
      \inprod{u}{v}_{\mathcal{H}_0} \mydefn \sum_{j = 1}^{\infty}
      \alpha_j^{-1} \cdot \inprod{u}{\psi_j} \cdot \inprod{v}{\psi_j}.
  \end{align*}
  Note that the summation shown above is actually finite, since both
  both sequences $(\inprod{u}{\psi_j})_{j = 1}^{+ \infty},
  ~(\inprod{v}{\psi_j})_{j = 1}^{+ \infty}$ only have a finite number
  of non-zero entries. We then define the inner product space
  $(\mathcal{H}, \inprod{\cdot}{\cdot}_{\mathcal{H}})$ as the
  completion of $(\mathcal{H}_0,
  \inprod{\cdot}{\cdot}_{\mathcal{H}_0})$. It is easy to see that
  $\mathcal{H}$ is a Hilbert space, and a linear subspace of
  $\Ltwospace (\statespace, \stationary)$.
  
For any $\valuestar \in \mathcal{H}$, the estimation error is at most
(in the worst-case)
\begin{align}
\label{eq:sieve-fullbd}
    \Exs \big[ \vecnorm{\valuehat_\numobs - \valuestar}{\Ltwospace
        (\statespace, \stationary)}^2 \big] \leq \tfrac{c}{1 -
      \discount} \cdot \alpha_{d_n}
    \vecnorm{\valuestar}{\mathcal{H}}^2 + \tfrac{c \varbound^2 d_n
      \mixingtime}{(1 - \discount)^2 n}.
\end{align}
For example, when the eigenvalues of Hilbert space $\mathcal{H}$ decay
as $\alpha_j \asymp j^{- 2 s}$ for some $s > 0$, the estimator
achieves a rate of $\myorder \big( (\mixingtime /
\numobs)^{\frac{2s}{2s + 1}} \big)$, which matches the minimax optimal
rate proved by~\cite{duan2021optimal} in the i.i.d. setting, but with
a multiplicative correction to the effective sample size by a factor
$\mixingtime$ to accommodate Markovian observations. Furthermore,
since one can estimate the quantities $(\SpecMat,
\SigStar_{\mathrm{Mkv}}, \SigStar_{\mathrm{MG}})$ in the
bound~\eqref{eq:td0-final-bound} using $\order{d}$ samples,
instance-dependent model selection can in principle be
conducted. Bounds of the form~\eqref{eq:sieve-fullbd} thus open the
door to asking important questions of this type.


\subsection{TD$(\lambda)$ methods}
\label{subsec:consequence-td-lambda}

Now we turn to stochastic approximation methods for the TD$(\lambda)$
projected fixed-point equation~\eqref{eq:population-td-lambda}, with
some given $\lambda \in [0, 1)$. With observation model $(\Lmap_{t +
    1} (\slidingwindow_t), \bmap_{t + 1} (\slidingwindow_t))$ given by
  equation~\eqref{eq:fit-tdlambda-into-our-setup}, the averaged SA
  procedure~\eqref{eq:lsa} can be written as
\begin{subequations}
\label{eq:td-lambda}
\begin{align}
\label{eq:td-lambda-iterates}  
\theta_{t + 1} & = \theta_t - \stepsize \Big\{ \inprod{\phi (s_t) -
  \discount \phi (s_{t + 1})^\top }{ \theta_t} - R_t (s_t) \Big\} g_t,
\quad \mbox{where} \\
\label{eq:td-lambda-g-update}
g_t &= \discount \lambda g_{t - 1} + \phi (s_t) \quad \mbox{and,} \\
\label{eq:td-lambda-avg}
\thetahat_\numobs &= \tfrac{1}{\numobs - \numburn} \sum_{t =
  \numburn}^{\numobs - 1} \theta_t.
\end{align}
\end{subequations}
The update on $g_t$ is the so-called ``eligibility trace'' in the
TD$(\lambda)$ algorithm.  As before, we assume the two bounds in
equation~\eqref{eq:td0-bounded-feature-vec}, and assume that the
mixing time condition in Assumption~\ref{assume-markov-mixing} holds
true for the chain $(\state_t)_{t \geq 1}$, with discrete metric and
mixing time $\mixingtime$. We consider the augmented Markov chain
$\slidingwindow_t \mydefn \big( s_t, s_{t + 1}, \tfrac{1 - \discount
  \lambda}{\varsigma \sqrt{\beta d}} g_t \big) \in \statespace^2
\times \ball(0, 1)$ and begin by establishing mixing conditions on
this augmented chain.
\begin{proposition}
  \label{prop:mixing-td-lambda-augmented-chain}
Under the setup above, consider the metric
\begin{subequations}  
\begin{align}
  \metric \big( (s_1, s_2, h), (s_1', s_2', h') \big) \mydefn
  \tfrac{1}{4} \big( \bm{1}_{s_1 \neq s_1'} + \bm{1}_{s_2 \neq s_2'} +
  \vecnorm{h - h'}{2} \big).
\end{align}
Taking $\tau = 4 \big( \mixingtime + \frac{1}{1 - \discount \lambda}
\big)$, the augmented chain $\big\{\slidingwindow_t =(s_t, s_{t + 1},
\frac{1 - \discount \lambda}{\varsigma \sqrt{\beta d}} g_t) \big\}_{t
  \geq 0}$ satisfies the mixing bound
\begin{align}
      \Wass_{1, \metric} \big( \law (\slidingwindow_\tau), \law
      (\slidingwindow_\tau') \big) \leq \frac{1}{2} \metric
      \big(\slidingwindow_0, \slidingwindow_0' \big)
  \end{align}
for two chains $(\slidingwindow_t)_{t \geq 0}$ and
$(\slidingwindow_t')_{t \geq 0}$ starting from $\slidingwindow_0$ and
$\slidingwindow_0'$, respectively. In particular, the stationary
distribution $\widetilde{\stationary}$ of the chain
$(\slidingwindow_t)_{t \geq 0}$ exists and is unique.
\end{subequations}  
\end{proposition}
\noindent See
Appendix~\ref{subsubsec:proof-mixing-td-lambda-augmented-chain} for
the proof of this proposition.\\

Taking this proposition as given, we are now ready to present our main
corollary for TD$(\lambda)$ procedures. We consider the following
instantiation of quantities in Theorem~\ref{thm:markov-main}: \\

The projected linear operator $(1 - \lambda) \sum_{k = 0}^{+ \infty}
\lambda^k (\discount \projecttolin \transition)^{k + 1}$ in the
equation~\eqref{eq:population-td-lambda} can be represented in the
orthonormal basis of the subspace $\LinSpace$ as
\begin{align*}
    \SpecMat_\lambda &\mydefn I_d - \Bmat^{-1/2} \Exs_{(s, s^+,
      \frac{1 - \discount \lambda}{\varsigma \sqrt{\beta d}} g) \sim
      \widetilde{\stationary}} \big[ g \phi (s)^\top - \discount g
      \phi (s^+)^\top \big] \Bmat^{-1/2}\\ &=(1 - \lambda)
    \Bmat^{-1/2} \sum_{t = 0}^{\infty} \lambda^t \discount^{t + 1}
    \Exs \big[ \phi (s_{0}) \phi (s_{t + 1}) \big] \Bmat^{-1/2}.
\end{align*}
The Markovian and martingale part of the noise (in the low-dimensional
subspace $\LinSpace$) takes the form
\begin{align*}
    \varepsilon_{\mathrm{Mkv}, \lambda} &\big(s, s^+, \frac{1 -
      \discount \lambda}{\varsigma \sqrt{\beta d}} g \big) =
    \Bmat^{-1/2} \big(\phi (s)^\top \thetabar - \discount \phi (s^+)
    \thetabar - r(s) \big) g,\\ \varepsilon_{\mathrm{MG}, \lambda}
    &\big(s, s^+, \frac{1 - \discount \lambda}{\varsigma \sqrt{\beta
        d}} g \big) = \Bmat^{-1/2} (R_0 (s) - r(s)) g
\end{align*}
Finally, we define the covariance matrices $\SigStar_{\mathrm{Mkv},
  \lambda}$ and $\SigStar_{\mathrm{MG}, \lambda}$ according to
equation~\eqref{eq:def-sigstar}:
\begin{align*}
    \SigStar_{\mathrm{Mkv}, \lambda} &\mydefn \sum_{t = -
      \infty}^{\infty} \Exs \big[\varepsilon_{\mathrm{Mkv}, \lambda}
      \big(s_t, s_{t + 1}, \frac{1 - \discount \lambda}{\varsigma
        \sqrt{\beta d}} g_t \big) \varepsilon_{\mathrm{Mkv}, \lambda}
      \big(s_0, s_{1}, \frac{1 - \discount \lambda}{\varsigma
        \sqrt{\beta d}} g_0 \big)^\top \big]
    ,\\ \SigStar_{\mathrm{MG}, \lambda} &\mydefn \Exs_{s \sim
      \stationary} \big[ \Exs \big[ \varepsilon_{\mathrm{MG}, \lambda}
        (s) \varepsilon_{\mathrm{MG}, \lambda} (s)^\top \mid s \big]
      \big].
\end{align*}
As before, we let $\beta \mydefn \lammax (\Bmat)$, $\mu \mydefn
\lambda_{\min} (\Bmat)$ and $\kappa_\lambda \mydefn \frac{1}{2
}\lammax (\SpecMat_\lambda + \SpecMat_\lambda^\top)$, and define the
quantity $\varbound$ according to
equation~\eqref{eq:varbound-in-td}. Note that a straightforward
calculation reveals that $\kappa_\lambda \leq \frac{(1 - \lambda)
  \discount}{1 - \lambda \discount}< 1$. Assuming that $\mu > 0$, we
are then ready to state our main result for TD$(\lambda)$ methods.
\begin{corollary}
  \label{cor:TD-lambda-finite-state-markov}
Under the setup above, take the stepsize and burn-in period as
    \begin{subequations}
\begin{align}
\stepsize = \tfrac{(1 - \discount \lambda)^{2/3}}{c \beta \big(
  (\varsigma^4 + 1) d (1 - \kappa_\lambda) \numobs^2 \big( \mixingtime
  + \frac{1}{1 - \discount \lambda} \big) \big)^{1/3}}, \quad
\mbox{and} \quad \numburn = \tfrac{1}{2}
\numobs,\label{eq:td-lambda-stepsize-burnin-req}
\end{align}
and suppose that $\tfrac{\numobs}{\log^3 \numobs} \geq \tfrac{2 (
  \mixingtime + \tfrac{1}{1 - \discount \lambda} ) \: (\varsigma^4 d +
  1) \beta^2}{(1 - \kappa_\lambda)^2 (1 - \discount \lambda)^2
  \mu^2}$. Then the value function estimate
\mbox{$\valuehat_\numobs(\state) \mydefn
  \inprod{\thetahat_\numobs}{\phi(\state)}$} obtained from the
Polyak--Ruppert procedure~\eqref{eq:td-lambda} has MSE bounded as
\begin{multline}
    \Exs \big[ \vecnorm{\valuehat_\numobs -
        \valuebar^{(\lambda)}}{\Ltwospace (\statespace,
        \stationary)}^2 \big] \leq c \numobs^{-1} \mathrm{Tr} \big(
    (I_d - \SpecMat_\lambda)^{-1} (\SigStar_{\mathrm{Mkv}} +
    \SigStar_{\mathrm{MG}}) (I_d - \SpecMat_\lambda)^{- \top} \big)
    \\ + c \big( \tfrac{\beta^2 \varsigma^4 d \big(\mixingtime +
      \frac{1}{1 - \discount \lambda} \big)}{\mu^2 (1 -
      \kappa_\lambda)^2 (1 - \discount \lambda)^2 \numobs} \big)^{4/3}
    \varbound^2 \log^2 \numobs,\label{eq:td-lambda-final-bound}
\end{multline}
where $\valuebar^{(\lambda)}$ is the solution to the projected
fixed-point equation~\eqref{eq:projected-fixed-point}.
\end{subequations}
\end{corollary}
\noindent See
Appendix~\ref{subsec:proof-cor-td-lambda-finite-state-markov} for the
proof of this corollary.\\

A few remarks are in order. First, using the same argument as in
Corollaries~\ref{prop:td-0-continuous-statespace}
and~\ref{prop:sieve-td-results}, one can extend the results for
TD$(\lambda)$ to the cases of continuous state spaces with Wasserstein
mixing, as well as to nonparametric sieve estimators. As is
well-known, different choices of the tuning parameter $\lambda$
interpolate the ``temporal difference'' method, in which we aim at
solving the Bellman equation, and the ``Monte Carlo'' method, in which
the value function is estimated directly by averaging the rollout of a
Markovian trajectory. For example, on the one hand, letting $\lambda =
0$ recovers the instance-dependent upper bound for TD$(0)$ method in
Corollary~\ref{cor:TD-0-finite-state-markov}. On the other hand, by
taking $\lambda = \discount$, we have $\kappa_\lambda \leq
\frac{\discount}{1 + \discount} \leq \frac{1}{2}$, and the dependence
on the discount factor $\discount$ appears only through the variance
of the noise, instead of through the conditioning of the matrix
$\SpecMat_\lambda$. In the next section, we sketch a recipe for the
instance-dependent selection of $\lambda$ that also takes the
approximation error into account.


\subsubsection{Using instance-dependent results to select $\lambda$}

Recall that the TD$(\lambda)$ algorithm aims at estimating the
solution $\valuebar^{(\lambda)}$ to the projected fixed-point
equation~\eqref{eq:population-td-lambda}. The linear operator in the
unprojected fixed-point
equation~\eqref{eq:population-td-lambda-unprojected} satisfies the
norm bound
\begin{align*}
\matsnorm{(1 - \lambda) \sum_{k = 0}^{\infty} \lambda^k \discount^{k +
    1} \transition^{k + 1}}{\Ltwospace (\statespace, \stationary)
  \rightarrow \Ltwospace (\statespace, \stationary)} \leq (1 -
\lambda) \sum_{k = 0}^{\infty} \lambda^k \discount^{k + 1} = \tfrac{(1
  - \lambda) \discount}{1 - \lambda \discount}.
\end{align*}
Consequently, invoking Theorem 1 in the paper~\cite{mou2020optimal},
the approximation error satisfies the bound
\begin{align*}
    \vecnorm{\valuebar^{(\lambda)} - \valuestar}{\Ltwospace
      (\statespace, \stationary)}^2 \leq \alpha \big(
    \SpecMat_\lambda, \tfrac{(1 - \lambda) \discount}{1 - \lambda
      \discount} \big) \cdot \inf_{\valuefunc \in \LinSpace}
    \vecnorm{\valuefunc - \valuestar}{\Ltwospace (\statespace,
      \stationary)}^2,
\end{align*}
where $\alpha (\SpecMat, z) \mydefn 1 + \lammax \big( (I_d -
\SpecMat)^{-1} (z^2 I_d - \SpecMat \SpecMat^\top) (I_d - \SpecMat)^{-
  \top} \big)$ is the approximation factor.  Combining with
Corollary~\ref{cor:TD-lambda-finite-state-markov}, we obtain the
following bound on the distance to the true value function:
\begin{multline}
\label{eq:td-lambda-complete-bound}
  \Exs \big[ \vecnorm{\valuehat_\numobs - \valuestar}{\Ltwospace
      (\statespace, \stationary)}^2 \big] \\
      \leq c \alpha \big(
  \SpecMat_\lambda, \tfrac{(1 - \lambda) \discount}{1 - \lambda
    \discount} \big) \cdot \inf_{\valuefunc \in \LinSpace}
  \vecnorm{\valuefunc - \valuestar}{\Ltwospace (\statespace,
    \stationary)}^2 + \tfrac{c}{\numobs} \mathrm{Tr} \big(
  (I_d - \SpecMat_\lambda)^{-1} (\SigStar_{\mathrm{Mkv}} +
  \SigStar_{\mathrm{MG}}) (I_d - \SpecMat_\lambda)^{- \top} \big)\\ + c \big( \tfrac{\beta^2 \varsigma^4 d
    \big(\mixingtime + \frac{1}{1 - \discount \lambda} \big)}{\mu^2 (1
    - \kappa_\lambda)^2 (1 - \discount \lambda)^2 \numobs} \big)^{4/3}
  \varbound^2 \log^2 \numobs
\end{multline}
for a universal constant $c > 0$.

It can be seen that $\alpha \big( \SpecMat_\lambda, \tfrac{(1 -
  \lambda) \discount}{1 - \lambda \discount} \big) \leq c' \tfrac{1 -
  \lambda \discount}{1 - \discount} $ for a universal constant.  We
also recall that $\kappa_\lambda \leq \tfrac{(1 - \lambda)
  \discount}{1 - \lambda \discount}$. If we take the parameters $(\mu,
\beta, \varsigma)$ to be of constant order, in the worst case, the
upper bound~\eqref{eq:td-lambda-complete-bound} takes the simplified
form~\footnote{Here we use a coarse bound by taking $\vecnorm{\thetabar}{2} = O ( \tfrac{1}{1 - \discount})$. Note that one could potentially improve the worst-case dependence on the effective horizon by a careful upper bound on the variance of $\varepsilon_{\mathrm{Mkv}, \lambda}$. We omit the details for simplicity.}
\begin{align*}
\Exs \big[ \vecnorm{\valuehat_\numobs - \valuestar}{\Ltwospace
    (\statespace, \stationary)}^2 \big] \leq c \frac{1 - \lambda
  \discount}{1 - \discount} \inf_{\valuefunc \in \LinSpace}
\vecnorm{\valuefunc - \valuestar}{\Ltwospace (\statespace,
  \stationary)}^2 + c \frac{\big( \mixingtime + \frac{1}{1 - \discount
    \lambda} \big)\usedim}{(1 - \discount)^4\numobs}.
\end{align*}
From such an upper bound, it may appear that the optimal choice of
$\lambda$ is always \sloppy \mbox{$\lambda = \discount \wedge (1 - 1 /
  \mixingtime)$}, so that the approximation factor is minimized and
the variance remains controlled. However, this choice could be overly
conservative, since the actual variance with small $\lambda$ can be
significantly smaller, with the feature vectors still having bounded
one-step cross-correlation. Choosing the parameter $\lambda$ close to
$1$ cannot take advantage of small one-step correlation. On the other
hand, a fine-grained bound of the
form~\eqref{eq:td-lambda-complete-bound} can be used to perform
instance-dependent model selection, as follows:

\bcar
\item Construct a uniform finite grid $0 = \lambda_1 < \lambda_2 <
  \cdots < \lambda_{m} = \discount$ for possible values of $\lambda$.
\item For each $\ell \in [m]$, compute the TD$(\lambda_\ell)$
  estimator, and construct empirical plug-in estimates $\big(
  \widehat{\SpecMat}_{\lambda, \numobs},
  \widehat{\SigStar}_{\mathrm{Mkv}, \lambda, \numobs},
  \widehat{\SigStar}_{\mathrm{MG}, \lambda, \numobs} \big)$ for the
  matrices $\big( \SpecMat_{\lambda}, \SigStar_{\mathrm{Mkv},
    \lambda}, \SigStar_{\mathrm{MG}, \lambda} \big)$ by replacing the
  expectations by empirical averages. Similarly replace
  $\thetabar^{(\lambda)}$ by $\thetahat_\numobs$.
\item Estimate the approximation factor $\alpha \big(
  \SpecMat_\lambda, \frac{(1 - \lambda) \discount}{1 - \lambda
    \discount} \big)$ and the covariance $(I_d -
  \SpecMat_\lambda)^{-1} (\SigStar_{\mathrm{Mkv}} +
  \SigStar_{\mathrm{MG}}) (I_d - \SpecMat_\lambda)^{- \top}$ by
  plugging in the estimated matrices described above, for each
  $\lambda = \lambda_\ell$ with $\ell \in [m]$. Based on prior
  knowledge about the scale of the optimal approximation error
  $\inf_{\valuefunc \in \LinSpace} \vecnorm{\valuefunc -
    \valuestar}{\Ltwospace (\statespace, \stationary)}^2$, select
  $\lambda_\ell$ in the grid that minimizes our estimate of the total
  error according to equation~\eqref{eq:td-lambda-complete-bound}.
  \ecar

Note that the procedure above is simply a sketch; a formal proof of
correctness would show bounds that are uniform over all $m$
estimators.  It is an important direction of future work to provide
sharp non-asymptotic analysis of such a model selection procedure.

\subsection{Autoregressive models}

Next, we turn to Example~\ref{example:autoregressive}, the
multivariate auto-regressive model.  We study the stochastic
approximation procedure in which, for any $i \in [k]$, we have
\begin{align*}
    A^{(i)}_{t + 1} = A^{(i)}_t - \stepsize \big( \sum_{j = 0}^{k - 1}
    A^{(j)}_t X_{t - j} X_{t + 1 - i}^\top - X_{t + 1} X_{t + 1-
      i}^\top \big), \quad \mbox{and} \quad \widehat{A}^{(i)}_\numobs
    &= \frac{1}{\numobs - \numburn} \sum_{t = \numburn}^{\numobs - 1}
    A^{(i)}_t.
\end{align*}
The first step in our analysis is to establish necessary and
sufficient conditions for the existence and uniqueness of the
stationary distribution of the
process~\eqref{eq:autoregressive-process}.  The following $k m \times
k m$ matrix plays a crucial role in this context:
\begin{align*}
   \Rstar = \begin{bmatrix} \AstarInd{1} & \AstarInd{2} & \cdots & &
     \AstarInd{k} \\ I_m & 0 &\cdots && 0\\ 0 & I_m & 0&\cdots & 0\\ 0
     && \ddots && 0\\ 0 & \cdots & 0& I_m &0
    \end{bmatrix}.
\end{align*}
In the noiseless case, the stability of the linear dynamical system is
equivalent to the following \emph{Lyapunov stability condition} (see
e.g.~\cite{nemirovski2001lectures}, Section 3.3):
\begin{align}
    \exists \Pstar \succ 0, \Qstar \succ 0, \quad \mbox{such that }
    \Rstar^\top \Pstar \Rstar = \Pstar -
    \Qstar. \label{eq:lyapunov-stability}
\end{align}
Clearly we have $\Pstar \succ \Qstar$.  We let $\beta \mydefn
\lammax (\Pstar)$ and $\mu \mydefn \lambda_{\min} (\Qstar)$. Based on
stability theory for discrete-time linear
systems~\cite{brockwell2009time},
condition~\eqref{eq:lyapunov-stability} is necessary for the stationary
distribution to exist. In the following proposition, we show that this
condition is also sufficient, with a concrete mixing time bound.
\begin{proposition}
\label{prop:lyap-stable}
Under the Lyapunov stability condition~\eqref{eq:lyapunov-stability}
and assuming that the noise has bounded first moment $\Exs \big[
  \vecnorm{\varepsilon_t}{2} \big] < \infty$, the stationary
distribution $\widetilde{\stationary}$ for the sliding window
$\slidingwindow_t = (X_{t + 1}, X_t, \cdots, X_{t - k + 1})$ of the
auto-regressive process~\eqref{eq:autoregressive-process} exists and
is unique. Furthermore, the mixing
assumption~\ref{assume-markov-mixing} is satisfied with Wasserstein
distance in $\real^{(k + 1)m}$ and a mixing time bound $\mixingtime =
ck + c \frac{\beta}{\mu} \big(1 + \log \frac{\beta}{\mu} \big)$.
\end{proposition}
\noindent See Section~\ref{subsec:proof-prop-lyap-stable} for the
proof of this claim.

In addition to this mixing guarantee, we also make the following
assumptions on the noise:
\begin{align}
\Exs \big[ \varepsilon_t \big] = 0, \quad \sup_{u \in \sphere^{d - 1}}
\Exs \big[ \inprod{u}{\varepsilon_t}^4 \big] \leq \varsigma^4, \quad
\mbox{and}\quad \vecnorm{\varepsilon_t}{2} \leq \varsigma \sqrt{m},
~\mathrm{a.s.}
\end{align}
We are now in a position to consider the problem of parameter
estimation using stochastic approximation. Consider the vectorized
version of the parameter $\theta = \vectorize \big( \begin{bmatrix}
  A^{(1)}; A^{(2)}; \cdots; A^{(k)} \end{bmatrix} \big) \in \real^{k
  m^2}$. The population-level Yule--Walker estimation
equation~\eqref{eq:yule-walker-pop} can be written as
\begin{align}
  \big( \underbrace{\begin{bmatrix} \Gamma_{j - i} \end{bmatrix}_{i, j
      \in [k]}}_{H^*} \otimes I_m \big) \theta = \vectorize \big(
  \big[\Gamma_1; \Gamma_2; \cdots; \Gamma_k \big]
  \big),\label{eq:yule-walker-vectorized}
\end{align}
where $\Gamma_i \mydefn \Exs \big[ X_i X_0^\top \big] \in \real^{m
  \times m}$, for $i \in \mathbb{Z}$. We assume that
\begin{align*}
    \frac{1}{2} \big(H^* + (H^*)^\top \big) \succeq h^* I_{km}, \quad
    \mbox{for some $h^* > 0$}.
\end{align*}

In order to state the main corollary of Theorem~\ref{thm:markov-main}
to auto-regressive models, the following quantities are relevant:
\begin{align*}
\varepsilon_{\mathrm{Mkv}} (\slidingwindow_t) &\mydefn \vectorize
\big( \big( \sum_{j = 0}^{k - 1} A^{(j)}_* X_{t - j} - X_{t + 1} \big)
\cdot \begin{bmatrix} X_{t - 1}^\top & X_{t - 2}^\top & \cdots X_{t -
    k}^\top \end{bmatrix} \big)\\ \SigStar_{\mathrm{Mkv}} &\mydefn
\sum_{t = - \infty}^{\infty} \Exs \big[ \varepsilon_{\mathrm{Mkv}}
  (\slidingwindow_t) \varepsilon_{\mathrm{Mkv}}
  (\slidingwindow_0)^\top \big].
\end{align*}

Let $\varbound$ be defined according to
equation~\eqref{eq:defn-effective-noise}. We have the following
corollary for autoregressive models.
\begin{subequations}
\begin{corollary}\label{cor:auto-regressive}
Under the setup above, take the stepsize and burn-in period as
\begin{align}
  \stepsize = \tfrac{1}{c \big( n^2 \big( \frac{\beta}{\mu} \log
    \frac{\beta}{\mu} \big) (h^*)^2 \varsigma^4 k^3 m^2 \beta^8 /
    \mu^8 \big)^{1/3}}, \quad \mbox{and} \quad \numburn = \tfrac{1}{2}
  \numobs,
\end{align}
and suppose that $\frac{\numobs}{\log^3 \numobs} \geq \big( k +
\frac{\beta}{\mu} \log \frac{\beta}{\mu} \big) \varsigma^4 k^3 m^2
\frac{\beta^8}{\mu^8 (h^*)^2}$. Then the Polyak--Ruppert estimator
$(\widehat{A}^{(j)}_n)_{j \in [k]}$ satisfies
\begin{multline}
\label{eq:auto-regressive-final-bound}  
  \sum_{j = 1}^k \Exs \big[ \matsnorm{\widehat{A}_n^{(j)} -
      \AstarInd{j}}{F}^2 \big] \leq \tfrac{c}{n} \mathrm{Tr} \big(
  \big( H^* \otimes I_m \big)^{-1} \mathrm{\Sigma}_{\mathrm{Mkv}}
  \big( H^* \otimes I_m \big)^{-1} \big) \\ + \Big\{ \tfrac{k
    \vardim^2 \varsigma^2 \cdot \lammax \big( \Exs \big[
      \varepsilon_{\mathrm{Mkv}} (s_0) \varepsilon_{\mathrm{Mkv}}
      (s_0)^\top \big] \big)}{(h^*)^2 n} \big( k + \frac{\beta}{\mu}
  \log \frac{\beta}{\mu} \big) \Big\}^{4/3} \varbound^2 \log^2 n.
    \end{multline}
\end{corollary}
\end{subequations}
A few remarks are in order. First, the leading-order term in the
bound~\eqref{eq:auto-regressive-final-bound} matches the variance of
asymptotic efficient estimators for AR$(m)$ models, up to a constant
factor (see~\cite{brockwell2009time}, Section 8). This simply follows
from the fact that the plug-in Yule-Walker estimator is asymptotically
efficient for auto-regressive models. On the other hand,
Corollary~\ref{cor:auto-regressive} is completely non-asymptotic,
holding true for any reasonably large sample size.  Note that the
sample complexity lower bound exhibits an $\order{\beta^9 / \mu^9}$
dependency on the conditioning $\beta / \mu$ of the Lyapunov stability
certificate $(\Pstar, \Qstar)$.  The contributions arise from a term
linear in $\beta / \mu$ arises from the mixing time $\frac{\beta}{\mu}
\log \frac{\beta}{\mu}$, and all other factors are from the
almost-sure bounds on $\vecnorm{X_t}{2}$ and moment bound $\sup_{u \in
  \sphere^{m - 1}} \inprod{u}{ X_t}^4$.
If we had made other assumptions on these uniform or moment bounds as
in some past work~\cite{jain2021streaming}, these would have reflected
in our result instead of the factor $\beta^8 \varsigma^4 k^2 / \mu^8$.


\section{Proof of Proposition~\ref{prop:lsa-markov-iterate-bound}}
\label{subsec:proof-prop-lsa-iterate}

This section is devoted to proving the bound on the last iterate
claimed in Proposition~\ref{prop:lsa-markov-iterate-bound}.  We begin
in Section~\ref{SecInitialRecursion} by deriving a key recursion that
underlies the analysis.  In Section~\ref{SecIntuition}, we provide a
high-level overview of the proof structure, and the remaining
subsections deal with the technical arguments.

\subsection{An initial recursion}
\label{SecInitialRecursion}

Define the error term $\Delta_t \mydefn \theta_t - \thetabar$, as well
as the noise terms
\begin{subequations}
\begin{align}
\multnoiseMG_{t + 1} \mydefn \Lmat_{t + 1} - \Lmap (\state_t), & \quad
\qquad \addnoiseMG_{t + 1} \mydefn (\Lmat_{t + 1} - \Lmap (\state_t))
\thetabar + (\bvec_{t + 1} - \bmap
(\state_t)), \label{eq:defs-mg-noise} \\
\multnoiseMarkov_{t} \mydefn \Lmap (\state_t) - \Lbar, &\quad \qquad
\addnoiseMarkov_t \mydefn (\Lmap (\state_t) - \Lbar) \thetabar +
(\bmap (\state_t) - \bbar).\label{eq:defs-markov-noise}
\end{align}
\end{subequations}
Using this notation, we have the recursion
\begin{align}
\label{eq:recursive-relation-delta}  
\Delta_{t + 1} = (I - \stepsize (I - \Lbar)) \Delta_t + \stepsize
\big( \multnoiseMarkov_t + \multnoiseMG_{t + 1} \big) \Delta_t +
\stepsize (\addnoiseMarkov_t + \addnoiseMG_{t + 1}).
\end{align}
Taking squared norms on both sides yields the bound
$\vecnorm{\Delta_{t + 1}}{2}^2 \leq \sum_{i=1}^4 T_i$, where
\begin{align*}
  T_1 &\mydefn \vecnorm{(I - \stepsize (I - \Lbar)) \Delta_t}{2}^2,\\
T_2 &\mydefn 2 \stepsize \inprod{(I - \stepsize (I - \Lbar)) \Delta_t}{
  \multnoiseMarkov_t \Delta_t + \addnoiseMarkov_t},\\
   T_3 &\mydefn 2 \stepsize \inprod{(I - \stepsize (I
    - \Lbar)) \Delta_t}{\big( \multnoiseMG_{t + 1} \Delta_t +
    \addnoiseMG_{t + 1} \big)}, \\
T_4 &\mydefn 4 \stepsize^2 \big( \vecnorm{\multnoiseMarkov_t
  \Delta_t}{2}^2 + \vecnorm{\multnoiseMG_{t + 1} \Delta_t}{2}^2 +
\vecnorm{\addnoiseMG_{t + 1}}{2}^2 + \vecnorm{\addnoiseMarkov_t}{2}^2
\big).
\end{align*}

Beginning with the term $T_1$, expanding the square and then invoking
the condition~\eqref{eq:kappa-opnorm} yields
\begin{align*}
\Term_1 &= \statnorm{\Delta_t}^2 - 2 \stepsize \statinprod{\Delta_t}{(I
  - \Lbar) \Delta_t} + \stepsize^2 \statnorm{(I - \Lbar) \Delta_t}^2 \\
  &
\leq \big(1 - 2 \stepsize (1 - \kappa) + 2 \stepsize^2 (1 + \conmax^2)
\big) \statnorm{\Delta_t}^2.
\end{align*}
As for the cross terms involved in $T_2$ and $T_3$, we note that
\begin{align*}
2 \inprod{(I - \Lbar) \Delta_t}{ \multnoiseMarkov_t \Delta_t } &\leq
\vecnorm{(I - \Lbar) \Delta_t}{2}^2 + \vecnorm{ \multnoiseMarkov_t
  \Delta_t}{2}^2 \; \leq \; 2 (1 + \conmax^2) \vecnorm{\Delta_t}{2}^2
+ \vecnorm{ \multnoiseMarkov_t \Delta_t}{2}^2,\\
2 \inprod{(I - \Lbar) \Delta_t}{ \addnoiseMarkov_t} &\leq \vecnorm{(I
  - \Lbar) \Delta_t}{2}^2 + \vecnorm{ \addnoiseMarkov_t}{2}^2 \; \leq
\; 2 (1 + \conmax^2) \vecnorm{\Delta_t}{2}^2 +
\vecnorm{\addnoiseMarkov_t}{2}^2, \\
2 \inprod{(I - \Lbar) \Delta_t}{ \multnoiseMG_{t + 1} \Delta_t } &\leq
\vecnorm{(I - \Lbar) \Delta_t}{2}^2 + \vecnorm{ \multnoiseMG_{t + 1}
  \Delta_t}{2}^2 \; \leq \; 2 (1 + \conmax^2) \vecnorm{\Delta_t}{2}^2
+ \vecnorm{ \multnoiseMG_{t + 1} \Delta_t}{2}^2, \\
2 \inprod{(I - \Lbar) \Delta_t}{ \addnoiseMG_{t + 1}} & \leq
\vecnorm{(I - \Lbar) \Delta_t}{2}^2 + \vecnorm{ \addnoiseMG_{t + 1}
}{2}^2 \; \leq \; 2 (1 + \conmax^2) \vecnorm{\Delta_t}{2}^2 +
\vecnorm{\addnoiseMG_{t + 1} }{2}^2.
\end{align*}

We collect the above bounds on the sum $\sum_{i=1}^4 T_i$ and use the
stepsize bound $\stepsize \leq \frac{1 - \kappa}{12 (1 + \conmax^2)}$,
which results in the recursive inequality
\begin{multline} 
\vecnorm{\Delta_{t + 1}}{2}^2  \leq \big(1 - \stepsize (1 - \kappa)
\big) \vecnorm{\Delta_t}{2}^2 + 2 \stepsize \underbrace{\big(
  \inprod{\Delta_t}{ \multnoiseMarkov_t \Delta_t } +
  \inprod{\Delta_t}{\addnoiseMarkov_t} \big)}_{\mydefn H_1 (t)}\\  \quad \quad + 2 \stepsize\underbrace{ \big(
  \inprod{\Delta_t}{ \multnoiseMG_{t + 1} \Delta_t } +
  \inprod{\Delta_t}{\addnoiseMG_{t + 1} } \big)}_{\mydefn H_2 (t)} + 8
\stepsize^2 \underbrace{\big( \vecnorm{\multnoiseMarkov_t
    \Delta_t}{2}^2 + \vecnorm{\multnoiseMG_{t + 1} \Delta_t}{2}^2 +
  \vecnorm{\addnoiseMG_{t + 1}}{2}^2 +
  \vecnorm{\addnoiseMarkov_t}{2}^2 \big)}_{\mydefn H_3 (t)}.\label{eq:one-step-recursion-error}
\end{multline}
Multiplying both sides by $e^{\stepsize (1 - \kappa) (t+1)}$ and using
the fact that $\big(1 - \stepsize (1 - \kappa) \big) \leq e^{-
  \stepsize (1- \kappa)}$, we have
\begin{multline*}
  e^{\stepsize (1 - \kappa) (t + 1)} \vecnorm{\Delta_{t + 1}}{2}^2
  \leq e^{\stepsize (1 - \kappa) t} \vecnorm{\Delta_{t}}{2}^2 \\+ 2
  \stepsize e^{\stepsize (1 - \kappa) (t + 1)} \big( H_1 (t) + H_2 (t)
  \big) + 8 \stepsize^2 e^{\stepsize (1 - \kappa) (t + 1)} H_3 (t).
\end{multline*}
Unrolling this expression yields
\begin{multline}
e^{\stepsize (1 - \kappa) \numobs} \vecnorm{\Delta_{\numobs}}{2}^2
\leq \vecnorm{\Delta_0}{2}^2 + 2 \stepsize \sum_{t = 0}^{\numobs - 1}
e^{\stepsize (1 - \kappa) (t + 1)} \big( H_1 (t) + H_2 (t) \big) \\
+ 8
\stepsize^2 \sum_{t = 0}^{\numobs - 1} e^{\stepsize (1 - \kappa) (t +
  1)} H_3 (t),
\label{eq:error-expansion-markov}
\end{multline}
which is the key recursion underlying our analysis.

\subsection{High-level overview of the proof strategy}
\label{SecIntuition}

Before diving into the remainder of the proof, let us provide a brief
overview of our strategy, highlighting the key technical challenges
and our solutions to them.

For simplicity, let us give intuition for the analysis under
mean-squared error. In order to analyze the recursive error
expansion~\eqref{eq:one-step-recursion-error}, we need to bound the
terms $\Exs [H_1 (t)]$, $\Exs [H_2 (t)]$, and $\Exs [H_3 (t)]$,
respectively. For the martingale noise part, we have $\Exs [H_2 (t)] =
0$. As for the term $H_3 (t)$, following
Assumption~\ref{assume:noise-moments}, we have that
\begin{align*}
    \Exs [\vecnorm{\addnoiseMarkov_t}{2}^2 + \vecnorm{\addnoiseMG_{t +
          1}}{2}^2] \lesssim d \quad \mbox{and} \quad \Exs
         [\vecnorm{\multnoiseMG_{t + 1} \Delta_t}{2}^2] \lesssim
         \usedim \cdot \Exs [\vecnorm{\Delta_t}{2}^2].
\end{align*}
These bounds are similar to the analysis under $\mathrm{i.i.d.}$ setup
(see the paper~\cite{mou2020optimal}). If other terms were not
present, we could unroll this recursion and obtain a last-iterate
error bound of $O (\stepsize \usedim)$, as long as $\stepsize \ll
\usedim^{-1}$. The technical challenges arise, however, with the
interaction between Markovian noises and the error $\Delta_t$. In
particular, we observe the following facts:
\begin{itemize}
    \item Since $\Delta_t$ and $(\addnoiseMarkov_t,
      \multnoiseMarkov_t)$ are inter-dependent, the term $H_1 (t)$
      does not have zero expectation. If we simply bound it using
      Assumption~\ref{assume-lip-mapping}, for any stepsize $\stepsize
      > 0$, the error recursion will diverge as $t$ grows.
    \item Assumption~\ref{assume-lip-mapping} implies that $\Exs
      [\vecnorm{\multnoiseMarkov_t \Delta_t}{2}^2] \lesssim \usedim^2
      \cdot \Exs [\vecnorm{\Delta_t}{2}^2]$. In order to unroll
      recursion using this bound and obtain convergent result, we need
      the stepsize $\stepsize \lesssim \usedim^{-2}$. This will lead
      to a sub-optimal sample complexity, since we need at least
      $\numobs \gtrsim \stepsize^{-1}$ steps.
\end{itemize}
In tackling the aforementioned difficulties, our first proof technique
makes use of the rapid mixing nature of the underlying Markov
chain---the Markov chain state after $O (\mixingtime)$ steps is nearly
independent with the current iterate. We elaborate on the key ideas as
follows.

\paragraph{Multi-step looking-back for the cross terms:}
Let $\tau \asymp \mixingtime \log (d / \stepsize)$. The dependence
between $\Delta_{t - \tau}$ and $\state_t$ is weak, and consequently,
we can show that
\begin{align*}
    \big| \Exs \big[ \inprod{\multnoiseMarkov_t \Delta_{t -
          \tau}}{\Delta_{t - \tau}} \big] \big| &\lesssim \stepsize
    \usedim \cdot \Exs [\vecnorm{\Delta_{t - \tau}}{2}^2],\\ \big|
    \Exs \big[ \inprod{\addnoiseMarkov_t}{\Delta_{t - \tau}} \big]
    \big| &\lesssim \stepsize \usedim + \stepsize \usedim \cdot \Exs
          [\vecnorm{\Delta_{t - \tau}}{2}^2], \quad \mbox{and}\\ \Exs
          \big[ \vecnorm{\multnoiseMarkov_t \Delta_{t - \tau}}{2}^2
            \big] &\lesssim \usedim \cdot \Exs [\vecnorm{\Delta_{t -
                \tau}}{2}^2].
\end{align*}
In showing these bounds, we construct an auxiliary process
$(\widetilde{\state}_{t - \tau + \ell})_{ \ell \geq 0}$, which starts
from $\widetilde\state_{t - \tau} \sim \stationary$ independent with
the data, and moves according to the optimal coupling that achieves
the Wasserstein mixing. With the value $\tau$ given above, we can
ensure that $\Wass_{1, \metric} (\state_t, \widetilde{\state}_t)
\lesssim \stepsize / \usedim$. We can then apply bounds under
independent $\widetilde{\state}_t$ and $\Delta_{t - \tau}$, and bound
the residual using Wasserstein distance and the Lipschitz
assumption~\ref{assume-lip-mapping}. See the proof of
Lemma~\ref{LemNewH1} for details.

However, this does not complete the analysis, as we originally need to
bound the cross terms between $(\addnoiseMarkov_t,
\multnoiseMarkov_t)$ and $\Delta_t$, instead of the $\tau$-step
looking-back version $\Delta_{t - \tau}$. In order to convert above
estimates to useful bound for analyzing the
recursion~\eqref{eq:recursion-for-mse-bound}, we need a stability
estimate, i.e., an upper bound on $\Exs \big[ \vecnorm{\theta_t -
    \theta_{t - \tau}}{2}^2 \big]$. This is the major technical
challenge we face in order to obtain the sharp dimension
dependence. In tackling this challenge, we introduce a novel
bootstrapping argument, which may be of independent interest.

\paragraph{Bootstrapping arguments for stability bounds:}
Expanding the recursion~\eqref{eq:lsa} yields
\begin{align*}
 \theta_{t + \tau} - \theta_{t} = \stepsize \cdot \sum_{\ell =
   1}^{\tau} \Big\{ (\Lmat_{t + \ell} - I_d) \theta_{t + \ell - 1} +
 \bvec_{t + \ell} \Big\}.
\end{align*}
If we use triangle inequality and
Assumptions~\ref{assume:noise-moments},~\ref{assume-lip-mapping} to
bound the difference, some calculations will lead to a coarse bound
(see Lemma~\ref{lemma:norm-not-blowup})
\begin{align}
 \sqrt{\Exs \big[ \vecnorm{\theta_{t + \tau} - \theta_t}{2}^2 \big]}
 \lesssim \stepsize \tau \usedim \sqrt{\Exs [\vecnorm{\Delta_t}{2}^2]}
 + \stepsize \tau
 \sqrt{\usedim}.\label{eq:coarse-bound-in-proof-overview}
\end{align}
However, if we directly substitute this bound into the arguments
above, we will need the stepsize $\stepsize$ to satisfy $\stepsize
\lesssim \usedim^{-2}$ in order to make the iterates stable. As we
have discussed above, this will cost us a sub-optimal sample
complexity of $\numobs \gtrsim \usedim^2$. In order to make the
arguments work with a larger stepsize, we need the pre-factor in the
first term of the right-hand-side of
equation~\eqref{eq:coarse-bound-in-proof-overview} to be scaling as $O
(\tau \stepsize \sqrt{\usedim})$. To achieve this goal, we start with
the bound~\eqref{eq:coarse-bound-in-proof-overview}, and gradually
improve it using a bootstrapping lemma. In
Lemma~\ref{lemma:local-moves} to follow, we shows a bootstrapping
result: as long as we have the bound
\begin{align*}
  \sqrt{\Exs \big[ \vecnorm{\theta_{t + \tau} - \theta_t}{2}^2 \big]} \lesssim \stepsize \tau \omega \sqrt{\Exs [\vecnorm{\Delta_t}{2}^2]} + \stepsize \tau \beta,
\end{align*}
we can establish the improved bound
\begin{align*}
  \sqrt{\Exs \big[ \vecnorm{\theta_{t + \tau} - \theta_t}{2}^2
            \big]} \lesssim \stepsize \tau \Big( \frac{\omega}{2} +
        \sqrt{\usedim} \Big) \sqrt{\Exs [\vecnorm{\Delta_t}{2}^2]} +
        \stepsize \tau \Big( \frac{\beta}{2} + \stepsize
        \sqrt{\usedim} \cdot \omega + \sqrt{\usedim} \Big).
\end{align*}
Once again, the proof of this lemma relies on the multi-step
looking-back arguments explained above: when analyzing the
iterate~\eqref{eq:lsa}, we can gain the near-independent by replacing
$\theta_t$ with $\theta_{t - \tau}$, at an additional cost depending
on the stability bound $\Exs [\vecnorm{\theta_t - \theta_{t -
      \tau}}{2}^2]$. By repeatedly applying this lemma, we obtain a
sequence of pairs $(\omega, \beta)$, which converges to the fixed
point
\begin{align*}
  \omega \asymp \sqrt{\usedim}, \quad \mbox{and} \quad \beta
        \asymp \sqrt{\usedim},
\end{align*}
which yield the desirable stability bound.

\paragraph{Completing the proof by solving the recursion:}
The improved stability bound allows us to establish sharp bounds on
the terms $\Exs [H_1 (t)]$ and $\Exs [H_3 (t)]$. These bounds involves
not only the current iterate error $\Delta_t$, but also the
looking-back iterate error $\Delta_{t - \tau}$. In order to analyze
this type of recursion, we multiply the inequality with an
exponentially growing factor $e^{\stepsize (1 - \kappa) t}$ and
telescope the summation. Solving it directly yields the MSE bound. As
for higher-order moments, we apply martingale concentration
inequalities to the martingale noise $H_2 (t)$ and the martingales
created from the auxiliary processes in analyzing $H_1 (t)$. The
recursive inequalities in this case can be solved using techniques
similar to our prior work~\cite{mou2020linear}.


\subsection{Analyzing the recursion~\eqref{eq:error-expansion-markov}}

Note that the running sum $M_2(\numobs) \mydefn \sum_{t = 0}^{\numobs
  - 1} e^{\stepsize (1 - \kappa) t} H_2 (t)$ is, by construction, a
martingale adapted to the filtration $(F_t)_{t \geq 0}$.  In contrast,
the analogous quantity defined in terms of the process $H_1$ is
\emph{not} an adapted martingale.  In order to circumvent this
obstacle, our proof is based on introducing a \emph{surrogate version}
$\Htil_1$ of the process $H_1$, such that the running sum
\begin{align*}
\widetilde{M}_1 (\numobs) \mydefn \sum_{t = 0}^{\numobs - 1}
e^{\stepsize (1 - \kappa) (t + \tau)} \Htil_1 (t + \tau)
\end{align*}
can be decomposed as a sum of $\tau$ martingales.  See the proof of
Lemma~\ref{LemNewH1} for the details of the construction of $\Htil_1$.
This decomposition allows us to apply standard maximal inequalities
for martingales.  Of course, we also need the bound the moments of the
differences $\Htil_1(t) - H_1(t)$; see Lemma~\ref{LemNewH1} for the
bound that we provide on this difference.

We prove the MSE bounds and higher-moment bounds using slightly
different analysis tools. In order to study the mean-squared error (the case
$p = 1$), we note that both $\widetilde{M}_1 (t)$ and $H_2 (t)$ have
zero expectation for any $t \geq 0$. Taking expectations on both sides
of equation~\eqref{eq:error-expansion-markov}, we obtain the bound
\begin{multline}
    e^{\stepsize (1 - \kappa) \numobs} \Exs
    \big[\vecnorm{\Delta_{\numobs}}{2}^2 \big] \leq
    \vecnorm{\Delta_0}{2}^2 + 2 \stepsize \sum_{t = 0}^{\numobs - 1}
    e^{\stepsize (1 - \kappa) (t + 1)} \Exs \big[ \abss{H_1 (t) -
        \Htil_1 (t)} \big] \\
        + 8 \stepsize^2 \sum_{t = 0}^{\numobs - 1}
    e^{\stepsize (1 - \kappa) (t + 1)} \Exs \big[ H_3 (t) \big].
\label{eq:last-iterate-proof-lyap-decomp-mse}  
\end{multline}

For higher moments, our analysis of the
recursion~\eqref{eq:error-expansion-markov} is based on a Lyapunov
function $\lyap_\numobs$ and auxiliary function $\Lambda_\numobs$
given by
\begin{align*}
  \lyap_\numobs \mydefn \big( \Exs \big[ \sup_{0 \leq t \leq \numobs}
    e^{\stepsize (1 - \kappa) t p} \vecnorm{\Delta_t}{2}^{2p} \big]
  \big)^{1/p}, \quad \mbox{and} \quad \Lambda_\numobs = \max
  \limits_{t \in \{0, 1, \ldots, \numobs\}} e^{- \frac{\stepsize (1 -
      \kappa) t}{2}} \lyap_t.
\end{align*}
By applying Minkowski's inequality to the
recursion~\eqref{eq:error-expansion-markov}, we obtain the upper bound
\begin{multline}
\label{eq:last-iterate-proof-lyap-decomp}  
  \lyap_\numobs \leq \lyap_0 + 4 \stepsize \big( \Exs \sup_{0 \leq t
    \leq \numobs} |\widetilde{M}_1 (t)|^p \big)^{1/p} + 4 \stepsize
  \big( \Exs \big(\sum_{t = 0}^{\numobs - 1} e^{\stepsize (1 - \kappa)
    t} |H_1 (t) - \widetilde{H}_1 (t)| \big)^{p}\big)^{1/p} \\
+ 4 \stepsize \big( \Exs \sup_{0 \leq t \leq \numobs} |M_2 (t)|^p
\big)^{1/p} + 16 \stepsize^2 \big( \Exs \big(\sum_{t = 0}^{\numobs -
  1} e^{\stepsize (1 - \kappa) t} H_3 (t) \big)^{p}\big)^{1/p}.
\end{multline}

In order to complete the proof, we need to control each of the terms
on the right-hand side.  The following auxiliary results provide the
needed control; in all cases, the quantities $(c, c_0)$ etc. denote
universal constants; the number $\numobs$ in the following lemmas is
seen as a general iteration index, instead of the total sample size in
the final statement of the theorem.

Our first auxiliary result guarantees the existence of the surrogate
variables $\Htil_1(t)$ with desirable properties:
\begin{lemma}
  \label{LemNewH1}
There is a surrogate version $\{\Htil_1(t) \}_{t \geq 0}$ of the
process $\{H_1(t) \}_{t \geq 0}$ such that $\Exs \big[ \Htil (t) \big]
= 0$ for any $t \geq 0$, and for any integer $p \in [1, \pmax / 2]$,
scalar $\tau \geq c p \mixingtime \log (c_0 \mixingtime d / \stepsize)$ and
stepsize $\stepsize \leq \frac{1}{c \mixingtime (\conmax + p \sigmaA
  d)}$, we have the following bounds for any $\numobs > 0$:
\begin{subequations}
  \begin{align}
\label{EqnDifferenceBound}    
        \big( \Exs \big[ \abss{H_1 (\numobs) - \widetilde{H}_1
            (\numobs)}^p \big]\big)^{1/p} & \leq c \stepsize p^2 \tau
        \big( ( d \sigmaA^2 + \conmax^2) \cdot \big( \Exs
        \vecnorm{\Delta_{\numobs - \tau \vee 0}}{2}^{2p}
        \big)^{\frac{1}{p}} + \varbound^2 d \big),
        \end{align}
        and for any $p \geq 2$, we have that
        \begin{align}
\label{EqnMartingaleBound}        
        \big( \Exs \sup_{0 \leq t \leq \numobs} |\widetilde{M}_1
        (t)|^p \big)^{1/p} & \leq \frac{c p^{3/2}}{\sqrt{\stepsize (1
            - \kappa)}} \big( \sigmaA \sqrt{d} \lyap_\numobs +
        \varbound \sqrt{e^{\stepsize (1 - \kappa)
            \numobs}\lyap_\numobs d} \big).
    \end{align}
\end{subequations}
\end{lemma}
\noindent See Section~\ref{SecProofNewH1} for the proof of this claim.
We note that it is especially challenging to prove the
bound~\eqref{EqnDifferenceBound}. \\

\noindent Our second auxiliary result is a more straightforward bound
on a martingale supremum:
\begin{lemma}
\label{LemNewH2}
The process $M_2$ is a martingale adapted to the filtration
$(\filtration_t)_{t \geq 0}$.  Furthermore, for each $p \in [1, \pmax
  / 2]$, $\tau \geq 2 p \mixingtime \log (c_0 d)$ and $\stepsize \leq
\frac{1}{c (\conmax + \sigmaA d) \tau}$, for any $\numobs > 0$, we
have that
\begin{align}
  \big( \Exs \sup_{0 \leq t \leq \numobs} |M_2 (t)|^p \big)^{1/p}
  &\leq \frac{c p^{3/2} \tau^{1/2}}{\sqrt{\stepsize (1 - \kappa)}}
  \big( \sigmaA \sqrt{d} \lyap_\numobs + \varbound \sqrt{e^{\stepsize
      (1 - \kappa) \numobs} \lyap_\numobs d} \big).
\end{align}
\end{lemma}
\noindent See Section~\ref{SecProofNewH2} for the proof of this
claim.\\

\noindent Finally, our third auxiliary result provides control on the
process $H_3(t)$:
\begin{lemma}
  \label{LemNewH3}
   There is a universal constant $c$ such that given $\tau \geq c p
   \mixingtime \log (c_0 \mixingtime d / \stepsize)$ and stepsize $\stepsize \leq
   \frac{1}{c \mixingtime (\conmax + \sigmaA d)}$, for any $p \in [1,
     \pmax / 2]$, we have
\begin{align}
\label{EqnAuxH3}  
  \big( \Exs \big[ H_3 (t)^p \big] \big)^{1/p} \leq c \big(p^2
  \sigmaA^2 d + \conmax^2 \big) \big( \Exs \big[ \vecnorm{\Delta_{t -
        \tau \vee 0}}{2}^{2p} \big]\big)^{1/p} + c p^2 \varbound^2 d.
\end{align}
\end{lemma}
\noindent See Section~\ref{SecProofNewH3} for the proof of this
claim.\\

We now use these three lemmas to complete the proof of
Proposition~\ref{prop:lsa-markov-iterate-bound}. We prove the case of
$\pmax = 2$ and $\pmax \geq \log \numobs$ separately.

\paragraph{Proof in the case of $\pmax = 2$:} By Lemma~\ref{LemNewH1} with $\tau = c \mixingtime \log (c_0 \mixingtime d / \stepsize)$ and Cauchy--Schwarz inequality, we have that
\begin{align*}
    &\Exs \big[\sum_{t = 0}^{\numobs - 1} e^{\stepsize (1 - \kappa) t}
      |\widetilde{H}_1 (t) - H_1 (t)| \big]\\
       &\leq c \stepsize \tau
    \sum_{t = 0}^{\numobs - 1} e^{\stepsize (1 - \kappa) t} \big(
    (\sigmaA^2 d + \conmax^2) \Exs \big[ \vecnorm{\Delta_{t - \tau
          \vee 0}}{2}^2 \big] + \varbound^2 d \big)\\ &\leq \frac{c
      \tau \varbound^2 d}{1 - \kappa} e^{ \stepsize (1 - \kappa)
      \numobs} + c e \stepsize \tau (\sigmaA^2 d + \conmax^2) \sum_{t
      = 0}^{\numobs - 1} e^{\stepsize (1 - \kappa) t} \Exs \big[
      \vecnorm{\Delta_t}{2}^2 \big].
\end{align*}

Similarly, by applying Lemma~\ref{LemNewH3} to the last term of
equation~\eqref{eq:last-iterate-proof-lyap-decomp-mse}, we obtain the
bound
\begin{align*}
    \sum_{t = 0}^{\numobs - 1} e^{\stepsize (1 - \kappa) (t + 1)} \Exs
    \big[ H_3 (t) \big] \leq \frac{c \varbound^2 d}{(1 - \kappa)
      \stepsize} e^{ \stepsize (1 - \kappa) \numobs} + c e (\sigmaA^2
    d + \conmax^2) \sum_{t = 0}^{\numobs - 1} e^{\stepsize (1 -
      \kappa) t} \Exs \big[ \vecnorm{\Delta_t}{2}^2 \big].
\end{align*}
Combining them with the
decomposition~\eqref{eq:last-iterate-proof-lyap-decomp-mse}, for any $n = 1,2, \cdots$, we find
that $e^{\stepsize (1 - \kappa) \numobs} \Exs
\big[\vecnorm{\Delta_{\numobs}}{2}^2 \big]$ is upper bounded by
\begin{align}
  \label{eq:recursion-for-mse-bound}  
  \vecnorm{\Delta_0}{2}^2 + c \frac{\stepsize \tau \varbound^2 d}{1 -
    \kappa} e^{ \stepsize (1 - \kappa) \numobs} + c \stepsize^2 \tau
  (\sigmaA^2 d + \conmax^2) \sum_{t = 0}^{\numobs - 1} e^{\stepsize (1
    - \kappa) t} \Exs \big[ \vecnorm{\Delta_t}{2}^2 \big].
\end{align}

In order to exploit this recursive upper bound, we define the partial
sum sequence \mbox{$S_\numobs \mydefn \sum_{t = 0}^{n} e^{\stepsize (1
    - \kappa) t} \Exs \big[ \vecnorm{\Delta_t}{2}^2 \big]$.}
Equation~\eqref{eq:last-iterate-proof-lyap-decomp-mse} implies that
\begin{align*}
S_\numobs & \leq S_0 + c \frac{\stepsize \tau \varbound^2 d}{1 -
  \kappa} e^{ \stepsize (1 - \kappa) \numobs} + \big(1 + c \stepsize^2
\tau (\sigmaA^2 d + \conmax^2) \big) S_{n - 1}\\ &\leq S_0 \cdot
\sum_{t = 0}^{\numobs} e^{c \stepsize^2 \tau (\sigmaA^2 d + \conmax^2)
  t} + c \frac{\stepsize \tau \varbound^2 d}{1 - \kappa} \cdot \sum_{t
  = 0}^{\numobs} e^{c \stepsize^2 \tau (\sigmaA^2 d + \conmax^2) t +
  \stepsize (1 - \kappa) (\numobs - t)} \\
& \leq \frac{3}{(1 - \kappa) \stepsize} e^{\stepsize (1 - \kappa)
  \numobs / 3} S_0 + \frac{3 c\tau \varbound^2 d}{(1 - \kappa)^2}
e^{\stepsize (1 - \kappa) \numobs}.
\end{align*}
Substituting back into the
recursion~\eqref{eq:recursion-for-mse-bound} yields
\begin{align*}
\Exs \big[\vecnorm{\Delta_{\numobs}}{2}^2 \big] &\leq \frac{6}{(1 -
  \kappa) \stepsize} e^{- \stepsize (1 - \kappa) \numobs / 3}
\vecnorm{\Delta_0}{2}^2 + c \frac{\stepsize \tau \varbound^2 d}{1 -
  \kappa} + c \stepsize^2 \tau (\sigmaA^2 d + \conmax^2) \cdot \frac{2
  c\tau \varbound^2 d}{(1 - \kappa)^2}\\ &\leq e^{- \stepsize (1 -
  \kappa) \numobs / 2} \vecnorm{\Delta_0}{2}^2 + c' \frac{\stepsize
  \tau \varbound^2 d}{1 - \kappa},
\end{align*}
which completes the proof of the MSE bound.

\paragraph{Proof in the case of $\pmax \geq \log \numobs$:}

Now we turn to prove the $p$-th moment bound under
Assumption~\ref{assume:noise-moments} with $\pmax \geq \log
\numobs$. Recall that we analyze the growth of the Lyapunov function
$\lyap_\numobs$, and we start from the
decomposition~\eqref{eq:last-iterate-proof-lyap-decomp}.

The first term in equation~\eqref{eq:last-iterate-proof-lyap-decomp}
is simply $\vecnorm{\Delta_0}{2}^2$, and the second term is controlled
using equation~\eqref{EqnMartingaleBound} in Lemma~\ref{LemNewH1}. In
order to bound the third term, we apply H\"{o}lder's inequality, and
obtain the bound
\begin{align*}
\Exs \big(\sum_{t = 0}^{\numobs - 1} e^{\stepsize (1 - \kappa) t}
\abss{H_1 (t) - \Htil_1 (t)} \big)^{p} & \leq \big( \sum_{t =
  0}^{\numobs - 1} e^{\frac{\stepsize (1 - \kappa) p t}{2 (p - 1)}}
\big)^{p - 1} \cdot \sum_{t = 0}^{\numobs - 1} e^{\frac{\stepsize p (1
    - \kappa )t}{2} } \Exs \abss{H_1 (t) - \Htil_1 (t)}^p .
\end{align*}
By equation~\eqref{EqnDifferenceBound} in Lemma~\ref{LemNewH1}, this
quantity is at most
\begin{align*}
\frac{e^{\frac{\stepsize (1 - \kappa) p
      \numobs}{2}} }{(\stepsize (1 - \kappa))^{1 - p}} \sum_{t = 0}^{\numobs - 1} e^{\frac{\stepsize p (1
      - \kappa) t}{2}} \big(c \tau \big(p^2 \sigmaA^2 d + \conmax^2
  \big) \big( \Exs \big[ \vecnorm{\Delta_{t - \tau \vee 0}}{2}^{2p}
    \big]\big)^{1/p} + c \tau p^2 \varbound^2 d \big)^p.
\end{align*}
We then obtain the inequality:
\begin{align*}
&\big( \Exs \big(\sum_{t = 0}^{\numobs - 1} e^{\stepsize (1 - \kappa)
    t} \abss{H_1 (t) - \Htil_1 (t)} \big)^{p} \big)^{1/p}\\ &\leq c
  p^2 \frac{e^{\stepsize (1 - \kappa) \numobs}}{\stepsize (1 -
    \kappa)} \varbound^2 \tau d + c \big(p^2 \sigmaA^2 d + \conmax^2
  \big) \tau \frac{e^{\frac{1}{2} \stepsize (1 - \kappa) \numobs}
  }{\stepsize (1 - \kappa)} \big( \sum_{t = 0}^{\numobs - 1}
  e^{\frac{1}{2}\stepsize p (1 - \kappa) t} \Exs \big[
    \vecnorm{\Delta_{t}}{2}^{2p} \big] \big)^{1/p}\\ &\leq c p^2
  \frac{e^{\stepsize (1 - \kappa) \numobs}}{\stepsize (1 - \kappa)}
  \varbound^2 \tau d + c \big(p^2 \sigmaA^2 d + \conmax^2 \big) \tau
  \frac{e^{\frac{1}{2} \stepsize (1 - \kappa) \numobs} }{\stepsize (1
    - \kappa)} \big( \sum_{t = 0}^{\numobs - 1} e^{-
    \frac{1}{2}\stepsize p (1 - \kappa) t} \lyap_t^p
  \big)^{1/p}\\ &\leq c p^2 \frac{e^{\stepsize (1 - \kappa)
      \numobs}}{\stepsize (1 - \kappa)} \varbound^2 \tau d + c
  \big(p^2 \sigmaA^2 d + \conmax^2 \big) \tau \frac{e^{\frac{1}{2}
      \stepsize (1 - \kappa) \numobs} }{\stepsize (1 - \kappa)}
  \numobs^{1/p} \Lambda_\numobs.
\end{align*}
Similarly, the fourth term on the right hand side is controlled using
Lemma~\ref{LemNewH2}, and the bounds for the last term are based on
Lemma~\ref{LemNewH3} and the same strategy as above. Concretely,
combining H\"{o}lder's inequality with the bound~\eqref{EqnAuxH3}
yields
\begin{align*}
\Exs \big(\sum_{t = 0}^{\numobs - 1} e^{\stepsize (1 - \kappa) t} H_3
(t) \big)^{p} & \leq \big( \sum_{t = 0}^{\numobs - 1}
e^{\frac{\stepsize (1 - \kappa) p t}{2 (p - 1)}} \big)^{p - 1} \cdot
\sum_{t = 0}^{\numobs - 1} e^{\frac{\stepsize p (1 - \kappa )t}{2} }
\Exs [H_3 (t)^p].
\end{align*}
This quantity is at most
\begin{align*}
(\stepsize (1 - \kappa))^{1 - p} e^{\frac{\stepsize (1 - \kappa) p
      \numobs}{2}} \sum_{t = 0}^{\numobs - 1} e^{\frac{\stepsize p (1
      - \kappa) t}{2}} \big(c \big(p^2 \sigmaA^2 d + \conmax^2 \big)
  \big( \Exs \big[ \vecnorm{\Delta_{t - \tau \vee 0}}{2}^{2p}
    \big]\big)^{1/p} + c p^2 \varbound^2 d \big)^p.
\end{align*}
Noting that each term satisfies the inequality $ e^{\frac{\stepsize p
    (1 - \kappa) t}{2}} \big( \Exs \big[ \vecnorm{\Delta_{t - \tau
      \vee 0}}{2}^{2p} \big]\big)^{1/p} \leq \Lambda_\numobs$ for $t
\in [0, \numobs]$.  We conclude that the moment $\big( \Exs
\big(\sum_{t = 0}^{\numobs - 1} e^{\stepsize (1 - \kappa) t} H_3 (t)
\big)^{p} \big)^{1/p}$ is upper bounded by
\begin{align*}
  c p^2 \frac{e^{\stepsize (1 - \kappa) \numobs}}{\stepsize (1 -
    \kappa)} \varbound^2 d + c \big(p^2 \sigmaA^2 d + \conmax^2 \big)
  \frac{e^{\frac{1}{2} \stepsize (1 - \kappa) \numobs} }{\stepsize (1
    - \kappa)} n^{1/p} \Lambda_\numobs.
\end{align*}
  
Collecting the above bounds and substituting into the
decomposition~\eqref{eq:last-iterate-proof-lyap-decomp}, we note that
\begin{align*}
    \lyap_\numobs &\leq \lyap_0 + c \sqrt{\frac{p^3 \stepsize}{1 -
        \kappa}} \big( \sigmaA \sqrt{d} \lyap_\numobs + \varbound
    \sqrt{e^{\stepsize (1 - \kappa) \numobs}\lyap_\numobs d} \big)
    \\ &\qquad + c p^2 \frac{e^{\stepsize (1 - \kappa)
        \numobs}}{\stepsize (1 - \kappa)} \varbound^2 \tau d +
    \big(p^2 \sigmaA^2 d + \conmax^2 \big) \frac{e^{\frac{1}{2}
        \stepsize (1 - \kappa) \numobs} }{\stepsize (1 - \kappa)} \tau
    \numobs^{1/p} \Lambda_\numobs\\ &\leq \lyap_0 + 4c \sigmaA
    \sqrt{\frac{p^3 \tau \stepsize d}{1 - \kappa}} \lyap_\numobs +
    \frac{1}{4} \lyap_\numobs + c \stepsize \frac{\varbound^2 p^3 d
      \tau}{1 - \kappa} \cdot e^{\stepsize (1 - \kappa) \numobs} \\
& \qquad + c p^2 \stepsize \frac{e^{\stepsize (1 - \kappa) \numobs}}{1
      - \kappa} \varbound^2 \tau d + c \stepsize \big(p^2 \sigmaA^2 d
    + \conmax^2 \big) \frac{e^{\frac{1}{2} \stepsize (1 - \kappa)
        \numobs} }{1 - \kappa} \tau \Lambda_\numobs
\end{align*}
In the last step, we apply Young's inequality to the term
$\sqrt{e^{\stepsize (1 - \kappa) \numobs}\lyap_\numobs d}$, and use
the condition $p \geq \log \numobs$ to the last term so that
$\numobs^{1/p} \leq e$.
 
Taking the stepsize $\stepsize \leq \frac{1 - \kappa}{64 c^2 \sigmaA^2
  \tau d p^3}$, we arrive at the following bound valid for any
$\numobs \in [1, e^p]$:
\begin{align*}
e^{- \frac{\stepsize (1 - \kappa) \numobs}{2}} \lyap_\numobs \leq 2
\lyap_0 + c p^3 \stepsize \frac{e^{\frac{1}{2} \stepsize (1 - \kappa)
    \numobs}}{1 - \kappa} \varbound^2 \tau d + c \stepsize \frac{ p^2
  \sigmaA^2 d + \conmax^2 }{1 - \kappa} \tau \Lambda_\numobs.
\end{align*}
Note that the right-hand-side of above expression is monotonic
increasing in the index $\numobs$. For any integer pair $(t, n)$ such
that $0 < t \leq \numobs \leq e^p$, we have the inequality:
\begin{align*}
    e^{- \frac{\stepsize (1 - \kappa) t}{2}} \lyap_t &\leq 2 \lyap_0 +
    c p^3 \stepsize \frac{e^{\frac{1}{2} \stepsize (1 - \kappa) t}}{1
      - \kappa} \varbound^2 \tau d + c \stepsize \frac{ p^2 \sigmaA^2
      d + \conmax^2 }{1 - \kappa} \tau \Lambda_t\\ &\leq 2 \lyap_0 + c
    p^3 \stepsize \frac{e^{\frac{1}{2} \stepsize (1 - \kappa)
        \numobs}}{1 - \kappa} \varbound^2 \tau d + c \stepsize \frac{
      p^2 \sigmaA^2 d + \conmax^2 }{1 - \kappa} \tau \Lambda_\numobs.
\end{align*}
Given the value of $\numobs$ fixed and taking supremum over $t \in
\{0, 1, 2, \cdots, \numobs\}$ in the left-hand-side, we arrive at the
conclusion:
\begin{align*}
   \Lambda_\numobs = \sup_{t \in \{0, 1, \cdots, \numobs\}} e^{-
     \frac{\stepsize (1 - \kappa) t}{2}} \lyap_t \leq 2 \lyap_0 + c
   p^3 \stepsize \frac{e^{\frac{1}{2} \stepsize (1 - \kappa)
       \numobs}}{1 - \kappa} \varbound^2 \tau d + c \stepsize \frac{
     p^2 \sigmaA^2 d + \conmax^2 }{1 - \kappa} \tau \Lambda_\numobs.
\end{align*}
 Given the stepsize $\stepsize \leq
\frac{1 - \kappa}{2 c ( p^3 \sigmaA^2 d + \conmax^2 ) \tau}$, we
arrive at the bound
\begin{align*}
  \big( \Exs \vecnorm{\Delta_t}{2}^p \big)^{1/p} \leq e^{- \frac{1}{2}
    \stepsize (1 - \kappa) \numobs} \Lambda_\numobs \leq e^{-
    \frac{1}{2} \stepsize (1 - \kappa) \numobs} \big( \Exs
  \vecnorm{\Delta_0}{2}^p \big)^{1/p} + \frac{c p^3 \stepsize}{1 -
    \kappa} \varbound^2 \tau d,
\end{align*}
which completes the proof of the theorem. \\

\vspace*{0.05in}

\noindent It remains to prove our three auxiliary lemmas.

\subsection{Proof of Lemma~\ref{LemNewH1}}
\label{SecProofNewH1}
We break the proof into three steps.  In the first step, given in
Section~\ref{SecConstruct}, we construct the surrogate process,
whereas the remaining two steps are devoted to the proving the
bounds~\eqref{EqnMartingaleBound} and~\eqref{EqnDifferenceBound}, as
detailed in Sections~\ref{SecMartingaleBound}
and~\ref{SecDifferenceBound} respectively.


\subsubsection{Construction of the surrogate process}
\label{SecConstruct}
We first claim that for any $t = 1, 2, \ldots$ and any $\tau \in \{0,
\ldots, t\}$, there is a random variable $\widetilde{\state}_t
\in \statespace$ such that $\widetilde{\state}_t \mid \filtration_{t -
  \tau} \sim \stationary$, and
\begin{align}
  \big( \Exs \big[ \metric (\state_t, \widetilde{\state}_t)^p \mid
    \filtration_{t - \tau} \big] \big)^{1/p} \leq c_0 \exp \big ( -
  \tfrac{\tau}{2 \mixingtime p} \big ) \quad \mbox{for each $p \geq
    2$.}
\end{align}
Here $c_0$ is a universal constant.

\noindent Our construction is based on the following bound on the
Wasserstein distance:
\begin{lemma}
\label{LemWassDecay}
Under Assumptions~\ref{assume-markov-mixing}
and~\ref{assume:stationary-tail}, the Wasserstein distance is upper
bounded as
\begin{align*}
\Wass_{1, \metric} \big(\delta_x \transition^{\tau}, \stationary \big)
\leq c_0 2^{- \lfloor \tfrac{\tau}{\mixingtime} \rfloor},
\end{align*}
valid for any $x \in \statespace$ and $\tau \geq 0$.
\end{lemma}
\noindent See Appendix~\ref{AppLemWassDecay} for the proof of this
claim. \\

We now use Lemma~\ref{LemWassDecay} to construct the desired process.
We begin by constructing a coupling conditionally on the
$\sigma$-field $\filtration_{t - \tau}$: let $\widetilde{\state}_t$ be
a state whose conditional law is $\stationary$, satisfying the
identity:
\begin{align}
    \Exs \big[ \metric (\state_t, \widetilde{\state}_t) \mid
      \filtration_{t - \tau} \big] = \Wass_{1, \metric} \big( \law
    (\state_t \mid \filtration_{t - \tau}), \stationary
    \big).\label{eq:coupling-wass-in-step-back-proof}
\end{align}
The existence of such $\widetilde{\state}_t$ is guaranteed by the
definition of Wasserstein distance. We now bound the relevant
quantities based on this construction.

Combining the identity~\eqref{eq:coupling-wass-in-step-back-proof}
with Lemma~\ref{LemWassDecay} yields \mbox{$\Exs \big[ \metric
    (\state_t, \widetilde{\state}_t) \mid \filtration_{t - \tau}\big]
  \leq c_0 \cdot 2^{- \lfloor \frac{\tau}{\mixingtime} \rfloor}$.}
Applying Cauchy--Schwarz inequality and invoking
Assumption~\ref{assume:stationary-tail}, we find that
\begin{align}
\big(\Exs \big[ \metric (s_t, \widetilde{s}_t)^p \mid \filtration_{t -
    \tau} \big] \big)^{1/p} & \leq \big( \Exs \big[ \metric (s_t,
  \widetilde{s}_t) \mid \filtration_{t - \tau}
  \big]\big)^{\frac{1}{2p}} \cdot \big( \Exs \big[ \metric (s_t,
  \widetilde{s}_t)^{2p - 1} \mid \filtration_{t - \tau} \big]
\big)^{\frac{1}{2p}} \nonumber \\
& \leq \big( \Exs \big[ \metric (s_t, \widetilde{s}_t) \mid
  \filtration_{t - \tau} \big]\big)^{\frac{1}{2p}} \nonumber \\
\label{eq:high-moment-coupling-estimate-in-step-back-proof}
& \leq c_0 \cdot 2^{1 - \frac{\tau}{2 \mixingtime p}},
\end{align}
which establishes the claim.

We now use the sequence of random variables $\widetilde{\state}_t$ just
constructed to define the extended filtration
$\widetilde{\filtration}_t \mydefn \sigma \big( (\state_k)_{0 \leq k
  \leq t}, (\widetilde{\state}_k)_{0 \leq k \leq t}, \big( (\Lmat_k,
\bvec_k) \big)_{0\leq k \leq t} \big)$, as well as the surrogate
quantities
\begin{align*}
\widetilde{\addnoiseMarkov}_t &\mydefn \big( \Lmap (\widetilde{s}_t) -
\Lbar \big) \thetabar + \big( \bmap (\widetilde{s}_t) - \bbar \big),
\quad \mbox{and}\\
 \widetilde{H}_1(t) &\mydefn \inprod{\Delta_{(t -
    \tau) \vee 0}}{\widetilde{\addnoiseMarkov}_t} + \inprod{\Delta_{(t
    - \tau) \vee 0}}{\big(\Lmap (\widetilde{s}_t) -
  \Lbar\big)\Delta_{(t - \tau) \vee 0}}.
\end{align*}
Note that by definition, we have $\Exs \big[ \widetilde{H}_1 (t) \mid
  \widetilde{\filtration}_{(t - \tau) \vee 0} \big] = 0$ for each $t =
0, 1, 2, \ldots$.


\subsubsection{Proof of the bound~\eqref{EqnMartingaleBound}}
\label{SecMartingaleBound}

We first perform a decomposition on the process $\widetilde{M}_1$. In
particular, for $\ell \in \{0,1, \cdots, \tau - 1\}$, we define the
stochastic process $\widetilde{M}_1^{(\ell)}(n) \mydefn \sum_{t =
  0}^{n - 1} e^{\stepsize (1 - \kappa) (t + \tau)} \widetilde{H}_1 (t
+ \tau) \bm{1}_{\{t \operatorname{mod} \tau = \ell\}}$.  Clearly, we
have $\widetilde{M}_1 (\numobs) = \sum_{\ell = 0}^{\tau - 1}
\widetilde{M}_1^{(\ell)} (\numobs)$ for any $\numobs \geq
0$. Furthermore, we note for any $t \geq 0$, we have the relations:
\begin{align*}
\Exs \big[ \widetilde{H}_1 (t + \tau) \mid
  \widetilde{\filtration}_{t}\big] = 0, \quad \mbox{and} \quad
\widetilde{H}_1 (t) \in \widetilde{\filtration}_t.
\end{align*}
So for each $\ell \in [0, \tau - 1]$, the process
$\widetilde{M}_1^{(\ell)}$ is a martingale adapted to the filtration
$\big(\widetilde{\filtration}_{t}\big)_{t \geq 0}$.

By the BDG inequality, we have the maximal inequality $\big(\Exs
\sup_{0 \leq t \leq \numobs} |\widetilde{M}_1^{(\ell)} (t)|^p
\big)^{1/p} \leq c p \big( \Exs \big(
\quadvar{\widetilde{M}_1^{(\ell)}}{\numobs} \big)^{p/2} \big)^{1/p}$,
valid for all $\ell = 0, 1, \ldots, \tau-1$.  Similarly, for the
quadratic variation term
$\quadvar{\widetilde{M}_1^{(\ell)}}{\numobs}$, we have that
\begin{align*}
&\Exs \big[ \big( \quadvar{\widetilde{M}_1^{(\ell)}}{\numobs}
  \big)^{p/2} \big]\\
   & = \Exs \big[ \big( \sum_{k = 0}^{\lfloor
    \frac{\numobs - 1}{\tau} \rfloor} e^{\stepsize (1 - \kappa) (k
    \tau + \tau + \ell)} \vecnorm{\widetilde{H}_1 (k \tau +
    \ell)}{2}^2 \big)^{p/2} \big] \\
& \leq \big(\sum_{k = 0}^{\lfloor \frac{\numobs - 1}{\tau} \rfloor}
e^{\stepsize (1 - \kappa) p (k \tau + \tau + \ell)} \Exs \big[
  \vecnorm{\widetilde{H}_1 (k \tau + \ell)}{2}^p\big] \big) \cdot
\big( \sum_{t = 0}^{\numobs - 1} e^{- \frac{p^2}{2p - 4} \tau
  \stepsize (1 - \kappa) t } \big)^{\frac{p - 2}{2}},
\end{align*}
which is at most
\begin{align*}
  \sum_{t
    = \tau}^{\numobs - 1} \tfrac{e^{\stepsize (1 - \kappa) t p}}{{( \stepsize \tau (1 - \kappa))^{p/2 - 1}} } \big( \Exs
  \big[ \abss{2 \inprod{\Delta_{t - \tau}}{(\Lbar
        (\widetilde{\state}_t) - \Lbar) \Delta_{t - \tau} }}^{p} \big]
  + \Exs \big[ \abss{2 \inprod{\widetilde{\addnoiseMarkov}_t}{
        \Delta_{t - \tau} }}^{p} \big] \big) \bm{1}_{\{t
    \operatorname{mod} \tau = \ell\}}.
\end{align*}
Invoking the tail condition in Assumption~\ref{assume:noise-moments} under the
stationary distribution, we have that
\begin{align*}
     \Exs \big[\abss{2 \inprod{\Delta_{t - \tau}}{(\Lbar
           (\widetilde{\state}_t) - \Lbar) \Delta_{t - \tau} }}^{p}
       \mid \filtration_{t - \tau} \big] &\leq \big( p \sigmaA
     \sqrt{d} \cdot \vecnorm{\Delta_{t - \tau}}{2}^2 \big)^p,
     \quad\mbox{and}\\ \Exs \big[ \abss{
         \inprod{\widetilde{\addnoiseMarkov}_t}{ \Delta_{t - \tau}
       }}^{p} \mid \filtration_{t - \tau} \big] &\leq \big( p
     \varbound \sqrt{d} \cdot \vecnorm{\Delta_{t - \tau}}{2} \big)^p.
\end{align*}
Substituting into the moment bounds for
$\quadvar{\widetilde{M}_1^{(\ell)}}{\numobs}$ and combining the
results for $\ell = 0,1, \cdots, \tau - 1$ using Minkowski's
inequality, we arrive at the bound
\begin{align*}
     &\big( \Exs \sup_{0 \leq t \leq \numobs} |\widetilde{M}_1 (t)|^p
  \big)^{1/p}\\ &\leq \sum_{\ell = 0}^{\tau - 1} \big( \Exs \sup_{0
    \leq t \leq \numobs} |\widetilde{M}_1^{(\ell)} (t)|^p
  \big)^{1/p}\\ &\leq \frac{ \tau \cdot \numobs^{\frac{1}{p}}
    \sqrt{p}}{\big(\stepsize \tau (1 - \kappa)\big)^{\frac{1}{2} +
      \frac{1}{p}}} \big\{ p \sigmaA \sqrt{d} \cdot \max_{0 \leq t
    \leq \numobs} \big[ e^{\stepsize (1 - \kappa) t} \big( \Exs
    \vecnorm{\Delta_t}{2}^{2p}\big)^{1/p} \big]]\\
    &\qquad + e^{\frac{\stepsize
      (1 - \kappa)\numobs}{2}} p \varbound \sqrt{d} \max_{0 \leq t
    \leq \numobs}\big[ e^{\stepsize (1 - \kappa) t / 2} \big( \Exs
    \vecnorm{\Delta_t}{2}^{p}\big)^{1/p} \big] \big\}\\
    \ &\leq
  \sqrt{\frac{\tau p}{\stepsize (1 - \kappa)}} \big( p \sigmaA
  \sqrt{d} \lyap_\numobs + p \varbound \sqrt{e^{\stepsize (1 - \kappa)
      \numobs}\lyap_\numobs d} \big),
\end{align*}
which completes the proof of this lemma.

\subsubsection{Proof of the bound~\eqref{EqnDifferenceBound}}
\label{SecDifferenceBound}

By Minkowski's inequality, we can upper bound the error as
\mbox{$\big( \Exs \big[ (H_1(t) - \widetilde{H}_1 (t))^p \big]
  \big)^{1/p} \leq \sum_{k = 1}^6 J_k$,} where
\begin{align*}
J_1 &\mydefn \big( \Exs \big[ \inprod{\Delta_{t -
      \tau}}{\addnoiseMarkov_t - \widetilde{\addnoiseMarkov}_t}^p
  \big]\big)^{1/p},\\
   J_2 &\mydefn \big( \Exs \big[
  \inprod{\Delta_{t} - \Delta_{t - \tau}}{\addnoiseMarkov_t}^p
  \big]\big)^{1/p} \\
J_3 &\mydefn \big( \Exs \big[ \inprod{\Delta_{t - \tau}}{\big( \Lmap
    (\widetilde{s}_t) - \Lmap (s_t) \big) \Delta_{t - \tau}}^p
  \big]\big)^{1/p}\\
   J_4 &\mydefn \big( \Exs \big[
  \inprod{\Delta_t - \Delta_{t - \tau}}{ \multnoiseMarkov_t \Delta_{t
      - \tau}}^p \big]\big)^{1/p} \\
J_5 &\mydefn \big( \Exs \big[ \inprod{\Delta_t}{ \multnoiseMarkov_t
    (\Delta_t - \Delta_{t - \tau})}^p \big]\big)^{1/p} \\
     J_6
&\mydefn \big( \Exs \big[ \inprod{\Delta_t - \Delta_{t - \tau}}{
    \multnoiseMarkov_t (\Delta_t - \Delta_{t - \tau})}^p
  \big]\big)^{1/p}
\end{align*}
The terms $J_1$ and $J_3$ can be controlled using the bound on
$\metric (\state_t, \widetilde{\state}_t)$ and the Lipschitz
condition~\eqref{assume-lip-mapping}; doing so yields the bound
\begin{align*}
J_1 & \leq \varbound d \big( \Exs \big[ \vecnorm{\Delta_{t -
      \tau}}{2}^p \cdot \Exs \big[ \metric (\state_t,
    \widetilde{\state}_{t})^p \mid \filtration_{t - \tau} \big]
  \big]\big)^{1/p} \leq 2^{1 - \frac{\tau}{2
    p \mixingtime}} c_0 \varbound d \big( \Exs
\vecnorm{\Delta_{t - \tau}}{2}^p \big)^{1/p}, \\ J_3 &\leq \matLip d \big( \Exs
\big[ \vecnorm{\Delta_{t - \tau}}{2}^{2p} \cdot \Exs \big[ \metric
    (\state_t, \widetilde{\state}_{t})^p \mid \filtration_{t - \tau}
    \big] \big]\big)^{1/p} \leq 2^{1 -
  \frac{\tau}{2 p \mixingtime}}c_0 \matLip d \big( \Exs
\vecnorm{\Delta_{t - \tau}}{2}^{2p} \big)^{1/p}.
\end{align*}
Given the time lag parameter $\tau \geq c p \mixingtime \log (c_0
\mixingtime d) \geq 2 p \mixingtime \log
\big(\tfrac{d}{\stepsize}\big)$, we have the bound
\begin{align}
    J_1 \leq \stepsize \varbound \sqrt{d} \big( \Exs
    \vecnorm{\Delta_{t - \tau}}{2}^p \big)^{1/p}, \quad \mbox{and}
    \quad J_3 \leq \stepsize \sigmaA \sqrt{d} \big( \Exs
    \vecnorm{\Delta_{t - \tau}}{2}^{2p}
    \big)^{1/p}.\label{eq:lemma-H1-j1-j3-bound}
\end{align}
Turning to the $J_2$ term, applying the Cauchy--Schwarz inequality
yields
\begin{align}
    J_2 \leq \big( \Exs \vecnorm{\Delta_t - \Delta_{t - \tau}}{2}^{2p}
    \big)^{\frac{1}{2p}} \cdot \big( \Exs
    \vecnorm{\addnoiseMarkov_t}{2}^{2p} \big)^{\frac{1}{2p}} &
    \stackrel{(i)}{\leq} \big( \Exs \vecnorm{\Delta_t - \Delta_{t -
        \tau}}{2}^{2p} \big)^{\frac{1}{2p}} \cdot p \varbound
    \sqrt{\usedim}.\label{eq:lemma-H1-j2-bound}
\end{align}
where step (i) follows from Assumption~\ref{assume:noise-moments}.

The terms $J_4$ and $J_5$ can be controlled via once again replacing
$\state_t$ with its surrogate $\widetilde{\state}_t$. First, by
Cauchy--Schwarz inequality, we note that
\begin{align*}
    J_4 &\leq \big( \Exs\vecnorm{\Delta_t - \Delta_{t - \tau}}{2}^{2p}
    \big)^{\frac{1}{2p}} \cdot \big( \Exs\vecnorm{\multnoiseMarkov_t
      \Delta_{t - \tau}}{2}^{2p} \big)^{\frac{1}{2p}},\\
       J_5 &\leq
    \big( \Exs\vecnorm{\Delta_t - \Delta_{t - \tau}}{2}^{2p}
    \big)^{\frac{1}{2p}} \cdot \big(
    \Exs\vecnorm{\multnoiseMarkov_t^\top \Delta_{t - \tau}}{2}^{2p}
    \big)^{\frac{1}{2p}}.
\end{align*}
Using the decomposition $\multnoiseMarkov_t = (\Lmap (\widetilde{s}_t)
- \Lbar) + (\Lmap (s_t) - \Lmap (\widetilde{s}_t))$, we note that
\begin{align*}
    \big( \Exs\vecnorm{\multnoiseMarkov_t \Delta_{t - \tau}}{2}^{2p}
    \big)^{\frac{1}{2p}} \leq \big( \Exs\vecnorm{(\Lmap
      (\widetilde{s}_t) - \Lbar) \Delta_{t - \tau}}{2}^{2p}
    \big)^{\frac{1}{2p}} + \big( \Exs\vecnorm{(\Lmap (s_t) - \Lmap
      (\widetilde{s}_t)) \Delta_{t - \tau}}{2}^{2p}
    \big)^{\frac{1}{2p}}.
\end{align*}
We bound the conditional expectations of the quantities above. The
first term can be controlled via
Assumption~\ref{assume:noise-moments}:
\begin{align*}
    \Exs \big[ \vecnorm{(\Lmap (\widetilde{s}_t) - \Lbar ) \Delta_{t -
          \tau}}{2}^{2p} \mid \filtration_{t - \tau} \big] \leq
    (\sigmaA p \sqrt{d})^{2p} \vecnorm{\Delta_{t - \tau}}{2}^{2p},
\end{align*}
and the second term is controlled using the Lipschitz
condition~\ref{assume-lip-mapping}:
\begin{align*}
    \Exs \big[ \vecnorm{(\Lmap (s_t) - \Lmap (\widetilde{s}_t))
        \Delta_{t - \tau}}{2}^{2p} \mid \filtration_{t - \tau} \big]
    &\leq (\matLip d)^{2p} \cdot \Exs \big[ \metric (s_t,
      \widetilde{s}_t)^{2p} \mid \filtration_{t - \tau} \big] \cdot
    \vecnorm{\Delta_{t - \tau}}{2}^{2p}\\ &
    \leq (\matLip d)^{2p} \cdot
    c_0 \cdot 2^{1 - \frac{\tau}{\mixingtime}} \cdot
    \vecnorm{\Delta_{t - \tau}}{2}^{2p}.
\end{align*}
Consequently, taking $\tau \geq 2 \mixingtime p \log (c_0 d)$, we have
the bounds
\begin{align*}
    \big( \Exs\vecnorm{\multnoiseMarkov_t \Delta_{t - \tau}}{2}^{2p}
    \big)^{\frac{1}{2p}} &\leq \sigmaA p \sqrt{d} \cdot \big(
    \Exs\vecnorm{\Delta_{t - \tau}}{2}^{2p} \big)^{\frac{1}{2p}},
    \quad \mbox{and} \\
\big( \Exs\vecnorm{\multnoiseMarkov_t^\top \Delta_{t - \tau}}{2}^{2p}
\big)^{\frac{1}{2p}} &\leq \sigmaA p \sqrt{d} \cdot \big(
\Exs\vecnorm{\Delta_{t - \tau}}{2}^{2p} \big)^{\frac{1}{2p}}.
\end{align*}
Putting together the pieces, we arrive at the bound
\begin{align}
J_4 + J_5 \leq 2 \big( \Exs\vecnorm{\Delta_t - \Delta_{t -
    \tau}}{2}^{2p} \big)^{\frac{1}{2p}} \cdot \sigmaA p \sqrt{d}
\cdot \big( \Exs \vecnorm{\Delta_{t - \tau}}{2}^{2p}
\big)^{\frac{1}{2p}}.\label{eq:lemma-H1-j4-j5-bound}
\end{align}
By the Lipschitz condition~\eqref{assume-lip-mapping} and the assumed
boundedness~\eqref{assume:stationary-tail} of the metric space, the
term $J_6$ admits the simple upper bound
\begin{align}
 J_6 & \leq \big(\Exs \big[ \opnorm{\multnoiseMarkov_t}^{p}
   \vecnorm{\Delta_t - \Delta_{t - \tau}}{2}^{2p}
   \big]\big)^{\frac{1}{p}} \leq \sigmaA d \big(\Exs
 \vecnorm{\Delta_t - \Delta_{t - \tau}}{2}^{2p} \big)^{\frac{1}{p}}\label{eq:lemma-H1-j6-bound}
    \end{align}
From all of these bounds, we see that the remaining crucial piece is
to bound $\Exs \vecnorm{\Delta_t - \Delta_{t - \tau}}{2}^{2p}$. In
order to do so, we require the following two helper lemmas
\begin{lemma}
\label{lemma:norm-not-blowup}
Given $p \geq 2$ and $\ell > 0$, the iterates~\eqref{eq:lsa-iterates}
with stepsize $\stepsize \leq \big( 6 (\conmax + \sigmaA d) \ell
\big)^{-1}$ satisfy the bound
\begin{subequations}    
  \begin{align}
    \big(\Exs \big[ \vecnorm{\Delta_{t + \ell} - \Delta_t}{2}^p \big]
    \big)^{1/p} \leq e \stepsize \ell (\conmax + \sigmaA d) \big(
    \Exs \big[ \vecnorm{\Delta_t}{2}^p \big] \big)^{1/p} + 3
    \stepsize p \ell \sqrt{d} \varbound,
    \end{align}
and consequently,
\begin{align}
\frac{1}{2} \big( \Exs \big[ \vecnorm{\Delta_{t }}{2}^p \big]
\big)^{1/p} - 6 \stepsize p \ell \sqrt{d} \varbound \leq \big( \Exs
\big[ \vecnorm{\Delta_{t + \ell}}{2}^p \big] \big)^{1/p} \leq e
\big( \Exs \big[ \vecnorm{\Delta_{t }}{2}^p \big] \big)^{1/p} + 6
\stepsize p \ell \sqrt{d}\varbound.
\end{align}
\end{subequations}
\end{lemma}
\noindent See Appendix~\ref{App:proof-norm-not-blowup} for the proof of this claim. \\

Our second auxiliary result is of a bootstrap nature: it is based on
assuming that for some given an integer $p \geq 2$, fix any integer
$\tau \geq 2 \mixingtime p \log (c_0 d)$, there exist positive scalars
$\omega_p, \beta_p > 0$ such that
  \begin{align}
\label{eq:local-move-coarse-bound}    
\big(\Exs \big[ \vecnorm{\Delta_{t + \ell} - \Delta_t}{2}^p \big]
\big)^{1/p} \leq \stepsize \omega_p \cdot \big( \Exs \big[
  \vecnorm{\Delta_t}{2}^p \big] \big)^{1/p} + \stepsize \beta_p
\varbound
  \end{align}
  for any $t \geq 0$, $\stepsize \leq \frac{1}{48 (\conmax + \sigmaA
    d) \tau}$ and $\ell \in [0, \tau]$.  We then have the following
  guarantee:
\begin{lemma}
  \label{lemma:local-moves}
When the condition~\eqref{eq:local-move-coarse-bound} holds, then, for
any $t \geq 0$, $\stepsize \leq \frac{1}{48(\conmax + \sigmaA d)
  \tau}$, and $\ell \in [0, \tau]$, we have
\begin{multline}
\label{eq:local-move-bootstrap-bound}  
\big(\Exs \big[ \vecnorm{\Delta_{t + \ell} - \Delta_t}{2}^p \big]
\big)^{1/p} \leq \stepsize \big( 12 \big( p \sqrt{d} \sigmaA +
\conmax \big) \ell 
+ \frac{\omega_p}{2} \big) \big( \big(\Exs
\vecnorm{\Delta_{t{}}}{2}^p\big)^{1/p} + \stepsize p (\tau + \ell)
\sqrt{d} \varbound \big) \\
+ \stepsize \big( 2 p \ell \sqrt{d} + \frac{1}{2} \beta_p\big)
\varbound.
\end{multline}
\end{lemma}
\noindent See Appendix~\ref{App:proof-local-moves} for the proof of
this claim. \\

We now complete the proof of the bound~\eqref{EqnDifferenceBound} by
using a bootstrapping argument in order to obtain a sharp bound on
$\Exs \vecnorm{\Delta_t - \Delta_{t - \tau}}{2}^{p}$.  Let
$\omega_p^{(0)} \mydefn e \tau (\conmax + \sigmaA d)$ and
$\beta_p^{(0)} \mydefn p \tau \sqrt{d}$, and define the following
recursion:
\begin{align*}
  \begin{cases}
        \omega_p^{(i + 1)} = \frac{1}{2} \omega_p^{(i)} + 12 \big( p
        \sqrt{d} \sigmaA + \conmax \big) \tau, \\ \beta_p^{(i + 1)} =
        \frac{1}{2} \beta_p^{(i)} + 2 p \tau \sqrt{d} + 2 \stepsize
        \big( 12 \big( p \sqrt{d} \sigmaA + \conmax \big) \tau +
        \frac{1}{2} \omega_p^{(i)} \big) p \tau \sqrt{d}.
         \end{cases}
\end{align*}
It can be seen that as $i \rightarrow \infty$, the sequence
$(\omega_p^{(i)}, \beta_p^{(i)})$ converges to a unique limit
$(\omega_p^*, \beta_p^*)$; this limit is the unique fixed point of the
iterates defined above.

By Lemma~\ref{lemma:local-moves}, if the iterates satisfy the
bound~\eqref{eq:local-move-coarse-bound} with constants $\big(
\omega_p^{(i)}, \beta_p^{(i)}\big)$, then it also satisfy the bound
with constants $\big( \omega_p^{(i + 1)}, \beta_p^{(i + 1)}\big)$. By
Lemma~\ref{lemma:norm-not-blowup}, the iterates satisfy bound with
constants $\big( \omega_p^{(0)}, \beta_p^{(0)}\big)$. An induction
argument then yields the bound for any $\big( \omega_p^{(i)},
\beta_p^{(i)}\big)$. In particular, the bound is satisfied by the
fixed point $\big( \omega_p^*, \beta_p^*)$.

Solving directly for the fixed-point equation, we find that
\begin{align*}
  \omega_p^* = 24 \big( p \sqrt{d} \sigmaA + \conmax \big) \tau, \quad
  \mbox{and} \quad \beta_p^* = 4 p \tau \sqrt{d} + 96 \stepsize \big(
  p \sqrt{d} \sigmaA + \conmax \big) p \tau^2\sqrt{d}.
\end{align*}
Taking the stepsize $\stepsize \leq \frac{1}{48(\conmax + p\sigmaA d)
  \tau}$, we arrive at the bound
\begin{align}
    \big(\Exs \big[ \vecnorm{\Delta_{t + \ell} - \Delta_t}{2}^p \big]
    \big)^{1/p} \leq 24 \stepsize \tau \big( p \sqrt{d} \sigmaA +
    \conmax \big) \big(\Exs \vecnorm{\Delta_{t}}{2}^p\big)^{1/p} + 6
    \stepsize p \tau \sqrt{d}
    \varbound, \label{eq:local-move-final-bound}
\end{align}
for any $t \geq 0$ and $\ell \in [0, \tau]$.

Collecting the
bounds~\eqref{eq:lemma-H1-j1-j3-bound},~\eqref{eq:lemma-H1-j2-bound},~\eqref{eq:lemma-H1-j4-j5-bound},~\eqref{eq:lemma-H1-j6-bound}
and~\eqref{eq:local-move-final-bound} and taking the stepsize
$\stepsize \leq \frac{1}{c (\conmax + p \sigmaA d) \tau}$, we arrive
at the bound
\begin{align*}
    \big( \Exs \big[ (H_1 (t) - \widetilde{H}_1 (t))^p
      \big]\big)^{1/p} \leq c \stepsize p^2 \tau \big( ( d
    \sigmaA^2 + \conmax^2) \cdot \big( \Exs \vecnorm{\Delta_{t -
        \tau}}{2}^{2p} \big)^{\frac{1}{p}} + \varbound^2 d \big),
\end{align*}
thereby completing the proof of the bound~\eqref{EqnDifferenceBound}.


\subsection{Proof of Lemma~\ref{LemNewH2}}
\label{SecProofNewH2}

By the BDG inequality, we have the bound $\big( \Exs \sup_{0 \leq t
  \leq \numobs} |M_2 (t)|^p \big)^{1/p} \leq c p \big( \Exs \big(
\quadvar{M_2}{\numobs} \big)^{p/2} \big)^{1/p}$, valid for all $\ell =
0, 1, \ldots, \tau-1$.

As for the quadratic variation $\quadvar{M_2}{\numobs}$, applying
H\"{o}lder's inequality yields
\begin{align*}
 &\Exs \big[ \big( \quadvar{M_2}{\numobs} \big)^{p/2} \big]\\
  &=
 \Exs \big[ \big( \sum_{t = 0}^{\numobs - 1} e^{\stepsize (1 -
     \kappa) t} \vecnorm{H_2 (t)}{2}^2 \big)^{p/2} \big] \\
& \leq \big( \sum_{t = 0}^{\numobs - 1} e^{\stepsize (1 - \kappa) t
      p} \Exs \big[ \vecnorm{H_2 (t)}{2}^p\big] \big) \cdot
    \big( \sum_{t = 0}^{\numobs - 1} e^{- \frac{p^2}{2p - 4}
      \stepsize (1 - \kappa) t } \big)^{\frac{p - 2}{2}} \\
& \leq \big( \stepsize (1 - \kappa) \big)^{- \frac{p}{2} + 1 } \sum_{t
      = 0}^{\numobs - 1} e^{\stepsize (1 - \kappa) t p} \big( \Exs
    \big[ \abss{2 \inprod{\Delta_t}{\multnoiseMG_{t + 1} \Delta_t
      }}^{p} \big] + \Exs \big[ \abss{2 \inprod{\addnoiseMG_{t +
            1}}{ \Delta_t }}^{p} \big] \big).
\end{align*}

For the moment terms above, we invoke
Assumption~\ref{assume:noise-moments}, and obtain the following
bounds:
\begin{align*}
\Exs \big[ \abss{\inprod{\Delta_t}{\multnoiseMG_{t + 1} \Delta_t
  }}^{p} \mid \filtration_t \big] &\leq \vecnorm{\Delta_t}{2}^p
\cdot \Exs \big[ \big( \sum_{j = 1}^\usedim
  \inprod{\coordinate_j}{\multnoiseMG_{t + 1} \Delta_t
  }^2\big)^{p/2} \mid \filtration_t \big] \leq \big( p \sigmaA
\sqrt{d} \cdot \vecnorm{\Delta_t}{2}^2 \big)^p,\\ \Exs \big[ \abss{
    \inprod{\addnoiseMG_{t + 1}}{ \Delta_t }}^{p} \mid \filtration_t
  \big] &\leq \vecnorm{\Delta_t}{2}^p \cdot \Exs \big[ \big(
  \sum_{j = 1}^\usedim \inprod{\coordinate_j}{\addnoiseMG_{t +
      1}}^2\big)^{p/2} \mid \filtration_t \big] \leq \big( p
\varbound \sqrt{d} \cdot \vecnorm{\Delta_t}{2} \big)^p.
\end{align*}
Substituting into the bound above, we find that
\begin{align*}
     &\big( \Exs \big[ \big( \quadvar{M_2}{\numobs} \big)^{p/2} \big]
  \big)^{1/p}\\ & \leq \tfrac{ (\stepsize (1 - \kappa))^{-\frac{1}{p}}
    \cdot \numobs^{\frac{1}{p}}}{\sqrt{\stepsize (1 - \kappa)}} \big\{
  p \sigmaA \sqrt{d} \cdot \max_{0 \leq t \leq \numobs} \big[
    e^{\stepsize (1 - \kappa) t} \big( \Exs
    \vecnorm{\Delta_t}{2}^{2p}\big)^{1/p} \big]\\
    &\qquad \qquad + e^{\frac{\stepsize
      (1 - \kappa)\numobs}{2}} p \varbound \sqrt{d} \max_{0 \leq t
    \leq \numobs}\big[ e^{\stepsize (1 - \kappa) t / 2} \big( \Exs
    \vecnorm{\Delta_t}{2}^{p}\big)^{1/p} \big] \big\}\\ &\leq
  \tfrac{1}{\sqrt{\stepsize (1 - \kappa)}} \big( p \sigmaA \sqrt{d}
  \lyap_\numobs + p \varbound \sqrt{e^{\stepsize (1 - \kappa)
      \numobs}\lyap_\numobs d} \big).
\end{align*}


\subsection{Proof of Lemma~\ref{LemNewH3}}
\label{SecProofNewH3}

Recall the definitions~\eqref{eq:defs-mg-noise}
and~\eqref{eq:defs-markov-noise}.  By Minkowski's inequality, we have
the upper bound
\begin{align}
  \label{EqnInitialMin}
&\big( \Exs \big[ H_3 (t)^{p} \big]\big)^{1/p} \nonumber\\
& \leq \big( \Exs
\vecnorm{\multnoiseMarkov_t \Delta_t}{2}^{2p} \big)^{1/p} + \big(
\Exs \vecnorm{\multnoiseMG_{t + 1} \Delta_t}{2}^{2p} \big)^{1/p} +
\big( \Exs \vecnorm{\addnoiseMG_{t + 1}}{2}^{2p} \big)^{1/p} +
\big( \Exs \vecnorm{\addnoiseMarkov_t}{2}^{2p} \big)^{1/p},
\end{align}
For the martingale part of the noise, we note that
Assumption~\ref{assume:noise-moments} implies that
\begin{align*}
\big( \Exs \vecnorm{\multnoiseMG_{t + 1} \Delta_t}{2}^{2p} \mid
\filtration_t \big)^{1/p} \leq p^2 \sigmaA^2 \usedim \cdot
\vecnorm{\Delta_t}{2}^{2}, \quad \mbox{and} \quad \big( \Exs
\vecnorm{\addnoiseMG_{t + 1} }{2}^{2p} \big)^{1/p} \leq p^2
\varbound^2 \usedim.
\end{align*}
For the additive Markov noise, applying
Assumption~\ref{assume:noise-moments} yields the bound $\big( \Exs
\vecnorm{\addnoiseMarkov_t}{2}^{2p} \big)^{1/p} \leq p^2 \varbound^2
\usedim$.

For the Markov part of the multiplicative noise, we make use of the
construction given in Section~\ref{SecConstruct}, where we showed that
for a given $\tau > 0$, there exists a random variable
$\widetilde{\state}_t$ such that $\widetilde{\state}_t \mid
\filtration_{t - \tau} \sim \stationary$, and $\Exs \big[ \metric^p
  (\state_t, \widetilde{\state}_t) \mid \filtration_{t - \tau} \big]
\leq c_0 \cdot 2^{1 - \frac{\tau}{\mixingtime}}$.  Observe the
decomposition
\begin{align*}
  \multnoiseMarkov_t \Delta_t = \big(\Lmap (\state_t) - \Lmat
  (\widetilde{\state}_t) \big) \Delta_{t - \tau} + \big(\Lmat
  (\widetilde{\state}_t) - \Lbar \big) \Delta_{t - \tau} +
  \multnoiseMarkov_t \big(\Delta_t - \Delta_{t - \tau} \big).
\end{align*}
Using the Lipschitz condition~\eqref{assume-lip-mapping}, we have that
\begin{align*}
  \Exs \big[ \vecnorm{\big(\Lmap (\state_t) - \Lmap
      (\widetilde{\state}_t) \big) \Delta_{t - \tau}}{2}^{2p} \mid
    \filtration_{t - \tau} \big] \leq c_0 \cdot 2^{1 -
    \frac{\tau}{\mixingtime}} \big(\sigmaA d \vecnorm{\Delta_{t -
      \tau}}{2}\big)^{2p}.
\end{align*}
For any $\tau \geq 2 p \mixingtime \log d$, we have the bound
\begin{align*}
\big( \Exs \big[ \vecnorm{\big(\Lmap (\state_t) - \Lmap
    (\widetilde{\state}_t) \big) \Delta_{t - \tau}}{2}^{2p}\big]
\big)^{1/p} \leq p^2 \sigmaA^2 d \cdot \big( \Exs
\vecnorm{\Delta_t}{2}^{2p} \big)^{1/p}.
\end{align*}
By the moment bounds~\eqref{assume:noise-moments} on the stationary
distribution, we have
\begin{align*}
 \Exs \big[ \vecnorm{ \big(\Lmap (\widetilde{\state}_t) - \Lbar \big)
     \Delta_{t - \tau}}{2}^{2p} \mid \filtration_{t - \tau}
   \big]\leq \big( 2 p \sigmaA \sqrt{d} \vecnorm{\Delta_{t -
     \tau}}{2} \big)^{2p}.
\end{align*}
For the last term, we use the Lipschitz
condition~\ref{assume-lip-mapping} as well as the boundedness
condition~\ref{assume:stationary-tail} of metric space. In conjunction
with the inequality~\eqref{eq:local-move-final-bound}, for $\tau \geq
2 p \mixingtime \log (c_0 d)$ and stepsize $\stepsize \leq \frac{1}{48
  \tau (\sigmaA d + \conmax)}$, we arrive at the bound
\begin{align*}
 &\big(\Exs \big[ \vecnorm{\multnoiseMarkov_t (\Delta_t - \Delta_{t -
       \tau})}{2}^{2p} \big]\big)^{1/p}\\
        &\leq \sigmaA^2 d^2 \cdot
 \big( \Exs \big[\vecnorm{\Delta_t - \Delta_{t -
       \tau}}{2}^{2p}\big]\big)^{1/p} \\
& \leq c \stepsize^2 \sigmaA^2 d^2 \tau^2 \big(p^2 \sigmaA^2 d +
 \conmax^2 \big) \big( \Exs \big[ \vecnorm{\Delta_{t -
       \tau}}{2}^{2p} \big]\big)^{1/p} + c\stepsize^2 p^2
 \sigmaA^2 \varbound^2 d^3 \tau^2 \\
& \leq c \big(p^2 \sigmaA^2 d + \conmax^2 \big) \big( \Exs \big[
   \vecnorm{\Delta_{t - \tau}}{2}^{2p} \big]\big)^{1/p} + c p^2
 \varbound^2 d,
\end{align*}
for a universal constant $c > 0$.

Collecting the bounds above and substituting into our initial
bound~\eqref{EqnInitialMin}, we find that
\begin{align*}
  \big( \Exs \big[ H_3 (t)^{p} \big]\big)^{1/p} & \leq c
  \big(p^2 \sigmaA^2 d + \conmax^2 \big) \big( \Exs \big[
    \vecnorm{\Delta_{t - \tau}}{2}^{2p} \big]\big)^{1/p} + c p^2
  \varbound^2 d,
\end{align*}
as claimed.


\section{Proof of Theorem~\ref{thm:markov-main}}
\label{subsec:proof-main-prj}

From the defining equations~\eqref{eq:lsa-iterates}
and~\eqref{eq:lsa-average}, we have the telescoping relation
\begin{align}
\label{eq:telescope-in-proof-of-markov-case}  
 \tfrac{\theta_\numobs - \theta_\numburn}{\stepsize (\numobs -
   \numburn)} & = \tfrac{1}{\numobs - \numburn} \sum_{t =
   \numburn}^{{\numobs} - 1} \big( \theta_t - \Lmat_{t + 1} \theta_t -
 \bvec_{t + 1} \big) = (I - \Lbar) (\thetahat_\numobs - \thetabar) +
 \tfrac{1}{\numobs - \numburn}\martingale_{\numburn, \numobs}
 +\tfrac{1}{\numobs-\numburn} \Upsilon_{\numburn, \numobs}
\end{align}
where $\martingale_{\numburn, \numobs} = \sum_{t =
  \numburn}^{{\numobs} - 1} \big(\Lmat_{t + 1} \theta_t + \bvec_{t +
  1} - \Exs \big[ \Lmat_{t + 1} \theta_t + \bvec_{t + 1} |
  \filtration_t \big] \big)$ and $\Upsilon_{\numburn, \numobs} \defn
\sum_{t = \numburn}^{\numobs - 1} \big(
\Lmap (s_t) \theta_t + \bmap (s_t) - \Lbar \theta_t - \bbar \big)$.
Some algebra yields
\begin{align}
\label{eq:lsa-markov-decomposition-3}
\thetahat_\numobs - \thetabar & = \tfrac{(I - \Lbar)^{-1}
  \big(\theta_\numobs - \theta_\numburn \big)}{\stepsize (\numobs -
  \numburn)} - \tfrac{ (I - \Lbar)^{-1} \martingale_{\numburn,
    \numobs} }{\numobs - \numburn} - \tfrac{(I - \Lbar)^{-1}
  \Upsilon_{\numburn, \numobs} }{\numobs - \numburn} =: I_1 + I_2 +
I_3
\end{align}
From the triangle inequality, it suffices to bound the norms of $I_1$,
$I_2$ and $I_3$.

In the following, we prove a slightly stronger claim, which gives
bounds on an arbitrary quadratic loss functional. In particular, given
a matrix $\testMat \succ 0$, we seek bounds on the $\testMat$-norm
$\vecnorm{\thetahat_\numobs - \thetabar}{\testMat} \mydefn
\sqrt{(\thetahat_\numobs - \thetabar)^\top \testMat (\thetahat_\numobs
  - \thetabar)}$.

\subsection{Bounding the three terms}

We now bound each term in the
decomposition~\eqref{eq:lsa-markov-decomposition-3} in turn. \\

\subsubsection{Bounding the term $I_1$}

The bound for term $I_1$ follows directly from
Proposition~\ref{prop:lsa-markov-iterate-bound}.  In particular, given
a sample size $\numobs \geq \frac{8}{\stepsize (1 - \kappa)} \log
\big( \vecnorm{\theta_0 - \thetabar}{2} d / \stepsize \big)$ and
burn-in period $\numburn = \numobs / 2$, we have
\begin{align*}
  \Exs \big[ \vecnorm{\theta_\numobs -
    \thetabar}{2}^{2} \big] \leq \frac{c \stepsize}{1 -
    \kappa} \varbound^2 \tau d, \quad \mbox{and} \quad  \Exs
  \big[ \vecnorm{\theta_\numburn - \thetabar}{2}^{2}\big]\leq
  \frac{c \stepsize}{1 - \kappa} \varbound^2 \tau d.
\end{align*}
Noting that $\opnorm{(I - \Lbar)^{-1}} \leq (1 - \kappa)^{-1}$, we
conclude that
\begin{align}
\label{eq:prj-thm-term-i1-bound}  
\Exs \big[\vecnorm{I_1}{\testMat}^2 \big] \leq \lammax (\testMat) \Exs
\big[\vecnorm{I_1}{2}^2 \big] \leq \lammax (\testMat) \cdot \tfrac{c
  \varbound^2 \tau d}{\stepsize (1 - \kappa)^3 \numobs^2}.
\end{align}


\subsubsection{Bounding the term $I_2$}

For the term $I_2$, note that the process $(\Psi_t)_{t \geq \numburn}$
is a martingale adapted to the natural filtration. Its second moment
equals the quadratic variation:
\begin{align*}
    \Exs \big[\vecnorm{I_2}{\testMat}^2 \big] &= \frac{4}{\numobs^2} \Exs
    \big[ \quadvar{\testMat^{1/2} (I - \Lbar)^{-1} \Psi}{\numburn, \numobs} \big]\\
    &= \frac{4}{\numobs^2} \sum_{t = \numburn}^{\numobs - 1} \Exs
    \big[ \vecnorm{(I - \Lbar)^{-1} \big( (\Lmat_{t + 1} - \Lmap (
        \state_t) ) \theta_t + \bvec_{t + 1} - \bmap (\state_t)
        \big)}{\testMat}^2 \big].
\end{align*}
By the Cauchy--Schwarz inequality, we have the bound
\begin{align}
    \Exs \big[ \vecnorm{I_2}{\testMat}^2 \big] & \leq
    \tfrac{8}{\numobs^2} \sum_{t = \numburn}^{\numobs - 1} \Exs \big[
      \vecnorm{(I - \Lbar)^{-1} \addnoiseMG_{t + 1}}{\testMat}^2 \big]
    + \tfrac{8}{\numobs^2} \sum_{t = \numburn}^{\numobs - 1} \Exs
    \big[ \vecnorm{(I - \Lbar)^{-1} \multnoiseMG_{t + 1} \Delta_t
      }{\testMat}^2 \big] \nonumber \\
    & \leq \tfrac{16}{\numobs} \mathrm{Tr} \big( \testMat (I -
    \Lbar)^{-1} \SigStar_{\mathrm{MG}} (I - \Lbar)^{-\top} \big) +
    \tfrac{16 \sigmaA^2 \lammax (\testMat) \usedim}{(1 - \kappa)^2
      \numobs^2} \sum_{t = \numburn}^{\numobs - 1} \Exs \big[
      \vecnorm{\Delta_t}{2}^2 \big]\nonumber\\
    \label{eq:prj-thm-i2-bound}    
    & \leq \tfrac{16}{\numobs} \mathrm{Tr} \big( (I - \Lbar)^{-1}
    \SigStar_{\mathrm{MG}} (I - \Lbar)^{-\top} \big) + \lammax
    (\testMat) \cdot \tfrac{16 \sigmaA^2 \usedim}{(1 -
      \kappa)^2\numobs} \cdot \tfrac{c \stepsize \usedim \tau}{1 -
      \kappa} \varbound^2.
\end{align}


\subsubsection{Bounding the term $I_3$}

Applying the Cauchy-Schwarz inequality yields
\begin{align}
\label{eq:markov-noise-decomposition-for-avg-iterates}  
  \Exs \big[ \vecnorm{(I - \Lbar)^{-1} \Upsilon_{\numburn,
        \numobs}}{2}^2 \big] \leq 2 \Exs \big[ \vecnorm{\sum_{t =
        \numburn}^{\numobs - 1} (I - \Lbar)^{-1}
      \addnoiseMarkov_t}{2}^2 \big] + 2 \Exs \big[ \vecnorm{\sum_{t
        = \numburn}^{\numobs - 1} (I - \Lbar)^{-1} \multnoiseMarkov_t
      \Delta_t}{2}^2 \big].
\end{align}
We make use of the two auxiliary lemmas in order to control the terms
in the
decomposition~\eqref{eq:markov-noise-decomposition-for-avg-iterates}.
\begin{lemma}
\label{lemma:relate-finite-sum-to-sigstar-mkv}
Under the setup above, for a sample size $\numobs$ satisfying the bound
$\frac{\numobs}{\log \numobs} \geq 2 \mixingtime \log (c_0 d)$, there
exists a universal constant $c > 0$ such that
\begin{multline*}
  \Exs \big[ \vecnorm{\sum_{t = \numburn}^{\numobs - 1} (I -
      \Lbar)^{-1} \addnoiseMarkov_t}{\testMat}^2 \big] \leq (\numobs -
  \numburn) \cdot \mathrm{Tr} \big( \testMat (I - \Lbar)^{-1}
  \SigStar_{\mathrm{Mkv}} (I - \Lbar)^{-\top} \big) \\
  + \lammax
  (\testMat) \cdot\frac{c \mixingtime^2 \varbound^2 \usedim}{(1 -
    \kappa)^2} \log^2 (c_0 d).
\end{multline*}
\end{lemma}
\noindent See Section~\ref{AppRelated} for the proof of this claim. \\

\begin{lemma}
\label{lemma:bound-sum-markov-mult-noise}
Under the above conditions, there exists a universal constant $c > 0$
such that for any scalar $\tau \geq 3 \mixingtime \log^2 (c_0 d
\numobs)$, stepsize $\stepsize \in \big(0, \frac{1 - \kappa}{c \tau
  (\sigmaA^2 d + \conmax^2)}\big]$ and burn-in time $\numburn \geq
  \tau + \frac{2}{(1 - \kappa) \stepsize} \log (\numobs d)$, we have
  $\Exs \big[ \vecnorm{\sum_{t = \numburn}^{\numobs - 1}
      \multnoiseMarkov_t \Delta_t}{2}^2 \big] \leq c \stepsize^2
  \numobs^2 \tau^2 \usedim^2 \sigmaA^2 \varbound^2$.
\end{lemma}
\noindent See Section~\ref{AppMultNoise} for the proof of this
claim. \\

\noindent We now exploit the preceding two lemmas to upper bound the
term $I_3$.  We have
\begin{align}
\Exs \big[ \vecnorm{I_3}{\testMat}^2 \big] & \leq \tfrac{2}{(\numobs -
  \numburn)^2} \Exs \big[ \vecnorm{\sum_{t = \numburn}^{\numobs - 1}
    (I - \Lbar)^{-1} \addnoiseMarkov_t}{\testMat}^2 \big] +
\tfrac{2}{(\numobs - \numburn)^2} \Exs \big[ \vecnorm{\sum_{t =
      \numburn}^{\numobs - 1} (I - \Lbar)^{-1} \multnoiseMarkov_t
    \Delta_t}{\testMat}^2 \big] \nonumber\\
& \leq \tfrac{8\mathrm{Tr} \big( \testMat (I - \Lbar)^{-1}
  \SigStar_{\mathrm{Mkv}} (I - \Lbar)^{-\top} \big)}{\numobs} +
\lammax (\testMat) \Big\{ \tfrac{c \mixingtime^2 \varbound^2
  \usedim}{(1 - \kappa)^2 \numobs^2} \log^2 (c_0 d) +  \tfrac{ c \stepsize^2 \tau^2 \usedim^2 \sigmaA^2
  \varbound^2 }{(1 - \kappa)^2} \Big\}.\label{eq:prj-thm-i3-bound}
\end{align}
Collecting the
bounds~\eqref{eq:prj-thm-term-i1-bound},~\eqref{eq:prj-thm-i2-bound},
and~\eqref{eq:prj-thm-i3-bound}, we find that
\begin{align*}
    \Exs \big[ \vecnorm{\thetahat_\numobs - \thetabar}{\testMat}^2
      \big] & \leq \tfrac{c}{\numobs} \mathrm{Tr} \big( \testMat (I -
    \Lbar)^{-1} (\SigStar_{\mathrm{MG}} + \SigStar_{\mathrm{Mkv}})(I -
    \Lbar)^{-\top} \big) \\ & \qquad+ \lammax (\testMat) \cdot \big[
      \frac{c \varbound^2 \mixingtime d }{\stepsize (1 - \kappa)^3
        \numobs^2} + \frac{16 \sigmaA^2 \usedim}{(1 -
        \kappa)^2\numobs} \cdot \frac{c \stepsize \usedim \mixingtime
      }{1 - \kappa} \varbound^2 \big]\\ &\qquad\qquad+ \lammax
    (\testMat) \cdot \big[ \frac{c \mixingtime^2 \varbound^2
        \usedim}{(1 - \kappa)^2 \numobs^2} \log^2 (c_0 d \numobs) +
      \frac{ c \stepsize^2 \mixingtime^2 \usedim^2 \sigmaA^2
        \varbound^2}{(1 - \kappa)^2} \big].
\end{align*}
For a sample size $\numobs$ lower bounded as $\frac{\numobs}{\log^2
  \numobs} \geq \frac{2 \mixingtime (\sigmaA^2 \usedim +
  \conmax^2)}{(1 - \kappa)^2} \log (c_0 d)$, we can take the optimal
stepsize $\stepsize = \big[c \big( (1 - \kappa) \numobs^2 \mixingtime
  (\sigmaA^2 d + \conmax^2) \big) \big]^{-1/3}$.  With this choice, we
have
\begin{multline}
\label{eq:estimation-error-bound-testmat-general}  
  \Exs \big[ \vecnorm{\thetahat_\numobs - \thetabar}{\testMat}^2 \big]
   \leq \tfrac{c}{\numobs} \mathrm{Tr} \big( \testMat (I -
  \Lbar)^{-1} (\SigStar_{\mathrm{MG}} + \SigStar_{\mathrm{Mkv}})(I -
  \Lbar)^{-\top} \big)\\
   + c \lammax (\testMat) \cdot \big(
  \frac{\sigmaA^2 \usedim \mixingtime}{(1 - \kappa)^2 \numobs}
  \big)^{4/3} \varbound^2 \log^2 \numobs.
\end{multline}
Setting $\testMat \mydefn I_d$ completes the proof.


\subsection{Proof of auxiliary results}

In this section, we prove the two auxiliary results used in the proof
of Theorem~\ref{thm:markov-main}: namely,
Lemma~\ref{lemma:relate-finite-sum-to-sigstar-mkv} and
Lemma~\ref{lemma:bound-sum-markov-mult-noise}.


\subsubsection{Proof of Lemma~\ref{lemma:relate-finite-sum-to-sigstar-mkv}}
\label{AppRelated}
Given an integer $k \geq 0$, we define the $k$-step correlation under
the stationary Markov chain as
\begin{align*}
\crosscorrelation_k \mydefn \Exs_{\state \sim \stationary, \state'
  \sim \transition^k \delta_\state} \big[ \inprod{ \testMat^{1/2} (I - \Lbar)^{-1}
    \addnoiseMarkov (\state)}{ \testMat^{1/2} (I - \Lbar)^{-1} \addnoiseMarkov
    (\state')} \big].
\end{align*}
Clearly, we have $\crosscorrelation_0 \geq 0$, and by Cauchy--Schwarz
inequality, for any $k \geq 0$, there is:
\begin{align*}
     \abss{\crosscorrelation_k} \leq \sqrt{\Exs_{\state \sim
         \stationary} \vecnorm{(I - \Lbar)^{-1} \addnoiseMarkov
         (\state)}{\testMat}^2} \cdot \sqrt{\Exs_{\state' \sim \stationary}
       \vecnorm{(I - \Lbar)^{-1} \addnoiseMarkov (\state')}{\testMat}^2} =
     \crosscorrelation_0.
\end{align*}
The desired quantity can be written as $\mathrm{Tr} \big(
\testMat^{1/2} (I - \Lbar)^{-1} \SigStar_{\mathrm{Mkv}} (I -
\Lbar)^{-\top} \testMat^{1/2} \big) = \crosscorrelation_0 + 2 \sum_{k
  = 1}^{+ \infty} \crosscorrelation_k$.  Expanding the squared norm
yields
\begin{align*}
  &\Exs \big[ \vecnorm{\sum_{t = \numburn}^{\numobs - 1} \testMat^{1/2} (I -
      \Lbar)^{-1} \addnoiseMarkov_t}{2}^2 \big]\\
       &= \sum_{\numburn
    \leq t_1, t_2 \leq \numobs - 1} \Exs \big[\inprod{\testMat^{1/2} (I -
      \Lbar)^{-1} \addnoiseMarkov (\state_{t_1})}{\testMat^{1/2} (I - \Lbar)^{-1}
      \addnoiseMarkov (\state_{t_2})} \big]\\ &= (\numobs -
  \numburn) \crosscorrelation_0 + 2\sum_{k = 1}^{\numobs - \numburn -
    1} (\numobs - \numburn - k) \crosscorrelation_k.
\end{align*}

We claim that the cross-correlations $\crosscorrelation_k$ satisfy the
bound
\begin{align}
\label{eq:cross-correlation-decay}  
|\crosscorrelation_k| \leq c_0 \frac{\varbound^2 \opnorm{\testMat}
  \usedim^2}{(1 - \kappa)^2} \cdot 2^{1 - \frac{k}{2\mixingtime}}.
\end{align}
We return to prove this fact momentarily.  Taking it as given, this
inequality, in conjunction with the bound $|\crosscorrelation_k| \leq
\crosscorrelation_0$, can be employed to bound the tail sums needed
for the proof. We have
\begin{align*}
\abss{\sum_{k = 1}^{\numobs - \numburn - 1} k \crosscorrelation_k}
&\leq \sum_{k = 1}^\tau \tau |\crosscorrelation_k| + \sum_{k = \tau +
  1}^{\infty} k |\crosscorrelation_k| \leq \tau^2 \crosscorrelation_0
+ 2 c_0 \frac{\varbound^2 \opnorm{\testMat} \usedim^2}{(1 - \kappa)^2}
\sum_{k = \tau + 1}^{\infty} k \cdot 2^{- \frac{k}{2 \mixingtime}}.
\end{align*}
With the choice $\tau \mydefn 2\mixingtime \log (c_0 d)$, simplifying
yields
\begin{align*}
\abss{\sum_{k = 1}^{\numobs - \numburn - 1} k \crosscorrelation_k} &
\leq \frac{\tau^2 \varbound^2 \usedim \opnorm{\testMat}}{(1 -
  \kappa)^2} + 2 c_0 \frac{\varbound^2 \usedim^2 \opnorm{\testMat}}{(1
  - \kappa)^2} \cdot 2 \mixingtime \big( \tau + 1 + 2 \mixingtime
\big) \cdot 2^{- \frac{\tau + 1}{\mixingtime}}\\ &\leq \frac{2 \tau^2
  \varbound^2 \usedim}{(1 - \kappa)^2} \opnorm{\testMat},
\end{align*}
and for $\numobs$ satisfying $\frac{\numobs}{\log \numobs} \geq 2 \log
(c_0 \usedim \mixingtime)$, we have:
\begin{multline*}
\sum_{k = \numobs - \numburn}^{\infty} |\crosscorrelation_k| \leq 2
c_0 \frac{\varbound^2 \usedim^2 \opnorm{\testMat}}{(1 - \kappa)^2} \sum_{k =
  \frac{1}{2}\numobs}^{\infty} \cdot 2^{- \frac{k}{2 \mixingtime}}\\
\leq 2 c_0 \frac{\varbound^2 \usedim^2 \opnorm{\testMat}}{(1 - \kappa)^2} \cdot 2^{-
  \frac{\numobs}{2 \mixingtime}} \leq 2 c_0 \frac{\varbound^2
  \usedim}{(1 - \kappa)^2 \numobs^2}\opnorm{\testMat}.
\end{multline*}
Putting together these bounds yields
\begin{align*}
  &\Exs \big[ \vecnorm{\sum_{t = \numburn}^{\numobs - 1} (I -
      \Lbar)^{-1} \addnoiseMarkov_t}{\testMat}^2 \big]\\
       &= (\numobs -
  \numburn) \big( \crosscorrelation_0 + 2 \sum_{k = 1}^{\infty}
  \crosscorrelation_k \big) - 2 (\numobs - \numburn) \sum_{k =
    \numobs - \numburn}^{\infty} \crosscorrelation_k - 2 \sum_{k =
    1}^{\numobs - \numburn - 1} k \mu_k \\
  & \leq (\numobs - \numburn) \cdot \mathrm{Tr} \big( (I -
  \Lbar)^{-1} \SigStar_{\mathrm{Mkv}} (I - \Lbar)^{-1} \big) +
  \frac{3 \tau^2 \varbound^2 \usedim}{(1 - \kappa)^2}\opnorm{\testMat},
\end{align*}
which completes the proof of the lemma.

\paragraph{Proof of equation~\eqref{eq:cross-correlation-decay}}

Let $\state_0 \sim \stationary$ and $(\state_t)_{t \geq 0}$ be a
stationary Markov chain starting from $\state_0$.  By the construction
given in Section~\ref{SecConstruct}, there exists a random variable
$\widetilde{\state}_k$, such that $\widetilde{\state}_k$ is
independent of $\state_0$, $\widetilde{\state}_k \sim \stationary$,
and such that $\Exs \big[ \metric( \state_k, \widetilde{\state}_k)
  \mid \state_0 \big] \leq c_0 \cdot 2^{1 - \frac{k}{\mixingtime}}$.
We then obtain the bound
\begin{align}
 |\crosscorrelation_k| & = \abss{\Exs \big[ \inprod{\testMat^{1/2} (I
       - \Lbar)^{-1} \addnoiseMarkov (\state_0)}{ \testMat^{1/2} (I -
       \Lbar)^{-1} \addnoiseMarkov (\state_k)} \big]} \nonumber\\
 & \leq \abss{\Exs \big[ \inprod{\testMat^{1/2} (I - \Lbar)^{-1}
       \addnoiseMarkov (\state_0)}{\Exs \big[\testMat^{1/2} (I -
         \Lbar)^{-1} \addnoiseMarkov (\widetilde{\state}_k) \mid
         \state_0 \big]} \big]} \nonumber\\
 & \qquad + \abss{\Exs \big[\testMat^{1/2} \inprod{(I - \Lbar)^{-1}
       \addnoiseMarkov (\state_0)}{\Exs \big[\testMat^{1/2} (I -
         \Lbar)^{-1} \big( \addnoiseMarkov ({\state}_k) -
         \addnoiseMarkov (\widetilde{\state}_k) \big) \mid \state_0
         \big]} \big]} \nonumber \\
 & \leq 0 + \sqrt{\Exs \big[ \vecnorm{\testMat^{1/2} (I - \Lbar)^{-1}
       \addnoiseMarkov (\state_0)}{2}^2 \big]} \cdot \sqrt{ \Exs \big[
     \vecnorm{\testMat^{1/2} (I - \Lbar)^{-1} \big( \addnoiseMarkov
       ({\state}_k) - \addnoiseMarkov (\widetilde{\state}_k)
       \big)}{2}^2 \big]} \nonumber\\
 & \leq \sqrt{\crosscorrelation_0 \opnorm{Q}} \cdot \frac{\sqrt{\opnorm{Q}}}{1 - \kappa}
 \sqrt{\Exs \big[ \metric (\state_k, \widetilde{\state}_k)^2 \cdot
    \varbound^2 d^2 \big]}
 \nonumber\\
\label{EqnCoconutBall}
& \leq c_0 \frac{\varbound \usedim}{1 - \kappa}
     \sqrt{\crosscorrelation_0} \cdot 2^{1 - \frac{k}{2\mixingtime}}.
\end{align}
On the other hand, applying the moment
condition~\eqref{assume:noise-moments} yields $\crosscorrelation_0
\leq \frac{1}{(1 - \kappa)^2} \cdot \Exs \big[
  \vecnorm{\addnoiseMarkov (s_0)}{\testMat}^2 \big] \leq
\frac{\varbound^2 d}{(1 - \kappa)^2} \opnorm{\testMat}$.  Substituting
this bound into our previous inequality~\eqref{EqnCoconutBall}
completes the proof.


\subsubsection{Proof of Lemma~\ref{lemma:bound-sum-markov-mult-noise}}
\label{AppMultNoise}

The proof of this claim relies on a bootstrap argument: we bound the
summation of interest by a more complicated summation that involves
products of noise matrices. Recursively applying the result for $m =
\log \usedim$ times yields the desired bound.

\begin{lemma}
\label{lemma:bound-sum-markov-mult-noise-bootstrap}
 Given any integer $m \geq 0$, deterministic sequence $0 = k_0 < k_1 <
 \cdots < k_m < \numobs_0$, and scalar $\tau \geq 3 m \mixingtime p
 \log (c_0 d \numobs)$, we have the second moment bound
\begin{subequations}    
  \begin{align}
    &\Exs \big[ \vecnorm{\sum_{t = \numburn}^{\numobs - 1}
        \big(\prod_{j = 0}^{m} \multnoiseMarkov_{t - k_j} \big)
        \Delta_{t - k_{m}} }{2}^2 \big] \nonumber\\ &\leq 2 \numobs^2
    d^{2m} \sigmaA^{2m + 2} \cdot \frac{c \stepsize }{1 - \kappa} d
    \mixingtime \varbound^2 + 4 \stepsize^2 \tau \sum_{k_{m + 1} = k_m
      + 1}^{k_m + \tau} \Exs \big[ \vecnorm{ \sum_{t =
          \numburn}^{\numobs} \big\{\prod_{j = 0}^{m + 1}
        \multnoiseMarkov_{t - k_j} \Delta_{t - k_{m + 1}} \big\}}{2}^2
      \big] \nonumber\\ &\qquad + 4 \stepsize^2 \tau \sum_{k_{m + 1} =
      k_m + 1}^{k_m + \tau} \Exs \big[ \vecnorm{ \sum_{t =
          \numburn}^{\numobs} \big\{\prod_{j = 0}^{m}
        \multnoiseMarkov_{t - k_j} \big( \addnoiseMarkov_{t - k_{m +
            1}} + \addnoiseMG_{t - k_{m + 1} + 1} \big) \big\}}{2}^2
      \big],\label{eq:bound-sum-markov-mult-noise-bootstrap-induction}
    \end{align}
    and in the special case $m = 0$, we have
    \begin{align}
    &\Exs \big[ \vecnorm{\sum_{t = \numburn}^{\numobs - 1}
        \multnoiseMarkov_{t} \Delta_{t}}{2}^2 \big] \nonumber\\
         &\leq c \sigmaA^2
    d \cdot \big( \numobs \tau + \numobs^2 \stepsize^2 \sigmaA^2 d
    \tau^2 \big) \frac{c \stepsize }{1 - \kappa} d \mixingtime
    \varbound^2 + 4 \stepsize^2 \tau \sum_{k_1 = 1}^{\tau} \Exs \big[
      \vecnorm{ \sum_{t = \numburn}^{\numobs} \multnoiseMarkov_{t}
        \multnoiseMarkov_{t - k_1} \Delta_{t - k_{1}}}{2}^2 \big]
    \nonumber\\ &\qquad + 4 \stepsize^2 \tau \sum_{k_{1} = 1}^{ \tau}
    \Exs \big[ \vecnorm{ \sum_{t = \numburn}^{\numobs}
        \multnoiseMarkov_{t} \big( \addnoiseMarkov_{t - k_{1}} +
        \addnoiseMG_{t - k_{1} + 1} \big) }{2}^2
      \big].\label{eq:bound-sum-markov-mult-noise-bootstrap-base}
    \end{align}
\end{subequations}    
\end{lemma}
\noindent See
Appendix~\ref{subsubsec:proof-lemma-bound-sum-markov-mult-noise-bootstrap}
for the proof of this lemma.

The following lemma controls the last term of the
bound~\eqref{eq:bound-sum-markov-mult-noise-bootstrap-induction}:
\begin{lemma}\label{lemma:partial-sum-prod-intricate-structure}
 Under the setup above, there exists a universal constant $c > 0$,
 such that for any integer $m > 0$ and deterministic sequence $0 = k_0
 < k_1 < \cdots < k_m < \numburn$, we have:
    \begin{align*}
      \Exs \big[ \vecnorm{\sum_{t = \numburn}^{\numobs - 1} \big(
          \prod_{j = 0}^{m - 1} \multnoiseMarkov_{t - k_j} \big) \big(
          \addnoiseMarkov_{t - k_m} + \addnoiseMG_{t - k_m + 1}
          \big)}{2}^2 \big] \leq c \big( \numobs^2 + \numobs \usedim
      (k_m + \mixingtime \log (c_0 d)) \big) \sigmaA^{2m}
      d^{2m}\varbound^2.
    \end{align*}
\end{lemma}
\noindent See Appendix~\ref{subsubsec:proof-lemma-intricate-structure}
for the proof of this lemma. \\

\vspace*{0.05in}

Taking these lemmas as given, we now proceed with the proof of
Lemma~\ref{lemma:bound-sum-markov-mult-noise}.  Given the scalar $\tau
\mydefn 3 \mixingtime \log^2 (c_0 d n)$, we define
\begin{align*}
\bigquantity_m \mydefn \sup_{0 = k_0 < k_1 < \cdots < k_m \leq
  \tau}\Exs \big[ \vecnorm{\sum_{t = \numburn}^{\numobs - 1}
    \big(\prod_{j = 0}^{m} \multnoiseMarkov_{t - k_j} \big) \Delta_{t
      - k_{m}} }{2}^2 \big]
\end{align*}
for $m = 0, 1, 2, \cdots, \log d$.  By
equation~\eqref{eq:bound-sum-markov-mult-noise-bootstrap-base} and
Lemma~\ref{lemma:partial-sum-prod-intricate-structure}, we have the
bound
\begin{align*}
\bigquantity_0 & \leq c \sigmaA^2 d \cdot \big( \numobs \tau +
\numobs^2 \stepsize^2 \sigmaA^2 d \tau^2 \big) \frac{c \stepsize }{1 -
  \kappa} d \mixingtime \varbound^2\\
  &\qquad + 4 \stepsize^2 \tau^2
\bigquantity_1 + 4 c \stepsize^2 \tau^2 \big( \numobs^2 + \numobs
\usedim (\tau + \mixingtime \log (c_0 d)) \big) \sigmaA^2
d^2\varbound^2 \\
& \leq 4 \stepsize^2 \tau^2 \bigquantity_1 + c' \stepsize^2 \numobs^2
\tau^2 \usedim^2 \sigmaA^2 \varbound^2 .
\end{align*}
In deriving the last inequality, we used the inequalities $\stepsize
\leq \frac{1 - \kappa}{ \sigmaA^2 d \tau}$ and $\numobs \geq
\frac{1}{(1 - \kappa) \stepsize}$.

By equation~\eqref{eq:bound-sum-markov-mult-noise-bootstrap-induction}
and Lemma~\ref{lemma:partial-sum-prod-intricate-structure}, we have
the recursive relation
\begin{align*}
\bigquantity_m & \leq 4 \stepsize^2 \tau^2 \bigquantity_{m + 1} + c
\numobs^2 d^{2m + 1} \tau \sigmaA^{2m + 2} \cdot \frac{\stepsize
  \log^3 n}{1 - \kappa} \varbound^2 + c \stepsize^2 \tau^2 \numobs^2
\sigmaA^{2m + 2} d^{2m + 2} \varbound^2\\ &\leq 4 \stepsize^2 \tau^2
\bigquantity_{m + 1} + c\numobs^2 \sigmaA^{2m} d^{2m} \varbound^2
\cdot \log^3 \numobs.
\end{align*}
Recursively applying these bounds yields
\begin{align*}
\bigquantity_0 & \leq (4 \stepsize^2 \tau^2)^m \bigquantity_m + c
\stepsize^2 \numobs^2 \tau^2 \usedim^2 \sigmaA^2 \varbound^2 + c \cdot
\sum_{q = 1}^{m - 1} (4 \stepsize^2 \tau^2)^q \numobs^2 \sigmaA^{2q}
d^{2q} \varbound^2 \\
& \leq (4 \stepsize^2 \tau^2)^m \bigquantity_m + 3
c \stepsize^2 \numobs^2 \tau^2 \usedim^2 \sigmaA^2 \varbound^2 .
\end{align*}
In order to control the term $\bigquantity_m$, we employ the coarse
bound
\begin{align*}
\Exs \big[ \vecnorm{\sum_{t = \numburn}^{\numobs - 1} \big(\prod_{j =
      0}^{m} \multnoiseMarkov_{t - k_j} \big) \Delta_{t - k_{m}}
  }{2}^2 \big] &\leq \numobs \sum_{t = \numburn}^{\numobs - 1} \Exs
\big[ \vecnorm{\big(\prod_{j = 0}^{m} \multnoiseMarkov_{t - k_j} \big)
    \Delta_{t - k_{m}} }{2}^2 \big]\\
    & \leq \numobs^2 (\sigmaA d)^{2m +
  2} \cdot \frac{c \stepsize \mixingtime d\varbound^2 }{1 - \kappa} .
\end{align*}
Taking the supremum and noting that $\stepsize \leq \frac{1 -
  \kappa}{\sigmaA^2 \usedim \tau}$ leads to $\bigquantity_m \leq c
\numobs^2 \sigmaA^{2m} d^{2m + 2} \varbound^2$.  Consequently, we have
established that $\bigquantity_0 \leq 3 c \stepsize^2 \numobs^2 \tau^2
\usedim^2 \sigmaA^2 \varbound^2 \big[1 + \big(2 \stepsize \tau \sigmaA
  \usedim\big)^{\frac{2m + 2}{2 m}} \big]$.  Taking \mbox{$m = \lceil
  \log d \rceil$} and $\stepsize \leq \frac{1}{6 \tau \sigmaA
  \usedim}$, we have $(2 \stepsize \tau \sigmaA \usedim^{\frac{2m +
    2}{2 m}})^{2m} < 1$, and thus $\bigquantity_0 \leq 6 c \stepsize^2
\numobs^2 \tau^2 \usedim^2 \sigmaA^2 \varbound^2 \log^3 \numobs$,
which completes the proof of this lemma.


\section{Discussion}
\label{sec:discussion}

In this paper, we established sharp instance-optimal guarantees for
linear stochastic approximation (SA) procedures based on Markovian
data.  Under ergodicity along with natural tail conditions, we proved
non-asymptotic upper bounds on the squared error of both the last
iterate of a standard SA scheme, as well as the Polyak--Ruppert
averaged sequence. The results highlight two important aspects: an
optimal sample complexity of $O(\mixingtime \usedim)$ for problems in
dimension $\usedim$ with mixing time $\mixingtime$; and an
instance-dependent error upper bound for the averaged estimator with
carefully chosen stepsize. Complementary to the upper bound, we also
showed a non-asymptotic local minimax lower bound over a small
neighborhood of a given Markov chain instance, certifying the
statistical optimality of the proposed estimators. Our proof of the
upper bounds uses a bootstrapping argument of possibly independent
interest.

Throughout the paper, we have introduced novel techniques of analysis
and motivated several open questions. In the following, we collect a
few interesting future directions:
\bcar
\item \textbf{Nonlinear stochastic approximation and controlled
  dynamics:} Our paper focuses on linear $Z$-equations where the
  underlying Markov chain does not involve a control. Though this
  setting already covers many important examples (as described in
  Section~\ref{subsec:examples}), its applicability to practical
  problems is still relatively restricted. To set up a general
  framework, one could consider a \emph{controlled Markov chain}
  $(s_t)_{t \geq 0}$ where the transition is given by $s_{t + 1} \sim
  \transition (\cdot | s_t, \theta_t)$. For any $\theta \in \real^d$,
  let $\stationary_\theta$ be the stationary distribution of the
  Markov chain $\transition (\cdot | \cdot , \theta)$ induced by the
  control $\theta$. Given a non-linear operator $H: \statespace \times
  \real^d \rightarrow \real^d$, suppose that we wish to solve the
  equation $\Exs_{s \sim \stationary (\theta)} \big[ H(\theta; s)
    \big] = 0$; see the book~\cite{benveniste2012adaptive} for a
  summary of classical asymptotic theory for such problems. The
  analysis tools introduced in this paper provide an avenue by which
  one could obtain optimal sample complexity bounds (especially in
  terms of dimension dependency) and instance-dependent guarantees for
  such problems. In particular, the multi-step looking-back technique
  and bootstrapping stability bounds introduced in
  Proposition~\ref{prop:lsa-markov-iterate-bound} could be extended to
  nonlinear operators, and it would be very interesting to see how
  Markovian SA achieves optimal dependence on $(\mixingtime, \usedim)$
  in general. On the other hand, the proof of
  Theorem~\ref{thm:markov-main} is specialized to linear operators, as
  it explicitly involves bounding product of random matrices (see
  Lemma~\ref{lemma:bound-sum-markov-mult-noise-bootstrap}). Obtaining
  sharp and instance-dependent results for the non-linear SA may
  require novel proof techniques, and is an important direction of
  future work.\\
\item \textbf{Online statistical inference:} By carefully choosing the
  burn-in period, one can show that the Polyak--Ruppert estimator
  $\thetahat_\numobs$ is asymptotically normal and locally minimax
  optimal. In particular, under suitable conditions, we have the
  following limiting result (see the paper~\cite{fort2015central} for
  details):
    \begin{align}
        \sqrt{\numobs} (\thetahat_\numobs - \thetabar) \xrightarrow{d}
        \mathcal{N} \big( (I_d - \Lbar)^{-1} (\SigStar_{\mathrm{MG}} +
        \SigStar_{\mathrm{Mkv}})(I_d - \Lbar)^{-\top}
        \big).\label{eq:clt-for-sa-procedure}
    \end{align}
    In order to construct confidence intervals for the solution
    $\thetabar$ with streaming data, it suffices to estimate the
    asymptotic covariance in
    equation~\eqref{eq:clt-for-sa-procedure}. In the $\myiid$ setting,
    online procedures have been developed to estimate such
    covariances, with non-asymptotic error
    guarantees~\cite{chen2020statistical}. The problem becomes more
    subtle in the Markovian setting, as the matrix
    $\SigStar_{\mathrm{Mkv}}$ involves auto-correlations of the noise
    process. It is an important open direction to construct online
    estimators of this matrix to enable inference in a streaming
    fashion. \\
\item \textbf{Model selection and optimal methods for policy
  evaluation} The policy evaluation problem involves manual
  choice of two important parameters: the feature vector dimension
  $\usedim$ and the resolvent parameter $\lambda$ in TD$(\lambda)$. In
  Section~\ref{subsubsec:sieve-consequence}
  and~\ref{subsec:consequence-td-lambda}, we provide optimal
  instance-dependent guarantees on both the approximation factor and
  the estimation error, for a fixed choice of $d$ and $\lambda$. An
  important direction of future research is to select such parameters
  adaptively based on data, possibly under a streaming computational
  model. Ideally, we want the risk of such estimator to attain the
  infimum of the right hand side of
  equation~\eqref{eq:td-lambda-final-bound}, over $\lambda \in (0, 1)$
  and $d \in \mathbb{N}_+$. A possible candidate approach towards such
  a model selection problem is the celebrated Lepskii method for
  adaptive bandwidth selection~\cite{lepskii1991problem}.
  \ecar.



\section*{Acknowledgements}
We thank Yaqi Duan for helpful discussions.
AP was partially supported by the National Science Foundation through grants 2107455 and 2210734, and by awards/gifts from Adobe, Amazon, and Mathworks.
MJW and PLB were partially supported by the NSF through grants IIS-1909365 and DMS-2023505. PLB was partially supported by the ONR through MURI award N000142112431.
This work was partially supported by NSF grant CCF-1955450, ONR grant
N00014-21-1-2842, and NSF grant DMS-2311072 to MJW.

\bibliographystyle{alpha}
\bibliography{reference}

\newcommand{\etalchar}[1]{$^{#1}$}
\begin{thebibliography}{KMMW19}

\bibitem[BB96]{bradtke1996linear}
Steven~J Bradtke and Andrew~G Barto.
\newblock Linear least-squares algorithms for temporal difference learning.
\newblock {\em Machine Learning}, 22(1-3):33--57, 1996.

\bibitem[BCN18]{bottou2018optimization}
L{\'e}on Bottou, Frank~E Curtis, and Jorge Nocedal.
\newblock Optimization methods for large-scale machine learning.
\newblock {\em SIAM Review}, 60(2):223--311, 2018.

\bibitem[BD09]{brockwell2009time}
Peter~J Brockwell and Richard~A Davis.
\newblock {\em Time series: theory and methods}.
\newblock Springer Science \& Business Media, 2009.

\bibitem[Ber11]{bertsekas2011temporal}
Dimitri~P Bertsekas.
\newblock Temporal difference methods for general projected equations.
\newblock {\em IEEE Transactions on Automatic Control}, 56(9):2128--2139, 2011.

\bibitem[Ber19]{bertsekas2019reinforcement}
Dimitri~P Bertsekas.
\newblock {\em Reinforcement Learning and Optimal Control}.
\newblock Athena Scientific Belmont, MA, 2019.

\bibitem[Bil61]{billingsley1961statistical}
Patrick Billingsley.
\newblock Statistical methods in {M}arkov chains.
\newblock {\em The Annals of Mathematical Statistics}, pages 12--40, 1961.

\bibitem[BJN{\etalchar{+}}20]{bresler2020least}
Guy Bresler, Prateek Jain, Dheeraj Nagaraj, Praneeth Netrapalli, and Xian Wu.
\newblock Least squares regression with {M}arkovian data: Fundamental limits
  and algorithms.
\newblock {\em arXiv preprint arXiv:2006.08916}, 2020.

\bibitem[BMP12]{benveniste2012adaptive}
Albert Benveniste, Michel M{\'e}tivier, and Pierre Priouret.
\newblock {\em Adaptive Algorithms and Stochastic Approximations}, volume~22.
\newblock Springer Science \& Business Media, 2012.

\bibitem[Bor09]{borkar2009stochastic}
Vivek~S Borkar.
\newblock {\em Stochastic Approximation: a Dynamical Systems Viewpoint},
  volume~48.
\newblock Springer, 2009.

\bibitem[Bor21]{borkar2021concentration}
Vivek~S Borkar.
\newblock A concentration bound for contractive stochastic approximation.
\newblock {\em Systems \& Control Letters}, 153:104947, 2021.

\bibitem[Boy02]{boyan2002technical}
Justin~A Boyan.
\newblock Technical update: Least-squares temporal difference learning.
\newblock {\em Machine Learning}, 49(2-3):233--246, 2002.

\bibitem[BREZ20]{bou2020coupling}
Nawaf Bou-Rabee, Anreas Eberle, and Raphael Zimmer.
\newblock Coupling and convergence for {H}amiltonian {M}onte {C}arlo.
\newblock {\em The Annals of Applied Probability}, 30(3):1209--1250, 2020.

\bibitem[BRS18]{bhandari2018finite}
Jalaj Bhandari, Daniel Russo, and Raghav Singal.
\newblock A finite time analysis of temporal difference learning with linear
  function approximation.
\newblock {\em arXiv preprint arXiv:1806.02450}, 2018.

\bibitem[CDBM20]{chen2020explicit}
Shuhang Chen, Adithya Devraj, Ana Busic, and Sean Meyn.
\newblock Explicit mean-square error bounds for {M}onte-{C}arlo and linear
  stochastic approximation.
\newblock In {\em International Conference on Artificial Intelligence and
  Statistics}, pages 4173--4183. PMLR, 2020.

\bibitem[CLTZ20]{chen2020statistical}
Xi~Chen, Jason~D Lee, Xin~T Tong, and Yichen Zhang.
\newblock Statistical inference for model parameters in stochastic gradient
  descent.
\newblock {\em The Annals of Statistics}, 48(1):251--273, 2020.

\bibitem[CMSS21]{chen2021lyapunov}
Zaiwei Chen, Siva~Theja Maguluri, Sanjay Shakkottai, and Karthikeyan Shanmugam.
\newblock A {L}yapunov theory for finite-sample guarantees of asynchronous
  {Q}-learning and {TD}-learning variants.
\newblock {\em arXiv preprint arXiv:2102.01567}, 2021.

\bibitem[DDA21]{debavelaere2021convergence}
Vianney Debavelaere, Stanley Durrleman, and St{\'e}phanie Allassonni{\`e}re.
\newblock On the convergence of stochastic approximations under a subgeometric
  ergodic {M}arkov dynamic.
\newblock {\em Electronic Journal of Statistics}, 15(1):1583--1609, 2021.

\bibitem[DDB20]{dieuleveut2020bridging}
Aymeric Dieuleveut, Alain Durmus, and Francis Bach.
\newblock Bridging the gap between constant step size stochastic gradient
  descent and {M}arkov chains.
\newblock {\em The Annals of Statistics}, 48(3):1348--1382, 2020.

\bibitem[DDG18]{daskalakis2018testing}
Constantinos Daskalakis, Nishanth Dikkala, and Nick Gravin.
\newblock Testing symmetric {M}arkov chains from a single trajectory.
\newblock In {\em Conference On Learning Theory}, pages 385--409. PMLR, 2018.

\bibitem[DMN{\etalchar{+}}21]{durmus2021stability}
Alain Durmus, Eric Moulines, Alexey Naumov, Sergey Samsonov, and Hoi-To Wai.
\newblock On the stability of random matrix product with {M}arkovian noise:
  Application to linear stochastic approximation and {TD} learning.
\newblock {\em arXiv preprint arXiv:2102.00185}, 2021.

\bibitem[DNPR20]{doan2020finite}
Thinh~T Doan, Lam~M Nguyen, Nhan~H Pham, and Justin Romberg.
\newblock Finite-time analysis of stochastic gradient descent under {M}arkov
  randomness.
\newblock {\em arXiv preprint arXiv:2003.10973}, 2020.

\bibitem[DS94]{dayan1994td}
Peter Dayan and Terrence~J Sejnowski.
\newblock {TD}($\lambda$) converges with probability 1.
\newblock {\em Machine Learning}, 14(3):295--301, 1994.

\bibitem[DWW21]{duan2021optimal}
Yaqi Duan, Mengdi Wang, and Martin~J Wainwright.
\newblock Optimal policy evaluation using kernel-based temporal difference
  methods.
\newblock {\em arXiv preprint arXiv:2109.12002}, 2021.

\bibitem[EDM03]{EveMan03}
E.~Even-{D}ar and Y.~Mansour.
\newblock Learning rates for {$Q$}-learning.
\newblock {\em Journal of {M}achine {L}earning {R}esearch}, 5:1--25, 2003.

\bibitem[For15]{fort2015central}
Gersende Fort.
\newblock Central limit theorems for stochastic approximation with controlled
  {M}arkov chain dynamics.
\newblock {\em ESAIM: Probability and Statistics}, 19:60--80, 2015.

\bibitem[GL95]{gill1995applications}
Richard~D Gill and Boris~Y Levit.
\newblock Applications of the van {T}rees inequality: a {B}ayesian
  {C}ram{\'e}r-{R}ao bound.
\newblock {\em Bernoulli}, 1(1-2):59--79, 1995.

\bibitem[GL12]{ghadimi2012optimal}
Saeed Ghadimi and Guanghui Lan.
\newblock Optimal stochastic approximation algorithms for strongly convex
  stochastic composite optimization {I}: A generic algorithmic framework.
\newblock {\em SIAM Journal on Optimization}, 22(4):1469--1492, 2012.

\bibitem[GP17]{gadat2017optimal}
S{\'e}bastien Gadat and Fabien Panloup.
\newblock Optimal non-asymptotic bound of the {R}uppert--{P}olyak averaging
  without strong convexity.
\newblock {\em arXiv preprint arXiv:1709.03342}, 2017.

\bibitem[GSW80]{ghosh1980second}
JK~Ghosh, Bimal~K Sinha, and HS~Wieand.
\newblock Second order efficiency of the mle with respect to any bounded
  bowl-shape loss function.
\newblock {\em The Annals of Statistics}, pages 506--521, 1980.

\bibitem[GW95]{greenwood1995efficiency}
Priscilla~E Greenwood and Wolfgang Wefelmeyer.
\newblock Efficiency of empirical estimators for {M}arkov chains.
\newblock {\em The Annals of Statistics}, pages 132--143, 1995.

\bibitem[Ham20]{hamilton2020time}
James~Douglas Hamilton.
\newblock {\em Time series analysis}.
\newblock Princeton university press, 2020.

\bibitem[HKL{\etalchar{+}}19]{hsu2019mixing}
Daniel Hsu, Aryeh Kontorovich, David~A Levin, Yuval Peres, Csaba
  Szepesv{\'a}ri, and Geoffrey Wolfer.
\newblock Mixing time estimation in reversible {M}arkov chains from a single
  sample path.
\newblock {\em The Annals of Applied Probability}, 29(4):2439--2480, 2019.

\bibitem[JKNN21]{jain2021streaming}
Prateek Jain, Suhas~S Kowshik, Dheeraj Nagaraj, and Praneeth Netrapalli.
\newblock Streaming linear system identification with reverse experience
  replay.
\newblock {\em arXiv preprint arXiv:2103.05896}, 2021.

\bibitem[JZ13]{johnson2013accelerating}
Rie Johnson and Tong Zhang.
\newblock Accelerating stochastic gradient descent using predictive variance
  reduction.
\newblock {\em Advances in neural information processing systems}, 26:315--323,
  2013.

\bibitem[KB18]{karmakar2018two}
Prasenjit Karmakar and Shalabh Bhatnagar.
\newblock Two time-scale stochastic approximation with controlled {M}arkov
  noise and off-policy temporal-difference learning.
\newblock {\em Mathematics of Operations Research}, 43(1):130--151, 2018.

\bibitem[KC78]{kushnerclarks}
Harold~J. Kushner and Dean~S. Clark.
\newblock {\em Stochastic approximation methods for constrained and
  unconstrained systems}, volume~26 of {\em Applied Mathematical Sciences}.
\newblock Springer-Verlag, New York-Berlin, 1978.

\bibitem[KLL20]{kotsalis2020simple}
Georgios Kotsalis, Guanghui Lan, and Tianjiao Li.
\newblock Simple and optimal methods for stochastic variational inequalities,
  {II}: {M}arkovian noise and policy evaluation in reinforcement learning.
\newblock {\em arXiv preprint arXiv:2011.08434}, 2020.

\bibitem[KMMW19]{karimi2019non}
Belhal Karimi, Blazej Miasojedow, Eric Moulines, and Hoi-To Wai.
\newblock Non-asymptotic analysis of biased stochastic approximation scheme.
\newblock In {\em Conference on Learning Theory}, pages 1944--1974. PMLR, 2019.

\bibitem[KMN{\etalchar{+}}20]{kaledin2020finite}
Maxim Kaledin, Eric Moulines, Alexey Naumov, Vladislav Tadic, and Hoi-To Wai.
\newblock Finite time analysis of linear two-timescale stochastic approximation
  with {M}arkovian noise.
\newblock In {\em Conference on Learning Theory}, pages 2144--2203. PMLR, 2020.

\bibitem[KPR{\etalchar{+}}21]{khamaru2020temporal}
Koulik Khamaru, Ashwin Pananjady, Feng Ruan, Martin~J Wainwright, and Michael~I
  Jordan.
\newblock Is temporal difference learning optimal? {A}n instance-dependent
  analysis.
\newblock {\em SIAM Journal on Mathematics of Data Science}, 3(4):1013--1040,
  2021.

\bibitem[KT00]{konda2000actor}
Vijay~R Konda and John~N Tsitsiklis.
\newblock Actor-critic algorithms.
\newblock In {\em Advances in Neural Information Processing Systems}, pages
  1008--1014, 2000.

\bibitem[Kut97]{kutoyants1997efficiency}
Yury~A Kutoyants.
\newblock Efficiency of the empirical distribution for ergodic diffusion.
\newblock {\em Bernoulli}, pages 445--456, 1997.

\bibitem[KXWJ21]{khamaru2021instance}
Koulik Khamaru, Eric Xia, Martin~J Wainwright, and Michael~I Jordan.
\newblock Instance-optimality in optimal value estimation: Adaptivity via
  variance-reduced {Q}-learning.
\newblock {\em arXiv preprint arXiv:2106.14352}, 2021.

\bibitem[KY03]{kushner2003stochastic}
Harold Kushner and G~George Yin.
\newblock {\em Stochastic approximation and recursive algorithms and
  applications}, volume~35.
\newblock Springer Science \& Business Media, 2003.

\bibitem[L\"05]{Lut05}
H.~L\"{u}tkepohl.
\newblock {\em New introduction to multiple time series analysis}.
\newblock Springer, New York, 2005.

\bibitem[Lep91]{lepskii1991problem}
OV~Lepskii.
\newblock On a problem of adaptive estimation in {G}aussian white noise.
\newblock {\em Theory of Probability \& Its Applications}, 35(3):454--466,
  1991.

\bibitem[Lju77a]{ljung1977analysis}
Lennart Ljung.
\newblock Analysis of recursive stochastic algorithms.
\newblock {\em IEEE Transactions on Automatic Control}, 22(4):551--575, 1977.

\bibitem[Lju77b]{ljung1977positive}
Lennart Ljung.
\newblock On positive real transfer functions and the convergence of some
  recursive schemes.
\newblock {\em IEEE Transactions on Automatic Control}, 22(4):539--551, 1977.

\bibitem[LLP21]{li2021accelerated}
Tianjiao Li, Guanghui Lan, and Ashwin Pananjady.
\newblock Accelerated and instance-optimal policy evaluation with linear
  function approximation.
\newblock {\em preprint}, 2021.

\bibitem[LMWJ20]{li2020root}
Chris~Junchi Li, Wenlong Mou, Martin~J Wainwright, and Michael~I Jordan.
\newblock {ROOT-SGD}: Sharp nonasymptotics and asymptotic efficiency in a
  single algorithm.
\newblock {\em arXiv preprint arXiv:2008.12690}, 2020.

\bibitem[LP17]{levin2017markov}
David~A Levin and Yuval Peres.
\newblock {\em Markov Chains and Mixing Times}, volume 107.
\newblock American Mathematical Soc., 2017.

\bibitem[LS18]{lakshminarayanan2018linear}
Chandrashekar Lakshminarayanan and Csaba Szepesv\'{a}ri.
\newblock Linear stochastic approximation: {H}ow far does constant step-size
  and iterate averaging go?
\newblock In {\em International Conference on Artificial Intelligence and
  Statistics}, pages 1347--1355, 2018.

\bibitem[LWC{\etalchar{+}}20]{li2020breaking}
Gen Li, Yuting Wei, Yuejie Chi, Yuantao Gu, and Yuxin Chen.
\newblock Breaking the sample size barrier in model-based reinforcement
  learning with a generative model.
\newblock In {\em Advances in Neural Information Processing Systems},
  volume~33, pages 12861--12872. Curran Associates, Inc., 2020.

\bibitem[LWZ18]{li2018estimation}
Xudong Li, Mengdi Wang, and Anru Zhang.
\newblock Estimation of {M}arkov chain via rank-constrained likelihood.
\newblock In {\em International Conference on Machine Learning}, pages
  3033--3042. PMLR, 2018.

\bibitem[MB11]{moulines2011non}
\'{E}ric Moulines and Francis~R Bach.
\newblock Non-asymptotic analysis of stochastic approximation algorithms for
  machine learning.
\newblock In {\em Advances in Neural Information Processing Systems}, pages
  451--459, 2011.

\bibitem[MLW{\etalchar{+}}20]{mou2020linear}
Wenlong Mou, Chris~Junchi Li, Martin~J Wainwright, Peter~L Bartlett, and
  Michael~I Jordan.
\newblock On linear stochastic approximation: {F}ine-grained {P}olyak-{R}uppert
  and non-asymptotic concentration.
\newblock In {\em Proceedings of Thirty Third Conference on Learning Theory},
  volume 125, pages 2947--2997, 2020.

\bibitem[MP84]{metivier1984applications}
Michel Metivier and Pierre Priouret.
\newblock Applications of a {K}ushner and {C}lark lemma to general classes of
  stochastic algorithms.
\newblock {\em IEEE Transactions on Information Theory}, 30(2):140--151, 1984.

\bibitem[MPW20]{mou2020optimal}
Wenlong Mou, Ashwin Pananjady, and Martin~J Wainwright.
\newblock Optimal oracle inequalities for solving projected fixed-point
  equations.
\newblock {\em arXiv preprint arXiv:2012.05299}, 2020.

\bibitem[MS08]{munos2008finite}
R{\'e}mi Munos and Csaba Szepesv{\'a}ri.
\newblock Finite-time bounds for fitted value iteration.
\newblock {\em Journal of Machine Learning Research}, 9(May):815--857, 2008.

\bibitem[MT12]{meyn2012markov}
Sean~P Meyn and Richard~L Tweedie.
\newblock {\em Markov chains and stochastic stability}.
\newblock Springer Science \& Business Media, 2012.

\bibitem[Nem01]{nemirovski2001lectures}
Arkadi Nemirovski.
\newblock Lectures on modern convex optimization.
\newblock In {\em Society for Industrial and Applied Mathematics (SIAM}.
  Citeseer, 2001.

\bibitem[NJLS09]{nemirovski2009robust}
Arkadi Nemirovski, Anatoli Juditsky, Guanghui Lan, and Alexander Shapiro.
\newblock Robust stochastic approximation approach to stochastic programming.
\newblock {\em SIAM Journal on Optimization}, 19(4):1574--1609, 2009.

\bibitem[Pen91]{penev1991efficient}
Spiridon Penev.
\newblock Efficient estimation of the stationary distribution for exponentially
  ergodic {M}arkov chains.
\newblock {\em Journal of Statistical Planning and Inference}, 27(1):105--123,
  1991.

\bibitem[PJ92]{polyak1992acceleration}
Boris~T Polyak and Anatoli~B Juditsky.
\newblock Acceleration of stochastic approximation by averaging.
\newblock {\em SIAM Journal on Control and Optimization}, 30(4):838--855, 1992.

\bibitem[PW21]{pananjady2020instance}
Ashwin Pananjady and Martin~J. Wainwright.
\newblock Instance-dependent $\ell_\infty$-bounds for policy evaluation in
  tabular reinforcement learning.
\newblock {\em IEEE Transactions on Information Theory}, 67(1):566--585, 2021.

\bibitem[Rao62]{rao1962first}
C.~R. Rao.
\newblock First and second order asymptotic efficiencies of estimators.
\newblock {\em Annales scientifiques de l'Universit{\'e} de Clermont.
  Math{\'e}matiques}, 8(2):33--40, 1962.

\bibitem[RM51]{robbins1951stochastic}
Herbert Robbins and Sutton Monro.
\newblock A stochastic approximation method.
\newblock {\em The Annals of Mathematical Statistics}, pages 400--407, 1951.

\bibitem[RN94]{Rummery1994line}
Gavin~A Rummery and Mahesan Niranjan.
\newblock On-line {Q}-learning using connectionist systems.
\newblock Technical report, Cambridge University Engineering Department, 1994.

\bibitem[Rup88]{ruppert1988efficient}
David Ruppert.
\newblock Efficient estimations from a slowly convergent {R}obbins-{M}onro
  process.
\newblock Technical report, Cornell University Operations Research and
  Industrial Engineering, 1988.

\bibitem[Sut88]{sutton1988learning}
Richard~S Sutton.
\newblock Learning to predict by the methods of temporal differences.
\newblock {\em Machine Learning}, 3(1):9--44, 1988.

\bibitem[SWW{\etalchar{+}}18]{sidford2018near}
Aaron Sidford, Mengdi Wang, Xian Wu, Lin~F Yang, and Yinyu Ye.
\newblock Near-optimal time and sample complexities for solving {M}arkov
  decision processes with a generative model.
\newblock In {\em Proceedings of the 32nd International Conference on Neural
  Information Processing Systems}, pages 5192--5202, 2018.

\bibitem[SY19]{srikant2019finite}
Rayadurgam Srikant and Lei Ying.
\newblock Finite-time error bounds for linear stochastic approximation and {TD}
  learning.
\newblock In {\em Conference on Learning Theory}, pages 2803--2830. PMLR, 2019.

\bibitem[Sze98]{szepesvari1998asymptotic}
Csaba Szepesv{\'a}ri.
\newblock The asymptotic convergence-rate of {Q}-learning.
\newblock {\em Advances in neural information processing systems}, pages
  1064--1070, 1998.

\bibitem[Sze10]{szepesvari2010algorithms}
Csaba Szepesv\'{a}ri.
\newblock {\em Algorithms for Reinforcement Learning}.
\newblock Morgan \& Claypool Publishers, 2010.

\bibitem[Tsi94]{Tsi94}
J.~N. Tsitsiklis.
\newblock Asynchronous stochastic approximation and {$Q$}-learning.
\newblock {\em Machine Learning}, 16:185--202, 1994.

\bibitem[Tsy08]{tsybakov2008introduction}
Alexandre~B Tsybakov.
\newblock {\em Introduction to Nonparametric Estimation}.
\newblock Springer Science \& Business Media, 2008.

\bibitem[TVR97]{tsitsiklis1997analysis}
John~N Tsitsiklis and Benjamin Van~Roy.
\newblock Analysis of temporal-diffference learning with function
  approximation.
\newblock In {\em Advances in Neural Information Processing Systems}, pages
  1075--1081, 1997.

\bibitem[TVR99]{tsitsiklis1999optimal}
John~N Tsitsiklis and Benjamin Van~Roy.
\newblock Optimal stopping of {M}arkov processes: {H}ilbert space theory,
  approximation algorithms, and an application to pricing high-dimensional
  financial derivatives.
\newblock {\em IEEE Transactions on Automatic Control}, 44(10):1840--1851,
  1999.

\bibitem[vdV00]{van2000asymptotic}
Aad~W {v}an~der Vaart.
\newblock {\em Asymptotic Statistics}, volume~3.
\newblock Cambridge university press, 2000.

\bibitem[Wai19a]{wainwright2019high}
Martin~J Wainwright.
\newblock {\em High-dimensional Statistics: A Non-asymptotic Viewpoint},
  volume~48.
\newblock Cambridge University Press, 2019.

\bibitem[Wai19b]{wainwright2019stochastic}
Martin~J Wainwright.
\newblock Stochastic approximation with cone-contractive operators: Sharper
  $\ell_\infty $-bounds for {Q}-learning.
\newblock {\em arXiv preprint arXiv:1905.06265}, 2019.

\bibitem[Wai19c]{wainwright2019variance}
Martin~J Wainwright.
\newblock Variance-reduced {Q}-learning is minimax optimal.
\newblock {\em arXiv preprint arXiv:1906.04697}, 2019.

\bibitem[WD92]{watkins1992q}
Christopher~JCH Watkins and Peter Dayan.
\newblock Q-learning.
\newblock {\em Machine Learning}, 8(3-4):279--292, 1992.

\bibitem[WK21]{wolfer2021statistical}
Geoffrey Wolfer and Aryeh Kontorovich.
\newblock Statistical estimation of ergodic {M}arkov chain kernel over discrete
  state space.
\newblock {\em Bernoulli}, 27(1):532--553, 2021.

\bibitem[YB10]{yu2010error}
Huizhen Yu and Dimitri~P Bertsekas.
\newblock Error bounds for approximations from projected linear equations.
\newblock {\em Mathematics of Operations Research}, 35(2):306--329, 2010.

\bibitem[YBVE20]{yu2020analysis}
Lu~Yu, Krishnakumar Balasubramanian, Stanislav Volgushev, and Murat~A Erdogdu.
\newblock An analysis of constant step size {SGD} in the non-convex regime:
  Asymptotic normality and bias.
\newblock {\em arXiv preprint arXiv:2006.07904}, 2020.

\bibitem[YBW17]{yang2017statistical}
Fanny Yang, Sivaraman Balakrishnan, and Martin~J Wainwright.
\newblock Statistical and computational guarantees for the {B}aum--{W}elch
  algorithm.
\newblock {\em The Journal of Machine Learning Research}, 18(1):4528--4580,
  2017.

\end{thebibliography}

\appendix
\vspace{0.4cm}
\setcounter{section}{0}
\def\thesection{\Alph{section}}

\section{Additional related work}
\label{subsec:additional-related-works}

This paper analyzes stochastic approximation algorithms based on
Markov data, and has consequences for reinforcement learning.  So as
to put our results into context, we now provide more background on
past work in these areas.

\subsection{Statistical estimation based on Markov data}

There is a large body of past work on statistical estimation based on
observing a single trajectory of a Markov chain; for example, see~\cite{billingsley1961statistical} for an overview of some
classical results. For the problem of functional estimation under the
stationary distribution, the asymptotic efficiency of plug-in
estimators\footnote{These papers refer to such methods as
``empirical'' estimators.} has been established for discrete-state
Markov chains~\cite{penev1991efficient,greenwood1995efficiency} and
It\^{o} diffusion processes~\cite{kutoyants1997efficiency}. In
this paper, we provide non-asymptotic bounds, both upper and lower,
that depend on a certain instance-dependent functional that also
appears in an asymptotic analysis.  More recent work has seen
non-asymptotic results for statistical estimation with Markovian data,
including the estimation of transition
kernels~\cite{wolfer2021statistical,li2018estimation}, mixing
times~\cite{hsu2019mixing}, the parameters of Gaussian hidden Markov
models~\cite{yang2017statistical}, as well for certain testing
problems~\cite{daskalakis2018testing}.  These papers can be roughly
divided into two categories. Papers in the first category focus on
estimating parameters for each individual state of the Markov chain
(e.g., transition kernels), and thus require sample sizes that scale
with the complexity of the state space (e.g., its cardinality in the
discrete case).  By contrast, papers in the second category are
concerned with estimating properties of the Markov chain (e.g., the
expectation of a functional under the stationary distribution), and
the sample complexity of such problems need not depend on the size of
the state space. Our paper falls within the second category.

\subsection{Stochastic approximation methods}

The use of recursive stochastic procedures for solving fixed point
equations dates back to the seminal work of Robbins and Monro~\cite{robbins1951stochastic}; see the reference
books~\cite{borkar2009stochastic,benveniste2012adaptive,kushner2003stochastic}
for more background.  By averaging the iterates of the SA procedure,
it is known that one can obtain both an improved convergence rate and
central limit
behavior~\cite{polyak1992acceleration,ruppert1988efficient}.  A
variety of stochastic approximation procedures now serve as the
workhorse for modern large-scale machine learning and statistical
inference~\cite{nemirovski2009robust,bottou2018optimization}, and many
algorithmic techniques are known to accelerate their
convergence~\cite{ghadimi2012optimal,johnson2013accelerating,li2020root}. In
particular, non-asymptotic bounds matching the optimal Gaussian limit
have been established in a variety of
settings~\cite{moulines2011non,gadat2017optimal,dieuleveut2020bridging,mou2020linear,mou2020optimal}.

While the instance-dependent nature of this line of investigation
aligns with the objective of our work, prior work either assumes an
$\myiid$ observation model or imposes a martingale difference
assumption on the noise.\footnote{In the linear equation setup, the
martingale difference noise assumes that \mbox{$\Exs[ \Lmat_{t + 1}
    \mid \filtration_t ] = \Lbar$} and \mbox{$\Exs [ \bvec_{t + 1}
    \mid \filtration_t ] = \bbar$,} which does not cover the Markov
case.} The first study of SA procedures without a martingale
difference assumption was initiated by~\cite{kushnerclarks}, who give a general criteria for
convergence, as well as
\cite{ljung1977analysis,ljung1977positive}, who analyzed linear
problems motivated by control and filtering. The paper~\cite{metivier1984applications} analyzed general SA problems for
controlled Markov processes by applying the Kushner--Clark lemma. In
addition to this classical work, stochastic approximation in the
Markov setting has attracted much recent attention. The paper~\cite{chen2020explicit} provides finite-sample error bounds on the averaged iterate of Markovian linear stochastic approximation, with an optimal leading-order term. Central limit
theorems~\cite{fort2015central} and non-asymptotic convergence
rates~\cite{karimi2019non} have been established for controlled Markov
processes. In addition to the papers discussed in
Section~\ref{sec:intro}, several recent works have considered
particular aspects of SA with Markov data, including two-time-scale
variants~\cite{doan2020finite,karmakar2018two}, observation skipping schemes
for bias reduction~\cite{kotsalis2020simple}, Lyapunov function-based analysis
under general norms~\cite{chen2021lyapunov}, and
proving guarantees under weaker ergodicity
conditions~\cite{debavelaere2021convergence}.


\subsection{Application to RL problems}

Markovian observations arise naturally in the context of stochastic
control and reinforcement learning (RL).  See~\cite{benveniste2012adaptive} for a historical survey of
algorithms for stochastic control and filtering with Markovian
stochastic approximation, and the books~\cite{bertsekas2019reinforcement,szepesvari2010algorithms}
for more background on the RL setting. In RL problems, SA
algorithms are typically used to solve Bellman equations, a class of
linear or non-linear fixed-point equations. In policy evaluation
problems, temporal difference (TD) methods~\cite{sutton1988learning}
use linear stochastic approximation to estimate the value function of
a given policy, with asymptotic convergence
guarantees~\cite{dayan1994td,tsitsiklis1997analysis,boyan2002technical}
and non-asymptotic
bounds~\cite{bhandari2018finite,khamaru2020temporal,mou2020optimal}. In
the non-linear case, the Q-learning algorithm~\cite{watkins1992q} is a
stochastic approximation method that estimates the Q-function of a
Markov decision process from data. There is a long line of past work
on this algorithm, including convergence
guarantees~\cite{Tsi94,szepesvari1998asymptotic,EveMan03}, results on
linear function approximation for optimal stopping
problems~\cite{tsitsiklis1999optimal,bhandari2018finite}, and
non-asymptotic rates under general norms in both the $\myiid$
setting~\cite{wainwright2019stochastic,borkar2021concentration} as
well as the Markovian
setting~\cite{chen2021lyapunov}. A class of
variants of TD and Q-learning are also studied in literature,
including actor-critic methods~\cite{konda2000actor},
SARSA~\cite{Rummery1994line}, and methods that employ
variance-reduction~\cite{sidford2018near,khamaru2020temporal,wainwright2019variance,khamaru2021instance}. A
concurrent preprint to this manuscript~\cite{li2021accelerated} proves lower bounds on the oracle complexity of policy evaluation
with access to temporal difference operators, and develops an acceleration 
scheme with variance reduction to achieve these lower bounds while retaining the optimal sample complexity.

It should be noted that an important feature of reinforcement learning
is function approximation, i.e., using a given function class (e.g. a
linear subspace) to approximate the solution to the Bellman equation
of interest. This method enables estimation with a sample size
depending on the intrinsic complexity of the function class, instead
of the cardinality of state-action space. On the other hand, an
approximation error is induced by projecting the Bellman equation onto
this function class. This trade-off is central to the class of TD
algorithms, as studied in a line of past
work~\cite{tsitsiklis1997analysis,yu2010error,bertsekas2011temporal,munos2008finite,mou2020optimal}.
Prior work by a subset of the current authors~\cite{mou2020optimal} focuses on the $\myiid$ setting,
and shows that projected linear equations have a non-standard tradeoff
between approximation and estimation errors.  The current paper is
complementary in nature, building on this work by analyzing the
more challenging setting of Markov observations.  Among the concrete
consequences of this paper are an instance-optimal analysis of TD
algorithms in the Markov setting with linear function approximation.
This analysis provides the basis for a principled choice of the
parameter $\lambda$ in the broader class of TD($\lambda$) algorithms.

\section{Auxiliary truncation results related to the assumptions}\label{SecTruncate}

In this section, we present two auxiliary results on the relations between assumptions~\ref{assume:noise-moments},~\ref{assume:stationary-tail}, and~\ref{assume-lip-mapping}. These results are based on truncation arguments.

\subsection{Assumption~\ref{assume:noise-moments} (almost) implies assumption~\ref{assume-lip-mapping} under discrete metric}
For the discrete metric $\metric (x, y) \mydefn \bm{1}_{x \neq y}$, the Lipschitz assumption~\ref{assume-lip-mapping} is equivalent to the following uniform upper bounds:
\begin{align*}
    \opnorm{\Lmap_{t + 1} (\state_t) - \Lbar} \leq \sigmaA \usedim \quad \mbox{and} \quad \vecnorm{\bmap_{t + 1} (\state_t) - \bbar}{2} \leq \sigmab \sqrt{\usedim}.
\end{align*}
The following proposition provides uniform high-probability upper bounds on such quantities based on the moment assumption:
\begin{proposition}\label{prop:moment-bound-implies-uniform-bound}
    Under Assumption~\ref{assume:noise-moments} with $\pmax = + \infty$, there exists a universal constant $c > 0$, such that for any $\delta > 0$, the following bounds hold true uniformly over $t = 1,2, \cdots, \numobs$, with probability $1 - \delta$:
    \begin{align}
        \opnorm{\Lmap_{t + 1} (\state_t) - \Lbar} \leq c \usedim \cdot \sigmaA \log \frac{\numobs \usedim}{\delta} \quad \mbox{and} \quad
        \vecnorm{\bmap_{t + 1} (\state_t) - \bbar}{2} \leq c \sqrt{\usedim} \cdot \sigmab \log \frac{\numobs \usedim}{\delta}.\label{eq:moment-bound-implies-uniform-bound}
    \end{align}
\end{proposition}
\noindent We prove this proposition at the end of this section.

When the random observations $(\Lmat_{t + 1}, \bvec_{t + 1})$ are not
almost-surely bounded, but satisfies the moment
assumption~\ref{assume:noise-moments} with $\pmax = + \infty$, we can
apply our theorems on the event that
equation~\eqref{eq:moment-bound-implies-uniform-bound} holds true, and
the main theorems hold true conditionally on such an event, with
constants $(\sigmaA, \sigmab)$ inflated with a factor $\log (\numobs
\usedim / \delta)$.

\paragraph{Proof of Proposition~\ref{prop:moment-bound-implies-uniform-bound}:} For a given $t \in [\numobs]$, we note that:
\begin{align*}
  \opnorm{\Lmat_{t + 1} - \Lbar}^2 \leq \matsnorm{\Lmat_{t + 1} -
    \Lbar}{F}^2 = \sum_{j, \ell = 1}^d \Big[ \coordinate_j^\top \big(
    \Lmat_{t + 1} - \Lbar \big) \coordinate_\ell \Big]^2.
\end{align*}
For each pair $j , \ell \in [d]$,
Assumption~\ref{assume:noise-moments} implies that:
\begin{align*}
    \Prob \left( \abss{\coordinate_j^\top \big( \Lmat_{t + 1}
      (\state_t) - \Lbar \big) \coordinate_\ell } \geq c \sigmaA \log
    (\numobs \usedim/\delta) \right) \leq \frac{\delta}{2 d^2 \numobs}
\end{align*}
Taking union bound over all the coordinate pairs $(j, \ell)$ and
substituting into above expansion, we have that:
\begin{align*}
    \Prob \left( \opnorm{\Lmat_{t + 1} - \Lbar} \geq c \usedim \cdot
    \sigmaA \log (\numobs \usedim/\delta) \right) \leq \delta /
    (2\numobs).
\end{align*}
Similarly, for the vector-valued observations $\bvec_{t + 1}$, we have
the following bounds with probability $1 - \delta / n$:
\begin{align*}
    \vecnorm{\bvec_{t + 1} - \bbar}{2}^2 \leq \sum_{j = 1}^\usedim
    \big( \coordinate_j^\top (\bvec_{t + 1} - \bbar) \big)^2 \leq c
    \sigmab^2 \usedim \cdot \log^2 (\numobs \usedim/\delta) .
\end{align*}
Taking union bound over $t = 1,2, \cdots, \numobs$, we complete the
proof of this proposition.

\subsection{On the stationary tail and boundedness assumption~\ref{assume:stationary-tail}}

Note that in many applications, the Markov chain $(\state_t)_{t \geq
  0}$ lives in an unbounded state space. However, as long as the
stationary distribution $\stationary$ of $\transition$ is sufficiently
light-tailed, a simple truncation argument applies, which we
illustrate for completeness. Concretely, suppose that there exists a
constant $\sigma_\metric > 0$, such that the following bound holds
true for any $p \geq 2$:
\begin{align}
\label{eq:unbounded-space-tail-assumption}  
\Exs_{s \sim \stationary} \big[ \metric(s, s_0)^p \big] \leq p!  \cdot
\sigma_\metric^p.
\end{align}
Given a stationary Markovian trajectory $\{ \state_t
\}_{t=1}^\numobs$, consider the event
\begin{align*}
 \Event_{n, \delta} = \big\{ \forall t \in [1, \numobs], \metric
 (\state_0, \state_t) \leq 2 \sigma_\metric \log
 \tfrac{\numobs}{\delta} \big\}.
\end{align*}
By the tail assumption~\eqref{eq:unbounded-space-tail-assumption} and
a union bound, it directly follows that $\Prob \big( \Event_{n,
  \delta} \big) \geq 1 - \delta$. Consider a truncated Markov
transition kernel $\transition'$ defined as
\begin{align*}
    \transition' (x, Z) \mydefn \transition \big(x, Z \cap \ball
    \big(0, 2 \sigma_\metric \log (n / \delta) \big) \big) +
    \transition \big(x, \ball \big(0, 2 \sigma_\metric \log (n /
    \delta) \big)^c \big) \bm{1}_{\state_0 \in Z},
\end{align*}
for any $x \in \statespace$ and $Z \subseteq \statespace$.

In words, the Markov chain $\transition'$ attempts to make the transition from $\state_t$ to $\state_{t + 1}$ according to the transition kernel $\transition'$. If the state $\state_{t + 1}$ lies in the ball $\ball \big(0, 2 \sigma_\metric \log (n / \delta) \big)^c$, we keep it as is; otherwise, we let the next-step transition be deterministically $\state_0$.

Given a trajectory $\{\state'_t \}_{t =1}^\numobs$ of the Markov chain
$\transition'$, there exists a coupling such that
\begin{align*}
    \Prob \big( \{\state_t\}_{t=1}^\numobs \neq
    \{\state_t'\}_{t=1}^\numobs \big) \leq \Prob \big( \Event_{n,
      \delta}^c \big) \leq \delta.
\end{align*}
One can then proceed by working on the high probability event
$\Event_{n, \delta}$, where the Markov chain has a effective diameter
of $O \big( \sigma_\metric \log \frac{\numobs}{\delta} \big)$.

\subsection{Proof of Corollary~\ref{cor:non-stationary-initial}}\label{app:subsec-proof-non-stationary-initial}
Suppose that $\state_0 \sim \pi_0$, by Lemma~\ref{LemWassDecay} and convexity of the Wasserstein distance, we have
\begin{align*}
  \Wass_{1, \metric} \big( \pi_0 \transition^\numcold, \stationary \big) \leq 2^{ - \lfloor \numcold / \mixingtime \rfloor} \leq 2 \exp \Big( - \frac{\numobs}{8 \mixingtime} \Big).
\end{align*}
Let $\big( \widetilde{\state}_t \big)_{t \geq 0}$ be a stationary chain with $\widetilde{\state}_0 \sim \stationary$ independent of $\state_0$. There exists a coupling between the paths, such that
\begin{align*}
  \Exs \big[ \metric (\state_{\numcold}, \widetilde{\state}_{\numcold}) \big] \leq  2 \exp \Big( - \frac{\numobs}{8 \mixingtime} \Big).
\end{align*}
Applying Assumption~\ref{assume-markov-mixing} (b) conditionally on $(\state_\numcold, \widetilde{\state}_\numcold)$, since $c_0 = 1$, there exists a coupling between the next-step transitions, such that
\begin{align*}
  \Exs \Big[ \metric (\state_{\numcold + 1}, \widetilde{\state}_{\numcold + 1}) \mid (\state_\numcold, \widetilde{\state}_\numcold) \Big] \leq \metric (\state_{\numcold}, \widetilde{\state}_{\numcold}).
\end{align*}
Similarly, we can inductively construct the coupling between $\state_{\numcold + i + 1}$ and $\widetilde{\state}_{\numcold + i + 1}$ conditionally on the pair $(\state_{\numcold + i}, \widetilde{\state}_{\numcold + i})$, for $i = 1,2, \cdots$. Putting them together, we obtain a coupling between the two paths, such that $\big(\metric (\state_{\numcold + i + 1}, \widetilde{\state}_{\numcold + i + 1})\big)_{i \geq 0}$ is a super-martingale. By Markov inequality, for each $t \geq \numcold$, we have
\begin{align*}
  \Prob \left( \metric (\state_{t}, \widetilde{\state}_{t}) \geq e^{- \frac{\numobs}{16 \mixingtime}} \right) \leq 2 \exp \Big( - \frac{\numobs}{8 \mixingtime} \Big).
\end{align*}
Define the event
\begin{align*}
  \Event \mydefn \Big\{\metric (\state_t, \widetilde{\state}_t) \leq e^{- \frac{\numobs}{16 \mixingtime}} ~: ~ t = \numcold, \numcold + 1, \cdots, \numobs\Big\}.
\end{align*}
By union bound, we have
\begin{align*}
  \Prob \big( \Event \big) \geq 1 - 2 \numobs \exp \Big( - \frac{\numobs}{8 \mixingtime} \Big) \geq 1 - \exp \Big( - \frac{\numobs}{16 \mixingtime} \Big).
\end{align*}
Define the error scalar $\delta_\numobs \mydefn e^{- \frac{\numobs}{16 \mixingtime}}$, and let $(\theta_t)_{t \geq \numcold}, ~ (\widetilde{\theta}_t)_{t \geq \numcold}$ be the iterate sequences generated by the Markov chains $(\state_t)_{t \geq \numcold}$ and $(\widetilde{\state}_t)_{t \geq \numcold}$, respectively. For any $t \geq \numcold$, we note that
\begin{align*}
  \theta_{t + 1} &= \big( (1 - \stepsize) I_d + \stepsize\Lmap_{t + 1} (\state_t) \big) \theta_t + \stepsize \bmap_{t + 1} (\state_t),\\
  \widetilde{\theta}_{t + 1} &= \big( (1 - \stepsize) I_d + \stepsize\Lmap_{t + 1} (\widetilde{\state}_t) \big) \widetilde{\theta}_t + \stepsize \bmap_{t + 1} (\widetilde{\state}_t).
\end{align*}
Taking their difference and applying triangle inequality, on the event $\Event$, we have the almost-sure upper bound on the one-step error
\begin{align*}
  &\vecnorm{\widetilde{\theta}_{t + 1} - \theta_{t + 1}}{2}\\
  & \leq \opnorm{ (1 - \stepsize) I_d + \stepsize\Lmap_{t + 1} ({\state}_t) \big)} \cdot  \vecnorm{\widetilde{\theta}_{t} - \theta_{t}}{2} + \stepsize \opnorm{ \Lmap_{t + 1} (\widetilde{\state}_t) \big) - \Lmap_{t + 1} (\state_t) \big)} \cdot  \vecnorm{\widetilde{\theta}_{t}}{2} \\
  &\qquad \qquad+ \stepsize \vecnorm{\bmap_{t + 1} (\widetilde{\state}_t) - \bmap_{t + 1} (\state_t)}{2}\\
  &\leq \big(1 + \stepsize (\conmax + \usedim \sigmaA) \big) \vecnorm{\widetilde{\theta}_{t} - \theta_{t}}{2} + \stepsize \delta_\numobs \cdot \Big\{ \sigmaA \usedim \vecnorm{\widetilde{\theta}_t}{2} + \sigmab \sqrt{\usedim} \Big\},
\end{align*}
where in the last step, we use the Lipschitz assumption~\ref{assume-lip-mapping}.

Solving this recursion yields the uniform upper bound for $t \in \{\numcold, \numcold + 1, \cdots, \numobs\}$
\begin{align*}
  \vecnorm{\widetilde{\theta}_{t} - \theta_{t}}{2} \leq \stepsize \delta_\numobs \exp \Big(  \stepsize (\conmax + \usedim \sigmaA) \numobs \Big) \cdot \sum_{t = \numcold}^\numobs \Big\{ \sigmaA \usedim \vecnorm{\widetilde{\theta}_t}{2} + \sigmab \sqrt{\usedim} \Big\},
\end{align*}
holding with probability $1$ on the event $\Event$.

Given a stepsize satisfying $\stepsize \leq \frac{1}{32 \mixingtime (\conmax + \usedim \sigmaA)}$, we have
\begin{align*}
   \delta_\numobs \exp \Big(  \stepsize (\conmax + \usedim \sigmaA) \numobs \Big) \leq \exp \Big( - \frac{\numobs}{32 \mixingtime} \Big).
\end{align*}
For the summation term, we apply Cauchy--Schwarz inequality, and obtain the MSE bound
\begin{align*}
  \Exs \Big\{ \sum_{t = \numcold}^\numobs \sigmaA \usedim \vecnorm{\theta_t}{2} + \sigmab \sqrt{\usedim} \Big\}^2 &\leq 2 \numobs^2 \usedim^2 \big( \sigmab^2 + \sigmaA^2 \vecnorm{\thetabar}{2}^2 \big)  + 4 \numobs \sigmaA^2 \usedim^2 \sum_{t = \numcold}^\numobs \Exs \big[ \vecnorm{\widetilde{\theta}_t - \thetabar}{2}^2\big]\\
  &\overset{(i)}{\leq} \numobs^2 \usedim^2 \Big( 2\sigmab^2 + 6 \sigmaA^2 \vecnorm{\thetabar}{2}^2 + \frac{4 c \stepsize \mixingtime \stepsize}{1 - \kappa} \varbound^2 \log \numobs \Big)\\
  &\leq 12 \numobs^3 \usedim^2 \big(\sigmab^2 + \sigmaA^2 \vecnorm{\thetabar}{2}^2 \big),
\end{align*}
where in step $(i)$, we apply Proposition~\ref{prop:lsa-markov-iterate-bound} to the iterate sequence $(\widetilde{\theta}_t)_{t \geq \numcold}$.

Putting them together, we conclude that
\begin{align*}
  \Exs \big[ \vecnorm{\widetilde{\theta}_{t} - \theta_{t}}{2}^2 \bm{1}_{\Event} \big] \leq 12 \numobs^3 \usedim^2 \big(\sigmab^2 + \sigmaA^2 \vecnorm{\thetabar}{2}^2 \big) \exp \Big( - \frac{\numobs}{16 \mixingtime} \Big). 
\end{align*}
Let $\thetahat_\numobs' \mydefn \frac{1}{\numobs - \numburn} \sum_{t = \numburn}^{\numobs - 1} \widetilde{\theta}_t$, applying Cauchy--Schwarz inequality, we have
\begin{align*}
  \Exs \big[ \vecnorm{\thetahat_\numobs' - \thetahat_\numobs}{2}^2 \bm{1}_{\Event}  \big] &\leq \frac{4}{\numobs} \sum_{t = \numburn}^\numobs \Exs \big[ \vecnorm{\widetilde{\theta}_{t} - \theta_{t}}{2}^2 \bm{1}_{\Event} \big]\\
  &\leq  12 \numobs^3 \usedim^2 \big(\sigmab^2 + \sigmaA^2 \vecnorm{\thetabar}{2}^2 \big) \exp \Big( - \frac{\numobs}{16 \mixingtime} \Big)\\
  &\leq e^{- \frac{\numobs}{32 \mixingtime}}\big(\sigmab^2 + \sigmaA^2 \vecnorm{\thetabar}{2}^2 \big),
\end{align*}
for a sample size satisfying $\tfrac{\numobs}{\log \numobs} \geq 2400 \mixingtime \log \usedim$. Invoking Theorem~\ref{thm:markov-main} on the estimator $\thetahat_\numobs'$ completes the proof of this corollary.

\section{Auxiliary results underlying Proposition~\ref{prop:lsa-markov-iterate-bound}}

This appendix is devoted to the proofs of auxiliary lemmas that are
used in the proof of Proposition~\ref{prop:lsa-markov-iterate-bound}.


\subsection{Proof of Lemma~\ref{LemWassDecay}}
\label{AppLemWassDecay}

Throughout the proof, we let $x \in \statespace$ be an arbitrary but
fixed state.  Note that any positive integer $\tau$ can be represented
as $\tau = k \mixingtime + q$ with $k \in \mathbb{N}_+$ and $0 \leq q
\leq \mixingtime - 1$. We show the desired claim by induction over $k
\geq 0$.

\paragraph{Base case:}
When $k = 0$, Assumption~\ref{assume:stationary-tail} implies that
\begin{align*}
\Wass_{1, \metric} (\delta_x \transition^\tau, \stationary) \leq
\sup_{s, s' } \metric (s,s') \leq 1 \leq c_0,
\end{align*}
so that the base case ($k = 0$) holds for our induction proof.

\paragraph{Induction step:}
At step $k$ of the argument, the induction hypothesis ensures that
\begin{align}
\label{eq:induction-hypo-in-lemma-wass-decay}  
\Wass_{1, \metric} \big( \delta_x \transition^{k \mixingtime + q},
\stationary \big) \leq c_0 \cdot 2^{- k}, \quad \mbox{for } q = 0,1,
\cdots, \mixingtime - 1.
\end{align}
We now need to show that the result holds for any $\tau = (k + 1)
\mixingtime + q$, where $q \in \{0, 1, \ldots, \mixingtime -1 \}$ is
arbitrary.  We do so via a coupling argument. Take a random initial
state $y \sim \stationary$, and consider two processes $\{s_t \}_{t
  \geq 0}$ and $\{s_t'\}_{t \geq 0}$ starting from $x$ and $y$,
respectively. Their joint distribution is defined as follows: choose
the coupling between the law of $\state_{k \mixingtime + q}$ and
$\state_{k \mixingtime + q}'$ to satisfy the identity $\Exs \big[
  \metric (\state_{k \mixingtime + q}, \state_{k \mixingtime + q}')
  \big] = \Wass_{1, \metric} \big( \delta_x \transition^{k \mixingtime
  + q}, \stationary \big)$.  Conditionally on $(\state_{k \mixingtime
  + q}, \state_{k \mixingtime + q}')$,
Assumption~\ref{assume-markov-mixing} guarantees the existence of a
coupling between $\delta_{\state_{k \mixingtime + q}}
\transition^{\mixingtime}$ and $\state_{k \mixingtime + q}'
\transition^{\mixingtime}$ such that
\begin{align*}
  \Exs \big[ \metric \big( \state_{(k + 1) \mixingtime + q},
    \state_{(k + 1) \mixingtime + q}' \big) \mid (\state_{k
      \mixingtime + q}, \state_{k \mixingtime + q}') \big] \leq
  \tfrac{1}{2} \metric (\state_{k \mixingtime + q}, \state_{k
    \mixingtime + q}').
\end{align*}
Taking expectation on both sides and substituting with equation~\eqref{eq:induction-hypo-in-lemma-wass-decay}, we
find that
\begin{align*}
    \Wass_{1, \metric} \big( \delta_x \transition^{(k + 1)
      \mixingtime + q}, \stationary \big) \leq \Exs \big[ \metric
      \big( \state_{(k + 1) \mixingtime + q}, \state_{(k + 1)
        \mixingtime + q}' \big) \big] \leq c_0 \cdot 2^{- (k + 1)},
\end{align*}
which completes the proof of the induction step.

  
\subsection{Proof of Lemma~\ref{lemma:norm-not-blowup}}
\label{App:proof-norm-not-blowup}

Our proof is based on the following intermediate claim
\begin{align}
\label{eq:delta-norm-growth-bound}  
\big( \Exs \big[ \vecnorm{\Delta_{t + \ell }}{2}^p \big]
\big)^{1/p} \leq e \big( \Exs \big[ \vecnorm{\Delta_{t }}{2}^p
  \big] \big)^{1/p} + 6 \stepsize p \ell \varbound \sqrt{d} .
\end{align}
This bound, which we return to prove at the end of this section, is a
weaker form of the claim in the lemma.

We now use the bound~\eqref{eq:delta-norm-growth-bound} to prove the
lemma.  Applying Minkowski's inequality to the recursive
relation~\eqref{eq:recursive-relation-delta}, we find that for any $p
\geq 2$, the $p$-th moment is upper bounded as
\begin{multline*}
    \big( \Exs \big[ \vecnorm{\Delta_{t + \ell + 1} - \Delta_t}{2}^p
      \big] \big)^{1/p} \leq \big( \Exs \big[ \vecnorm{ \Delta_{t +
          \ell} - \Delta_t }{2}^p \big] \big)^{1/p} + \stepsize
    \big( \Exs \big[ \vecnorm{ \Lmat_{t + \ell + 1} \Delta_{t + \ell}
      }{2}^p \big] \big)^{1/p} \\
      + \stepsize \big( \Exs \big[
      \vecnorm{\addnoiseMarkov_{t + \ell}
       + \addnoiseMG_{t + \ell +
          1}}{2}^p \big] \big)^{1/p}.
\end{multline*}

For the martingale part of the noise, we take the decomposition
$\Lmat_{t + \ell + 1} = \Lmap (\state_{t + \ell}) + \multnoiseMG_{t +
  \ell + 1}$.  By Assumption~\ref{assume:noise-moments} and
H\"{o}lder's inequality, we have the bounds
\begin{align*}
 \Exs \big[ \vecnorm{ \multnoiseMG_{t + \ell + 1} \Delta_{t + \ell}
   }{2}^p \mid \filtration_{t} \big] &\leq \usedim^{\tfrac{p}{2}}
 \sum_{j = 1}^{\usedim} \Exs \big[
   \inprod{\coordinate_j}{\multnoiseMG_{t + \ell + 1} \Delta_{t +
       \ell} }^p \mid \filtration_t \big] \\
       &\leq \big( p \sigmaA
 \sqrt{d} \big)^p \Exs \big[ \vecnorm{\Delta_{t + \ell}}{2}^p \mid
   \filtration_t \big], \quad \mbox{and} \\
\Exs \big[ \vecnorm{ \addnoiseMG_{t + \ell + 1}}{2}^p  \big] & \leq \usedim^{\tfrac{p}{2}} \sum_{j =
  1}^{\usedim} \Exs \big[ \inprod{\coordinate_j}{ \addnoiseMG_{t +
      \ell + 1}}^p  \big] \leq (p
\sqrt{\usedim})^p \cdot \varbound^p.
\end{align*}
Similarly, for the Markov part of the noise, we have:
\begin{align*}
    \Exs \big[ \vecnorm{\addnoiseMarkov_{t + \ell + 1}}{2}^p \big]
    \leq (p \sqrt{d})^p \cdot \varbound^p.
\end{align*}

On the other hand, the Lipschitz condition~\eqref{assume-lip-mapping}
and the boundedness condition~\eqref{assume:stationary-tail} of the
metric space imply that
\begin{align*}
  \opnorm{\Lmat_{t + \ell + 1} (\state) - \Lbar} \leq \sigmaA \usedim,\quad  \mbox{for all $\state
    \in \statespace$.}
\end{align*}
Substituting into the decomposition above,
we arrive at the bounds
\begin{align*}
 \big( \Exs \big[ \vecnorm{ \Lmat_{t + \ell + 1} \Delta_{t + \ell}
   }{2}^p \big] \big)^{1/p} &\leq (\conmax + \sigmaA p \sqrt{d} +
 \sigmaA d) \big( \Exs \big[ \vecnorm{ \Delta_{t + \ell} }{2}^p \big]
 \big)^{1/p}, \quad \mbox{and} \\
 \big( \Exs \big[ \vecnorm{\addnoiseMarkov_{t + \ell} +
     \addnoiseMG_{t + \ell + 1}}{2}^p \big] \big)^{1/p} &\leq 2 p \varbound
 \sqrt{d} .
\end{align*}
Applying equation~\eqref{eq:delta-norm-growth-bound} yields
\begin{multline*}
  \big( \Exs \big[ \vecnorm{\Delta_{t + \ell + 1} - \Delta_t}{2}^p
    \big] \big)^{1/p} \leq \big( \Exs \big[ \vecnorm{ \Delta_{t
        + \ell} - \Delta_t }{2}^p \big] \big)^{1/p} + e \stepsize
  (\conmax + \sigmaA d) \big( \Exs \big[ \vecnorm{\Delta_t}{2}^p
    \big] \big)^{1/p} \\
  + 2 (1 + 6 \stepsize \ell) \stepsize p \varbound \sqrt{d},
\end{multline*}
where the second inequality comes from the definition of $\varbound$.

Solving this recursion leads to the bound
\begin{align*}
  \big( \Exs \big[ \vecnorm{\Delta_{t + \ell} - \Delta_t}{2}^p \big]
  \big)^{1/p} & \leq e \stepsize \ell (\conmax + \sigmaA d) \big(
  \Exs \big[ \vecnorm{\Delta_t}{2}^p \big] \big)^{1/p} + 3 \stepsize
  p \ell \varbound \sqrt{d} ,
\end{align*}
which establishes the first claim.

Since the stepsize is upper bounded as $\stepsize \leq \big(2 e
\stepsize \ell (\conmax + \sigmaA d) \big)^{-1}$, we have the lower
bound
\begin{align*}
  \big( \Exs \big[ \vecnorm{\Delta_{t + \ell}}{2}^p \big]
  \big)^{1/p} & \geq \big( \Exs \big[ \vecnorm{\Delta_{t}}{2}^p
    \big] \big)^{1/p} - \big( \Exs \big[ \vecnorm{\Delta_{t + \ell}
      - \Delta_t}{2}^p \big] \big)^{1/p} \\
& \geq \tfrac{1}{2} \big( \Exs \big[ \vecnorm{\Delta_{t}}{2}^p \big]
  \big)^{1/p} - 3 \stepsize p \ell \varbound \sqrt{d},
\end{align*}
which, in conjunction with the
bound~\eqref{eq:delta-norm-growth-bound}, establishes the second
claim.


\paragraph{Proof of equation~\eqref{eq:delta-norm-growth-bound}:}

Applying Minkowski's inequality to the recursive
relation~\eqref{eq:recursive-relation-delta} yields (for any $p \geq
2$) a bound on the $p^{th}$ conditional moment:
\begin{align}
\big( \Exs \big[ \vecnorm{\Delta_{t + \ell + 1}}{2}^p \big]
\big)^{1/p} \leq \big( \Exs \big[ \vecnorm{ (I - \stepsize \Lmat_{t
      + \ell + 1}) \Delta_{t + \ell} }{2}^p \big] \big)^{1/p} +
\stepsize \big( \Exs \big[ \vecnorm{\addnoiseMarkov_{t + \ell} +
    \addnoiseMG_{t + \ell + 1}}{2}^p \big]
\big)^{1/p}.\label{eq:recursive-relation-in-norm-not-blowup-proof}
\end{align}
Our next step is to bound the two terms above.

Substituting into the recursive
relation~\eqref{eq:recursive-relation-in-norm-not-blowup-proof}, and
applying Minkowski's inequality, we find that the moment $\big( \Exs
\big[ \vecnorm{\Delta_{t + \ell + 1}}{2}^p \big] \big)^{1/p}$ is upper
bounded by
\begin{align*}
  (1 + \stepsize \conmax)\big( \Exs \big[ \vecnorm{\Delta_{t +
        \ell}}{2}^p \big] \big)^{1/p} + \stepsize \sigmaA \usedim
  \big( \Exs \big[ \vecnorm{ \Delta_{t + \ell} }{2}^p \big]
  \big)^{1/p} + 2 \stepsize p \sqrt{d} \varbound.
\end{align*}
Solving this recursive inequality leads to
\begin{align*}
  \big( \Exs \big[ \vecnorm{\Delta_{t + \ell }}{2}^p \big]
  \big)^{1/p} \leq \exp \big( \stepsize \ell (\conmax + \sigmaA d)
  \big) \big( \big( \Exs \big[ \vecnorm{\Delta_{t }}{2}^p \big]
  \big)^{1/p} + 2 \stepsize p \ell \sqrt{d} \varbound.
\end{align*}
For any stepsize $\stepsize \in \big(0, \tfrac{1}{(\conmax + \sigmaA
  d) \ell} \big]$, we have
\begin{align*}
  \big( \Exs \big[ \vecnorm{\Delta_{t + \ell }}{2}^p \big]
  \big)^{1/p} \leq e \big( \Exs \big[ \vecnorm{\Delta_{t }}{2}^p
    \big] \big)^{1/p} + 6 \stepsize p \ell \sqrt{d} \varbound,
\end{align*}
which establishes the claim.


\subsection{Proof of Lemma~\ref{lemma:local-moves}}
\label{App:proof-local-moves}

For notational simplicity, we extend the process $(\Delta_t)_{t \geq
  0}$ to the entire set $\mathbb{Z}$ of integers, in particular by
defining $\Delta_t \defn \Delta_0$ for negative integer $t$. Note that
under our assumption, Lemma~\ref{lemma:norm-not-blowup} and the
assumed bound~\eqref{eq:local-move-coarse-bound} both hold true for
the extended process, with index set $t \in \mathbb{Z}$.  Moreover, as
in the proof of Lemma~\ref{lemma:norm-not-blowup}, for each $p \geq
2$, we have the moment bound
\begin{multline*}
  \big( \Exs \big[ \vecnorm{\Delta_{t + \ell + 1} - \Delta_t}{2}^p
    \big] \big)^{1/p} \leq \big( \Exs \big[ \vecnorm{ \Delta_{t +
        \ell} - \Delta_t }{2}^p \big] \big)^{1/p} + \stepsize \big(
  \Exs \big[ \vecnorm{ \Lmat_{t + \ell + 1} \Delta_{t + \ell} }{2}^p
    \big] \big)^{1/p} \\
  + \stepsize \big( \Exs \big[
    \vecnorm{\addnoiseMarkov_{t + \ell} + \addnoiseMG_{t + \ell +
        1}}{2}^p \big] \big)^{1/p}.
\end{multline*}

Our next step is to exploit the coarse
bound~\eqref{eq:local-move-coarse-bound} so as to obtain upper bounds
on the second term $\big( \Exs \big[ \vecnorm{ \Lmat_{t + \ell + 1}
    \Delta_{t + \ell} }{2}^p \big] \big)^{1/p}$.  Given the time lag
$\tau > 0$, we take the decomposition $\Delta_{t + \ell} = \Delta_{t +
  \ell - \tau} + (\Delta_{t + \ell} - \Delta_{t + \ell - \tau})$, and
by Minkowski's inequality, we have that
\begin{multline}
\label{eq:local-move-step-back-decomposition}  
  \big( \Exs \big[ \vecnorm{ \Lmat_{t + \ell + 1} \Delta_{t + \ell}
    }{2}^p \big] \big)^{1/p} \\
    \leq \big( \Exs \big[ \vecnorm{
      \Lmat_{t + \ell + 1} \Delta_{t + \ell - \tau} }{2}^p \big]
  \big)^{1/p} + \big( \Exs \big[ \vecnorm{ \Lmat_{t + \ell + 1}
      (\Delta_{t + \ell} - \Delta_{t + \ell - \tau}) }{2}^p \big]
  \big)^{1/p}.
\end{multline}
The latter term of the
bound~\eqref{eq:local-move-step-back-decomposition} can be controlled
through Assumption~\ref{assume-lip-mapping}:
\begin{align*}
\vecnorm{\Lmat_{t + \ell + 1} (s_{t + \ell}) (\Delta_{t + \ell} -
  \Delta_{t + \ell - \tau})}{2} & \leq (\conmax + \sigmaA d )
\vecnorm{\Delta_{t + \ell} - \Delta_{t + \ell - \tau}}{2}, \quad
\mbox{a.s.}
\end{align*}
The distance $\vecnorm{\Delta_{t + \ell} - \Delta_{t + \ell -
    \tau}}{2}$ is controlled via the coarse
bound~\eqref{eq:local-move-coarse-bound}. Putting together the pieces,
we find that
\begin{align}
 \label{eq:local-move-proof-part1-bound}  
  \big( \Exs \big[ \vecnorm{ \Lmat_{t + \ell + 1} (\Delta_{t + \ell}
      - \Delta_{t + \ell - \tau})}{2}^p \big] \big)^{1/p} \leq
  \stepsize \big( \conmax + \sigmaA d \big) \cdot \big(\omega_p
  \big( \Exs \big[ \vecnorm{\Delta_{t + \ell - \tau}}{2}^p \big]
  \big)^{1/p} + \beta_p \varbound \big).
\end{align}
In order to bound the former term $\big( \Exs \big[ \vecnorm{ \Lmat_{t
      + \ell + 1} \Delta_{t + \ell - \tau} }{2}^p \big] \big)^{1/p}$
in the bound~\eqref{eq:local-move-step-back-decomposition}, we invoke
Lemma~\ref{LemWassDecay}, and obtain a random variable
$\widetilde{\state}_{t + \ell}$, such that
\begin{align}
\label{eq:mixing-estimate-for-replaced-state-in-local-move-proof}  
  \widetilde{\state}_{t + \ell} \mid \filtration_{t + \ell - \tau}
  \sim \stationary, \quad \mbox{and} \quad \big( \Exs \big[ \metric
    (\state_{t + \ell}, \widetilde{\state}_{t + \ell - \tau})^p \mid
    \filtration_{t + \ell - \tau} \big] \big)^{1/p} \leq c_0 \cdot
  2^{1 - \tfrac{\tau}{2 \mixingtime p}}.
\end{align}
By Assumption~\ref{assume:noise-moments}, we have the bounds
\begin{subequations}
  \label{eq:step-back-estimates-in-local-move-proof-abcd}
\begin{align}
    \Exs \big[ \vecnorm{\multnoiseMG_{t + \ell + 1} \Delta_{t + \ell
          - \tau}}{2}^p \mid \filtration_{t + \ell - \tau} \big]
    &\leq (p \sqrt{d} \sigmaA)^p \vecnorm{\Delta_{t + \ell -
        \tau}}{2}^p, \quad \mbox{and}\\ \Exs \big[ \vecnorm{ \big(
        \Lmap (\widetilde{\state}_{t + \ell - \tau}) - \Lbar \big)
        \cdot \Delta_{t + \ell - \tau}}{2}^p \mid \filtration_{t +
        \ell - \tau} \big]&\leq (p \sqrt{d} \sigmaA)^p
    \vecnorm{\Delta_{t + \ell - \tau}}{2}^p.
\end{align}
Invoking the moment
bound~\eqref{eq:mixing-estimate-for-replaced-state-in-local-move-proof}
and using the Lipschitz condition~\eqref{assume-lip-mapping}, we find
that
\begin{align}
&\Exs \big[ \vecnorm{ \big( \Lmap (\widetilde{\state}_{t + \ell -
      \tau}) - \Lmap ({\state}_{t + \ell - \tau}) \big) \cdot
    \Delta_{t + \ell - \tau}}{2}^p \mid \filtration_{t + \ell - \tau}
  \big] \nonumber \\
  & \leq \Exs \big[ \opnorm{ \Lmap (\widetilde{\state}_{t + \ell
      - \tau}) - \Lmap ({\state}_{t + \ell - \tau}) }^p \mid
  \filtration_{t + \ell - \tau} \big] \cdot \vecnorm{\Delta_{t + \ell
    - \tau}}{2}^p \nonumber\\ &\leq \big(\sigmaA c_0 d \cdot 2^{1 -
  \tfrac{\tau}{2 \mixingtime p}} \vecnorm{\Delta_{t + \ell -
    \tau}}{2}\big)^p.
\end{align}
Finally, we have the operator norm bound
\begin{align}
    \vecnorm{\Lbar \Delta_{t + \ell - \tau}}{2} \leq \conmax \vecnorm{\Delta_{t + \ell - \tau}}{2}.
\end{align}
\end{subequations}
Collecting the results from
equations~\eqref{eq:step-back-estimates-in-local-move-proof-abcd}(a)---(d),
we arrive at the bound
\begin{align}
\label{eq:local-move-proof-part2-bound}  
   \big( \Exs \big[ \vecnorm{ \Lmat_{t + \ell + 1} \Delta_{t + \ell -
         \tau} }{2}^p \mid \filtration_{t + \ell - \tau} \big]
   \big)^{1/p} \leq \big( 2 p \sqrt{d} \sigmaA + \conmax + \sigmaA
   c_0 d \cdot 2^{1 - \tfrac{\tau}{2 \mixingtime p}} \big)
   \vecnorm{\Delta_{t + \ell - \tau} }{2}.
\end{align}
According to Lemma~\ref{lemma:norm-not-blowup}, given a stepsize
bounded as $\stepsize \leq \big( 6 (\conmax + \sigmaA d) \tau
\big)^{-1}$, we have
\begin{align*}
  \big( \Exs \vecnorm{\Delta_{t + \ell - \tau}}{2}^p \big)^{1/p}
  \leq 2 \big( \Exs \vecnorm{\Delta_{t + \ell}}{2}^p \big)^{1/p} +
  12 \stepsize p \tau \varbound \sqrt{d} .
\end{align*}

Collecting the bounds~\eqref{eq:local-move-proof-part1-bound}
and~\eqref{eq:local-move-proof-part2-bound}, and substituting into the
decomposition~\eqref{eq:local-move-step-back-decomposition}, for $\tau
\geq 2 \mixingtime p \log (c_0 d)$, we arrive at the inequality:
\begin{multline*}
 \big( \Exs \big[\vecnorm{\Lmat_{t + \ell + 1} \Delta_{t + \ell}}{2}^p
   \big] \big)^{1/p}\\
    \leq 2 \big( \big( p \sqrt{d} \sigmaA + \conmax
 \big) + \stepsize \omega_p \big( \conmax + \sigmaA d \big) \big)\cdot
 \big( \big(\Exs \vecnorm{\Delta_{t + \ell}}{2}^p\big)^{1/p} +
 \stepsize p \tau \sqrt{d} \varbound \big)\\
  + \stepsize \big( \conmax
 + \sigmaA d \big) \beta_p \varbound.
\end{multline*}

By following the derivation in the proof of
Lemma~\ref{lemma:norm-not-blowup}, we can show that the third term is
upper bounded as
\begin{align*}
\big( \Exs \big[ \vecnorm{\addnoiseMarkov_{t + \ell} + \addnoiseMG_{t
      + \ell + 1}}{2}^p \big] \big)^{1/p} \leq 2 p \varbound \sqrt{d} .
\end{align*}
Substituting back into the original decomposition, we find that the
difference in moments \mbox{$D \defn \big( \Exs \big[
    \vecnorm{\Delta_{t + \ell + 1} - \Delta_t}{2}^p \big] \big)^{1/p}
  - \big( \Exs \big[ \vecnorm{ \Delta_{t + \ell} - \Delta_t }{2}^p
    \big] \big)^{1/p}$} is bounded as
\begin{multline*}
  D \leq 2 \stepsize \Big\{ \big( p \sqrt{d} \sigmaA + \conmax \big) +
  \stepsize \omega_p \big( \conmax + \sigmaA d \big) \Big \} \cdot
  \big( \big(\Exs \vecnorm{\Delta_{t + \ell}}{2}^p\big)^{1/p} +
  \stepsize p \tau \sqrt{d} \varbound \big) \\
  + \big(2 \stepsize p
  \sqrt{d} + \stepsize^2 \big( \conmax + \sigmaA d \big) \beta_p \big)
\end{multline*}
Lemma~\ref{lemma:norm-not-blowup} implies that $\big( \Exs \big[
  \vecnorm{\Delta_{t + \ell}}{2}^p \big] \big)^{1/p} \leq e \big( \Exs
\big[ \vecnorm{\Delta_{t }}{2}^p \big] \big)^{1/p} + 6 \stepsize p
\ell \sqrt{d}\varbound$ and solving the recursion, we arrive at the
bound
\begin{align*}
     &\big( \Exs \big[ \vecnorm{\Delta_{t + \ell} - \Delta_t}{2}^p
    \big] \big)^{1/p}\\ &\leq 12 \stepsize \ell \big( \big( p
  \sqrt{d} \sigmaA + \conmax \big) + \stepsize \omega_p \big( \conmax
  + \sigmaA d \big) \big)\cdot \big( \big(\Exs
  \vecnorm{\Delta_{t}}{2}^p\big)^{1/p} + \stepsize p (\tau + \ell)
  \sqrt{d} \varbound \big) \\ &\qquad+ \big(2 \stepsize p \sqrt{d} +
  \stepsize^2 \big( \conmax + \sigmaA d \big) \beta_p \big) \ell
  \varbound\\ &\leq \stepsize \big( 12 \big( p \sqrt{d} \sigmaA +
  \conmax \big) \ell + \tfrac{\omega_p}{2} \big) \big( \big(\Exs
  \vecnorm{\Delta_{t}}{2}^p\big)^{1/p} + \stepsize p (\tau + \ell)
  \sqrt{d} \varbound \big) + \stepsize \big( 2 p \ell \sqrt{d} +
  \tfrac{1}{2} \beta_p\big) \varbound,
\end{align*}
for any $\tau \geq 2 \mixingtime p \log (c_0 d)$ and stepsize choice
$\stepsize \leq \tfrac{c}{48(\conmax + \sigmaA d)}$.


\section{Auxiliary results underlying Theorem~\ref{thm:markov-main}}
\label{AppTheoremUpper}

In this appendix, we prove two auxiliary lemmas that were used in the
proof of Theorem~\ref{thm:markov-main}.


\subsection{Proof of Lemma~\ref{lemma:bound-sum-markov-mult-noise-bootstrap}}
\label{subsubsec:proof-lemma-bound-sum-markov-mult-noise-bootstrap}

According to Lemma~\ref{LemWassDecay}, given $\tau > 0$ fixed, for any
$t \geq \tau + k_m$, there exists a random variable
$\widetilde{\state}_{t - k_m}$ such that $\widetilde{\state}_{t - k_m}
\mid \filtration_{t - k_m - \tau} \sim \stationary$, and $\Exs \big[
  \metric (\state_{t - k_m}, \widetilde{\state}_{t - k_m}) \mid
  \filtration_{t - \tau - k_m} \big] \leq c_0 \cdot 2^{1 -
  \tfrac{\tau}{\mixingtime}}$.  By
Assumption~\ref{assume-markov-mixing}, conditionally on the pair of
states $(\state_{t - k_m}, \widetilde{\state}_{t - k_m} )$, we have
the following bound for $j \in [m]$:
\begin{align*}
\Wass_{\metric, 1} \big( \transition^{k_j - k_{j - 1}}
\delta_{\state_{t - k_j}}, \transition^{k_j - k_{j - 1}}
\delta_{\widetilde{\state}_{t - k_j}} \big) \leq c_0 \cdot \metric
\big( \state_{t - k_j}, \widetilde{\state}_{t - k_j} \big),
\quad\mathrm{a.s.}
\end{align*}
Consequently, there exists a sequence of random variables
$(\widetilde{\state}_{t - k_j})_{0 \leq j \leq m - 1}$, such that the
following relations hold true for $j = 1,2, \cdots, m$:
\begin{align*}
    &\widetilde{\state}_{t - k_{j - 1}} \mid \filtration_{t - k_m}
  \sim \transition^{k_j - k_{j - 1}} \delta_{\widetilde{\state}_{t -
      k_j}}, \quad \mbox{and}\\ &\Exs \big[ \metric \big(
    \widetilde{\state}_{t - k_{j - 1}}, \state_{t - k_{j - 1}} \big)
    \mid \filtration_{t + k - \ell } \big] \leq c_0^{m + 1 - j}
  \cdot \metric \big( \state_{t - k_m}, \widetilde{\state}_{t - k_m}
  \big).
\end{align*}

Based on above construction, we consider the following decomposition:
\begin{align}
    \big( \prod_{j = 0}^m \multnoiseMarkov_{t - k_j} \big) \Delta_{t - k_m} &= \big(\prod_{j = 0}^m \multnoiseMarkov (\state_{t - k_j}) - \prod_{j = 0}^m \multnoiseMarkov (\widetilde{\state}_{t - k_j})  \big) \Delta_{t - k_m - \tau} + \big(\prod_{j = 0}^m \multnoiseMarkov (\widetilde{\state}_{t - k_j}) \big) \cdot \Delta_{t - k_m- \tau} \nonumber\\
    &\qquad+ \big(\prod_{j = 0}^m \multnoiseMarkov (\state_{t - k_j})\big) \cdot \big(\Delta_{t - k_m} - \Delta_{t - \tau - k_m} \big) \mydefn Q_1 (t) + Q_2 (t) + Q_3 (t).\label{eq:cross-term-decomp-for-optimal-cov}
\end{align}
In the following, we bound the moments for the summation of the three terms above, respectively.
For the first term, we note the telescoping equation:
\begin{multline*}
    \prod_{j = 0}^m \multnoiseMarkov (\state_{t - k_j}) - \prod_{j = 0}^m \multnoiseMarkov (\widetilde{\state}_{t - k_j})\\
     =  \sum_{q = 0}^m \big(\prod_{j = 0}^{q - 1} \multnoiseMarkov (\state_{t - k_j})\big) \cdot \big( \Lmap (\state_{t - k_q}) - \Lmap (\widetilde{\state}_{t - k_q}) \big) \cdot \big(\prod_{j = q + 1}^{m} \multnoiseMarkov (\widetilde{\state}_{t - k_j})\big).
\end{multline*}
Note that each matrix in the product has operator norm uniformly bounded by $\sigmaA d$.
We can then use the Lipschitz condition~\ref{assume-lip-mapping} as well as the bound on the distance $\metric (\state_{t - k_q}, \widetilde{\state}_{t - k_q})$, and obtain the bound
\begin{align*}
& \Exs \big[ \opnorm{\prod_{j = 0}^m \multnoiseMarkov (\state_{t -
        k_j}) - \prod_{j = 0}^m \multnoiseMarkov
      (\widetilde{\state}_{t - k_j})}^2 \mid \filtration_{t - k_m -
      \tau} \big] \\
& \leq (m + 1) \cdot (\sigmaA d)^m \sum_{q = 0}^m \Exs \big[
    \opnorm{\Lmap (\state_{t - k_q}) - \Lmap (\widetilde{\state}_{t -
        k_q})}^2 \mid \filtration_{t - k_m - \tau} \big] \\
    & \leq  (m
  + 1)^2 (c_0 \sigmaA d)^{m + 1} \cdot 2^{- \tfrac{\tau}{\mixingtime}}.
\end{align*}
Applying the bound on $\vecnorm{\Delta_{t - \tau}}{2}$ in
Proposition~\ref{prop:lsa-markov-iterate-bound} and taking $\tau \geq
3 m \mixingtime p \log (c_0 d \numobs)$, we find that
\begin{align}
\Exs \big[ \vecnorm{Q_1 (t)}{2}^2 \big] &\leq \Exs \big[ \Exs \big[
    \opnorm{\prod_{j = 0}^m \multnoiseMarkov (\state_{t - k_j}) -
      \prod_{j = 0}^m \multnoiseMarkov (\widetilde{\state}_{t -
        k_j})}^2 \mid \filtration_{t - k_m - \tau} \big] \cdot
  \vecnorm{\Delta_{t - \tau - k_m}}{2}^2 \big] \nonumber\\ &\leq (m +
1)^2 (c_0 \sigmaA d)^{m + 1} \cdot 2^{- \tfrac{\tau}{\mixingtime}} c
\varbound^2 \tfrac{\stepsize \tau \usedim \log^2 n}{1 - \kappa} \; \leq
\; \tfrac{\sigmaA^{m + 1}}{\numobs^2}
\varbound^2. \label{eq:bound-sum-markov-mult-part-1}
\end{align}
Now we turn to bounding the term $Q_2(t)$. First, we note that
\begin{multline*}
    \Exs \big[ \vecnorm{Q_2 (t)}{2}^2 \big] \leq \Exs \big[ \opnorm{\prod_{j = 0}^{m - 1} \multnoiseMarkov (\widetilde{\state}_{t - k_j})}^2 \cdot \vecnorm{ \multnoiseMarkov (\widetilde{\state}_{t - k_m}) \Delta_{t - k_m - \tau}}{2}^2 \big]\\
    \leq (\sigmaA \usedim)^{2 m} \Exs \big[\vecnorm{ \multnoiseMarkov (\widetilde{\state}_{t - k_m}) \Delta_{t - k_m - \tau}}{2}^2 \big] \leq (\sigmaA \usedim)^{2 m} \cdot \sigmaA^2 d \cdot \Exs \big[ \vecnorm{\Delta_{t - k_m - \tau}}{2}^2 \big].
\end{multline*}
By Proposition~\ref{prop:lsa-markov-iterate-bound}, for $t \geq
\numburn$ and $\numburn \geq 2 (\tau + k_m)$, we have: $\Exs \big[
  \vecnorm{\Delta_{t - k_m - \tau}}{2}^2 \big] \leq \tfrac{c \stepsize
}{1 - \kappa} \mixingtime d \varbound^2$.  If $m = 0$, we have that
$\Exs \big[ \multnoiseMarkov (\widetilde{\state}_{t + \tau }) \mid
  \filtration_t \big] = 0$ almost surely for each $t \geq
\numburn$. For $m \geq 1$, the conditional unbiasedness does not hold
true, but we still have the following upper bound on the bias
\begin{align*}
&\opnorm{\Exs \big[\prod_{j = 0}^m \multnoiseMarkov
    (\widetilde{\state}_{t + k_m + \tau - k_j}) \mid \filtration_t
    \big]} \\
    & = \sup_{u, v \in \sphere^{d - 1}} \Exs \big[
  \inprod{u}{\prod_{j = 0}^m \multnoiseMarkov (\widetilde{\state}_{t +
      k_m + \tau - k_j}) v} \big] \\
& \leq \sup_{u, v \in \sphere^{d - 1}} \Exs \big[
  \vecnorm{\multnoiseMarkov (\widetilde{\state}_{t + k_m + \tau})^\top
    u}{2} \cdot \opnorm{\prod_{j = 1}^{m - 1} \multnoiseMarkov
    (\widetilde{\state}_{t + k_m + \tau - k_j}) } \cdot
  \vecnorm{\multnoiseMarkov (\widetilde{\state}_{t + \tau}) v}{2}
  \big] \\
& \leq (\sigmaA d)^{m - 1} \sup_{u, v \in \sphere^{d - 1}} \sqrt{\Exs
  \vecnorm{\multnoiseMarkov (\widetilde{\state}_{t + k_m + \tau})^\top
    u}{2}^2 \cdot \Exs \vecnorm{\multnoiseMarkov
    (\widetilde{\state}_{t + \tau}) v}{2}^2}\\
& \leq (\sigmaA d)^{m - 1} \cdot \sigmaA^2 d.
\end{align*}
Denote $Y_t \mydefn \prod_{j = 0}^m \multnoiseMarkov (\state_{t -
  k_j})$ and $\widetilde{Y}_t \mydefn \prod_{j = 0}^m \multnoiseMarkov
(\widetilde{\state}_{t - k_j})$ for any $t \geq k_m$. We have the
expansion:
\begin{align*}
 &\Exs \big[ \vecnorm{\sum_{t = \numburn}^{\numobs - 1} Q_2 (t)}{2}^2
   \big] \\
   & \leq 2 \Exs \big[ \vecnorm{\sum_{t = \numburn}^{\numobs
       - 1} \Exs [\widetilde{Y}_t] \cdot \Delta_{t - k_m - \tau}}{2}^2
   \big] + 2 \Exs \big[ \vecnorm{\sum_{t = \numburn}^{\numobs - 1}
     (\widetilde{Y}_t - \Exs [\widetilde{Y}_t]) \cdot \Delta_{t - k_m
       - \tau}}{2}^2 \big] \\
 & \leq 2 \numobs \big( d^m \sigmaA^{m + 1} \big)^2 \sum_{t =
   \numburn}^{\numobs} \Exs \vecnorm{ \Delta_{t - k_m - \tau}}{2}^2 \\
& \qquad \quad + 2 \sum_{\numburn \leq s, t \leq \numobs - 1} \Exs
\big[ \inprod{(\widetilde{Y}_t - \Exs [\widetilde{Y}_t]) \cdot
    \Delta_{t - k_m - \tau}}{(\widetilde{Y}_s - \Exs
    [\widetilde{Y}_s]) \cdot \Delta_{s - k_m - \tau}} \big].
\end{align*}
Note that in the special case of $m = 0$, we have $\Exs
[\widetilde{Y}_t] = 0$ so that the bound holds without the first term
on the RHS. \\

\noindent For $t > s + \tau + k_m$, we have the relations
\begin{align*}
    \Exs \big[(\widetilde{Y}_t - \Exs [\widetilde{Y}_t]) \cdot
      \Delta_{t - k_m - \tau} \mid \widetilde{\filtration}_{t - k_m -
        \tau} \big] = 0, \quad \mbox{and} \quad (\widetilde{Y}_s -
    \Exs [\widetilde{Y}_s]) \cdot \Delta_{s - k_m - \tau} \in
    \widetilde{\filtration}_{t - k_m - \tau},
\end{align*}
meaning that the product term vanishes when $|s - t| > \tau +
k_m$. Therefore, we arrive at the bound
\begin{multline}
\Exs \big[ \vecnorm{\sum_{t = \numburn}^{\numobs - 1} Q_2 (t)}{2}^2
  \big] \\
   \leq \begin{cases} \big( 2n^2 \big(d^m \sigmaA^{m +
    1}\big)^2 + 4 \numobs (k_m + \tau) \cdot (\sigmaA d)^{2m} \cdot
  \sigmaA^2 d \big) \cdot \tfrac{c \stepsize}{1 - \kappa} d
  \mixingtime \varbound^2 & m \geq 1,\\ 4 \numobs \tau \sigmaA^2
  \usedim \cdot \tfrac{c \stepsize }{1 - \kappa} d \mixingtime
  \varbound^2 & m = 0.
     \end{cases}\label{eq:bound-sum-markov-mult-part-2}
\end{multline}

\noindent Now we turn to the last term in the
decomposition~\eqref{eq:cross-term-decomp-for-optimal-cov}.  We start
with the decomposition:
\begin{align*}
  \Delta_t - \Delta_{t - \tau} = \stepsize \sum_{\ell = 1}^{\tau}
  \big( \Lmat_{t - \ell + 1} (\state_{t - \ell}) \Delta_{t - \ell} +
  \addnoiseMarkov_{t - \ell} + \addnoiseMG_{t - \ell + 1} \big).
\end{align*}
We therefore have the following decomposition:
\begin{align*}
&\Exs \big[ \vecnorm{\sum_{t = \numburn }^{\numobs - 1} Q_3 (t)}{2}^2
  \big]\\
   & \leq 4\stepsize^2 \Exs \big[ \vecnorm{ \sum_{t =
      \numburn}^{\numobs} \big\{Y_t \cdot \big(\sum_{\ell = 1}^\tau
    \multnoiseMG_{t - k_m - \ell + 1} \Delta_{t - k_m - \ell} \big)
    \big\}}{2}^2 \big] + 4\stepsize^2 \Exs \big[ \vecnorm{
    \sum_{t = \numburn}^{\numobs} \big\{Y_t \cdot \big(\Lbar
    \sum_{\ell = 1}^\tau \Delta_{t - k_m - \ell} \big) \big\}}{2}^2
  \big]\\ &\qquad + 4\stepsize^2 \Exs \big[ \vecnorm{ \sum_{t =
      \numburn}^{\numobs} \big\{Y_t \cdot \big(\sum_{\ell = 1}^\tau
    \multnoiseMarkov_{t - k_m - \ell } \Delta_{t - k_m - \ell} \big)
    \big\}}{2}^2 \big] \\ &\qquad \qquad+ 4\stepsize^2 \Exs \big[
  \vecnorm{ \sum_{t = \numburn}^{\numobs} \big\{Y_t \cdot
    \big(\sum_{\ell = 1}^\tau (\addnoiseMarkov_{t - k_m - \ell} +
    \addnoiseMG_{t - k_m - \ell + 1}) \big) \big\}}{2}^2 \big]
\end{align*}
For the martingale component of the noise, note that each term $\prod_{j = 0}^m \multnoiseMarkov (
    \state_{t - k_j}) \cdot \multnoiseMG_{t - \ell + 1} (\state_{t - \ell})$ has zero conditional mean conditioned on $\filtration_{t - \ell}$. We have that
\begin{multline*}
    \Exs \big[ \vecnorm{  \sum_{t = \numburn}^{\numobs}Y_t  \multnoiseMG_{t - k_m - \ell + 1} (\state_{t - k_m - \ell}) \Delta_{t - k_m - \ell} }{2}^2 \big]
    = \sum_{t = \numburn}^{\numobs - 1} \Exs \big[ \vecnorm{ Y_t  \multnoiseMG_{t - k_m - \ell + 1} (\state_{t - k_m - \ell}) \Delta_{t - k_m - \ell}}{2}^2 \big]\\
    \leq (\sigmaA d)^{2 (m + 1)} \sum_{t = \numburn}^{\numobs - 1} \Exs \big[ \vecnorm{\multnoiseMG_{t - k_m - \ell + 1} (\state_{t - k_m - \ell}) \Delta_{t - k_m- \ell}}{2}^2 \big]
    \leq \sigmaA^{2 m + 4} d^{2m + 3} \numobs \cdot  \tfrac{c \stepsize}{1 - \kappa} d \mixingtime \varbound^2.
\end{multline*}
From the Lipschitz condition~\eqref{assume-lip-mapping} and the
boundedness condition~\eqref{assume:stationary-tail} on the metric
space, it follows that $\opnorm{Y_t} \leq (\sigmaA d)^{m + 1}$ almost
surely.  Using this fact, the second term can be bounded as
\begin{align*}
\Exs \big[ \vecnorm{ \sum_{t = \numburn}^{\numobs} \big\{Y_t \cdot
    \big(\Lbar \sum_{\ell = 1}^\tau \Delta_{t - k_m - \ell} \big)
    \big\}}{2}^2 \big] &  \leq \numobs \tau (\sigmaA d)^{2 m + 2}
\conmax^2 \sum_{t = \numburn}^{\numobs - 1} \sum_{\ell = 1}^\tau \Exs
\vecnorm{\Delta_{t - k_m - \ell}}{2}^2 \\
& \leq \numobs^2 \tau^2 (\sigmaA d)^{2 m + 2} \conmax^2 \cdot \tfrac{c
  \stepsize}{1 - \kappa} d \mixingtime \varbound^2.
\end{align*}
Collecting equations~\eqref{eq:bound-sum-markov-mult-part-1}
and~\eqref{eq:bound-sum-markov-mult-part-2} as well as the above
bounds for $Q_3$, we arrive at the upper bound \mbox{$\Exs \big[
    \vecnorm{\sum_{t = \numburn}^{\numobs - 1} \big( \prod_{j = 0}^m
      \multnoiseMarkov_{t - k_j} \big) \Delta_{t - k_m}}{2}^2 \big]
  \leq \sum_{j=1}^3 T_j$,} where
\begin{align*}
T_1 & \mydefn \numobs^2 d^{2m} \sigmaA^{2m + 2} \big(1 + \stepsize^2
\tau^2 \conmax^2 d^2 \sigmaA^2 + \stepsize^2 \tau^2 d^3 \sigmaA^2
/\numobs \big) \cdot \tfrac{c \stepsize }{1 - \kappa} d
\mixingtime \varbound^2 \\
T_2 & \mydefn 4 \stepsize^2 \Exs \big[ \vecnorm{ \sum_{t =
      \numburn}^{\numobs} \big\{Y_t \big(\sum_{\ell = 1}^\tau
    \multnoiseMarkov_{t - k_m - \ell } \Delta_{t - k_m - \ell} \big)
    \big\}}{2}^2 \big], \quad \mbox{and} \\ T_3 & \mydefn 4\stepsize^2
\Exs \big[ \vecnorm{ \sum_{t = \numburn}^{\numobs} \big\{Y_t
    \big(\sum_{\ell = 1}^\tau (\addnoiseMarkov_{t - k_m - \ell} +
    \addnoiseMG_{t - k_m - \ell + 1}) \big) \big\}}{2}^2 \big].
\end{align*}
In the special case of $m = 0$, we have:
\begin{align*}
    &\Exs \big[ \vecnorm{\sum_{t = \numburn}^{\numobs - 1}
        \multnoiseMarkov_{t} \Delta_{t}}{2}^2 \big]\\
         &\leq c
    \sigmaA^2 d \cdot \big( \numobs \tau + \numobs^2 \stepsize^2
    \sigmaA^2 d \tau^2 \big) \frac{c \stepsize }{1 - \kappa}
    d \mixingtime \varbound^2 + 4 \stepsize^2 \tau \sum_{k_1 =
      1}^{\tau} \Exs \big[ \vecnorm{ \sum_{t = \numburn}^{\numobs}
        \multnoiseMarkov_{t} \multnoiseMarkov_{t - k_1} \Delta_{t -
          k_{1}}}{2}^2 \big]\\ &\qquad + 4 \stepsize^2 \tau
    \sum_{k_{1} = 1}^{ \tau} \Exs \big[ \vecnorm{ \sum_{t =
          \numburn}^{\numobs} \multnoiseMarkov_{t} \big(
        \addnoiseMarkov_{t - k_{1}} + \addnoiseMG_{t - k_{1} + 1}
        \big) }{2}^2 \big].
\end{align*}
which completes the proof of this lemma.


\subsection{Proof of Lemma~\ref{lemma:partial-sum-prod-intricate-structure}}
\label{subsubsec:proof-lemma-intricate-structure}

We study the bias and variance of the summation separately. For the
bias term, we have:
\begin{align}
&\vecnorm{\Exs \big[\big( \prod_{j = 0}^{m - 1} \multnoiseMarkov_{t -
      k_j} \big) \big( \addnoiseMarkov_{t - k_m} + \addnoiseMG_{t -
      k_m + 1} \big) \big]}{2} \nonumber\\
       & = \sup_{z \in \sphere^{d -
    1}} \Exs \big[ \inprod{ \big( \prod_{j = 0}^{m - 1}
    \multnoiseMarkov_{t - k_j} \big) \big( \addnoiseMarkov_{t - k_m} +
    \addnoiseMG_{t - k_m + 1} \big)}{z} \big] \nonumber \\
\label{eq:bias-bound-intricate-structure}
& \overset{(i)}{\leq} \sup_{z \in \sphere^{d - 1}} \sqrt{\Exs
  \vecnorm{\multnoiseMarkov_t^\top z}{2}^2} \cdot \big[ \Exs
  \vecnorm{\big( \prod_{j = 1}^{m - 1} \multnoiseMarkov_{t - k_j}
    \big) \big(\addnoiseMarkov_{t - k} + \addnoiseMG_{t - k +
      1}\big)}{2}^2 \big]^{1/2} \nonumber \\
& \overset{(ii)}{\leq} \sigmaA \sqrt{d} \cdot (\sigmaA d)^{m - 1}
\cdot 2 \varbound \sqrt{d} = 2 (\sigmaA d)^m \varbound,
\end{align}
where step (i) uses the Cauchy--Schwarz inequality, and step (ii)
follows by invoking the moment assumption~\ref{assume:noise-moments}
as well as the Lipschitz assumption~\ref{assume-lip-mapping}.

For $t \in [k_m, \numobs]$, we define
\begin{align*}
\crosstermdebias_t \mydefn \big( \prod_{j = 0}^{m - 1}
\multnoiseMarkov_{t - k_j} \big) \big( \addnoiseMarkov_{t - k_m} +
\addnoiseMG_{t - k_m + 1} \big) - \Exs \big[\big( \prod_{j = 0}^{m -
    1} \multnoiseMarkov_{t - k_j} \big) \big( \addnoiseMarkov_{t -
    k_m} + \addnoiseMG_{t - k_m + 1} \big) \big].
\end{align*}
We have
\begin{multline*}
\Exs \big[ \vecnorm{\crosstermdebias_t}{2}^2 \big] \; \leq \; \Exs
\big[ \big( \prod_{j = 0}^{m - 1} \opnorm{\multnoiseMarkov_{t -
      k_j}}^2 \big) \cdot \vecnorm{\addnoiseMarkov_{t - {k_m}} +
    \addnoiseMG_{t - {k_m} + 1}}{2}^2 \big] \\
    \leq (\sigmaA d)^{2m}
\cdot \Exs \big[ \vecnorm{\addnoiseMarkov_{t - k} + \addnoiseMG_{t - k
      + 1}}{2}^2 \big] 
     \leq d^{2 m + 1} \sigmaA^{2m} \varbound^2.
\end{multline*}
For integers $t \geq 0$ and $\ell \geq k_m$, by
Lemma~\ref{LemWassDecay}, there exists a random variable
$\widetilde{s}_{t + \ell - k_m}$, such that $\widetilde{s}_{t + \ell -
  k_m} \mid \filtration_t \sim \stationary$, and that $\Exs \big[
  \metric (\state_{t + \ell - k_m}, \widetilde{\state}_{t + \ell -
    k_m}) \mid \filtration_t \big] \leq c_0 \cdot 2^{1 - \frac{\ell -
    k_m}{\mixingtime}}$.  By Assumption~\ref{assume-markov-mixing},
conditionally on the pair of states $(\state_{t + \ell - k_m},
\widetilde{\state}_{t + \ell - k_m} )$, we have the following bound
for $j \in [m]$:
\begin{align*}
    \Wass_{\metric, 1} \big( \transition^{k_j - k_{j - 1}}
    \delta_{\state_{t + \ell - k_j}}, \transition^{k_j - k_{j - 1}}
    \delta_{\widetilde{\state}_{t + \ell - k_j}} \big) \leq c_0
    \cdot \metric \big( \state_{t + \ell - k_j},
    \widetilde{\state}_{t + \ell - k_j} \big), \quad\mathrm{a.s.}
\end{align*}
Consequently, there exists a sequence of random variables
$(\widetilde{\state}_{t + \ell - k_j})_{0 \leq j \leq m - 1}$, such
that the following relations hold true for $j = 1,2, \cdots, m$:
\begin{align*}
    &\widetilde{\state}_{t + \ell - k_{j - 1}} \mid \filtration_{t +
    \ell - k_m} \sim \transition^{k_j - k_{j - 1}}
  \delta_{\widetilde{\state}_{t + \ell - k_j}}, \quad
  \mbox{and}\\ &\Exs \big[ \metric \big( \widetilde{\state}_{t + \ell
      - k_{j - 1}}, \state_{t + \ell - k_{j - 1}} \big) \mid
    \filtration_{t + \ell - k_m } \big] \leq c_0^{m + 1 - j} \cdot
  \metric \big( \state_{t + \ell - k_m}, \widetilde{\state}_{t + \ell
    - k_m} \big).
\end{align*}
Given the random variables constructed above, we can then construct
the proxy random variable for $\crosstermdebias_{t + \ell}$:
\begin{multline*}
    \widetilde{\crosstermdebias}_{t + \ell} \mydefn \big( \prod_{j =
      0}^{m - 1} \multnoiseMarkov (\widetilde{\state}_{t + \ell -
      k_j}) \big) \big( \addnoiseMarkov ( \widetilde{\state}_{t + \ell
      - k_m}) + \addnoiseMG_{t + \ell - k_m + 1} (
    \widetilde{\state}_{t + \ell - k_m}) \big) \\
    - \Exs \big[\big(
      \prod_{j = 0}^{m - 1} \multnoiseMarkov_{t - k_j} \big) \big(
      \addnoiseMarkov_{t - k_m} + \addnoiseMG_{t - k_m + 1} \big)
      \big].
\end{multline*}
By stationarity, we have $\Exs \big[ \widetilde{\crosstermdebias}_{t +
    \ell} \mid \filtration_t \big] = 0$ almost surely.  In order to
bound the difference, we note the telescope relation:
$\widetilde{\crosstermdebias}_{t + \ell} - \crosstermdebias_{t + \ell}
= \sum_{q = 0}^{m - 1} E^{(mix)}_q + \bar{E}^{(mix)}$, where
\begin{multline*}
E^{(mix)}_q \defn \big(\prod_{j = 0}^{q - 1} \multnoiseMarkov
(\state_{t + \ell - k_j})\big) \big( \Lbar (\widetilde{\state}_{t +
  \ell - k_q}) - \Lmap (\state_{t + \ell - k_q}) \big) \\
  \times \big(\prod_{j =
  q + 1}^{m - 1} \multnoiseMarkov (\widetilde{\state}_{t + \ell -
  k_j})\big) \big( \addnoiseMarkov ( \widetilde{\state}_{t + \ell -
  k_m}) + \addnoiseMG_{t + \ell - k_m + 1} ( \widetilde{\state}_{t +
  \ell - k_m}) \big),
\end{multline*}
and $\bar{E}^{(mix)} \defn \prod_{j = 0}^{m - 1} \multnoiseMarkov
(\state_{t + \ell - k_j}) \cdot \big( \addnoiseMarkov (
\widetilde{\state}_{t + \ell - k_m}) + \addnoiseMG_{t + \ell - k_m +
  1} ( \widetilde{\state}_{t + \ell - k_m}) - \addnoiseMarkov
({\state}_{t + \ell - k_m}) + \addnoiseMG_{t + \ell - k_m + 1} (
          {\state}_{t + \ell - k_m}) \big)$.

Using the Wasserstein distance bounds and Lipschitz
condition~\ref{assume-lip-mapping}, we find the conditional
expectation $A = \Exs \big[ \vecnorm{E^{(mix)}_q}{2} \mid \filtration_t
  \big]$ is bounded as
\begin{align*}
A & \leq (\sigmaA \usedim)^{m - 1} \Exs \big[ \opnorm{\Lmap (\state_{t
      + \ell - k_q}) - \Lmap (\widetilde{\state}_{t + \ell - k_q})}
  \cdot \vecnorm{\addnoiseMarkov (\widetilde{\state}_{t + \ell - k}) +
    \addnoiseMG_{t + \ell - k + 1} (\widetilde{\state}_{t + \ell -
      k})}{2} \mid \widetilde{\filtration}_t \big] \\
& \leq (\sigmaA d)^m \sqrt{\Exs [\metric (\state_{t + \ell - k_q},
    \widetilde{\state}_{t + \ell - k_q})^2 \mid
    \widetilde{\filtration}_t]} \cdot \sqrt{\Exs [
    \vecnorm{\addnoiseMarkov (\widetilde{\state}_{t + \ell - k}) +
      \addnoiseMG_{t + \ell - k + 1} (\widetilde{\state}_{t + \ell -
        k})}{2}^2 \mid \widetilde{\filtration}_t]}\\ &\leq (\sigmaA
d)^m c_0 \cdot 2^{1 - \frac{\ell - k_q}{2\mixingtime}} \cdot 2 d
\varbound,
\end{align*}
and the conditional expectation
$B = \Exs \big[ \vecnorm{\bar{E}^{(mix)}}{2} \mid \filtration_t
  \big]$ is bounded as
\begin{align*}
B & \leq (\sigmaA \usedim)^{m} \big(\sqrt{\Exs
  \big[\vecnorm{\addnoiseMG_{t + \ell - k + 1} (\state_{t + \ell - k})
      - \addnoiseMG_{t + \ell - k + 1} (\widetilde{\state}_{t + \ell -
        k}) }{2}^2 \mid \filtration_t \big]} \\
        &\qquad+ \sqrt{\Exs \big[
    \vecnorm{\addnoiseMarkov (\state_{t + \ell - k}) - \addnoiseMarkov
      (\widetilde{\state}_{t + \ell - k}) }{2}^2 \mid \filtration_t
    \big]} \big) \\
& \leq (\sigmaA d)^m d \varbound c_0 \cdot 2^{1 - \frac{\ell - k_m}{2
      \mixingtime}}.
\end{align*}
Consequently, we can bound the cross term as
\begin{align*}
\Exs \big[ \inprod{\crosstermdebias_t}{\crosstermdebias_{t + \ell}}
  \big] &= \Exs \big[ \inprod{\crosstermdebias_t}{\Exs
    \big[\widetilde{\crosstermdebias}_{t + \ell} \mid \filtration_t
      \big]} \big] + \Exs \big[ \inprod{\crosstermdebias_t}{\Exs
    \big[\crosstermdebias_{t + \ell} - \widetilde{\crosstermdebias}_{t
        + \ell} \mid \filtration_t \big]} \big]\\ &\leq 0 + \Exs
\big[ \vecnorm{\crosstermdebias_t}{2} \cdot \Exs
  \big[\vecnorm{\crosstermdebias_{t + \ell} -
      \widetilde{\crosstermdebias}_{t + \ell}}{2} \mid \filtration_t
    \big] \big] \\
& \leq 12 c_0 d^{m + 1} \sigmaA^m \varbound \cdot 2^{- \frac{\ell -
    k}{2 \mixingtime}} \cdot \sqrt{\Exs
  \vecnorm{\crosstermdebias_t}{2}^2} \\
& \leq 12 c_0 d^{2 m + 2} \sigmaA^{2 m} \varbound^2 \cdot 2^{-
  \frac{\ell - k}{2 \mixingtime}}.
\end{align*}
Taking $\tau = 16 \mixingtime \log (c_0 d)$, we can control the cross
terms in two different ways:
\begin{align*}
  \Exs \big[ \inprod{\crosstermdebias_t}{\crosstermdebias_{t + \ell}}
    \big] \leq \begin{cases} \sqrt{\Exs
      \vecnorm{\crosstermdebias_t}{2}^2} \cdot \sqrt{\Exs
      \vecnorm{\crosstermdebias_{t + \ell}}{2}^2} \leq d^{2m + 1}
    \sigmaA^{2m} \varbound^2, & 0 \leq \ell \leq k_m + \tau,\\ 12 c_0
    d^{2m + 2} \sigmaA^{2m} \varbound^2 \cdot 2^{- \frac{\ell - k}{2
        \mixingtime}} \leq d^{2m} \sigmaA^{2m} \varbound^2 & \ell \geq
    k_m + \tau.
     \end{cases}
\end{align*}
Summing up these terms yields
\begin{multline*}
    \Exs \big[ \vecnorm{\sum_{t = \numburn}^{\numobs - 1}
        \crosstermdebias_t}{2}^2 \big] = \sum_{t =
      \numburn}^{\numobs - 1} \Exs \vecnorm{\crosstermdebias_t}{2}^2 +
    2 \sum_{\numburn \leq t_1 < t_2 \leq \numobs - 1} \Exs \big[
      \inprod{\crosstermdebias_{t_1}}{\crosstermdebias_{t_2}} \big]\\
    \leq (k + \tau + 1) \numobs d^{2m + 1} \sigmaA^{2m} \varbound^2 +
    \numobs^2 d^{2m}\sigmaA^{2m} \varbound^2.
\end{multline*}
Combining with the bound~\eqref{eq:bias-bound-intricate-structure}, we
find that
\begin{align*}
 &\Exs \big[ \vecnorm{\sum_{t = \numburn}^{\numobs - 1} \big( \prod_{j
       = 0}^{m - 1} \multnoiseMarkov_{t - k_j} \big) \big(
     \addnoiseMarkov_{t - k_m} + \addnoiseMG_{t - k_m + 1} \big)
   }{2}^2 \big] \\
   & = \vecnorm{\sum_{t = \numburn}^{\numobs - 1}\Exs
   \big[\big( \prod_{j = 0}^{m - 1} \multnoiseMarkov_{t - k_j} \big)
     \big( \addnoiseMarkov_{t - k_m} + \addnoiseMG_{t - k_m + 1} \big)
     \big]}{2}^2 + \Exs \big[ \vecnorm{\sum_{t = \numburn}^{\numobs -
       1} \crosstermdebias_t}{2}^2 \big]\\ &\leq c \big( \numobs^2 +
 (k_m + \tau) \numobs \usedim \big) \sigmaA^{2m} d^{2m}\varbound^2,
\end{align*}
for a universal constant $c > 0$.


\section{Proof of Theorem~\ref{thm:local-minimax}}
\label{subsec:proof-lower-bound}

Our strategy is to prove a Bayes risk lower bound. We construct a
prior distribution over transition kernels by perturbing
the base matrix $\transition_0$ appropriately. We then apply the
Bayesian Cram\'{e}r--Rao lower bound to obtain our result.

Let us describe the construction in more detail. For each $\state
\in \statespace$, suppose we have a perturbation vector
$\perturb_\state \in \real^\statespace$. Use these to define the
perturbed transition kernel
\begin{align*}
 \transition_\perturb (x, y) \mydefn \tfrac{\transition_0 (x, y) \,
   e^{\perturb_x (y)}}{\sum_{z \in \statespace} \transition_0 (x, z)
   e^{\perturb_x (z)}} \qquad \text{ for each } x, y \in \statespace.
\end{align*}
Note that by construction, for any $x \in \statespace$ and any
$\perturb_x \in \real^{\statespace}$, we have \mbox{$\mathrm{supp}
  \big( \transition_\perturb (x, \cdot) \big) = \mathrm{supp} \big(
  \transition_0 (x, \cdot) \big)$.} Since $\transition_0$ is
irreducible and aperiodic, so is $\transition_\perturb$. Therefore,
the stationary distribution $\stationary_\perturb$ of
$\transition_\perturb$ exists and is unique. When the perturbation is
small enough, a quantitative perturbation principle can be obtained,
which we collect in Lemma~\ref{lemma:perturbation-principle-markov} below. 



It remains to specify how the perturbation vectors are generated. We
parameterize $\perturb$ with a linear transformation, writing $h = Qw$
for a linear operator $Q$ to be specified shortly, and a random vector
$w \in \mathbb{R}^d$ drawn from a distribution $\rho$. In particular,
given a collection of vectors $\{q_{x} (y)\}_{x , y \in \statespace}
\subseteq \real^d$, we consider the linear transformation $Q: \real^d
\rightarrow \real^{\statespace\times \statespace}$ given by $w \mapsto
\big[ \inprod{w}{q_x (y)} \big]_{x, y \in \statespace}$.

Next we specify the prior $\rho$, and along with some associated
notation.  Define the subspace $\statnullspace_\perturb \mydefn \big\{
f \in \real^\statespace: \Exs_{\stationary_\perturb}[f (\state)] = 0
\big\}$, and note that $\transition_\perturb$ maps
$\statnullspace_\perturb$ to itself.  Furthermore, since
$\transition_\perturb$ is irreducible and aperiodic, the mapping $(I -
\transition_\perturb)$ is invertible on $\statnullspace_\perturb$.
Consequently, for any function $f : \statespace \rightarrow \real$,
the following Green function operator is well-defined:
\begin{align*}
    \greenOp_\perturb f \mydefn (I - \transition_\perturb)^{-1}
    \big|_{\statnullspace_\perturb} \cdot \big(f -
    \Exs_{\stationary_\perturb} [f] \big) \in \real^{\statespace}.
\end{align*}

We also define an operator $\transitionOp_\perturb$ on the space of
real-valued functions on $\statespace$ as follows:
\begin{align*}
 \transitionOp_\perturb f(x) \mydefn \Exs_{Y \sim \transition_\perturb
   (x, \cdot)} [f(Y)].
\end{align*}
Importantly,
$\transitionOp_\perturb$ is an operator mapping functions to
functions, and distinct from the matrix $\transition_\perturb$. It is
straightforward to see that the operator $\transitionOp_\perturb$
commutes with the operator $\greenOp_\perturb$, for any perturbation
matrix $\perturb$. Indeed, if we denote $\mathcal{L}_\perturb \mydefn \transitionOp_\perturb - I$ as the generator. The green function $\greenOp_\perturb f$ solves the Poisson equation $- \mathcal{L}_\perturb u = f - \Exs_{\stationary_\perturb}[f (\state)]$.

Finally, for any $\perturb \in \real^{\statespace \times \statespace}$ and for all $x \in \statespace$, we define 
\begin{align}
\label{eq:def-poisson-eq-in-lower-bound-proof}  
\greenFunc_\perturb (x) = \big(I_d - \Exs_{\stationary_\perturb}
          [\Lmap (\state)] \big)^{-1} \big( \greenOp_\perturb \Lmap
          (x) \cdot \thetabar (\transition_\perturb) +
          \greenOp_\perturb \bmap (x) \big).
\end{align}

Since the proof works under the perturbed probability transition
kernel $\transition_\perturb$, it is useful to study the effect of
small perturbation on its stationary distribution. The following lemma
provides non-asymptotic bounds on the mixing time of perturbed Markov
chain and its stationary distribution $\stationary_\perturb$, which
will be useful throughout the proof.

\begin{lemma}\label{lemma:perturbation-principle-markov}
    Under the setup above, suppose that $\perturb_{\max} \mydefn \max_{x \in \statespace} \vecnorm{\perturb_x}{\infty} < \frac{1}{128 \mixingtime}$. Then the perturbed transition kernel satisfies the following.
    \begin{itemize}
        \item The Markov transition kernel $\transition_\perturb$
          satisfies the mixing condition
          (Assumption~\ref{assume-markov-mixing}) with the discrete
          metric and mixing time $4 \mixingtime$.
        \item The stationary distribution $\stationary_\perturb$
          satisfies the bound
\begin{align*}
  \max_{\state \in \statespace} \big\{ \log \tfrac{\stationary_0
    (\state)}{\stationary_\perturb (\state)},~ \log
  \tfrac{\stationary_\perturb (\state)}{\stationary_0 (\state)} \big\}
  \leq \mixingtime \big(2 + \log \perturb_{\max}^{-1} +
  \log\tfrac{1}{\min_x \stationary_0 (x)} \big)\perturb_{\max}.
\end{align*}
    \end{itemize}
\end{lemma}
\noindent See
Section~\ref{subsubsec:proof-perturbation-principle-markov} for the
proof of this lemma.

With this notation in hand, we are ready to construct the prior distribution on $w$.
We begin with the following one-dimensional density function, taken from~\cite{tsybakov2008introduction}:
\begin{subequations}\label{eqs:prior-density-in-lower-bound-proof}
\begin{align}
    \mu(t) \mydefn \cos^2 \big(
  \tfrac{\pi t}{2} \big) \cdot \bm{1}_{t \in [-1, 1]}.\label{eq:tsyabkov-density-in-lower-bound-proof}
\end{align}
Also, define the positive-definite matrix \mbox{$\Lambda \mydefn
  \Exs_{X \sim \stationary_0} \big[ \cov_{Y \sim \transition_0 (X,
      \cdot)} \big(\greenFunc_0 (Y) \mid X \big) \big]$,} and let
$\Lambda = U D U^\top$ denote its eigen-decomposition.  For a random
variable $\psi \sim \mu^{\otimes d}$, define the perturbation
parameter
\begin{align}
\label{eq:prior-density-construction}  
  w = \tfrac{1}{\sqrt{\numobs}} U D^{-1/2} \psi,
\end{align}
and let its density denote the prior distribution $\rho$. Note that
for any $w \in \mathrm{supp} (\rho)$, we have
\begin{align}
  \vecnorm{\Lambda w}{2} = \vecnorm{U D^{1/2} \psi}{2} =
  \vecnorm{D^{1/2} \psi}{2} \leq \sqrt{\trace(D) / n} =
  \sqrt{\trace (\Lambda) / n}.\label{eq:prior-support-size}
\end{align}
The final ingredient in our construction is to specify the linear
transformation $Q$. For each $x, y \in \statespace$, we set
\begin{align}
 q_x(y) \mydefn \greenFunc_0 (y) - \Exs_{\state' \sim \transition_0
   (x, \cdot)} \big[ \greenFunc_0 (\state') \big],
\end{align}
\end{subequations}
where the Green function $\greenFunc$ is defined in
equation~\eqref{eq:def-poisson-eq-in-lower-bound-proof}. Recall that
$h = Qw$ for $w \sim \rho$. This specifies our prior over transition
kernels, and concludes the construction.

Next, we state the version of the Bayesian Cram\'{e}r--Rao bound that
we use. Before stating the result, it is useful to introduce the
general setup and basic notation for parametric models. Given a
family $\mathcal{P}_{\Theta} = \big( \Prob_\eta: \eta \in \Theta
\big)$ of probability distributions of sample $X \in \Xspace$,
parameterized by $\eta \in \Theta$, where $\Theta$ is an open subset of
$\real^d$. Assume that each element in this family is absolute
continuous with respect to a base measure $\lambda$ over $\Xspace$,
and denote the Radon--Nikodym derivative by $p_\eta \mydefn \tfrac{d
  \Prob_\eta}{d \lambda}$. Assuming differentiability and
integrability of relevant quantities, for any $\eta \in \Theta$, we
define the Fisher information matrix $I(\eta)$ as
\begin{align*}
  I (\eta) \mydefn \Exs_{X \sim \Prob_\eta} \big[ \nabla_\eta \log p_\eta (X) \nabla_\eta \log p_\eta (X)^\top  \big] \in \real^{d \times d}.
\end{align*}
Now we are ready to state the Bayesian Cram\'{e}r--Rao lower bound.
\begin{proposition}[Theorem 1 of~\cite{gill1995applications}, special case]
  \label{prop:van-tree-functional}
  Under the setup above, given a prior distribution $\rho$ with
  continuously differentiable density and bounded support contained
  within $\Theta$, let $T: \support (\rho) \mapsto \real^d$ denote a
  locally continuously differentiable functional.  Then for any
  estimator $\widehat{T}$ based on observing $X$, we have
\begin{align}
\underset{\eta \sim \rho}{\Exs} \;\; \underset{X \sim p_\eta}{\Exs}
\vecnorm{\widehat{T} (X) - T (\eta)}{2}^2 \geq \tfrac{\big( \int
  \trace \big(\tfrac{\partial T}{\partial \eta} (\eta) \big) \rho
  (\eta) \usedim \eta \big)^2}{\int \trace \big( I (\eta) \big) \rho
  (\eta) d \eta + \int \vecnorm{\nabla \log \rho (\eta)}{2}^2 \rho
  (\eta) d \eta}.
\end{align}
\end{proposition}

In order to complete the proof, we provide non-asymptotic estimates on
the three quantities involved in the right-hand-side of
Proposition~\ref{prop:van-tree-functional}. These require a few
technical lemmas, whose proofs can be found at the end of the section.

\paragraph{Bounds on the term $\trace \big( \nabla_w \thetabar  \big)$:} We state two technical lemmas that are helpful in bounding this quantity.
The first computes the Jacobian matrix of the desired functional $\thetabar (\perturb)$ with respect to the parameter $w$.
\begin{lemma}
\label{lemma:jacobi-wrt-parameter}
Under the given set-up, for any $w \in \real^d$, we have
\begin{align}
\label{eq:jacobian-eval}
\nabla_w \thetabar (\transition_\perturb) = \Exs_{X \sim
  \stationary_\perturb} \Big[ \cov_{Y \sim \transition_\perturb (X,
    \cdot) } \big \{ \greenFunc_\perturb(Y) - \transitionOp_\perturb
  \greenFunc_\perturb(X), \greenFunc_0(Y) - \transitionOp_0
  \greenFunc_0(X) \mid X \big \} \Big].
    \end{align}
\end{lemma}
\noindent See Section~\ref{subsubsec:proof-lemma-jacobi-wrt-parameter} for the proof of this lemma.
Next, we control the RHS of equation~\eqref{eq:jacobian-eval} by replacing $\greenFunc_\perturb$ with $\greenFunc_0$.
\begin{lemma}
  \label{lemma:green-func-error-estimate}
Under the given set-up and for a sample size lower bounded as $\numobs
\geq \tfrac{c \mixingtime^2 \sigmaA^2 d^2 \log^2 d}{(1 - \kappa)^2}$
and $\max_{x \in \statespace} \vecnorm{\perturb_x}{\infty} \leq
\tfrac{1}{128 \mixingtime}$, we have
\begin{align*}
  \Exs_{Z \sim \stationary_\perturb} \big[
    \vecnorm{\greenFunc_\perturb (Z) - \greenFunc_0 (Z)}{2}^2 \big]
  \leq \tfrac{c (1 + \sigmaA^2) \varbound^2 \mixingtime^4 d^2}{(1 -
    \kappa)^4 \numobs} \log^6 \tfrac{d}{\min_x \stationary_0 (x)}.
 \end{align*}
Furthermore, for any $w$ in the support of $\rho$, we have
\begin{align*}
  \vecnorm{\thetabar (\transition_\perturb) - \thetabar
    (\transition_0)}{2} \leq \tfrac{3}{2} \sqrt{\mathrm{trace}
    (\Lambda) / \numobs} + \sqrt{\tfrac{c (1 + \sigmaA^2) \varbound^2
      \mixingtime^4 d^3}{(1 - \kappa)^4 \numobs^2} \log^6
    \tfrac{d}{\min_x \stationary_0 (x)}}.
\end{align*}
\end{lemma}
\noindent See
Section~\ref{subsubsec:proof-lemma-green-func-error-estimate} for the
proof of this lemma. \\

Combining these two lemmas
yields
\begin{align*}
    &\trace \big(\nabla_w \thetabar \big) \\ &\geq \Exs_{X \sim
    \stationary_\perturb} \big[ \var_{Y \sim \transition_\perturb (X,
      \cdot) } \big( \greenFunc_0 (Y) - \transitionOp_0 \greenFunc_0
    (X) \mid X \big) \big]\\ &\qquad- \Exs_{X \sim
    \stationary_\perturb} \big[ \sqrt{\var_{Y \sim
        \transition_\perturb (X, \cdot) } \big( \greenFunc_0 (Y) -
      \transitionOp_0 \greenFunc_0 (X) \mid X \big)} \big] \cdot
  \sqrt{\Exs_{Z \sim \stationary_\perturb} \big[
      \vecnorm{\greenFunc_\perturb (Z) - \greenFunc_0 (Z)}{2}^2 \big]}
  \\ &\geq \trace \big( \Lambda \big) - \sqrt{\trace \big( \Lambda
    \big)} \cdot \tfrac{c (1 + \sigmaA) \varbound \mixingtime^2 d}{(1 -
    \kappa)^2 \sqrt{\numobs} } \log^3 \tfrac{d}{\min_x \stationary_0
    (x)}.
\end{align*}
Now given a sample size lower bounded as $\numobs \geq \tfrac{c
  \mixingtime^2 \sigmaA^2 d^2 \log^2 d}{(1 - \kappa)^2} + \tfrac{2c (1
  + \sigmaA^2) \varbound^2 \mixingtime^4 d^2}{(1 - \kappa)^4 \trace
  (\Lambda)} \log^6 \tfrac{d}{\min_x \stationary_0 (x)}$, we can
conclude that
\begin{align}
  \label{eq:jacobian-trace}
\trace \big(\nabla_w \thetabar \big) \geq \tfrac{1}{2} \trace
(\Lambda) \qquad \mbox{for any $w$ in the support of $\rho$.}
\end{align}

\paragraph{Bounds on the Fisher information $I^{(\numobs)}(w)$:}

We now state an upper bound on the Fisher information of the observed
trajectory:
\begin{lemma}
  \label{lemma:fisher-info-for-obs}
Under the given set-up, for any $w \in \real^d$, if $\perturb_{\max}
\mydefn \max_{x} \vecnorm{\perturb}{\infty}$ satisfies the inequality
${\perturb}_{\max}^{-1} \geq c \mixingtime \big(\log
\perturb_{\max}^{-1} + \log (\min \stationary_0)^{-1} \big)$, we have
\begin{align*}
  I^{(n)}(w) \mydefn \Exs_{\perturb} \big[ \nabla_w \log
    \Prob_{\perturb} \big( \state_0^{\numobs} \big) \nabla_w \log
    \Prob_{\perturb} \big( \state_0^{\numobs} \big)^\top \big] \preceq
  \tfrac{3 \numobs}{2} \Exs_{X \sim \stationary_\perturb} \big[
    \cov_{Y \sim \transition_\perturb (X, \cdot)} \big( q_X (Y) \mid X
    \big) \big].
    \end{align*}
\end{lemma}
\noindent See Section~\ref{subsubsec:proof-lemma-fisher-info-for-obs}
for the proof of this lemma. \\

In order to apply the preceding lemma, we must verify the condition on
$h_{\max}$ for our setting.  Under our construction, we have $\max_{x
  \in \statespace} \vecnorm{\perturb_x}{\infty} = \max_{x, y
  \in \statespace} \inprod{\greenFunc_0 (y) - \transitionOp_0
  \greenFunc_0 (x)}{w}$.  Note that
Assumption~\ref{assume:noise-moments} and
Lemma~\ref{lemma:poisson-eq-estimate} in Section~\ref{appendix:useful-moment-bound} together imply
the following bound for any $\delta > 0$:
\begin{align*}
\stationary_0 \big( \state : \abss{\inprod{\greenFunc_0 (s)}{w}} \leq
\tfrac{c \varbound \mixingtime \vecnorm{w}{2}}{1 - \kappa} \cdot \log^3
\tfrac{d}{\delta} \big) > 1 - \delta.
\end{align*}
Taking $\delta \mydefn \tfrac{1}{2} \min_{s \in \statespace} \stationary_0 (s) > 0$, we have the uniform bound
\begin{align*}
    \max_{s \in \statespace} \abss{\inprod{\greenFunc_0 (s)}{w}} \leq \tfrac{c \varbound \mixingtime \vecnorm{w}{2}}{1 - \kappa} \log^3 \big(d / \min_{s} \stationary_0 (s) \big).
\end{align*}
Note that $\transitionOp_0$ is a probability transition kernel, for
any $s \in \statespace$, the vector $\transitionOp_0 \greenFunc_0 (s)$
lies in the convex hull of $\big( \greenFunc_0 (s')\big)_{s'
  \in \statespace}$. So we have the bound $\max_{s \in \statespace}
\abss{\inprod{\transitionOp_0 \greenFunc_0 (s)}{w}} \leq \max_{s
  \in \statespace} \abss{\inprod{\greenFunc_0 (s)}{w}} \leq \tfrac{c
  \varbound \mixingtime \vecnorm{w}{2}}{1 - \kappa} \log^3 \big(d /
\min_{s} \stationary_0 (s) \big)$.  Putting them together leads to the
bound
\begin{align*}
\max_{x \in \statespace} \vecnorm{\perturb_x}{\infty} \leq 2 c
\varbound \mixingtime \vecnorm{w}{2} \log^3 \big(d / \min_{s}
\stationary_0 (s) \big).
\end{align*}
Now given a sample size
\begin{align}
\label{eq:sample-size-requirement-for-stability-in-lower-bound}  
\numobs \geq  c \mixingtime^3 \varbound^2 \cdot  \trace (\Lambda)\cdot
\log^3 \tfrac{d}{\min_s \stationary_0 (s)},
\end{align}
we have that $\max_{x} \vecnorm{\perturb_x}{\infty} < \tfrac{1}{128
  \mixingtime}$. This satisfies the condition in
Lemma~\ref{lemma:perturbation-principle-markov} in the
appendix. Applying this lemma, we see that the condition
\begin{align*}
{\perturb}_{\max}^{-1} \geq c \mixingtime \big(\log
\perturb_{\max}^{-1} + \log (\min \stationary_0)^{-1} \
\end{align*}
is satisfied, so that Lemma~\ref{lemma:fisher-info-for-obs} guarantees
that
\begin{align}
\trace \big( I^{(\numobs)}(w) \big) &\preceq \tfrac{3 \numobs}{2}
\Exs_{X \sim \stationary_\perturb} \big[ \var_{Y \sim
    \transition_\perturb (X, \cdot) } \big( \greenFunc_0 (Y) -
  \transitionOp_0 \greenFunc_0 (X) \mid X \big) \big] \notag
\\ &\preceq \big( \tfrac{3}{2} \big)^3 \numobs \cdot \Exs_{X \sim
  \stationary_0} \big[ \var_{Y \sim \transition_0 (X, \cdot) } \big(
  \greenFunc_0 (Y) \mid X \big) \big] \notag \\ &= \tfrac{27
  \numobs}{8} \trace \big( \Lambda
\big). \label{eq:fisher-info-channel}
\end{align}
The last inequality follows because $\stationary_\perturb \preceq
\tfrac{3}{2} \stationary_0$ $\transition_\perturb (x, \cdot) \preceq
\tfrac{3}{2} \transition_0 (x, \cdot)$ for all $x \in \statespace$.

\paragraph{Bounds on the prior Fisher information:}

From Lemma 10 in~\cite{mou2020optimal}, the density $\rho$
of $w$ has Fisher information
\begin{align}
\label{eq:prior-fisher-info}  
  I(\rho) = U D^{1/2} I \big( \mu^{\otimes d} \big) D^{1/2} U^\top =
  n \pi \Lambda.
\end{align}
Consequently, we have $\int \vecnorm{\nabla \log \rho (w)}{2}^2
\rho(w) ~dw \trace \big( I(\rho) \big) = n \pi \cdot \trace
(\Lambda)$.

\paragraph{Putting together the pieces:}

Combining the bounds~\eqref{eq:jacobian-trace},
\eqref{eq:fisher-info-channel}, and \eqref{eq:prior-fisher-info} and
applying Proposition~\ref{prop:van-tree-functional}, we obtain the
lower bound
\begin{align}
\inf_{\thetahat_\numobs} ~\int_{\real^d} \Exs_{X_1^\numobs \sim
  \Prob_{Q w}} \big[ \vecnorm{\thetahat_\numobs - \thetabar
    (\transition_{Q w})}{2}^2 \big] \rho (d w) \geq \tfrac{1}{4 (5 +
  \pi) \numobs} \trace
(\Lambda).\label{eq:bayes-risk-lower-bound-final}
\end{align}
It remains to relate the matrix $\Lambda$ to the local complexity
$\varepsilon_\numobs$ in the theorem. In order to do so, we require
the following lemma.
\begin{lemma}
\label{lemma:poisson-eq-identity}
Under the setup above, for any function $f : \statespace \rightarrow
\real$ such that $\Exs_{\stationary_0} [f (s)] = 0$, we have $\Exs_{X
  \sim \stationary_0, Y \sim \transition_0 (X, \cdot)} \big[ \big(
  \greenOp_0 f(Y) - \transitionOp_0 \greenOp_0 f(X) \big)^2 \big] =
\sum_{k = - \infty}^{\infty} \Exs \big[ f (s_0) f (s_k) \big]$, where
$(s_k)_{k \in \mathbb{Z}}$ is a stationary Markov chain following
$\transition_0$.
\end{lemma}
\noindent See Section~\ref{subsubsec:proof-lemma-poisson-eq-identity}
for the proof of this lemma. \\

Applying Lemma~\ref{lemma:poisson-eq-identity} with $f_j(s) =
\inprod{(I_d - \Lbar^{(0)})^{-1} \big( \Lmap (s) \thetabar
  (\transition_0) + \bmap (s)\big)}{\coordinate_j}$ for $j = 1,2,
\cdots, d$ respectively, we arrive at the chain of equalities
\begin{align*}
  \trace (\Lambda)& = \sum_{j = 1}^d \Exs_{X \sim \stationary_0, Y
    \sim \transition_0 (X, \cdot)} \big[ \big( \greenOp_0 f_j (Y) -
    \transitionOp_0 \greenOp_0 f_j (X) \big)^2 \big] \\ & = \sum_{j =
    1}^d \sum_{k = - \infty}^{\infty} \Exs \big[ f_j (s_0) f_j (s_k)
    \big] = \mathrm{trace} \big( (I - \Lbar^{(0)})^{-1}
  \SigStar_{\mathrm{Mkv}} (I - \Lbar^{(0)})^{-\top} \big) = \numobs
  \varepsilon_\numobs^2.
\end{align*}
Thus, the right-hand-side of
equation~\eqref{eq:bayes-risk-lower-bound-final} is exactly $\tfrac{
  \varepsilon_\numobs^2}{4 (5 + \pi)}$.

It remains to bound the size of the neighborhood. Given a sample size
$\numobs$ satisfying the
bound~\eqref{eq:sample-size-requirement-for-stability-in-lower-bound},
Lemma~\ref{lemma:green-func-error-estimate} implies that
$\vecnorm{\thetabar (\transition_\perturb) - \thetabar
  (\transition_0)}{2} \leq \sqrt{\tfrac{\trace (\Lambda)}{\numobs}}$.
Consequently, for any $w$ on the support of $\rho$, we have
$\transition_{Q w} \in \neighborhood_{\mathrm{Est}} (\transition_0, 2
\varepsilon_\numobs)$.

On the other hand, for any $w \in \mathrm{supp} (\rho)$ and any $x
\in \statespace$ and perturbation $\perturb = Q w$, we have
\begin{align*}
  \chisqdiv{\transition_{\perturb} (x, \cdot)}{\transition_0 (x,
    \cdot)} &= \Exs_{Y \sim \transition_0 (x, \cdot)} \big[ \big(
    \tfrac{\transition_\perturb (x, Y)}{\transition_0 (x, Y)} - 1
    \big)^2 \big] \\
  & = \var_{Y \sim \transition_0 (x, \cdot)}  \Big( \tfrac{e^{\perturb_x
      (Y)}}{\sum_{z \in \statespace}
  \transition_0 (x, z) e^{\perturb_x (z)}} \Big)\\
  & \overset{(i)}{\leq} \var_{Y \sim \transition_0 (x, \cdot)} 
    \big( e^{\perturb_x (Y)} \big) \\
    &\leq \Exs_{Y \sim \transition_0 (x, \cdot)} 
    \big[ \big( e^{\perturb_x (Y)} - 1 \big)^2 \big]
  \\
& \overset{(ii)}{\leq} e \cdot \Exs_{Y \sim \transition_0 (x, \cdot)}
  \big[\perturb_x (Y)^2\big],
\end{align*}
where step (i) follows by using Jensen's inequality to assert that
\begin{align*}
\sum_{z \in \statespace} \transition_0 (x, z) e^{\perturb_x (z)} \geq
e^{\sum_{z \in \statespace} \transition_0 (x, z) \perturb_x (z)} = 1,
\end{align*}
and step (ii) follows from the inequality $\abss{e^{x} - 1} \leq e
\cdot |x|$, valid for $x \in [-1, 1]$.

Accordingly, the average $\chi^2$-divergence admits the bound
\begin{align*}
\sum_{x \in \statespace} \stationary_0 (x)
\chisqdiv{\transition_{\perturb} (x, \cdot)}{\transition_0 (x, \cdot)}
&\leq e \cdot \Exs_{X \sim \stationary_0 , Y \sim \transition_0 (X,
  \cdot)} \big[ \inprod{w}{\greenFunc_0 (Y) -\transitionOp_0
    \greenFunc_0 (X)}^2\big] \\ & \leq e \cdot w^\top \Lambda w \leq
\tfrac{e d}{\numobs}.
\end{align*}
For any $w$ on the support of $\rho$, we thus have $\transition_{Q w}
\in \neighborhood_{\mathrm{Prob}} (\transition_0, e
\sqrt{\tfrac{d}{n}})$, as claimed. The Bayes risk lower
bound~\eqref{eq:bayes-risk-lower-bound-final} then implies the desired
minimax lower bound.


\subsection{Proof of Lemma~\ref{lemma:perturbation-principle-markov}} \label{subsubsec:proof-perturbation-principle-markov}

The proof relies on a total variation distance bound on the transition
kernel. In particular, for each $\state \in \statespace$, we have
\begin{multline}
    \totalvariation \big( \transition_0 (x, \cdot),
    \transition_\perturb (x, \cdot) \big) \leq \sqrt{\tfrac{1}{2}
      \chisqdiv{\transition_0 (x, \cdot)}{\transition_\perturb (x,
        \cdot)}} = \sqrt{\tfrac{1}{2} \sum_{y \in \statespace}
      \transition_0 (x, y) \cdot \big( \tfrac{\transition_\perturb (x,
        y)}{\transition_0 (x, y)} - 1 \big)^2}\\ \overset{(i)}{\leq}
    \sqrt{\tfrac{1}{2} \big( e^{\vecnorm{\perturb_x}{\infty}} - 1
      \big)^2} \overset{(ii)}{\leq} e \cdot \max_{x
      \in \statespace}\vecnorm{\perturb_x}{\infty}.\label{eq:from-perturb-l-infty-to-tv-bound}
\end{multline}
In step $(i)$, we use the fact
\begin{align*}
    \frac{\transition_\perturb (x, y)}{\transition_0 (x, y)} = \frac{e^{h_x (y)}}{\sum_{z \in \statespace} \transition_0 (x, y) e^{h_x (z)}} \in [e^{- \vecnorm{h_x}{\infty}}, e^{\vecnorm{h_x}{\infty}}],
\end{align*}
and in step $(ii)$, we use the fact $\vecnorm{\perturb_x}{\infty} < 1$.

Next, we turn to proofs of the two claims. We first prove the mixing
time bound. Note that the non-expansive
condition~\eqref{eq:assume-mixing-tmix}(b) is automatically satisfied
with $c_0 = 1$ for total variation distance (by a na\"{i}ve
coupling). Given a fixed pair $x, y \in \statespace$, invoking
Lemma~\ref{LemWassDecay} with $\tau = 4 \mixingtime$ yields the
existence of a joint distribution over the random sequence $\{x_k\}_{0
  \leq k \leq \tau}$ and $\{y_k\}_{0 \leq k \leq \tau}$, such that
$\{x_k\}$ and $\{y_k\}$ follows the Markov chain $\transition_0$,
starting from $x_0 = x$ and $y_0 = y$, respectively. Furthermore, we
have the bound $\Prob \big(x_{\tau} \neq y_{\tau} \big) \leq
\tfrac{1}{4}$.

Now we construct a coupling between the original chain and perturbed
chain. Taking the initial point $\widetilde{x}_0 = x$, we iteratively
construct the sequence $\{\widetilde{x}_k\}_{0 \leq k \leq \tau}$ as
follows: given $\widetilde{x}_k$ and $x_k$, we construct the
conditional distribution of $\widetilde{x}_{k + 1}$ as follows:
\begin{itemize}
    \item If $x_{k} = \widetilde{x}_k$, we let $\Prob
      \big(\widetilde{x}_{k + 1} \neq x_{k + 1} \mid x_k,
      \widetilde{x}_k \big) = \totalvariation \big( \transition_0
      (x_k, \cdot), \transition_\perturb (x_k, \cdot) \big)$.
    \item If $x_k \neq \widetilde{x}_k$, we simply take
      $\widetilde{x}_{k + 1}$ and $x_{k + 1}$ to be conditionally
      independent, following their respective transition kernels.
\end{itemize}
We construct the sequence $\{\widetilde{y}_k\}_{0 \leq k \leq \tau}$
in a similar fashion.

By the union bound, it follows that
\begin{align*}
    \Prob \big( x_{\tau} \neq \widetilde{x}_{\tau} \big) \leq \sum_{k
      = 0}^{\tau - 1} \Exs \big[ \Prob \big( x_{k + 1} \neq
      \widetilde{x}_{k + 1} \mid x_{k} = \widetilde{x}_k \big) \big]
    &=\sum_{k = 0}^{\tau - 1} \Exs \big[ \totalvariation \big(
      \transition_0 (x_k, \cdot), \transition_\perturb (x_k, \cdot)
      \big) \big] \\ &\leq 4 e \mixingtime \cdot \max_{x
      \in \statespace} \vecnorm{\perturb_x}{\infty} < \tfrac{1}{8}.
\end{align*}
In the last step, we have used the total variation distance bound~\eqref{eq:from-perturb-l-infty-to-tv-bound}.

Similarly, the process $\{\widetilde{y}_k\}$ satisfies the bound
$\Prob \big( y_{\tau} \neq \widetilde{y}_{\tau} \big) < \tfrac{1}{8}$.
Putting together the pieces, we conclude that
\begin{align*}
\totalvariation \big( \delta_x \transition_\perturb^{\tau}, \delta_y
\transition_\perturb^{\tau} \big) \leq \Prob \big(
\widetilde{x}_{\tau} \neq \widetilde{y}_\tau \big) & \leq \Prob \big(
\widetilde{x}_{\tau} \neq x_\tau \big) + \Prob \big( x_{\tau} \neq
y_\tau \big) + \Prob \big( y_\tau \neq \widetilde{y}_{\tau} \big)\\
& < \tfrac{1}{8} + \tfrac{1}{4} + \tfrac{1}{8} = \tfrac{1}{2},
\end{align*}
which shows that the perturbed chain $\transition_\perturb$ satisfies
the condition~\eqref{eq:assume-mixing-tmix}(a) with mixing time $\tau
= 4 \mixingtime$.

Next, we prove the perturbation result for the stationary
distribution.
Given any fixed initial distribution $\pi_0$, note that
for any deterministic sequence $(x_0, x_2, \cdots, x_n)$, we have the
following expression for the Radon-Nikodym derivative:
\begin{align*}
    \tfrac{d \Prob_\perturb \big(x_0, x_1, \cdots, x_n \big)}{d \Prob_0
      \big(x_0, x_1, \cdots, x_n \big)} = \prod_{k = 0}^{n - 1}
    \tfrac{\transition_\perturb (x_k, x_{k + 1})}{\transition_0 (x_k,
      x_{k + 1})} = \prod_{k = 0}^{n - 1} \tfrac{e^{\perturb_{x_k}
        (x_{k + 1})}}{\sum_{y \in \statespace} e^{\perturb_{x_k} (y)}
      \transition (x_k, y)}.
\end{align*}
We then have the max-divergence bound
\begin{align*}
    \maxdiv{ \Prob_\perturb \big(x_0^n \big)}{ \Prob_0 \big(x_0^n \big)} := \sup_{x_0^n \in \statespace^n} \abss{\log \tfrac{d \Prob_\perturb \big(x_0, x_1, \cdots, x_n \big)}{d \Prob_0 \big(x_0, x_1, \cdots, x_n \big)}} \leq n \cdot \max_{x} \vecnorm{\perturb_x}{\infty}.
\end{align*}
Taking the marginal distribution, we see that the bound $\maxdiv{
  \pi_0 \transition_\perturb^n}{\pi_0 \transition_0^n} \leq n \cdot
\perturb_{\max}$ holds for any initial distribution $\pi_0$ and any $n
> 0$.

To obtain the desired claim, we take the initial distribution to be the stationary
distribution $\stationary_\perturb$ of the chain
$\transition_\perturb$, and let $n = \mixingtime \log \big( \tfrac{2}{
  \perturb_{\max} \cdot \min_x \stationary_0 (x)} \big)$. Note that
$\stationary_\perturb \transition_\perturb^n = \stationary_\perturb$ in such case. On
the other hand, by Lemma~\ref{LemWassDecay}, the total variation
distance can be upper bounded as $\totalvariation\big( \stationary_\perturb
\transition_0^n, \stationary_0 \big) \leq 2^{1 -
  \tfrac{n}{\mixingtime}} \leq \perturb_{\max} \cdot \min_{x
  \in \statespace} \stationary_0(x)$. Therefore, for any $x \in \statespace$, we have
\begin{align*}
\abss{\tfrac{\stationary_\perturb \transition_0^n (x)}{ \stationary_0 (x)} - 1} \leq
\tfrac{ \totalvariation \big( \stationary_\perturb \transition_0^n, \stationary_0
  \big) }{ \min_{x \in \statespace} \stationary_0 (x)} \leq
\perturb_{\max} < \frac{1}{2}.
\end{align*}
Invoking the inequality $|\log z| \leq 2 |z - 1|$ for $|z| \leq 1/2$, we can translate the bound into
a max-divergence bound
\begin{align*}
\maxdiv{ \stationary_\perturb \transition_0^n}{ \stationary_0} = \max_{x
  \in \statespace} \abss{\log \tfrac{\stationary_\perturb \transition_0^n (x)}{
    \stationary_0 (x)}} \leq
2 \perturb_{\max}.
\end{align*}
Finally, applying the triangle inequality yields
\begin{align*}
\maxdiv{ \stationary_\perturb}{ \stationary_0} &\leq \maxdiv{ \stationary_\perturb
  \transition_\perturb^n}{\stationary_\perturb \transition_0^n} + \maxdiv{ \stationary_\perturb
  \transition_0^n}{ \stationary_0} \\ &\leq (n + 2) \perturb_{\max} \;
\leq \; \mixingtime \big(2 + \log \perturb_{\max}^{-1} + \log
\tfrac{1}{\min_x \stationary_0 (x)} \big)\perturb_{\max},
\end{align*}
which proves the second claim.

\subsection{Proof of Lemma~\ref{lemma:jacobi-wrt-parameter}}\label{subsubsec:proof-lemma-jacobi-wrt-parameter}

We first consider the functional $\perturb \mapsto \thetabar (\transition_\perturb) \mydefn \big(I - \Exs_{\stationary_\perturb} [\Lmap (\state)] \big)^{-1} \Exs_{\stationary_\perturb} [\bmap (\state)]$. Note that the stationary distribution $\stationary_\perturb$ satisfies the identity $\stationary_\perturb \transition_\perturb = \stationary_\perturb$. Taking derivatives, we obtain the following equality for all $x, y \in \statespace$:
\begin{align*}
    \frac{\partial \stationary_\perturb}{\partial \perturb_x (y)}
    \cdot (I - \transition_\perturb) = \stationary_\perturb \cdot
    \frac{\partial \transition_\perturb}{\partial \perturb_x (y)} =
    \stationary_\perturb (x) \transition_\perturb (x, y) \cdot \big[
      \bm{1}_{z = y} - \transition_\perturb (x, z) \big]_{z
      \in \statespace}.
\end{align*}

Note that the linear operator $(I - \transition_\perturb)$ is invertible on the subspace $\statnullspace_\perturb$. For any $f \in \statnullspace_\perturb$, we have
\begin{align*}
    &\frac{\partial}{\partial \perturb_x (y)} \Exs_{\stationary_\perturb} \big[ f (\state) \big] = \sum_{z \in \statespace} \frac{\partial \stationary_\perturb (z)}{\partial \perturb_x (y)} \cdot f (\state) \\
    &= \stationary_\perturb (x) \transition_\perturb (x, y) \cdot \big[ \bm{1}_{z = y} - \transition_\perturb (x, z)  \big]_{z \in \statespace} \cdot \big( I - \transition_\perturb \big)^{-1} \big|_{\statnullspace_\perturb} \cdot f.
\end{align*}
In the above expression, the notation $\big( I - \transition_\perturb \big)^{-1} \big|_{\statnullspace_\perturb}$ denotes the inverse of the operator $I - \transition_\perturb$ within the subspace $\statnullspace_\perturb$, a bounded linear operator on this space.
Note that the derivative is invariant under translation. For any $f \in \real^\statespace$, define the auxiliary function $\widetilde{f} \mydefn f - \Exs_{\stationary_\perturb} [f]$, and write
\begin{align}
    &\frac{\partial}{\partial \perturb_x (y)} \Exs_{\stationary_\perturb} \big[ f (\state) \big] \nonumber \\
    &= \frac{\partial}{\partial \perturb_x (y)} \Exs_{\stationary_\perturb} \big[ \widetilde{f} (\state) \big] = \stationary_\perturb (x) \transition_\perturb (x, y) \cdot \big[ \bm{1}_{z = y} - \transition_\perturb (x, z)  \big]_{z \in \statespace} \cdot \big( I - \transition_\perturb \big)^{-1} \big|_{\statnullspace_\perturb} \cdot \widetilde{f} \nonumber\\
    &= \stationary_\perturb (x) \transition_\perturb (x, y) \cdot \big[ \bm{1}_{z = y} - \transition_\perturb (x, z)  \big]_{z \in \statespace} \cdot \big( I - \transition_\perturb \big)^{-1} \big|_{\statnullspace_\perturb} \cdot \big(f - \Exs_{\stationary_\perturb} [f] \big) \nonumber\\
    &= \stationary_\perturb (x) \transition_\perturb (x, y) \cdot \big( \greenOp_\perturb f (y) - \sum_{z \in \statespace} \transition_\perturb (x, z) \greenOp_\perturb f (z) \big).\label{eq:expression-for-derivative-of-expectation-in-lower-bound-proof}
\end{align}

On the other hand, we can express the desired functional $\thetabar
(\transition_\perturb)$ in the form above. In particular, setting
$\Lbar^{(\perturb)} \mydefn \Exs_{\stationary_\perturb} \big[ \Lmap
  (\state) \big]$ and $\bbar^{(\perturb)} \mydefn
\Exs_{\stationary_\perturb} \big[ \bmap (\state) \big]$, we see that
for any $x, y \in \statespace$, we have
\begin{align*}
\frac{\partial \thetabar (\transition_\perturb)}{\partial \perturb_x
  (y)} &= \big(I - \Lbar^{(\perturb)} \big)^{-1} \frac{\partial
  \Lbar^{(\perturb)}}{\partial \perturb_x (y)}\big(I -
\Lbar^{(\perturb)} \big)^{-1} \bbar^{(\perturb)} +\big(I -
\Lbar^{(\perturb)} \big)^{-1} \frac{\partial
  \bbar^{(\perturb)}}{\partial \perturb_x (y)} \\ &=\big(I -
\Lbar^{(\perturb)} \big)^{-1} \big( \big( \frac{\partial}{\partial
  \perturb_x (y)} \Exs_{\stationary_\perturb} \big[ \Lmap (\state)
  \big]\big) \cdot \thetabar (\transition_\perturb) +
\frac{\partial}{\partial \perturb_x (y)} \Exs_{\stationary_\perturb}
\big[ \bmap (\state) \big] \big).
\end{align*}
Following the formula~\eqref{eq:expression-for-derivative-of-expectation-in-lower-bound-proof}, we conclude that
\begin{multline}
      \frac{\partial \thetabar (\transition_\perturb)}{\partial \perturb_x (y)} = \stationary_\perturb (x) \transition_\perturb (x, y) \big(I - \Lbar^{(\perturb)} \big)^{-1} \big[ \greenOp_\perturb \big( \Lmap (y) \thetabar (\transition_\perturb) + \bmap (y) \big) \big]\\
      -\stationary_\perturb (x) \transition_\perturb (x, y)  \sum_{z \in \statespace}  \transition_\perturb (x, z) \big(I - \Lbar^{(\perturb)} \big)^{-1} \big[ \greenOp_\perturb \big( \Lmap (z) \thetabar (\transition_\perturb) + \bmap (z) \big) \big].
\end{multline}
Recall the shorthand notation from before, where for each $\state \in \statespace$, we defined
\begin{align*}
    \greenFunc_\perturb (\state) = \big(I - \Lbar^{(\perturb)} \big)^{-1} \big[ \greenOp_\perturb \big( \Lmap (\state) \thetabar (\transition_\perturb) + \bmap (\state) \big) \big].
\end{align*}
Given $w \in \real^d$, if we parameterize the perturbation as $\perturb = Q w$, the chain rule yields
\begin{align*}
    &\nabla_w \thetabar (\transition_\perturb) = Q^\top \cdot \nabla_\perturb \thetabar (\transition_\perturb)\\
    &= \sum_{x \in \statespace} \stationary_\perturb (x) \big( \sum_{y \in \statespace} \transition_\perturb(x, y) \greenFunc (y) q_x (y)^\top \\
    &\qquad \qquad-   \big(\sum_{y \in \statespace} \transition_\perturb(x, y) \greenFunc (y)\big) \big(\sum_{y \in \statespace} \transition_\perturb(x, y) \greenFunc_\perturb (y) q_x (y)\big)^\top  \big)\\
    &= \Exs_{X \sim \stationary_\perturb} \big[ \cov_{Y \sim \transition_\perturb (X, \cdot)} \big( \greenFunc_\perturb (Y) - \transitionOp_\perturb \greenFunc_\perturb (X), q_X (Y) \mid X \big) \big],
\end{align*}
as claimed. \qed

\subsection{Proof of Lemma~\ref{lemma:green-func-error-estimate}} \label{subsubsec:proof-lemma-green-func-error-estimate}

The following technical lemma is used throughout the proof, and proved
in Section~\ref{sec:perturb-invert}.
\begin{lemma}
\label{lemma:perturbed-matrix-invertibility}
Given a perturbation vector $w$ satisfying $\vecnorm{w}{2} \leq
\frac{1 - \kappa}{2 c \mixingtime \sigmaA \sqrt{d \cdot
    \opnorm{\Lambda}} \log d} $, for $\perturb = Q w$, the matrix $I -
\Lbar^{(\perturb)}$ is invertible, with $\opnorm{\big( I -
  \Lbar^{(\perturb)} \big)^{-1}} \leq \frac{2}{1 - \kappa}$.
\end{lemma}

Before proceeding with the proof, we note two direct consequences of
Lemma~\ref{lemma:poisson-eq-estimate} from Section~\ref{appendix:useful-moment-bound}.
 First, by
taking $f(x) \mydefn \inprod{\coordinate_j}{\Lmap (x) u}$ and $f(x)
\mydefn \inprod{\coordinate_j}{\bmap (x)}$, applying the tail
assumption~\ref{assume:noise-moments} and the boundedness
assumption~\ref{assume-lip-mapping}, we have the following second
moment estimate for any $u \in \sphere^{d - 1}$ and $j \in [d]$:
\begin{subequations}\label{eq:poisson-eq-estimate-on-matrix-vec-observations}
\begin{align}
    \Exs_{X \sim \stationary_\perturb} \big[ \inprod{\coordinate_j}{\greenOp_\perturb \Lmap (X) u}^2 \big] &\leq c \mixingtime^2 \sigmaA^2 \log^2 d, \quad \mbox{and}\\
      \Exs_{X \sim \stationary_\perturb} \big[ \inprod{\coordinate_j}{\greenOp_\perturb \bmap (X)}^2 \big] &\leq c \mixingtime^2 \sigmab^2 \log^2 d.
\end{align}
\end{subequations}
Second, by taking $f_j (x) \mydefn \inprod{\coordinate_j}{\Lmap (x) \thetabar (\transition_\perturb) + \bmap (x)}$, for any integer $p \geq 1$ and $K > 0$, Markov's inequality yields the bound
\begin{align*}
    \Prob_{X \sim \stationary_\perturb} \Big[ \greenOp_\perturb f_j
      (X) \geq K \Big] & \leq K^{-2p} \Exs_{X \sim
      \stationary_\perturb} \big[ \greenOp_\perturb f_j (X)^{2p} \big]
    \leq \big( \tfrac{c p^2 \mixingtime \varbound \log d}{K} \big)^{2p}.
\end{align*}
By taking $K = 2 c p^2 \mixingtime \varbound \log d$ and $p = - 2 \log \min_{x \in \statespace}
\stationary_0 (x)$, we find that
\begin{align*}
\Prob_{X \sim \stationary_\perturb} \Big[ \greenOp_\perturb f_j (X)
  \geq 8 c \mixingtime \varbound
  \log^3 \big( \tfrac{d}{\min_{x \in \statespace} \stationary_0(x)}
  \big) \Big] < \frac{1}{2} \min_{x \in \statespace} \stationary_0 (x)
\leq \min_{x \in \statespace} \stationary_\perturb (x),
\end{align*}
Since $\stationary_\perturb$ is a discrete measure, this high-probability bound implies a deterministic bound
\begin{align*}
  \greenOp_\perturb f_j(x) \leq 8 c \mixingtime \varbound \log^3 \big(\tfrac{d}{\min_{x'
      \in \statespace} \stationary_0 (x')} \big) \qquad \mbox{for all
    $x \in \statespace$.}
\end{align*}
Combining the estimates for all $j$ coordinates yields the bound
\begin{align}
\label{eq:green-func-almost-sure-bounds}  
\max_{x \in \statespace} \vecnorm{\greenFunc_\perturb (x)}{2} \leq
\tfrac{1}{1 - \kappa} \max_{x \in \statespace} \vecnorm{
  \greenOp_\perturb \big[ f_j (x) \big]_{j \in [d]}}{2} \leq
\tfrac{c\mixingtime \varbound
  \sqrt{d}}{1 - \kappa} \log^3 \big( \tfrac{d}{\min_{x
    \in \statespace} \stationary_0 (x)} \big).
\end{align}

Given the two lemmas and facts derived above, we now proceed to the
proof of Lemma~\ref{lemma:green-func-error-estimate}.  Taking
derivatives on both sides of
equation~\eqref{eq:def-poisson-eq-in-lower-bound-proof}, we obtain
\begin{align*}
    \nabla_w \greenFunc_\perturb (z)  &= \big(I_d - \Lbar^{(\perturb)} \big)^{-1} \cdot \greenOp_\perturb  \Lmap (z) \cdot \nabla_w \thetabar (\transition_\perturb)\\
    & \qquad + \big(I_d - \Lbar^{(\perturb)} \big)^{-1} \cdot \big(\nabla_w \greenOp_\perturb\big)  \big(\Lmap (z) \thetabar (\transition_\perturb) + \bmap (z) \big)\\
    &\qquad \qquad -  \big(I_d - \Lbar^{(\perturb)} \big)^{-1} \nabla_w \big( \Lbar^{(\perturb)} \big) \big(I_d - \Lbar^{(\perturb)} \big)^{-1} (\greenOp_\perturb \Lmap (z) \cdot \thetabar (\transition_\perturb) + \greenOp_\perturb \bmap (z) )\\
    &=: J_1 (\perturb, z) + J_2 (\perturb, z) + J_3 (\perturb, z) .
\end{align*}
We then have the integral relation
\begin{align*}
    \greenFunc_\perturb (z) - \greenFunc_0 (z) = \int_0^1 \nabla_w
    \greenFunc_{s \perturb} (z) \cdot w ~ d s = \int_0^1 \big( J_1 (s
    \perturb, z) + J_2 ( s \perturb, z) + J_3 (s \perturb, z) \big)
    \cdot w ~ d s.
\end{align*}
It thus suffices to prove individual upper bounds on the terms $J_1 (s\perturb, z) \cdot w, J_2 (s\perturb, z) \cdot  w$ and $J_3 (s\perturb, z) \cdot w$.
\paragraph{Bounds on the term $J_1 (s\perturb, z) \cdot w$:}

Invoking Lemma~\ref{lemma:jacobi-wrt-parameter}, we have
\begin{align*}
   \nabla_w \thetabar (\transition_\perturb)
    = \Exs_{X \sim \stationary_{\perturb}, Y \sim \transition_{\perturb} (X, \cdot)} \big[ \big( \greenFunc_\perturb (Y) - \transitionOp_\perturb \greenFunc_\perturb (X) \big) \big(\greenFunc_0 (Y) - \transitionOp_0 \greenFunc_0 (X)  \big)^\top \big].
\end{align*}
Consequently, for $X \sim \stationary_{\perturb}$ and $Y \sim \transition_{\perturb} (X, \cdot)$, we have the error decomposition
\begin{multline*}
   \vecnorm{\nabla_w \thetabar (\transition_\perturb) w}{2} = \vecnorm{\Exs \big[ \cov \big( \greenFunc_\perturb (Y), \greenFunc_0 (Y)  \mid X\big) \big] w}{2} \\
   \leq  \vecnorm{\Exs \big[ \cov \big(\greenFunc_0 (Y)  \mid X\big) \big] w}{2} + \vecnorm{\Exs \big[ \cov \big( \greenFunc_\perturb (Y) - \greenFunc_0 (Y), \greenFunc_0 (Y)  \mid X\big) \big] w}{2}. 
\end{multline*}
For perturbation matrix $\perturb$ satisfying the condition $\max
\limits_{x \in \statespace} \vecnorm{\perturb_x}{\infty} \leq
\tfrac{1}{128 \mixingtime}$,
Lemma~\ref{lemma:perturbation-principle-markov} implies the sandwich
relations
\begin{align*}
\tfrac{1}{2} \stationary_0 \preceq \stationary_\perturb \preceq
\tfrac{3}{2} \stationary_0, \quad \mbox{and} \quad \tfrac{1}{2}
\transition_0 (x) \preceq \transition_\perturb (x, \cdot) \preceq
\tfrac{3}{2} \transition_0 (x), \quad \mbox{for all $x
  \in \statespace$.}
\end{align*}
For the first term in above decomposition, we have
\begin{align*}
   \vecnorm{\Exs \big[ \cov \big(\greenFunc_0 (Y)  \mid X\big) \big] w}{2} \leq \vecnorm{ \Exs \big[ (\greenFunc_0 (Y) - \transitionOp_0 \greenFunc_0 (X)) (\greenFunc_0 (Y) - \transitionOp_0 \greenFunc_0 (X))^\top \big] w}{2} \\
   \leq
     \tfrac{9}{4} \vecnorm{\cov_{X \sim \stationary_0, Y \sim
         \transition_0 (X, \cdot)} \big( \greenFunc_0 (Y) -
       \transitionOp_0 \greenFunc_0 (X) \big) \cdot w}{2} =
     \tfrac{9}{4} \vecnorm{\Lambda w}{2} \leq \tfrac{9}{4} \sqrt{\trace
       (\Lambda) / \numobs},
\end{align*}
where the last inequality is due to the bound~\eqref{eq:prior-support-size}.

For the second term in the decomposition, for $X \sim \stationary_{\perturb}$ and $Y \sim \transition_{\perturb} (X, \cdot)$, we have
\begin{align*}
    &\vecnorm{\Exs \big[ \cov \big( \greenFunc_\perturb (Y) - \greenFunc_0 (Y), \greenFunc_0 (Y)  \mid X\big) \big] w}{2}\\
    &= \sup_{v \in \sphere^{d - 1}} v^\top \Exs \big[ \cov \big( \greenFunc_\perturb (Y) - \greenFunc_0 (Y), \greenFunc_0 (Y)  \mid X\big) \big] w\\
    &\leq  \sup_{v \in \sphere^{d - 1}} \sqrt{\Exs \big[ \inprod{\big( \greenFunc_\perturb (Y) - \greenFunc_0 (Y) \big)}{v}^2 \big] } \cdot \sqrt{\Exs_{X \sim \stationary_{\perturb}, Y \sim \transition_{\perturb} (X, \cdot)} \big[ \big(\big(\greenFunc_0 (Y) - \transitionOp_0 \greenFunc_0 (X) \big)^\top w \big)^2 \big] }\\
    &\leq \tfrac{3}{2} \sqrt{w^\top \Lambda w} \sqrt{\Exs_{X \sim \stationary_\perturb} \vecnorm{\greenFunc_\perturb (X) - \greenFunc_0 (X)}{2}^2}.
\end{align*}
By equation~\eqref{eq:prior-density-construction}, on the support of
the prior density, we have the bound $w^\top \Lambda w =
  \numobs^{-1} \psi^\top D^{-1/2} U^\top \Lambda U D^{-1/2} \psi \leq
  \tfrac{d}{n}$.  Consequently, we have the upper bound
\begin{align}
\label{eq:perturbation-for-thetabar-in-lower-bound-proof}  
\vecnorm{\nabla_w \thetabar (\transition_\perturb) w}{2} \leq
\tfrac{9}{4} \sqrt{\tfrac{\trace(\Lambda)}{\numobs}} + \tfrac{3}{2}
\cdot \sqrt{\tfrac{d}{n} \cdot \Exs_{X \sim \stationary_\perturb}
  \vecnorm{\greenFunc_\perturb (X) - \greenFunc_0 (X)}{2}^2}.
\end{align}

Collecting the bounds above and invoking
equation~\eqref{eq:poisson-eq-estimate-on-matrix-vec-observations} and
Lemma~\ref{lemma:perturbed-matrix-invertibility}, we obtain the
following bound on the desired term:
\begin{align*}
    &\Exs_{Y \sim \stationary_\perturb} \big[ \vecnorm{J_1 (\ell
      \perturb, Y)w}{2}^2 \big]\\ &\leq \opnorm{\big( I_d -
    \Lbar^{(\ell \perturb)} \big)^{-1}}^2 \cdot \Exs_{Y \sim
    \stationary_\perturb} \big[ \vecnorm{\greenOp_{\ell \perturb}
      \Lmap (Y) \cdot \nabla_w \thetabar (\transition_{\ell \perturb})
      w}{2}^2 \big] \\
& \leq \tfrac{4}{(1 - \kappa)^2} \cdot \tfrac{3}{2} \Exs_{Y \sim
    \stationary_{\ell \perturb}} \big[ \vecnorm{\greenOp_{\ell
        \perturb} \Lmap (Y) \cdot \nabla_w \thetabar
      (\transition_{\ell \perturb}) w}{2}^2 \big]\\
& \leq \tfrac{6}{(1 - \kappa)^2} \cdot c \mixingtime^2 \sigmaA^2
  \usedim \log^2 \usedim \cdot \vecnorm{\nabla_w \thetabar
    (\transition_{\ell \perturb}) w}{2}^2 \\
& \leq \tfrac{ c \mixingtime^2 \sigmaA^2 \usedim \log^2 d}{(1 -
    \kappa)^2} \cdot \tfrac{\trace (\Lambda)}{n} + \tfrac{ c
    \mixingtime^2 \sigmaA^2 \usedim^2 \log^2 d}{(1 - \kappa)^2 n}
  \sup_{0 \leq \ell \leq 1} \Exs_{X \sim \stationary_{\ell \perturb}}
  \vecnorm{\greenFunc_{\ell \perturb} (X) - \greenFunc_0 (X)}{2}^2.
\end{align*}

\paragraph{Bounds on the term $J_2 (s\perturb, z) \cdot w$:} For any function
$\statespace \rightarrow \real^d$ and $x, y \in \statespace$, we note
that
\begin{align*}
\tfrac{\partial}{\partial \perturb_x (y)} \greenOp_\perturb f &= - (I
- \transitionOp_\perturb)^{-1} |_{\statnullspace_\perturb}\cdot
\tfrac{\partial \transitionOp_\perturb}{\partial \perturb_x (y)} \cdot
(I - \transitionOp_\perturb)^{-1} |_{\statnullspace_\perturb} f\\ &= -
\greenOp_\perturb \cdot \big[ \bm{1}_{s = x} \transition_\perturb (x,
  y) \cdot \big( \bm{1}_{s' = y} - \transition_\perturb (x, s') \big)
  \big]_{s, s' \in \statespace} \cdot \greenOp f\\ &= -
\greenOp_{\perturb} \cdot \big[ \bm{1}_{s = x} \transition_{\perturb}
  (x, y) \cdot \big( \greenOp_\perturb f (y) - \sum_{s'}
  \transition_\perturb (x, s') \greenOp_\perturb f (s') \big) \big]_{s
  \in \statespace}.
\end{align*}
We can then derive the formula for derivative with respect to the
parameter $w$, as
\begin{align*}
&\big( \nabla_w \greenOp_{\perturb} \big) f (z) = \sum_{x, y
    \in \statespace} \big(\frac{\partial}{\partial \perturb_x (y)}
  \greenOp_\perturb f (z) \big) \cdot q_x (y)^\top\\ &= - \sum_{x, y
    \in \statespace} \transition_{\perturb} (x, y) \greenOp_\perturb
  \bm{1}_x (z) \cdot \big( \greenOp_\perturb f (y) -
  \transitionOp_\perturb \greenOp_\perturb f (x) \big) \cdot \big(
  \greenFunc_0 (y) - \transitionOp_0 \greenFunc_0 (x) \big)^\top\\ &=
  - \sum_{x, y \in \statespace} \sum_{t = 0}^{\infty}\big(
  \transition_\perturb^t (z, x) - \stationary_\perturb (x) \big)
  \transition_\perturb (x, y) \big( \greenOp_\perturb f (y) -
  \transitionOp_\perturb \greenOp_\perturb f (x) \big) \big(
  \greenFunc_0 (y) - \transitionOp_0 \greenFunc_0 (x) \big)^\top.
\end{align*}
Substituting $f (z) = \Lmap (z) \thetabar (\transition_\perturb) + \bmap (z)$, we note that $\greenOp_\perturb f = \greenFunc_\perturb$, and consequently,
\begin{align*}
    &\big(\nabla_w \greenOp_\perturb\big)  \big(\Lmap (z) \thetabar (\transition_\perturb) + \bmap (z) \big)\\
    &= \sum_{t = 0}^{\infty} \big( \Exs_{X \sim \transition_\perturb^t (z, \cdot), Y \sim \transition_\perturb (X, \cdot)} \big[ \big( \greenFunc_\perturb (Y) - \transitionOp_\perturb \greenFunc_\perturb (X) \big)\big( \greenFunc_0 (Y) - \transitionOp_0 \greenFunc_0 (X) \big)^\top \big]\\
    &\qquad \qquad- \Exs_{X \sim \stationary_\perturb, Y \sim \transition_\perturb (X, \cdot)} \big[ \big( \greenFunc_\perturb (Y) - \transitionOp_\perturb \greenFunc_\perturb (X) \big)\big( \greenFunc_0 (Y) - \transitionOp_0 \greenFunc_0 (X) \big)^\top \big] \big)\\
    &=: \sum_{t = 0}^{\infty} D_t (z).
\end{align*}
Next, we estimate the difference term above in two different ways, depending on the value of $t$.
On the one hand, note that
\begin{align*}
    &\Exs_{Z \sim \stationary_\perturb} \vecnorm{\Exs_{X \sim \transition_\perturb^t (Z, \cdot), Y \sim \transition_\perturb (X, \cdot)} \big[ \big( \greenFunc_\perturb (Y) - \transitionOp_\perturb \greenFunc_\perturb (X) \big)\big( \greenFunc_0 (Y) - \transitionOp_0 \greenFunc_0 (X) \big)^\top \big] w}{2}^2\\
    &\leq \sup_{x, y \in \statespace} \vecnorm{\greenFunc_\perturb (y) - \transitionOp_\perturb \greenFunc_\perturb (x)}{2}^2 \cdot \Exs_{X \sim \stationary_\perturb, Y \sim \transition_\perturb (X, \cdot)} \big[ \inprod{w}{ \greenFunc_0 (Y) - \transitionOp_0 
    \greenFunc_0 (X)}^2 \big]\\
    &\leq 4 \sup_{x \in \statespace} \vecnorm{\greenFunc_\perturb (x)}{2}^2 \cdot \Exs_{X \sim \stationary_\perturb, Y \sim \transition_\perturb (X, \cdot)} \big[ \inprod{w}{ \greenFunc_0 (Y) - \transitionOp_0 
    \greenFunc_0 (X)}^2 \big],
\end{align*}
where the bound for the factor $\sup_{x \in \statespace}
\vecnorm{\greenFunc_\perturb (x)}{2}^2$ follows from
equation~\eqref{eq:green-func-almost-sure-bounds}. For the latter term
in the display above, we note that
\begin{multline*}
\Exs_{X \sim \stationary_\perturb, Y \sim \transition_\perturb (X,
  \cdot)} \big[ \inprod{w}{ \greenFunc_0 (Y) - \transitionOp_0
    \greenFunc_0 (X)}^2 \big] \\
     \leq 2 \Exs_{X \sim \stationary_0, Y
  \sim \transition_0 (X, \cdot)} \big[ \inprod{w}{ \greenFunc_0 (Y) -
    \transitionOp_0 \greenFunc_0 (X)}^2 \big]
 \leq 2 w^\top \Lambda w = \tfrac{2d}{n}.
\end{multline*}
Putting together the pieces yields the first estimate
\begin{align*}
    \Exs_{Z \sim \stationary_\perturb} \big[ \vecnorm{D_t (Z) w}{2}^2
      \big] \leq \tfrac{c \mixingtime^2 \varbound^2 d^2}{(1 -
      \kappa)^2 \numobs} \log^6 \big(\tfrac{d}{\min_{x
        \in \statespace} \stationary_0 (x)} \big).
\end{align*}

On the other hand, given $z \in \statespace$ and the Markov chain
$(\state_t)_{t \geq 0}$ starting from $s_0 = z$, for any $t > 0$,
there exists a random state $\widetilde{\state}_t$ such that
$\widetilde{\state}_t \sim \stationary_\perturb$, and we have
\mbox{$\Prob \big( \widetilde{\state}_t \neq \state_t \big) \leq
  2^{\lfloor \frac{t}{\mixingtime} \rfloor}$.}  Define a random
variable $\widetilde{\state}_{t + 1}$ by setting
$\widetilde{\state}_{t + 1} = \state_{t + 1}$ whenever $\state_t =
\widetilde{\state}_t$, and drawing $\widetilde{\state}_{t+1} \sim
\transition (\widetilde{\state}_t, \cdot)$ otherwise.  From this
construction, we have
\begin{align*}
&\vecnorm{D_t (z) w}{2}\\
 &\leq \sup_{u \in \sphere^{d - 1}} \bigg\{ \Exs
\big[ u^\top \big( \greenFunc_\perturb (\state_{t + 1}) -
  \transitionOp_\perturb \greenFunc_\perturb (\state_t) \big) \cdot
  w^\top \big( \greenFunc_0 (\state_{t + 1}) - \transitionOp_0
  \greenFunc_0 (\state_t) \big) \mid z \big]\\ &\qquad \qquad - \Exs
\big[ u^\top \big( \greenFunc_\perturb (\widetilde{\state}_{t + 1}) -
  \transitionOp_\perturb \greenFunc_\perturb (\widetilde{\state}_t)
  \big) \cdot w^\top \big( \greenFunc_0 (\widetilde{\state}_{t + 1}) -
  \transitionOp_0 \greenFunc_0 (\widetilde{\state_t}) \big) \mid z
  \big] \bigg\} \\ &\leq \sup_{u \in \sphere^{d - 1}} \Exs \big[
  u^\top \big( \greenFunc_\perturb (\state_{t + 1}) -
  \transitionOp_\perturb \greenFunc_\perturb (\state_t) \big) \cdot
  w^\top \big( \greenFunc_0 (\state_{t + 1}) - \transitionOp_0
  \greenFunc_0 (\state_t) \big) \bm{1}_{\state_t \neq
    \widetilde{\state}_t} \mid z \big]\\ &\qquad + \sup_{u \in
  \sphere^{d - 1}} \Exs \big[ u^\top \big( \greenFunc_\perturb
  (\widetilde{\state}_{t + 1}) - \transitionOp_\perturb
  \greenFunc_\perturb (\widetilde{\state}_t) \big) \cdot w^\top \big(
  \greenFunc_0 (\widetilde{\state}_{t + 1}) - \transitionOp_0
  \greenFunc_0 (\widetilde{\state_t}) \big) \bm{1}_{\state_t \neq
    \widetilde{\state}_t} \mid z \big].
\end{align*}
Applying the Cauchy--Schwarz inequality twice yields
\begin{align*}
    &\Exs_{Z \sim \stationary_\perturb} \big[ \vecnorm{D_t (Z) w}{2}^2
    \big]\\ &\leq \Exs \big[ \vecnorm{\greenFunc_\perturb (\state_{t +
        1}) - \transitionOp_\perturb \greenFunc_\perturb
      (\state_t)}{2}^4 \big]^{1/2} \cdot \Exs \big[ w^\top \big(
    \greenFunc_0 (\state_{t + 1}) - \transitionOp_0 \greenFunc_0
    (\state_t) \big)^8 \big]^{1/4} \cdot \Exs \big[ \bm{1}_{\state_t
      \neq \widetilde{\state}_t} \big]^{1/4} \\
& \qquad +\Exs \big[ \vecnorm{\greenFunc_\perturb
      (\widetilde{\state}_{t + 1}) - \transitionOp_\perturb
      \greenFunc_\perturb (\state_t)}{2}^4 \big]^{1/2} \cdot \Exs
  \big[ w^\top \big( \greenFunc_0 (\widetilde{\state}_{t + 1}) -
    \transitionOp_0 \greenFunc_0 (\widetilde{\state}_t) \big)^8
    \big]^{1/4} \cdot \Exs \big[ \bm{1}_{\state_t \neq
      \widetilde{\state}_t} \big]^{1/4} \\
& \leq \frac{c \mixingtime^4}{(1 - \kappa)^4} \varbound^4 d
  \vecnorm{w}{2}^2 \cdot \log^6 d \cdot 2^{1 - \frac{t}{4
      \mixingtime}},
\end{align*}
corresponding to the second estimate.

Finally, setting $\tau = c \mixingtime \log \tfrac{\mixingtime d}{1 -
  \kappa}$ yields
\begin{align*}
\Exs_{Z \sim \stationary_\perturb} \big[ \vecnorm{\sum_{t =
      0}^{\infty} D_t (Z) w}{2}^2 \big] & \leq \big( \sum_{t =
  0}^{\infty} e^{- \tfrac{t}{\tau}} \big) \cdot \big( \sum_{t =
  0}^{\infty} e^{\tfrac{t}{\tau}} \Exs_{Z \sim \stationary_\perturb}
\big[ \vecnorm{D_t (Z) w}{2}^2 \big] \big)\\ & \leq
\tfrac{c\mixingtime^4 \varbound^2 d^2}{(1 - \kappa)^2 \numobs} \log^6
\big( \tfrac{d}{\min_{x \in \statespace} \stationary_0 (x)} \big),
\end{align*}
so that
\begin{align*}
\Exs_{Z \sim \stationary_\perturb} \big[ \vecnorm{J_2 (\ell \perturb,
    Z) w}{2}^2 \big] \leq \tfrac{c\mixingtime^4 \varbound^2 d^2}{(1 -
  \kappa)^4 \numobs} \log^6 \big( \tfrac{d}{\min_{x \in \statespace}
  \stationary_0 (x)} \big).
\end{align*}

\paragraph{Bounds on the term $J_3 (s\perturb, z) \cdot w$:}

By
equation~\eqref{eq:expression-for-derivative-of-expectation-in-lower-bound-proof},
for any vector $u \in \sphere^{d - 1}$, we have
\begin{align*}
\nabla_w \big( \Lbar^{(\perturb)} u \big) = \sum_{x, y
  \in \statespace} \stationary_\perturb (x) \transition_\perturb (x,
y) \big( \greenOp_\perturb \Lbar^{(\perturb)} (y) - \sum_{z
  \in \statespace} \transition_\perturb (x, z) \greenOp_\perturb
\Lbar^{(\perturb)} (z) \big) u \cdot q_x (y)^\top.
\end{align*}
For any $z \in \statespace$, we obtain
\begin{align*}
   & \vecnorm{\nabla_w \big( \Lbar^{(\perturb)} \big) \greenFunc_\perturb (z) w}{2}\\
    &= \sup_{u \in \sphere^{d - 1}} \Exs_{X \sim \stationary_\perturb, Y \sim \transition_\perturb (X, \cdot)} \big[ u^\top \big( \greenOp_\perturb \Lbar^{(\perturb)} (Y) - \transitionOp_\perturb \greenOp_\perturb \Lbar^{(\perturb)} (X) \big) \greenFunc_\perturb (z) q_X (Y)^\top w \big]\\
    &\leq  \sup_{u \in \sphere^{d - 1}} \sqrt{\Exs \big(  u^\top \big( \greenOp_\perturb \Lbar^{(\perturb)} (Y) - \transitionOp_\perturb \greenOp_\perturb \Lbar^{(\perturb)} (X) \big) \greenFunc_\perturb (z) \big)^2} \cdot \sqrt{\Exs \big[ \big( q_X (Y)^\top w \big)^2 \big]}\\
    &\leq c \mixingtime \sigmaA \vecnorm{\greenFunc_\perturb (z)}{2} \log d \cdot \sqrt{\tfrac{d}{n}},
\end{align*}
where the final inequality is due to equation~\eqref{eq:poisson-eq-estimate-on-matrix-vec-observations}.
Combining with Lemma~\ref{lemma:perturbed-matrix-invertibility}, we have the bound
\begin{align*}
    \Exs_{Z \sim \stationary_\perturb} \big[ \vecnorm{J_3 (\ell \perturb, Z) w}{2}^2\big]
    &\leq \tfrac{c d^2}{(1 - \kappa)^2 \numobs} \cdot \mixingtime^2 \sigmaA^2 \log^2 d \cdot \Exs_{Z \sim \stationary_\perturb} \big[ \vecnorm{\greenFunc_\perturb (Z)}{2}^2 \big] \\
    &\leq \tfrac{c \sigmaA^2 \varbound^2 \mixingtime^4 d^2}{(1 - \kappa)^4 \numobs} \log^2 d.
\end{align*}

\paragraph{Finishing the proof.}
Collecting the bounds for $J_1,~ J_2$ and $J_3$ and for $\numobs \geq
\tfrac{c \mixingtime^2 \sigmaA^2 d^2 \log^2 d}{(1 - \kappa)^2}$, we
have
\begin{multline*}
\sup_{0 \leq \ell \leq 1} \Exs_{Z \sim \stationary_\perturb} \big[
  \vecnorm{\greenFunc_{\ell\perturb} (Z) - \greenFunc_0 (Z)}{2}^2 \big]\\
   \leq
\tfrac{c (1 + \sigmaA^2) \varbound^2 \mixingtime^4 d^2}{(1 - \kappa)^4
  \numobs} \log^6 \big( \tfrac{d}{\min_x \stationary_0 (x)} \big) +
\tfrac{1}{2} \sup_{0 \leq \ell \leq 1} \Exs_{Z \sim
  \stationary_\perturb} \big[ \vecnorm{\greenFunc_{\ell \perturb} (Z) -
    \greenFunc_0 (Z)}{2}^2 \big] ,
\end{multline*}
which completes the proof of the first claim of the lemma.

For the second claim, we combine the first claim with
equation~\eqref{eq:perturbation-for-thetabar-in-lower-bound-proof} and
obtain
\begin{align*}
\vecnorm{\nabla_w \thetabar (\transition_\perturb) w}{2} \leq
\tfrac{3}{2} \sqrt{\tfrac{\mathrm{trace}(\Lambda)}{\numobs}} +
\sqrt{\tfrac{c (1 + \sigmaA^2) \varbound^2 \mixingtime^4 d^3}{(1 -
    \kappa)^4 \numobs^2} \log^6 \big(\tfrac{d}{\min_x \stationary_0
    (x)} \big)}.
\end{align*}
Taking the integral yields
\begin{align*}
\vecnorm{\thetabar (\transition_\perturb) - \thetabar
  (\transition_0)}{2} &\leq \int_0^1 \vecnorm{\nabla_w \thetabar
  (\transition_{\ell \perturb}) w}{2} d \ell\\
  & \leq \tfrac{3}{2} \sqrt{
  \tfrac{\mathrm{trace} (\Lambda)}{\numobs}} + \sqrt{\tfrac{c (1 +
    \sigmaA^2) \varbound^2 \mixingtime^4 d^3}{(1 - \kappa)^4
    \numobs^2} \log^6 \big(\tfrac{d}{\min_x \stationary_0 (x)}\big)},
\end{align*}
which proves the second claim.

\subsection{Proof of Lemma~\ref{lemma:fisher-info-for-obs}}\label{subsubsec:proof-lemma-fisher-info-for-obs}

We first compute the Fisher information with respect to the
perturbation vector $\perturb$, and then transform this via chain rule into a formula
that holds with respect to the parameter $w$. We are interested in the matrix
$I^{(\numobs)} (\perturb) \mydefn \Exs_\perturb \big[ \nabla_\perturb
  \log \Prob_\perturb (\state_0^\numobs) \nabla_\perturb \log
  \Prob_\perturb (\state_0^\numobs)^\top \big]$.  When the Markov
chain $\transition_\perturb$ is run under the initial distribution
$\stationary_0$, the joint distribution of the observed trajectory
$(\state_t)_{t = 0}^\numobs$ can be factorized as $\Prob_\perturb
\big(\state_0, \state_1, \cdots, \state_\numobs \big) = \stationary_0
(\state_0) \cdot \prod_{t = 1}^\numobs \transition_\perturb (\state_{t
  - 1}, \state_t)$.

Let us now study the Fisher information matrix. For any pair $x, y
\in \statespace$ with $\transition (x, y) > 0$, performing some
algebra yields the expression
\begin{align*}
  \tfrac{\partial}{\partial \perturb_x (y)} \log \Prob_\perturb
  \big(\state_0, \state_1, \cdots, \state_\numobs \big) = \sum_{t =
    1}^\numobs \bm{1}_{\state_{t - 1} = x} \big( \bm{1}_{\state_{t} =
    y} - \transition_\perturb (x, y) \big).
\end{align*}
Consider the natural filtration $\filtration_t \mydefn \sigma
(\state_0, \state_1, \cdots, \state_t)$. Note that under the
transition kernel $\transition_\perturb$, we have the identity
\begin{align*}
    \Exs_\perturb \big[ \bm{1}_{\state_{t - 1} = x} \big( \bm{1}_{\state_{t} = y} - \transition_\perturb (x, y) \big) \mid \filtration_{t - 1} \big] = \bm{1}_{\state_{t - 1} = x} \cdot  \big( \Exs_\perturb \big[ \bm{1}_{\state_{t} = y}  \mid \state_{t - 1} = x\big] - \transition_\perturb (x, y) \big) = 0.
\end{align*}
Therefore, the process $\{ \nabla_\perturb \log \Prob_\perturb
(\state_0, \state_1, \cdots, \state_\numobs)\}_{\numobs \geq 0}$ is a
martingale adapted to the filtration $\{\filtration_t\}_{t \geq
  0}$. Its second moment is given by
\begin{align*}
S = \Exs \big[\nabla_\perturb \log \Prob_\perturb (\state_0^\numobs)
  \cdot \nabla_\perturb^\top \log \Prob_\perturb (\state_0^\numobs)
  \big] &= \sum_{t = 1}^\numobs \Exs \big[\nabla_\perturb \log
  \transition_\perturb (\state_{t - 1}, \state_t) \cdot
  \nabla_\perturb^\top \log \transition_\perturb (\state_{t - 1},
  \state_t) \big].
\end{align*}
We find that
\begin{align*}
S &= \big[\bm{1}_{x_1 = x_2} \cdot \sum_{t = 1}^\numobs \Exs \big[
    \bm{1}_{x_1 = \state_{t - 1}} \cdot \big( \bm{1}_{\state_{t} =
      y_1} - \transition_\perturb (x_1, y_1) \big) \big] \cdot \big(
  \bm{1}_{\state_{t} = y_2} - \transition_\perturb (x_2, y_2) \big)
  \big]_{(x_1, y_1), (x_2, y_2)}\\ &= \sum_{t = 1}^\numobs
\mathrm{diag} \big( \big\{ \Prob_\perturb \big(\state_{t - 1} = x
\big) \cdot \transition_\perturb (x, y) \big\}_{(x, y)} \big)\\
&\qquad \qquad -
\sum_{t = 1}^\numobs \big[ \Prob_\perturb (\state_{t - 1} = x) \cdot
  \transition_\perturb (x, y_1) \cdot \transition_\perturb (x, y_2)
  \big]_{(x, y_1), (x, y_2)}.
\end{align*}
Consequently, the Fisher information matrix is a block diagonal matrix
$I^{(\numobs)} (\perturb) = \mathrm{diag} \big( \big\{ I_x^{(\numobs)}
(\perturb) \big\}_{x \in \statespace} \big)$, where each block matrix
$I_x^{(\numobs)} (\perturb) \in \real^{\statespace
  \times \statespace}$ takes the form
\begin{align*}
I_x^{(\numobs)} (\perturb) = \sum_{t = 1}^\numobs \Prob_\perturb
\big(\state_{t - 1} = x \big) \cdot \big[ \mathrm{diag} \big( \big\{
  \transition_\perturb (x, y) \big\}_{ y \in \statespace} \big) -
  \big[ \transition_\perturb (x, y) \big]_{y \in \statespace} \big[
    \transition_\perturb (x, y) \big]_{y \in \statespace}^\top \big].
\end{align*}

By Lemma~\ref{lemma:perturbation-principle-markov},
for $\perturb_{\max}$ satisfying the inequality
${\perturb}_{\max}^{-1} \geq c \mixingtime \big(\log
\perturb_{\max}^{-1} + \log (\min \stationary_0)^{-1} \big)$ for some
constant $c > 0$, we have the bound $\tfrac{1}{2}
\stationary_{\perturb} \preceq \stationary_0 \preceq \tfrac{3}{2}
\stationary_{\perturb}$, and hence $\tfrac{1}{2}
\transition_\perturb^k \stationary_{\perturb} \preceq
\transition_\perturb^k \stationary_0 \preceq \tfrac{3}{2}
\transition_\perturb^k \stationary_{\perturb}$ \mbox{for each $k = 0,
  1, 2, \ldots$.}  From this sandwiching, we find that
\begin{align*}
  I_x^{(\numobs)} (\perturb) &\preceq \tfrac{3}{2} \sum_{t =
    1}^\numobs \transition_\perturb^{t - 1} \stationary_\perturb
  (x) \cdot \big[ \mathrm{diag} \big( \big\{ \transition_\perturb
    (x, y) \big\}_{ y \in \statespace} \big) - \big[
      \transition_\perturb (x, y) \big]_{y \in \statespace} \big[
      \transition_\perturb (x, y) \big]_{y \in \statespace}^\top
    \big]\\ &= \tfrac{3 \numobs}{2} \stationary_\perturb (x)\big[
    \mathrm{diag} \big( \big\{ \transition_\perturb (x, y) \big\}_{
      y \in \statespace} \big) - \big[ \transition_\perturb (x, y)
      \big]_{y \in \statespace} \big[ \transition_\perturb (x, y)
      \big]_{y \in \statespace}^\top \big].
\end{align*}
Turning to the Fisher information, we compute
\begin{align*}
    &I^{(\numobs)} (w) = Q^\top I^{(\numobs)} (\perturb) Q\\
     &\preceq
    \tfrac{3 \numobs}{2} \sum_{x \in \statespace} \stationary_\perturb
    (x) \big( \sum_{y \in \statespace} \transition_\perturb (x, y) q_x
    (y) q_x (y)^\top - \big( \sum_{y \in \statespace}
    \transition_\perturb (x, y) q_x (y) \big) \big( \sum_{y
      \in \statespace} \transition_\perturb (x, y) q_x (y) \big)^\top
    \big)\\
    & = \tfrac{3 \numobs}{2} \Exs_{X \sim \stationary_\perturb} \big[
      \Exs_{Y \sim \transition_\perturb (X, \cdot)} \big[ q_X (Y)
        q_X(Y)^\top \big] - \Exs_{Y \sim \transition_\perturb (X,
        \cdot)} \big[ q_X (Y) \big] \cdot \Exs_{Y \sim
        \transition_\perturb (X, \cdot)} \big[ q_X(Y) \big]^\top
      \big] \\
    & = \tfrac{3 \numobs}{2} \Exs_{X \sim \stationary_\perturb} \big[
      \cov_{\transition_\perturb (X, \cdot)} \big(q_X (Y) \mid X \big)
      \big].
\end{align*}

\subsection{Proof of Lemma~\ref{lemma:poisson-eq-identity}}
\label{subsubsec:proof-lemma-poisson-eq-identity}

For each $k \in \mathbb{Z}$, by the definition of the Green function,
we note that
\begin{align}
\label{eq:fA}
    f(s_k) = \greenOp_0 f(s_k) - \Exs \big[ \greenOp_0 f(s_{k + 1})
      \mid s_k \big] = \greenOp_0 f(s_k) - \transitionOp_0 \greenOp_0
    f(s_k).
\end{align}
By stationarity, we have
\begin{align*}
  &\sum_{k = - \infty}^{\infty} \Exs \big[ f(s_k) f (s_0) \big] = \Exs
      [f^2(s_0)] + 2 \sum_{k = 1}^{\infty} \Exs \big[ f(s_k) f(s_0)
        \big] \\
        & \overset{(i)}{=} - \Exs [f (s_0)^2] + 2 \Exs \big[
        f(s_0) \cdot \sum_{k = 0}^{\infty} \Exs \big[ f(s_k) \mid s_0
          \big] \big]
\end{align*}
where step (i) makes use of the dominated convergence theorem, in
particular by noting that $\abss{ \Exs \big[ f(s_k) \mid s_0 \big]}
\leq \vecnorm{f}{\infty} \cdot 2^{1 - k / \mixingtime}$ from
Lemma~\ref{LemWassDecay}.  Consequently, we can write
\begin{align*}
  &\sum_{k = - \infty}^{\infty} \Exs \big[ f(s_k) f(s_0) \big]\\
   & = -
  \Exs [f^2(s_0)] + 2 \Exs \big[ f(s_0) \cdot \greenOp_0 f(s_0) \big]
  \\
& \stackrel{(ii)}{=} - \Exs \big[ \big(\greenOp_0 f(s_0) -
    \transitionOp_0 \greenOp_0 f(s_0) \big)^2 \big] + 2 \Exs \big[
    \big(\greenOp_0 f(s_0) - \transitionOp_0 \greenOp_0 f (s_0) \big)
    \cdot \greenOp_0 f (s_0) \big] \\
  & = \Exs \big[ \big( \greenOp_0 f(s_0) \big)^2 \big] - \Exs \big[
    \big( \transitionOp_0 \greenOp_0 f(s_0) \big)^2 \big],
\end{align*}
where step (ii) follows from equation~\eqref{eq:fA}.

With $\Exs$ denoting expectation over
$X \sim \stationary_0, Y \sim \transition_0 (X, \cdot)$, we have
\begin{align*}
  &\Exs \big[ \big( \greenOp_0 f(Y) - \transitionOp_0 \greenOp_0 f(X)
    \big)^2 \big] \\
    &= \Exs \big[ \big( \greenOp_0 f(s_1) -
    \transitionOp_0 \greenOp_0 f(s_0) \big)^2 \big]\\ &= \Exs \big[
    \big( \greenOp_0 f(s_1) \big)^2 \big] + \Exs \big[ \big(
    \transitionOp_0 \greenOp_0 f(s_0) \big)^2 \big] - 2 \Exs \big[
    \big( \greenOp_0 f(s_1) \big) \cdot \big( \transitionOp_0
    \greenOp_0 f(s_0) \big) \big]\\ &= \Exs \big[ \big( \greenOp_0
    f(s_0) \big)^2 \big] + \Exs \big[ \big( \transitionOp_0 \greenOp_0
    f(s_0) \big)^2 \big] - 2 \Exs \big[ \Exs \big[ \greenOp_0 f(s_1)
      \mid s_0 \big] \cdot \big( \transitionOp_0 \greenOp_0 f(s_0)
    \big) \big]\\ &= \Exs \big[ \big( \greenOp_0 f(s_0) \big)^2 \big]
  - \Exs \big[ \big( \transitionOp_0 \greenOp_0 f(s_0) \big)^2 \big],
\end{align*}
and combining the pieces completes the proof of this lemma.
\qed

\subsection{Proof of Lemma~\ref{lemma:perturbed-matrix-invertibility}} \label{sec:perturb-invert}

By following the derivation of
equation~\eqref{eq:expression-for-derivative-of-expectation-in-lower-bound-proof},
we find that
\begin{align*}
\tfrac{\partial}{\partial \perturb_x (y)} \Lbar^{(\perturb)} =
\stationary_\perturb (x) \transition_\perturb (x, y) \Big \{
\greenOp_\perturb \Lmap (y) - \sum_{z \in \statespace}
\transition_\perturb (x, z) \greenOp \Lmap (z) \Big \}.
\end{align*}
Consequently, for any $u \in \sphere^{d - 1}$, we have the bound
\begin{align*}
&\opnorm{\nabla_w \big( \Lbar^{(\perturb)} u \big)}\\
 & \leq \sup_{z, v
  \in \sphere^{d - 1}} \sqrt{\Exs_{Y \sim \stationary_\perturb} \big[
    \big( z^\top \greenOp_\perturb \Lmap (Y) u \big)^2 \big]} \cdot
\sqrt{\Exs_{X \sim \stationary_\perturb, Y \sim \transition_\perturb
    (X, \cdot)} \big[ \big( ( \greenFunc_0 (Y) - \transitionOp_0
    \greenFunc_0 (X) )^\top v \big)^2 \big]} \\
& \leq \sup_{v \in \sphere^{d - 1}} \sqrt{\Exs_{Y \sim
    \stationary_\perturb} \big[ \vecnorm{ \greenOp_\perturb \Lmap (Y)
      u }{2}^2 \big]} \cdot \frac{3}{2} \sqrt{\Exs_{X \sim
    \stationary_0, Y \sim \transition_0 (X, \cdot)} \big[ \big( (
    \greenFunc_0 (Y) - \transitionOp_0 \greenFunc_0 (X) )^\top v
    \big)^2 \big]}\\ &\leq c \mixingtime \sigmaA \sqrt{d \cdot
  \opnorm{\Lambda}} \log d.
\end{align*}
We thus obtain 
\begin{align*}
    \opnorm{\Lbar^{(\perturb)} - \Lbar^{(0)}} &\leq \sup_{u \in
      \sphere^{d - 1}}\vecnorm{\big( \Lbar^{(\perturb)} - \Lbar^{(0)}
      \big) u}{2} \leq \int_0^1 \sup_{u \in \sphere^{d - 1}}
    \vecnorm{\nabla_w \big( \Lbar^{(s Q w)} u \big) \cdot w}{2} d s \\
& \leq c \mixingtime \sigmaA \sqrt{d \cdot \trace (\Lambda)} \log d
    \cdot \vecnorm{w}{2}.
\end{align*}
Now given a perturbation vector satisfying the bound $\vecnorm{w}{2} \leq \frac{1 - \kappa}{2 c \mixingtime \sigmaA \sqrt{d \cdot \opnorm{\Lambda}} \log d} $, we have the following bound for any $u \in \sphere^{d - 1}$:
\begin{multline*}
    \vecnorm{(I - \Lbar^{(\perturb)}) u}{2} \geq \vecnorm{(I -
      \Lbar^{(0)}) u}{2} - \vecnorm{(\Lbar^{(\perturb)} - \Lbar^{(0)})
      u}{2} \\
      \geq (1 - \kappa) - \opnorm{\Lbar^{(\perturb)} -
      \Lbar^{(0)}} \geq \tfrac{1 - \kappa}{2},
\end{multline*}
which implies that $\opnorm{I - \Lbar^{(\perturb)})^{-1}} \leq
\tfrac{2}{ 1- \kappa}$, as claimed.


\subsection{A useful moment bound}\label{appendix:useful-moment-bound}

Finally, we state and prove a moment bound that is useful in multiple
proofs. Recall that the operator $\transitionOp_\perturb$ is a the
perturbed probability transition kernel under perturbation matrix
$\perturb$, and the operator $\greenOp_\perturb$ is the Green function
operator associated with this transition kernel.

\begin{lemma}
\label{lemma:poisson-eq-estimate}
Consider a bounded function $f: \statespace \rightarrow \real$, and a
perturbation vector $\perturb$ satisfying the condition in
Lemma~\ref{lemma:perturbation-principle-markov}.  There there exists a
universal constant $c > 0$, such that for any integer $p \geq 1$
\begin{align*}
  \big( \Exs_{X \sim \stationary_\perturb} \big[
    \big(\greenOp_\perturb f (X)\big)^{2p} \big] \big)^{\frac{1}{2p}}
  \leq c \, p \, \mixingtime \Big[ \Exs_{X \sim \stationary_\perturb}
    \big[ f(X)^{2p} \big] \Big]^{\frac{1}{2p}} \log \Big \{
  \tfrac{\vecnorm{f}{\infty}^{2p}}{ \Exs_{X \sim \stationary_\perturb}
    \big[ f(X)^{2p} \big]} \Big \}
\end{align*}
\end{lemma}

The proof is similar to that of
Lemma~\ref{lemma:relate-finite-sum-to-sigstar-mkv}. For any function
$f: \statespace \rightarrow \real$ such that
$\Exs_{\stationary_\perturb}[f(X)] = 0$, we first observe that
$\greenOp_\perturb f(s) = \sum_{k = 0}^{\infty}
\transitionOp_\perturb^k f(s)$ for all $\state \in \statespace$.  Note
that Lemma~\ref{lemma:perturbation-principle-markov} guarantees that
the perturbed chain satisfies Assumption~\ref{assume-markov-mixing}
with mixing time $4 \mixingtime$.  By Lemma~\ref{LemWassDecay} and the
coupling definition of total variation distance, for each $t \geq 0$,
there exists a random variable $\widetilde{s}_{t}$ such that
$\widetilde{s}_{t} \mid s_0 \sim \stationary_\perturb$, and $\Prob
\big( \widetilde{s}_{t} \neq s_{t} \mid \state \big) \leq 2^{- \lfloor
  \tfrac{t}{4 \mixingtime} \rfloor}$.

By construction, the state $\widetilde{s}_{t}$ is independent of
$s$. Consequently, we have the equivalence \mbox{$\greenOp_\perturb
  f(s) = \sum_{k = 0}^{\infty} \Exs \big[ f(s_k) - f (\widetilde{s}_k)
    \mid s \big]$,} and for any $\alpha > 0$,
\begin{align*}
 \Exs_{s \sim \stationary_\perturb} \big[ \big( \greenOp_\perturb f(s)
   \big)^{2p} \big] & \leq \big( \sum_{k = 0}^{\infty} e^{2 p \alpha
   t} \Exs \big(\Exs \big[ f(s_k) - f (\widetilde{s}_k) \mid s \big]
 \big)^{2p} \big) \cdot \big( \sum_{k = 0}^{\infty} e^{- \frac{2p}{2p
     - 1} \alpha k} \big)^{2p - 1} \\
 & \leq \alpha^{1 - 2p } \sum_{k = 0}^{\infty} e^{2 p \alpha k} \Exs
 \big[ \abss{f(s_k) - f (\widetilde{s}_k)}^{2p} \big].
\end{align*}
We bound the moment of $f(s_k) - f (\widetilde{s}_k)$ for different
values of $k$ in two ways. On the one hand, Young's inequality
directly leads to the following naive bound
\begin{align*}
    \Exs \big[ \abss{f(s_k) - f (\widetilde{s}_k)}^{2p} \big] \leq
    2^{2p - 1} \big( \Exs \big[ f(s_k)^{2p} \big] + \Exs \big[
      f(\widetilde{s}_k)^{2p} \big] \big) = 2^{2p} \Exs_{s \sim
      \stationary_\perturb} \big[ f (s)^{2p} \big].
\end{align*}
On the other hand, for any bounded function $f$, we have
\begin{align*}
     \Exs \big[ \abss{f(s_k) - f (\widetilde{s}_k)}^{2p} \big] \leq
     \vecnorm{f}{\infty}^{2p}\cdot \Prob \big( s_k \neq
     \widetilde{s}_k \big) \leq \vecnorm{f}{\infty}^{2p} \cdot 2^{1 -
       \frac{k}{4\mixingtime}}.
\end{align*}
Combining the two estimates yields the bound
\begin{align*}
    \Exs \big[ (\greenOp_\perturb f (X))^{2p} \big] & \leq \alpha^{1 -
      2 p} \Big \{ 2^{2p} \cdot e^{2 p \alpha \tau} \tau \Exs_{s \sim
      \stationary_\perturb} \big[ f (s)^{2p} \big] +
    \vecnorm{f}{\infty}^{2p} \sum_{k = \tau + 1}^{\infty} e^{2p \alpha
      k} \cdot 2^{1 - \frac{k}{4\mixingtime}} \Big \},
\end{align*}
valid for any $\alpha > 0$ and $\tau > 0$.  Setting $\tau = c \,
\mixingtime \log \tfrac{\vecnorm{f}{\infty}^{2p}}{ \Exs[f(X)^{2p}]}$
and $\alpha = \tfrac{1}{16 \tau p}$ yields the claim.

\section{Proofs for the examples}
We collect the proofs of the consequences to specific examples in this section.

\subsection{Proofs for TD$(0)$}

We stated three corollaries applicable to this method, and in this
section, we prove each of them in turn.

\subsubsection{Proof of Corollary~\ref{cor:TD-0-finite-state-markov}}
\label{subsec:proof-cor-td-0-finite-state-markov}

The bulk of the proof involves verifying the conditions needed to
apply Proposition~\ref{prop:lsa-markov-iterate-bound} and
Theorem~\ref{thm:markov-main}, but some additional care is needed in
order to deal with non-orthonormal basis functions $(\phi_j)_{j \in
  [d]}$. First, we note that the SA procedure~\eqref{eq:td-0-iterates}
can be equivalently written as
\begin{align}
     \theta_{t + 1} & = (1 - \stepsize \beta) \theta_t + \stepsize
     \beta \Lmap_{t + 1}(\slidingwindow_t) \theta_t - \stepsize \beta
     \bmap_{t + 1}(\slidingwindow_t),\label{eq:td-0-iterates-rewritten}
\end{align}
where $\Lmap_{t + 1} (\slidingwindow_t) \defn \big(I_d - \beta^{-1}
\phi (s_t) \phi (s_t)^\top + \discount \beta^{-1} \phi (s_t) \phi
(s_{t + 1})^\top \big)$, and $\bmap_{t + 1}(\slidingwindow_t) \defn
\beta^{-1} R_t (s_t) \phi (s_t)$.  This is an SA scheme with
stepsize $\stepsize \beta$.

For any matrix $A \in \real^{d \times d}$, define $\kappa (A) \mydefn
\tfrac{1}{2} \lammax \big( A + A^\top \big)$.  We verify the
eigenvalue condition~\eqref{eq:kappa-opnorm} by noting that
\begin{align*}
  \tfrac{1}{2} \lammax \big( \Lbar + \Lbar^\top \big) &= 1 -
  \tfrac{1}{\beta} \, \kappa \big(\discount \Exs_{s \sim \stationary,
    s^+ \sim \transition (s, \cdot)} \big[ \phi(s) \phi (s^+)^\top
    \big] - \Exs_{\stationary} \big[ \phi (s) \phi(s)^\top \big]
  \big)\\
& = 1 - \tfrac{1}{ \beta} \lammax \big( \Bmat^{1/2} \big( I_d -
\tfrac{\SpecMat + \SpecMat^\top}{2} \big) \Bmat^{1/2} \big) \;= \; 1 -
\tfrac{\mu}{\beta} (1 - \kappa) < 1,
\end{align*}
and
\begin{align*}
\opnorm{\Lbar} \leq 1 + \tfrac{1}{\beta} \big(\opnorm{\Exs_{s \sim
    \stationary, s^+ \sim \transition (s, \cdot)} \big[ \phi (s)
    \phi(s^+)^\top \big]} + \opnorm{ \Exs_{\stationary} \big[ \phi (s)
    \phi (s)^\top \big]} \big) \leq 3.
\end{align*}
For the two-step sliding-window Markov chain $\slidingwindow_t =
(\state_t, \state_{t + 1})$, Assumption~\ref{assume-markov-mixing}
holds with mixing time $(\mixingtime + 1)$ in the discrete metric, and
the metric space has diameter at most $1$. It remains to verify the
boundedness and moment assumptions.

In order to verify Assumption~\ref{assume-lip-mapping}, we note that
the bounds~\eqref{eq:td0-bounded-feature-vec} imply that
\begin{align*}
\opnorm{\Lmap_{t + 1} (\state_t)} & \leq 1 + \tfrac{1}{\beta}
\big(\opnorm{ \phi (s_t) \phi (s_{t + 1})} + \opnorm{ \phi (s_t) \phi
  (s_t)^\top} \big) \leq (1 + \varsigma^2 ) d, \quad
\mbox{and}\\ \vecnorm{\bmap_{t + 1} (\state_t)}{2} & \leq
\tfrac{1}{\beta} |R_t (s_t)| \cdot \vecnorm{\phi (s_t)}{2} \leq
\varsigma^2 \sqrt{d / \beta}.
\end{align*}
Turning to the moment assumption, given any vector $u \in
\sphere^{\usedim - 1}$ and coordinate vector $\coordinate_j$, we have
the bounds
\begin{align*}
 \Exs_{s \sim \stationary, s^+ \sim \transition (s, \cdot)} \big[
   \big( \coordinate_j^\top \phi (s) \phi (s^+)^\top u \big)^2 \big] &
 \leq \sqrt{\Exs_{s \sim \stationary} \big[ \big( \coordinate_j^\top
     \phi (s) \big)^4 \big]} \cdot \sqrt{\Exs_{s \sim \stationary}
   \big[ \big( u^\top \phi (s) \big)^4 \big]} \leq \beta^2
 \varsigma^4, \\
 \Exs_{s \sim \stationary} \big[ \big( \coordinate_j^\top \phi (s)
   \phi (s)^\top u \big)^2 \big] & \leq \sqrt{\Exs_{s \sim
     \stationary} \big[ \big( \coordinate_j^\top \phi (s) \big)^4
     \big]} \cdot \sqrt{\Exs_{s \sim \stationary} \big[ \big( u^\top
     \phi (s) \big)^4 \big]} \leq \beta^2 \varsigma^4, \\
\Exs_{s \sim \stationary} \big[ \big( \coordinate_j^\top R_t (s) \phi
  (s) \big)^2 \big] &\leq \varsigma^2 \Exs_{s \sim \stationary} \big[
  \big( \coordinate_j^\top \phi (s) \big)^2 \big] \leq \beta
\varsigma^4.
\end{align*}
Finally, the quantity $\varbound$ from
equation~\eqref{eq:varbound-in-td} is bounded as
\begin{multline*}
\max_{j \in [d]} \Exs \big[ \inprod{\coordinate_j}{( \Lmap_{t + 1}
    (\slidingwindow_t) - \Lbar) \thetabar + (\bmap_{t + 1}
    (\slidingwindow_t) - \bvec) }^2 \big]  \\
\leq \max _{j \in [d]} \sqrt{\Exs \big[ \inprod{\coordinate_j}{\phi
      (s_t)}^4 \big]} \cdot \sqrt{\big( \Exs \big[ \phi (s_t)^\top
    \thetabar - \discount \phi (s_{t + 1})^\top \thetabar - R_t
    (s_t) \big)^4 \big]} \leq \varbound^2.
\end{multline*}
Invoking equation~\eqref{eq:estimation-error-bound-testmat-general}
with the test matrix $\testMat \mydefn \Bmat$ and substituting with
the representation $\valuefunc(\state) =
\inprod{\theta}{\phi(\state)}$ yields the claim.


\subsubsection{Proof of Corollary~\ref{prop:td-0-continuous-statespace}}
\label{subsubsec:proof-prop-td-0-cts-statespace}

We prove this corollary by verifying the assumptions used in our main
theorem. Assumption~\ref{assume:noise-moments} directly follows
from~\eqref{eq:continuous-space-td-fourth-moment} and the boundedness
of reward; Assumption~\ref{assume-markov-mixing} is exactly the
$\Wass_1$ mixing time bound imposed on the Markov chain. In order to
verify that \mbox{$\Lmap (\state, \state^+) = I_d - \beta^{-1} \big( \phi (\state) \phi
  (\state)^\top - \discount \phi (\state) \phi (\state^+)^\top\big)$}
satisfies Assumption~\ref{assume-lip-mapping}, we first note that
\begin{multline*}
\opnorm{\Lmap (\state_1, \state_1^+) - \Lmap (\state_2, \state_2^+)}\\
\leq \tfrac{1}{\beta}\opnorm{\phi(s_1) \phi (s_1)^\top - \phi (s_2) \phi (s_2)^\top} + \tfrac{\discount}{\beta} \opnorm{\phi (s_1) \phi (s_1^+)^\top - \phi (s_2) \phi (s_2^+)^\top}.
\end{multline*}
  By adding and
subtracting terms, we have the bound
\begin{align*}
  \opnorm{\phi(s_1) \phi(s_1)^\top - \phi(s_2) \phi(s_2)^\top} & \leq
  \big \{ \vecnorm{\phi(s_1)}{2} + \vecnorm{\phi(s_2)}{2} \big \} \;
  \vecnorm{\phi(s_1) - \phi(s_2)}{2} \\
  & \stackrel{(i)}{\leq} 2 \varsigma^2 \beta \usedim \vecnorm{s_1 - s_2}{2},
\end{align*}
The step $(i)$ follows from the Lipschitz condition~\eqref{eq:continuous-space-td-lip} and boundedness of the metric space $\statespace$. More precisely, we have $\vecnorm{\phi(s_1) - \phi(s_2)}{2} \leq \varsigma \sqrt{\beta \usedim}\vecnorm{s_1 - s_2}{2}$ and $\vecnorm{\phi (s_1)}{2} =\vecnorm{\phi (s_1) - \phi (0)}{2} \leq  \varsigma \sqrt{\beta \usedim}$. A
similar argument yields that
\begin{align*}
\opnorm{\phi (s_1) \phi (s_1^+)^\top - \phi (s_2) \phi (s_2^+)^\top} &
\leq  \varsigma^2 \usedim \big( \vecnorm{s^+_1 - s^+_2}{2} + \vecnorm{s_1 - s_2}{2} \big).
\end{align*}
Putting together the pieces, we have shown that the mapping
$\Lmap: \statespace \rightarrow \real^{d \times d}$ is $3 \varsigma^2
d$-Lipschitz with respect to the metric $\metric \big( (s_1, s_1^+),
(s_2, s_2^+) \big) = \vecnorm{s_1 - s_2}{2} + \vecnorm{s_1^+ -
  s_2^+}{2}$.

Similarly, for the vector observation $\bmap_t (\state) = R_t (\state)
\phi (\state)$, we note that for any $s_1, s_2 \in \statespace$,
\begin{align*}
    \vecnorm{\bmap_t (\state_1)  - \bmap_t (\state_2) }{2} &\leq \abss{R_t (s_1) - R_t (s_2)} \cdot \vecnorm{\phi (s_1)}{2} + \abss{R_t (s_2)} \cdot \vecnorm{\phi (s_1) - \phi (s_2)}{2} \\
    &\leq 2 \varsigma \sqrt{d / \beta} \vecnorm{\phi (s_1) - \phi (s_2)}{2},
\end{align*}
which shows that $\bvec: \statespace \rightarrow \real^{d / \beta}$ is $2
\varsigma^2 \sqrt{d}$-Lipschitz. Having verified the assumptions, we
complete the proof by following the same steps as in the proof as
Corollary~\ref{cor:TD-0-finite-state-markov}.


\subsubsection{Proof of Corollary~\ref{prop:sieve-td-results}}
\label{subsubsec:proof-prop-sieve-td-results}

In order to verify that Assumption~\ref{assume-lip-mapping} holds with
respect to the discrete metric, note that for any $d_n \geq 1$, we have
$\vecnorm{\bmap_t (s)}{2} \leq \tfrac{\varsigma}{\beta} \sqrt{\sum_{j
    = 1}^{d_n} \phi^2_j(s)} \leq \tfrac{\varsigma^2}{\beta}
\sqrt{d_n}$, and
\begin{align*}
\opnorm{\Lmap(s_1, s_2)} & \leq 1 + \tfrac{1}{\beta} \sum_{j =
  1}^{d_n} \phi_j^2(s_1) + \tfrac{1}{\beta} \sqrt{\sum_{j = 1}^{d_n}
  \phi_j^2(s_1)} \cdot \sqrt{\sum_{j = 1}^{d_n} \phi_j^2(s_2)} \leq
\tfrac{1 + \varsigma^2}{\beta} d_n.
\end{align*}

Turning to the moment condition, let $\Exs$ denote expectation over a
pair $s \sim \stationary$ and \mbox{$s^+ \sim \transition (s,
  \cdot)$.}  Then for any vector $u \in \sphere^{d_n - 1}$ and index
$j \in [d_n]$, we have
\begin{align*}
&\Exs \big[\inprod{\coordinate_j}{\Lmap (s, s^+) u}^2 \big]  \\
&\leq 3 +
\tfrac{3}{\beta^2} \Exs \Big[ \big( \inprod{\coordinate_j}{\phi(s)} \:
  \inprod{\phi(s^+)}{u} \big)^2 \Big] + \tfrac{3}{\beta^2} \Exs \Big[
  \big( \inprod{\coordinate_j}{\phi(s)} \: \inprod{\phi(s)}{u} \big)^2
  \Big]\\
  & \leq 3 + \tfrac{6}{\beta^2} \|\phi_j\|_\infty^2 \cdot \Exs \big[
  \inprod{\phi(s)}{u}^2 \big] \leq 3 + \tfrac{6}{\beta} \varsigma^2.
\end{align*}
For each $t = 1, 2, \ldots$, we also have $\Exs \big[
  \inprod{\coordinate_j}{\bmap_{t + 1} (s_t)}^2 \big] \leq
\tfrac{1}{\beta^2} \|R_t\|_\infty^2 \cdot \Exs_{s \sim \stationary}
\big[ \phi_j (s)^2 \big] \leq \tfrac{\varsigma^2}{\beta}$, which is an
order-one quantity.  Following the same steps as in the proof as
Corollary~\ref{cor:TD-0-finite-state-markov} then yields the claim.

 
\subsection{Proofs for TD$(\lambda)$}

We first prove
Proposition~\ref{prop:mixing-td-lambda-augmented-chain}---the mixing
time result---and then use it to establish
Corollary~\ref{cor:TD-lambda-finite-state-markov}.
 
\subsubsection{Proof of Proposition~\ref{prop:mixing-td-lambda-augmented-chain}}
\label{subsubsec:proof-mixing-td-lambda-augmented-chain}

We prove the claim via a coupling argument. Consider two initial
states $\slidingwindow_0 = (s_0, s_1, h_0)$ and $\slidingwindow_0' =
(s_0', s_1', h_1')$. By Assumption~\ref{assume-markov-mixing} (mixing
time) for the original chain in total variation distance, there exists
a coupling between a chains $(s_t)_{t \geq 1}$ and $(s_t')_{t \geq 1}$
starting from $s_1$ and $s_1'$ respectively, such that $\Prob \big(
s_{(k + 1) \mixingtime + 1} \neq s_{(k + 1) \mixingtime + 1}' \mid
\{s_t, s_t'\}_{t = 1}^{k \mixingtime +1} \big) \leq \tfrac{1}{2}$.
Furthermore, whenever $s_t = s_t'$ for some $t \geq 1$, the two
processes are always identical from then on. Let $(g_t)_{t \geq 0}$
and $(g_t')_{t \geq 0}$ be the eligibility trace
process~\eqref{eq:td-lambda-g-update} associated to $(s_t)_{t \geq 0}$
and $(s_t')_{t \geq 0}$, respectively, and let $h_t = \tfrac{1 -
  \lambda \discount}{\varsigma \sqrt{\beta d}} g_t$ and $h_t' =
\tfrac{1 - \lambda \discount}{\varsigma \sqrt{\beta d}} g_t'$.
 
Under this coupling, we note that $\Prob \big( s_{3 \mixingtime + 1}
\neq s_{3 \mixingtime + 1}' \big) \leq \tfrac{1}{8}$.  Conditioning on
the event \mbox{$\Event \mydefn \big\{ s_{3 \mixingtime + 1} = s_{3
    \mixingtime + 1}' \big\}$,} for any $t \geq 3 \mixingtime + 1$, we
have
\begin{align}
\label{eq:td-lambda-proof-h-contraction}  
\vecnorm{h_{t + 1} - h_{t + 1}'}{2} = \discount \lambda \vecnorm{h_t -
  h_t'}{2} = \cdots = (\discount \lambda)^{t - 3 \mixingtime - 1}
\vecnorm{h_{3 \mixingtime + 1} - h_{ 3 \mixingtime + 1}'}{2}.
\end{align}
We split the remainder of the proof into two cases.

\paragraph{Case I: $s_1 \neq s_1'$:}
The coupling bound implies that $\Prob (\Event) \geq \tfrac{7}{8}$. On
the event $\Event$, for $\tau \geq 3 \mixingtime + 1 + \tfrac{4}{1 -
  \discount \lambda}$, we have the bound $\vecnorm{h_{t + 1} - h_{t +
    1}'}{2} \leq \tfrac{1}{16} \vecnorm{h_{3 \mixingtime + 1} - h_{3
    \mixingtime + 1}'}{2} \leq \tfrac{1}{8}$ almost surely. Under this
coupling, we may write
\begin{align*}
  \Exs \big[ \metric \big( (\state_\tau, \state_{\tau + 1}, h_\tau),
    (\state_\tau', \state_{\tau + 1}' , h_\tau) \big) \big] &=
  \tfrac{1}{4} \big(\Prob \big(\state_\tau \neq \state_\tau' \big) +
  \Prob \big(\state_{\tau + 1} \neq \state_{\tau + 1}' \big) + \Exs
  \big[ \vecnorm{h_{\tau} - h_\tau'}{2} \big] \big)\\ &\leq
  \tfrac{3}{4} \Prob (\Event^c) + \tfrac{1}{4} \Exs \big[
    \vecnorm{h_{\tau} - h_\tau'}{2} \mid \Event \big]\\ &\leq
  \tfrac{1}{8} = \tfrac{1}{2} \cdot \tfrac{1}{4} \bm{1}_{s_1 \neq
    s_1'} \leq \tfrac{1}{2} \metric \big( (\state_0, \state_{1}, h_0),
  (\state_0', \state_{1}' , h_0) \big),
\end{align*}
which proves the Wasserstein contraction in this case.

\paragraph{Case II: $s_1 = s_1'$} In this case, the coupling
construction ensures that $s_t = s_t'$ for any $t \geq 1$. Invoking
the bound~\eqref{eq:td-lambda-proof-h-contraction} then yields
\begin{align*}
  \Exs \big[ \metric \big( (\state_\tau, \state_{\tau + 1}, h_\tau),
    (\state_\tau', \state_{\tau + 1}' , h_\tau) \big) \big] =
  \tfrac{1}{4} \Exs \big[ \vecnorm{h_{\tau} - h_\tau'}{2} \big] \leq
  \tfrac{1}{8} \vecnorm{h_{0} - h_0'}{2} \leq \tfrac{1}{2} \metric
  (\slidingwindow_0, \slidingwindow_0'),
\end{align*}
which establishes contraction in this case.  Combining the two cases
proves the proposition.


\subsubsection{Proof of Corollary~\ref{cor:TD-lambda-finite-state-markov}}
\label{subsec:proof-cor-td-lambda-finite-state-markov}

We note that the SA procedure~\eqref{eq:td-lambda-iterates} can be
written as
\begin{align*}
     \theta_{t + 1} &= (1 - \stepsize \beta) \theta_t + \stepsize\beta
     \Lmap_{t + 1}(\slidingwindow_t) \theta_t - \stepsize \beta
     \bmap_{t + 1}(\slidingwindow_t),
\end{align*}
where $\Lmap_{t+1}(\slidingwindow_t) = \big(I_d - \tfrac{1}{\beta}
g_t\phi (s_t)^\top + \discount \tfrac{1}{\beta} g_t \phi (s_{t +
  1})^\top \big)$ and $\bmap_{t+1}(\slidingwindow_t) =\tfrac{1}{\beta}
R_t (s_t) g_t$.

Recalling that $\SpecMat_\lambda = (1 - \lambda) \discount \sum_{t =
  0}^{\infty} \lambda^t \discount^{t + 1} \Bmat^{-1/2} \Exs \big[ \phi
  (s_0) \phi (s_{t + 1})^\top \big] \Bmat^{-1/2}$, we first study the
eigenvalues of the symmetrized version of $\SpecMat_\lambda$, and
relate these back to those of \mbox{$\Lbar =
  \Exs_{\widetilde{\stationary}} \big[ \Lmap_{t + 1}(\slidingwindow_t)
    \big]$.}  Note that by the Cauchy--Schwarz inequality, for any
vector $u \in \sphere^{d - 1}$, we have
\begin{align*}
    u^\top \Bmat^{-1/2} \Exs \big[ \phi(s_0) \phi (s_t)^\top \big]
    \Bmat^{-1/2} u \leq \sqrt{\Exs \big[\big( u^\top \Bmat^{-1/2} \phi
        (s_0) \big)^2 \big] \cdot \Exs \big[\big( u^\top \Bmat^{-1/2}
        \phi (s_t) \big)^2 \big]} = 1.
\end{align*}
We therefore have the bound $\tfrac{1}{2} \lambda_{\min}
(\SpecMat_\lambda + \SpecMat_\lambda^\top) \leq (1 - \lambda)
\discount \sum_{t = 0}^{\infty} (\discount \lambda)^{t } = \tfrac{(1 -
  \lambda) \discount}{1 - \lambda \discount}$.  As in the proof of
Corollary~\ref{cor:TD-0-finite-state-markov}, we can deduce that
\begin{align*}
    \tfrac{1}{2} \lammax \big( \Lbar + \Lbar^\top \big) = \tfrac{1}{
      \beta} \lammax \big( \Bmat^{1/2} \big( \tfrac{\SpecMat_\lambda +
      \SpecMat_\lambda^\top}{2} \big) \Bmat^{1/2} \big) \geq \tfrac{(1
      - \lambda) \discount}{1 - \lambda \discount}.
\end{align*}

Next, we verify Assumption~\ref{assume:noise-moments} on the noise
moments. By the update rule~\eqref{eq:td-lambda-g-update}, under a
stationary trajectory, we have the expression $g_t = \sum_{k =
  0}^{\infty} (\discount \lambda)^k \phi (s_{t - k})$.  For any $u \in
\sphere^{d-1}$, invoking H\"{o}lder's inequality yields
\begin{align*}
\Exs \big[ \inprod{g_t}{u}^4 \big] \leq \big( \sum_{k = 0}^{\infty}
(\discount \lambda)^{k} \big)^3 \cdot \sum_{k = 0}^{\infty} (\discount
\lambda)^{k} \Exs \big[ \inprod{u}{\phi (s_{t - k})}^4 \big] \leq
\beta^2 \big( \tfrac{\varsigma}{1 - \discount \lambda} \big)^4.
\end{align*}
In other words, for all standard basis vectors $e_j$, we have
\begin{align*}
  \Exs \big[ \inprod{\coordinate_j}{\Lmap_{t + 1} (\slidingwindow_t)
      u}^2 \big] & \leq 1 + \tfrac{2}{\beta^2} \sqrt{\Exs \big[
      \inprod{\coordinate_j}{\phi (\state_t)}^4 \big]} \cdot \sqrt{
    \Exs \big[ \inprod{g_t}{u}^4 \big] } \leq 1 + 2
  \tfrac{\varsigma^4}{(1 - \discount \lambda)^2},\\
\Exs \big[ \inprod{\coordinate_j}{\bmap_{t + 1} (\slidingwindow_t)
    u}^2\big] & \leq \tfrac{\varsigma^2}{\beta^2} \Exs \big[
  \inprod{g_t}{\coordinate_j}^2 \big] \leq \tfrac{\varsigma^4}{\beta(1
  - \discount \lambda)^2}.
\end{align*}
It remains to verify Assumption~\ref{assume-lip-mapping}. Note that
for any pair $\slidingwindow = (s, s_+, h)$ and $\slidingwindow'= (s',
s_+', h')$, the operator norm $T \defn \opnorm{\Lmap_{t + 1}
  (\slidingwindow) - \Lmap_{t + 1} (\slidingwindow')}$ is almost
surely upper bounded as
\begin{align*}
T & \leq \tfrac{\varsigma \sqrt{d / \beta}}{1 - \lambda \discount}
\cdot \big( \opnorm{h^\top \phi (s) - (h')^\top \phi (s')} +
\opnorm{h^\top \phi (s_+) - (h')^\top \phi (s_+')} \big)\\
  & \leq \tfrac{\varsigma \sqrt{d / \beta}}{1 - \lambda \discount}
\cdot \big( \opnorm{ (h - h')^\top \phi (s')} + \opnorm{ h^\top (\phi
  (s') - \phi (s))} \big) \\
  &\qquad \qquad + \tfrac{\varsigma \sqrt{d / \beta}}{1 - \lambda \discount}
\cdot \big( \opnorm{ (h - h')^\top \phi (s_+')} + \opnorm{
  h^\top (\phi (s_+') - \phi (s_+))} \big)\\
  & \leq \tfrac{2 \varsigma^2 d}{1 - \lambda \discount} \big( \bm{1}_{s
    \neq s'} + \bm{1}_{s_+ \neq s_+'} + \vecnorm{h - h'}{2} \big) =
  \tfrac{8 \varsigma^2 d}{1 - \lambda \discount} \metric
  (\slidingwindow, \slidingwindow').
\end{align*}

Finally, we note that the quantity $\varbound$ defined in
equation~\eqref{eq:varbound-in-td} satisfies the bound
\begin{align*}
    &\sup_{j \in [d]} \Exs \big[ \inprod{\coordinate_j}{( \Lmap_{t +
        1} (\slidingwindow_t) - \Lbar) \thetabar + (\bmap_{t + 1}
      (\slidingwindow_t) - \bvec) }^2 \big]\\ &\leq \sup_{j \in [d]}
  \sqrt{\Exs \big[ \inprod{\coordinate_j}{g_t}^4 \big]} \cdot
  \sqrt{\big( \Exs \big[ \phi (s_t)^\top \thetabar - \discount \phi
      (s_{t + 1})^\top \thetabar - R_t (s_t) \big)^4 \big]} \leq
  \tfrac{ \varbound^2}{(1 - \discount \lambda)^2}.
\end{align*}
Invoking equation~\eqref{eq:estimation-error-bound-testmat-general},
with the test matrix $\testMat \mydefn \Bmat$ and substituting the
expression \mbox{$\valuefunc(\state) = \inprod{\theta}{\phi(\state)}$}
yields the claim.


\subsection{Proofs for vector autoregressive estimation}

In this section, we present proofs of results on vector autoregressive
models, as introduced in Example~\ref{example:autoregressive}.

\subsubsection{Proof of Proposition~\ref{prop:lyap-stable}}
\label{subsec:proof-prop-lyap-stable}

We prove the claim by a direct construction of the coupling. Given two
initial points \mbox{$\slidingwindow_0 =  \begin{bmatrix} X_1^\top, X_0^\top, \cdots, X_{- k + 1}^\top \end{bmatrix} ^\top$} and \mbox{$\slidingwindow_0' =  \begin{bmatrix} X_1'^\top, X_0'^\top, \cdots, X_{- k + 1}'^\top \end{bmatrix}^\top$,} we consider a pair of stochastic processes $(X_t)_{t \geq
  1}$ and $(X_t')_{t \geq 1}$ starting from $\slidingwindow_0$ and
$\slidingwindow'$, respectively, driven by the same noise process
$(\varepsilon_t)_{t \geq 0}$.  Introduce the shorthand
$Y_{t+1}= \begin{bmatrix} X_{t+1} & \ldots &
  X_{t-k+2} \end{bmatrix}^\top$ (note that $Y_{t + 1}$ is a sliding window with length one unit shorter than $\slidingwindow_t$). We have:
\begin{align*}
\vecnorm{Y_{t+1} - Y'_{t+1}}{\Pstar}^2 \; = \; \vecnorm{ \Rstar \,
  (Y_t - Y'_t)}{\Pstar}^2 & = \vecnorm{Y_t - Y'_t}{\Pstar}^2 -
\vecnorm{Y_t - Y'_t}{\Qstar}^2 \\
& \leq \big( 1 - \tfrac{\mu}{\beta} \big) \vecnorm{Y_t -
  Y'_t}{\Pstar}^2.
\end{align*}
Consequently, the augmented processes $\slidingwindow_t = (X_{t + 1},
X_t, \cdots, X_{t - k + 1})$ and $\slidingwindow_t' = (X_{t + 1}',
X_t', \cdots, X_{t - k + 1}')$ satisfy the bound
\begin{multline*}
    \vecnorm{\slidingwindow_t - \slidingwindow_t'}{2} \leq \vecnorm{Y_{t + 1} - Y_{t + 1}'}{2} + \vecnorm{Y_{t} - Y_{t}'}{2}
    \leq
    \tfrac{1}{\sqrt{\lambda_{\min}(\Pstar)}} \big(\vecnorm{Y_{t+1} -
      Y'_{t+1}}{\Pstar} + \vecnorm{ Y_t - Y'_t}{\Pstar} \big) \\
      \leq 2\sqrt{\tfrac{\lammax(\Pstar)}{\lambda_{\min}(\Pstar)}} \big(1 -
    \tfrac{\mu}{2 \beta} \big)^{t} \vecnorm{\slidingwindow_0 -
      \slidingwindow_0'}{2}
\end{multline*}
Note that since $\Pstar \succeq \Qstar$, we have
$\lambda_{\min} (\Pstar) \geq \lambda_{\min} (\Qstar) = \mu$. Taking
$\mixingtime = c \tfrac{\beta}{\mu} \big(1 + \log \tfrac{\beta}{\mu}\big)$ yields
the contraction bound $\vecnorm{\slidingwindow_{\mixingtime} -
  \slidingwindow_{\mixingtime}'}{2} \leq \tfrac{1}{2} \vecnorm{\slidingwindow_0
  - \slidingwindow_0'}{2}$.  Taking expectations on both sides
completes the proof.


\subsubsection{Proof of Corollary~\ref{cor:auto-regressive}}

We begin by showing norm bounds and moment bounds on the process
$(X_t)_{t \geq 0}$.  By definition~\eqref{eq:autoregressive-process}
of the process and stability, the block vector $Y_t
\defn \begin{bmatrix} X_t & X_{t - 1} & \cdots & X_{t - k +
    1} \end{bmatrix}^\top$ satisfies the recursion $Y_t = \sum_{i =
  0}^{\infty} \Rstar^i \varepsilon_{t-i} e_1$, where $e_1$ is the
standard block basis vector equal to identify on the first block.  We
therefore have the bound
\begin{align*}
    \vecnorm{X_t}{2} \leq \tfrac{1}{\mu} \vecnorm{Y_t}{\Pstar} \leq
    \sum_{i = 0}^{\infty} \vecnorm{\Rstar^i \varepsilon_{t-i}
      e_1}{\Pstar} \leq \tfrac{1}{\mu} \sum_{i = 0}^{\infty} \big( 1 -
    \tfrac{\mu}{\beta} \big)^i \vecnorm{\varepsilon_{t -i} e_1}{\Pstar}
    \leq \tfrac{\beta^2}{\mu^2} \varsigma \sqrt{m}.
\end{align*}
Moreover, for each $u \in \sphere^{\vardim - 1}$, we have
\begin{align*}
\Exs \big[ \inprod{X_t}{u}^4 \big] & \leq \big( \sum_{i = 0}^{\infty}
e^{- \tfrac{i \mu}{6 \beta}} \big)^3 \cdot \sum_{i = 0}^{\infty}
e^{\tfrac{i \mu}{2 \beta}} \; \Exs \big[ \inprod{\Rstar^i
    \varepsilon_{t-i} e_1}{u e_1}^4 \big] \\
& \leq c \big( \beta / \mu \big)^3 \cdot \sum_{i = 0}^{\infty}
e^{\tfrac{i \mu}{2 \beta}} \cdot \tfrac{\beta^4}{\mu^4} \cdot e^{-
  \tfrac{i \mu}{ \beta}} \varsigma^4 \; \leq \; c' \big(
\tfrac{\beta^2 \varsigma}{\mu^2} \big)^4.
\end{align*}
Next, we proceed with verifying the assumptions used in
Theorem~\ref{thm:markov-main}. Letting $\nu \mydefn 1 / \opnorm{H^*}$,
the stochastic approximation procedure can be rewritten as
\begin{multline*}
 \theta_{t + 1} = (1 - \tfrac{\stepsize}{\nu}) \theta_t\\
  +
 \tfrac{\stepsize}{\nu} \big( \theta_t - \nu \big( \big[ X_{t - j}
   X_{t + 1 - i}^\top \big]_{i, j \in [\vardim]} \otimes I_\vardim
 \big) \theta_t + \nu \cdot \vectorize \big( \begin{bmatrix} X_{t + 1}
   X_{t}^\top & \cdots & X_{t + 1} X_{t - k + 1}^\top
    \end{bmatrix}\big) \big).
\end{multline*}
Observe that the matrix $\Lbar \defn I_{k \vardim^2} - \nu H^* \otimes
I_m$ satisfies the eigenvalue bound
\begin{align*}
  \tfrac{1}{2} \lammax (\Lbar + \Lbar^\top) \leq 1 -
  \tfrac{\nu}{2}\lambda_{\min}(H^* + (H^*)^\top) \leq 1 - \nu h^*.
\end{align*}
On the other hand, the empirical observations satisfy the almost-sure
bounds
\begin{align*}
    \opnorm{\Lmap_{t + 1} (\slidingwindow_t) - \Lbar} \leq \nu \cdot
    \opnorm{\big[ X_{t - j} X_{t + 1 - i}^\top \big]_{i, j \in
        [\vardim]}} \leq \nu \cdot \tfrac{\beta^4}{\mu^4} \varsigma^2
    \vardim k \text{ and } \\
\opnorm{\bmap_{t + 1} (\slidingwindow_t) - \bbar} \leq \nu \cdot
\matsnorm{ \begin{bmatrix} X_{t + 1} X_{t}^\top & \cdots & X_{t + 1}
    X_{t - k + 1}^\top
    \end{bmatrix}}{F} \leq \nu  \cdot \tfrac{\beta^4}{\mu^4} \varsigma^2 \vardim \sqrt{k}.
\end{align*}
For two collections of matrices $\mathcal{U} = \big( U^{(j)} \big)_{j
  = 1}^k$ and $\mathcal{V} = \big( V^{(j)} \big)_{j = 1}^k \subseteq
\real^{\vardim \times \vardim}$ such that $\sum_{j = 1}^k
\matsnorm{U^{(j)}}{F}^2 = \sum_{j = 1}^k \matsnorm{V^{(j)}}{F} = 1$,
the corresponding moment can be bounded as
\begin{align*}
 \Exs \big[ \inprod{\vectorize (\mathcal{U})}{\big( \Lmap_{t + 1}
     (\slidingwindow_t) - \Lbar \big) \vectorize (\mathcal{V})}^2
   \big] & \leq \nu^2 \Exs \big[ \big( \sum_{\ell = 0}^{k - 1}
   \inprod{U^{(\ell)}}{\sum_{j = 0}^{k - 1} V^{(j)} X_{t - j} X_{t -
       \ell}^\top}_F \big)^2 \big],
\end{align*}
which is in turn at most
\begin{multline*}
\nu^2 k^2 \sum_{\ell = 0}^{k - 1} \sum_{j = 0}^{k - 1} \sqrt{\Exs
  \big[ X_{t - \ell}^{\otimes 4} \big] \big[ (U^{(\ell)})^\top
    U^{(\ell)}, (U^{(\ell)})^\top U^{(\ell)} \big]} \\
    \times \sqrt{\Exs
  \big[ X_{t - j}^{\otimes 4} \big] \big[ (V^{(j)})^\top V^{(j)},
    (V^{(j)})^\top V^{(j)} \big]}.
\end{multline*}
In order to bound this last quantity, we let $(U^{(\ell)})^\top
U^{(\ell)} = \sum_{i = 1}^\vardim \lambda_i^2 u_i u_i^\top$ be its
singular value decomposition, and note that
\begin{multline*}
  \Exs \big[ X_{t - \ell}^{\otimes 4} \big] \big[ (U^{(\ell)})^\top
    U^{(\ell)}, (U^{(\ell)})^\top U^{(\ell)} \big] = \Exs \big[ X_{t -
      \ell}^{\otimes 4} \big] \big[ \sum_{i = 1}^\vardim \lambda_i^2
    u_i u_i^\top, \sum_{i = 1}^\vardim \lambda_i^2 u_i u_i^\top
    \big]\\ = \sum_{i, i'} \Exs \big[ X_{t - \ell}^{\otimes 4} \big]
       [u_i, u_i, u_{i'}, u_{i'}] \cdot \lambda_i^2 \lambda_{i'}^2
       \leq c' \big( \tfrac{\beta^2 \varsigma}{\mu^2} \big)^4
       \big(\sum_{i} \lambda_i^2 \big)^2 = c' \big( \tfrac{\beta^2
         \varsigma}{\mu^2} \big)^4 \matsnorm{U^{(\ell)}}{F}^2.
\end{multline*}
Putting together the pieces, we have
\begin{align*}
  &\Exs \big[ \inprod{\vectorize (\mathcal{U})}{\big( \Lmap_{t + 1}
      (\slidingwindow_t) - \Lbar \big) \vectorize (\mathcal{V})}^2
    \big] \\
    &\leq \nu^2 k^2 c' \big( \tfrac{\beta^2 \varsigma}{\mu^2}
  \big)^4 \cdot \sum_{\ell = 0}^{k - 1} \sum_{j = 0}^{k - 1}
  \matsnorm{U^{(\ell)}}{F}^2 \matsnorm{V^{(j)}}{F}^2 \leq c \big( \nu
  \cdot \tfrac{\beta^4 k \varsigma^2}{\mu^4} \big)^2.
\end{align*}
Similarly, we can prove analogous moment bounds on $\bmap_{t + 1}
(\slidingwindow_t)$. In particular, for indices $\ell \in [k]$ and
$i,j \in [\vardim]$, we consider the coordinate direction of the $(i,
j)$ entry in the $\ell$-th matrix to deduce that
\begin{multline*}
  \Exs \big[ \inprod{\coordinate_{\ell, i, j}}{ (\bmap_{t + 1}
      (\slidingwindow_t) - \bbar)}^2 \big] \leq \nu^2 \Exs \big[
    \inprod{e_i e_j^\top}{ X_{t + 1} X_{t - \ell + 1} }^2 \big] \\
  \leq \nu^2 \sqrt{\Exs \big[ \inprod{e_j^\top}{ X_{t + 1} }^4 \big]}
  \cdot \sqrt{\Exs \big[ \inprod{e_i^\top}{ X_{t - \ell + 1} }^4
      \big]} \leq c' \big( \nu \cdot \tfrac{\beta^2 \varsigma}{\mu^2}
  \big)^4.
\end{multline*}
Applying Theorem~\ref{thm:markov-main} completes the proof of this
corollary.

\end{document}